\documentclass[twoside, 10pt]{article}

\usepackage[hmargin= 0.9in, twosideshift=0in,height=9.3in]{geometry}
\usepackage{amssymb,amsthm, amsfonts,amsmath,mathrsfs, amsxtra, url, hyperref}

\usepackage{fancyhdr}
\fancypagestyle{plain}{
 \fancyhf{}

}

\renewcommand{\theequation}{\thesection.\arabic{equation}}
 \numberwithin{equation}{section}

\newtheorem {thm}{Theorem}[section]
\newtheorem {prop}{Proposition}[section]
\newtheorem {lemm}{Lemma}[section]
\newtheorem {deff}{Definition}[section]
\newtheorem {cor}{Corollary}[section]
\newtheorem {rem}{Remark}[section]

\newtheorem {eg}{Example}[section]

\def\ba{\begin{array}}
\def\ea{\end{array}}
\def\bea{\begin{eqnarray}}
\def\eea{\end{eqnarray}}
\def\beas{\begin{eqnarray*}}
\def\eeas{\end{eqnarray*}}
\def\bi{\begin{itemize}}
\def\ei{\end{itemize}}

\def\a{\alpha}
\def\g{\gamma}
\def\d{\delta}
\def\e{\varepsilon}
\def\z{\zeta}
\def\k{\kappa}
\def\l{\lambda}

\def\si{\sigma}

\def\t{\tau}
\def\th{\theta}
\def\o{\omega}
\def\f{\phi}
\def\vf{\varphi}
\def\p{\psi}

\def\D{\Delta}
\def\G{\Gamma}
\def\L{\Lambda}
\def\O{\Omega}
\def\F{\Phi}

\def\Th{\Theta}
\def\U{\Upsilon}

\def\bF{{\bf F}}

\def\cA{{\cal A}}

\def\cE{{\cal E}}
\def\cF{{\cal F}}
\def\cG{{\cal G}}

\def\cI{{\cal I}}

\def\cK{{\cal K}}
\def\cL{{\cal L}}

\def\cN{{\cal N}}
\def\cO{{\cal O}}
\def\cP{{\cal P}}

\def\cR{{\cal R}}
\def\cS{{\cal S}}

\def\cV{{\cal V}}

\def\cX{{\cal X}}
\def\cY{{\cal Y}}
\def\cZ{{\cal Z}}

\def\hB{\mathbb{B}}
\def\hC{\mathbb{C}}

\def\hE{\mathbb{E}}

\def\hH{\mathbb{H}}

\def\hK{\mathbb{K}}
\def\hL{\mathbb{L}}

\def\hN{\mathbb{N}}

\def\hQ{\mathbb{Q}}
\def\hR{\mathbb{R}}
\def\hS{\mathbb{S}}

\def\hV{\mathbb{V}}

\def\sB{\mathscr{B}}

\def\sD{\mathscr{D}}

\def\sK{\mathscr{K}}
\def\sL{\mathscr{L}}

\def\sN{\mathscr{N}}

\def\sP{\mathscr{P}}

\def\sY{\mathscr{Y}}

\def\fF{\mathfrak{F}}

\def\fL{\mathfrak{L}}
\def\fM{\mathfrak{M}}
\def\fN{\mathfrak{N}}

\def\fq{\mathfrak{q}}

\def\({\textnormal{(}}
\def\){\textnormal{)}}
\def\[{[\neg[}
\def\]{]\neg]}
\def\lan{\langle}
\def\ran{\rangle}

\def\no{\noindent}

\def\ss{\smallskip}
\def\ms{\medskip}
\def\bs{\bigskip}
\def\q{\quad}
\def\qq{\qquad}

\def\neg{\negthinspace}
\def\dneg{\neg \neg}
\def\tneg{\neg \neg \neg}

\def\ol{\overline}
\def\ul{\underline}
\def\ua{\mathop{\uparrow}}
\def\da{\mathop{\downarrow}}

\def\wt{\widetilde}
\def\wh{\widehat}

\def\dtp{{\hbox{$dt \otimes dP$-a.e.}}}
\def\pas{{\hbox{$P$-a.s.}}}

\def\hb{\hbox}
\def\dis{\displaystyle}
\def\cd{\cdot}
\def\cds{\cdots}

\def\fa{\,\forall \,}
\def\pa{\partial}
\def\es{\emptyset}

\def\dfnn{\stackrel{\triangle}{=}}
\def\b1{{\bf 1}}
\def\qed{\hfill $\Box$ \medskip}

\def\esssup{\mathop{\rm esssup}}
\def\liminf{\mathop{\ul{\rm lim}}}
\def\limsup{\mathop{\ol{\rm lim}}}
\newcommand{\esup}[1]{ \underset{#1}{\esssup}\,}

\newcommand{\lsup}[1]{ \underset{#1}{\limsup}}
\newcommand{\linf}[1]{ \underset{#1}{\liminf}}
\newcommand{\lmt}[1]{ \underset{#1}{\lim}}
\newcommand{\lmtu}[1]{ \underset{#1}{\lim} \neg \ua \,}
\newcommand{\lmtd}[1]{ \underset{#1}{\lim} \neg \da \,}

\begin{document}

\title{\bf  Quadratic Reflected BSDEs  with Unbounded Obstacles}
\author{
Erhan Bayraktar\thanks{ \noindent Department of
  Mathematics, University of Michigan, Ann Arbor, MI 48109; email:
{\tt erhan@umich.edu}. }  \thanks{E. Bayraktar is supported in part by the National Science Foundation under an applied mathematics research grant and a Career grant, DMS-0906257 and DMS-0955463, respectively, and in part by the Susan M. Smith Professorship.} $\,\,$,
$~~$Song Yao\thanks{
\noindent Department of
  Mathematics, University of Michigan, Ann Arbor, MI 48109; email: {\tt songyao@umich.edu}. } }
\date{ }

\maketitle

 \begin{abstract}

   In this paper,
       we analyze a real-valued reflected backward stochastic differential equation (RBSDE)
    with an unbounded obstacle and an unbounded terminal condition when its generator $f$ has quadratic growth in the $z$-variable. In particular,
    we  obtain existence, comparison, and stability results, and consider the optimal stopping for quadratic $g$-evaluations.   As an application of our results we analyze the obstacle problem for semi-linear parabolic PDEs in which the non-linearity appears as the square of the gradient. Finally, we prove a comparison theorem for these obstacle problems when the generator is convex or concave in the $z$-variable.

\end{abstract}

 \ms   {\bf Keywords: }\:Quadratic reflected backward stochastic differential equations,  convex/concave generator, $\th$-difference method,
 Legenre-Fenchel duality, optimal stopping problems for quadratic $g$-evaluations, stability,  obstacle problems for semi-linear parabolic PDEs,
 viscosity solutions.

\tableofcontents

\section{Introduction}
\setcounter{equation}{0}

We consider a reflected backward stochastic differential equation (RBSDE) with generator $f$, terminal condition $\xi$ and obstacle $L$
     \bea  \label{RBSDE}
  L_t \le Y_t= \xi  + \int_t^T f(s, Y_s, Z_s) \, ds +   K_T - K_t - \int_t^T  Z_s dB_s\,, \qq    t \in [0, T],
    \eea
    where the solution $(Y,Z,K)$     satisfies the so-called {\it flat-off}\, condition:
         \bea  \label{flat-off}
       \int_0^T   (  Y_t   -  L_t )     d K_t =0 ,
    \eea
    and $K$ is an increasing process.
    We will consider the case when $f$ is allowed to have quadratic growth in the $z$-variable. Moreover, we will allow $L$ and $\xi$ to be unbounded.

\ss
The theory of RBSDEs is closely related to the theory of optimal stopping in that the snell-envelope can be represented as a solution of an RBSDE. These equations were first introduced by \cite{EKPPQ-1997}.
The authors  provided the existence and uniqueness of  an adapted solution for a real-valued RBSDE with square-integrable terminal condition under the Lipschitz hypothesis on the generator. There has been a few developments after this seminal result. Some generalizations were obtained for backward stochastic differential equations (BSDEs) without an obstacle and later they were generalized to RBSDEs:

  \ss \no {\bf 1)}  \cite{Lep_San_97} showed  the existence of a  maximal and a minimal solution for  real-valued BSDEs, with square-integrable terminal condition when the generator $f$ is only continuous and has linear growth in variables $y$ and $z$. Then \cite{Matoussi_97} adapted this result to the case of RBSDEs.

    \ss \no {\bf 2)} \cite{Ko_2000} established the existence, comparison, and stability results for real-valued
quadratic BSDEs (when $f$ is allowed to have quadratic growth in the $z$-variable) with bounded terminal condition. In the spirit of \cite{PP-92}, the author gave a link between the solutions of BSDEs based on a diffusion
   and viscosity solutions of the corresponding semi-linear parabolic PDEs.   \cite{Lep_San_98} extended the existence result of quadratic BSDEs with bounded terminal condition to the case that the generator $f$ can have a superlinear growth in the $y$-variable.
          \cite{KLQT_RBSDE} made a counterpart study for  RBSDEs with bounded terminal condition and bounded obstacle when the generator $f$ has
        superlinear growth in $y$ and quadratic growth in $z$.

    \ss \no {\bf 3)} With   help of a localization procedure and a priori bounds, \cite{BH-06} showed that  the boundedness assumption on the terminal condition is not necessary for the existence of an adapted solution to a real-valued  quadratic BSDE:  One only needs to require the terminal condition has exponential moment of certain order. Correspondingly,  \cite{Lep_Xu_2007} derived the existence result for quadratic RBSDEs with such an unbounded terminal condition, but still with a bounded obstacle.

  Recently, \cite{BH-07}, under the assumption that the generator $f$ is additionally convex or concave in the $z$-variable, 
  used a so-called ``$\th$-difference" method to obtain comparison (thus uniqueness) and stability results for  quadratic BSDEs with solutions having  every exponential moment.
  Morever, \cite{DHR_2010} proved that uniqueness holds among solutions   having a given exponential moment by using a verification theorem that relies on     the Fenchel-Legendre dual of the generator.
   With these results they also showed that the solutions of BSDEs are viscosity solutions of PDEs which are quadratic in the gradient.
   On the other hand, \cite{Da_Lio_Ley_2010} showed that these PDEs have unique solutions.

  In the current paper, we extend the results of \cite{BH-07}, \cite{DHR_2010}, and \cite{Da_Lio_Ley_2010} to RBSDEs. Alternatively, our results can be seen as an extension of \cite{KLQT_RBSDE} and  \cite{Lep_Xu_2007} to the unbounded obstacles.
 We start by establishing   two a priori estimates which will serve as our basic tools; see Section~\ref{sec:a-priori-estimates}. The first one shows that any bounded $Y$ has an upper bound in term of  the terminal condition $\xi$ and the obstacle $L$. The second estimate is on the $\hL^p$ norms of $Z$ and $K$. With the help of these two estimates, we can establish a monotone stability result (see Theorem \ref{thm:monotone-stable}) in the spirit of  \cite{Ko_2000}.
Then the existence follows as a direct consequence; see Theorem \ref{thm:existence}.

\ss Next,  we apply the aforementioned $\th$-difference method to derive a comparison theorem (see Theorems \ref{thm:comparison} and \ref{thm:comparison2}) for quadratic RBSDEs with unbounded terminal conditions and unbound obstacles when the generator $f$ is additionally  convex or concave in the $z$-variable.
  Instead of estimating  the difference of two solutions $Y$ and $\wh{Y}$, we estimate  $Y-\th \wh{Y}$ for each
$\th \in (0,1)$, which allows us to utilize the convexity or the concavity of the generator $f$.
 In the concave-generator case,
    we prove a uniqueness result for RBSDEs using an argument that involves  the Fenchel-Legendre dual of the generator, see Theorem~\ref{thm:uniqueness}.  As opposed to \cite{DHR_2010} (or \cite{OS_CRM}), we are not relying on a verification argument but directly compare two solutions.
 Since it only requires a given exponential moment on solutions, this uniqueness result is   more general   than the one that would be implied by the above comparison theorem.
 We develop an alternative representation of the unique solution in Section \ref{sec:g-evaluation}, where we improve the results of Theorem 5.3 of \cite{OSNE2} on optimal stopping for quadratic $g$-evaluations. Moreover, the convexity/concavity assumption on generator $f$ in the $z$-variable as well as the $\th$-difference method are also used in deducing the stability result (see Theorem \ref{thm:stable}), which is crucial for the continuity property of the solutions of forward backward stochastic differential equations with respect to their initial conditions; see Proposition~\ref{prop:u_continuous}. This result together with the stability result gives a new proof of the flow property; see Proposition \ref{prop:u_markov}.      A Picard-iteration procedure was   introduced   to show this property for   BSDEs with Lipschitz generators,
       see e.g. Theorem 4.1 of \cite{EPQ-97}. However, it is not appropriate to apply such a Picard-iteration procedure to derive the flow property
       for quadratic RBSDEs.

  \ss Thanks to the flow property, the solution of the RBSDE is a viscosity solution of an associated obstacle problem for a semi-linear parabolic PDE, in which the non-linearity appears as the square of the gradient; see Theorem~\ref{thm:viscos_exist}.  It is worth pointing out that \cite{DHR_2010} shows  the existence of a viscosity solution to a similar PDE (with a quadratic gradient term) without obstacle by approximating the generator $f$  from below by a sequence of Lipschitz generators under a strong assumption that $f^-$ has a linear growth in variables $y$ and $z$.   However, such a strong assumption is not necessary if we directly use the flow property  to prove Theorem ~\ref{thm:viscos_exist}.
  Finally, we prove that in fact this obstacle problem has a unique solution, which is a direct consequence of Theorem~\ref{thm:viscos_comparison}, a comparison principle between  a viscosity subsolution and a viscosity supersolution. Although inspired by   Theorem 3.1 of \cite{Da_Lio_Ley_2010}, we
  prove Theorem~\ref{thm:viscos_comparison} in a quite different way   because  there are two gaps in the proof of Theorem 3.1 of \cite{Da_Lio_Ley_2010}, see Remark \ref{rem_gap}.


\subsection{Notation and Preliminaries} \label{sec:notation}

 \ms  Throughout this paper  we let $B$ be a $d$-dimensional standard Brownian
 Motion defined on a complete probability space $(\O,\cF, P)$,
 and consider the augmented filtration generated by it, i.e.,
 \beas
 \bF= \left\{\cF_t \dfnn \si \Big( \si\big(B_s; s\in [0,t]\big) \cup \cN \Big) \right\}_{t  \in  [0, \infty) }\,,
 \eeas
  where $\cN$ is the collection of all $P$-null sets in $\cF$.  We fix a finite time horizon $T>0$.
       Let $\cS_{0,T}$ be the collection of all $\bF$-stopping times $\nu$ such that $0 \le \nu \le T$,
   \pas~ For any $ \nu \in \cS_{0,T}$, we define $\cS_{\nu, T} \dfnn \{ \t \in \cS_{0,T}\, |\; \nu \le \t  \le T, ~\pas \}$.
  Moreover, we will use the convention $ \inf \{\es\} \dfnn \infty$.

 \ms  The following spaces of functions will be used in the sequel:

 \ss \no 1)\,   Let $ \hC [0,T]$  denote the set of all real-valued continuous functions on $[0,T]$,
  and let $ \hK [0,T]$ be the subset of $ \hC [0,T]$ that consists of all real-valued increasing and continuous  functions on $[0,T]$.
 For any $\{\ell_t\}_{t \in [0,T]} \in \hC [0,T]$,
  we define $\ell^\pm_* \dfnn \underset{t \in [0,T]}{\sup}  (\ell_t)^\pm$. Then
   \bea        \label{eq:a005}
  \ell_*  \dfnn    \underset{t \in [0,T]}{\sup} | \ell_t |  =  \underset{t \in [0,T]}{\sup} \big(   (\ell_t)^-  \vee (\ell_t)^+  \big)
    = \underset{t \in [0,T]}{\sup}   (\ell_t)^-    \vee \underset{t \in [0,T]}{\sup}   (\ell_t)^+   =\ell^-_* \vee \ell^+_*  .
   \eea

  \ss \no 2)\,   For any sub-$\si$-field  $\,\cG\,$  of $\,\cF\,$,  let

\ss   $\bullet$      $\mathbb{L}^0(\cG )$  be  the space of all real-valued,
$\,\cG$-measurable random variables;

\ss  $\bullet$  $\hL^p(\cG ) \dfnn \Big\{\xi \in \mathbb{L}^0(\cG ):\,\|\xi\|_{\hL^p(\cG )} \dfnn
 \Big\{ E \Big[   \,    |\xi |^p      \Big] \Big\}^{\frac{1}{p}}<\infty\Big\}$  for all $p \in [1, \infty)$;

\ss  $\bullet$  $\hL^\infty(\cG ) \dfnn \Big\{\xi \in \mathbb{L}^0(\cG ):\,\|\xi\|_{\hL^\infty(\cG )} \dfnn
 \esup{\o \in \O} \,|\xi(\o)| <\infty\Big\}$;

\ss  $\bullet$  $\hL^e(\cG ) \dfnn \Big\{\xi \in \mathbb{L}^0(\cG ):\, E \big[  e^{p    |\xi |  }   \big]  < \infty, ~\fa p \in  (1, \infty) \Big\}$.

\ss \no 3)\,   Let   $\hB$ be a generic  Banach  space with norm $|\cd|_{\hB}$. For any $p,q \in [1, \infty)$, we define three Banach spaces:

    \ss   $\bullet$   $ \hL^{p,q}_\bF([0,T]; \hB) $  denotes the space of all $\hB$-valued,
 $\bF$-adapted processes $X$ with
   \beas
   \|X\|_{\hL^{p,q}_\bF([0,T]; \hB)} \dfnn \left\{  E \left[ \left( \hb{$\int_0^T \neg |X_t|^p_\hB \,  dt$} \right)^{\frac{q}{p}} \right] \right\}^{\frac{1}{q}}<\infty ;
   \eeas

  \ss    $\bullet$    $ \hH^{p,q}_\bF([0,T]; \hB) $ \big(resp. $ \wh{\hH}^{p,q}_\bF([0,T]; \hB) $\big)
 $\dfnn \big\{X \in \hL^{p,q}_\bF([0,T]; \hB) :  X$ is $\bF$-predictable (resp. $\bF$-progressively measurable)\big\}.

 \ss \no When $p=q$, we simply write $\hL^p_\bF$, $\hH^p_\bF$ and $ \wh{\hH}^p_\bF$ for
  $\hL^{p,p}_\bF$, $\hH^{p,p}_\bF$ and $ \wh{\hH}^{p,p}_\bF$  respectively.
 Moreover we let

 \ss  $\bullet$ $\hH^{p,loc}_\bF([0,T];\hB)$ \big(resp.\;$\wh{\hH}^{p,loc}_\bF([0,T];\hB)$\big) denote the space of all $ \hB$-valued, $\bF$-predictable (resp. $\bF$-progressively measurable) processes $X$ with $\int_0^T |X_t|^p_{\hB} dt < \infty$, \pas ~ for any $p \in [1,\infty)$.

  \ms \no 4)\,    Let  $\hC^0_\bF[0,T]$ be the space of all real-valued, $\bF$-adapted continuous processes.  We need the following subspaces of $\hC^0_\bF[0,T]$.

 \ms   $\bullet$ $\hC^\infty_\bF[0,T]  \dfnn  \Big\{  X \in \hC^0_\bF[0,T] :\,  \|X\|_{\hC^\infty_\bF[0,T]} \dfnn \esup{\o \in \O} \Big( \underset{ t
\in [0,T] }{\sup}  \big| X_t(\o) \big| \Big) <\infty \Big\}$;

   $\bullet$ $\hC^p_\bF[0,T] \dfnn    \left\{  X \in \hC^0_\bF[0,T] :\,   \|X\|_{\hC^p_\bF[0,T]}
  \dfnn \bigg\{ E\Big[\,\underset{t \in [0,T]}{\sup}|X_t|^p\Big]  \bigg\}^{\frac{1}{p}}<\infty \right\}$
  for all $p \in [1, \infty)$;

  \ms  $\bullet$ $\hV_\bF[0,T] \dfnn \big\{  X \in \hC^0_\bF[0,T] :\,  X$ has finite variation\big\};

  \ss   $\bullet$  $  \hK_\bF[0,T]  \dfnn  \big\{  X \in \hC^0_\bF[0,T] :\,  X$ is an increasing process with $X_0=0 \big\}  \subset \hV_\bF[0,T]$;

   \ms   $\bullet$  $  \hK^p_\bF[0,T] \dfnn \big\{  X \in \hK_\bF[0,T] :\, X_T  \in \hL^p(\cF_T) 
   $\big\} for all $p \in [1, \infty)$;

\ms   $\bullet$  $\hE^{\l,\l'}_\bF[0,T] \dfnn \left\{ X \in \hC^0_\bF[0,T] :\,    E \left[  e^{\l     X^-_*  } + e^{\l'      X^+_*  } \right]  < \infty \right\}
   \subset \underset{p \in [1,  \infty)}{\cap}  \hC^p_\bF[0,T] $  for all $\l,\l' \in (0, \infty)$.  

 \ss \no For any $\l \in (0, \infty)$,  we set $\hE^\l_\bF[0,T] \dfnn \hE^{\l,\l}_\bF[0,T]$. For any $X \in \hC^0_\bF[0,T]$, one can deduce from
     \eqref{eq:a005} that
  \bea   \label{eq:a009}
  E  \neg \left[  e^{ \l    X_* }      \right]    =    E  \neg \left[  e^{\l  (X^-_* \vee X^+_* )  }   \right]
                   =  E  \neg \left[  e^{ \l  X^-_*    }  \vee   e^{ \l   X^+_*    }  \right] \le E  \neg \left[  e^{ \l  X^-_*    }  +   e^{ \l   X^+_*    }  \right]
                   \le   2 E  \neg \left[  e^{ \l    X_* }      \right]  ,
  \eea
  which implies that
    $\hE^\l_\bF[0,T]     =\left\{ X \in \hC^0_\bF[0,T] :\,    E \left[  e^{\l    X_*  }   \right]  < \infty \right\} $.
    Moreover,  for any $p \in [1, \infty)$, we set $ \hS^p_\bF[0,T] \dfnn \hE^p_\bF[0,T] \times  \hH^{2, 2p}_\bF([0,T];\hR^d)    \times \hK^p_\bF[0,T] $.

 \subsection{Reflected BSDEs}

 \ms Let $\sP$ denote the $\bF$-progressively measurable $\si$-field on $ [0,T] \times \O$.
 A {\it parameter set} $(\xi, f, L)$ consists of  a random variable $\xi \in \hL^0(\cF_T) $,
 a  function $f:  [0,T] \times \O \times \hR \times \hR^d \rightarrow  \hR $ and a process   $ L   \in   \hC^0_\bF[0,T]$ such  that
$f$ is $\sP \times \sB(\hR) \times   \sB(\hR^d)/\sB(\hR)$-measurable and that $L_T \le \xi$, \pas

\begin{deff}  \label{defn_the RBSDE}
   Given a parameter set $(\xi, f, L)$,
  a triplet $(Y, Z,K) \in \hC^0_\bF[0,T] \times \wh{\hH}^{2,loc}_\bF([0,T];\hR^d) \times \hK_\bF[0,T]$
  is called a solution of the {\rm reflected backward stochastic differential equation} with terminal condition $\xi$, generator $f$, and  obstacle $L$
 \neg  \big(RBSDE $(\xi, f, L)$ for short\big), if   \eqref{RBSDE} and \eqref{flat-off} hold \pas
     \end{deff}

A function   $f:  [0,T] \times \O \times \hR \times \hR^d \rightarrow  \hR $ is said to be {\it Lipschitz}  in $(y,z)$ if for some $\l>0$,
 it holds \dtp~ that
 \beas
  \;\;    \big| f  (t,\o, y_1, \th z_1   ) - f  (t,\o, y_2,   z_2 ) \big| \le \l \big( |y_1-y_2|+|z_1-z_2|\big),  \q
 \fa  y_1, y_2    \in    \hR, ~ \fa z_1, z_2 \in \hR^d.
 \eeas

\ss The theory of RBSDEs with Lipschitz generators was well developed in the seminal paper \cite{EKPPQ-1997}.  In this paper, we are interested in {\it quadratic} RBSDEs, i.e., the RBSDEs whose generators  have  quadratic growth in $z$ in the following sense:

 \ms \no {\bf (H1)}   For three constants $\a, \beta \ge 0$ and $\g>0$, it holds \dtp~ that
\beas
|f(t , \o, y,z)| \le \a  +\beta |y|+ \frac{\g}{2} |z|^2, \q \;\;  \fa  (y,z) \in \hR \times \hR^d .
\eeas
   In what follows,  for any  $\l \ge 0$ we let $c_\l$ denote a generic constant
  depending on $\l, \a , \beta,  \g  $ and $T$ (in particular, $c_0$ stands for a generic constant depending on $ \a , \beta,  \g  $ and $T$),
  whose form may vary from line to line.


 \section{Two A Priori Estimates}\label{sec:a-priori-estimates}

  \ms We first present an a priori estimate, which is an extension of Lemma 3.1 of \cite{Lep_Xu_2007}.

\begin{prop} \label{prop_estimate}
 Let $(\xi, f, L)$ be a parameter set such that $f$  satisfies  (H1).
 If $(Y,Z,K) $  is a solution of the quadratic RBSDE$(\xi,f,L)$ such that  $Y^+ \in \hC^\infty_\bF[0,T]$,
 then it holds \pas ~ that
     \bea        \label{a_priori}
          Y_t   \le     c_0  +  \frac{1}{\g} \ln  E\left[  e^{      \g  e^{\beta T  }   (\xi^+ \vee    L^+_*)  }\big|  \cF_t  \right]  , \q       t \in [0,T] .
    \eea
\end{prop}

\ss \no {\bf Proof:}   In light of It\^o's formula,   $(Y,Z,K) \in \hC^0_\bF[0,T] \times \wh{\hH}^{2,loc}_\bF([0,T];\hR^d) \times \hK_\bF[0,T] $ with $Y^+ \in \hC^\infty_\bF[0,T]$  is a solution of the  RBSDE$(\xi,f,L)$ if and only if
 \beas 
    \big(\wt{Y}, \wt{Z}, \wt{K} \big)\dfnn  \left( e^{\g Y},    \g e^{\g Y} Z ,
 \hb{$\g \int_0^\cd  e^{\g Y_s}  d K_s$}\right)  \in  \hC^\infty_\bF[0,T] \times \wh{\hH}^{2,loc}_\bF([0,T];\hR^d) \times \hK_\bF[0,T]
\eeas
 is a solution of  the RBSDE$(e^{\g \xi}, \wt{f}, e^{\g L}) $ with
 \beas
     \wt{f} (t, \o, y, z) \dfnn \b1_{\{y > 0\}}
   \left\{ \g y f\left(t, \o, \frac{\ln y}{\g}, \frac{z}{\g y}\right) -\frac12 \frac{|z|^2}{y}
\right\}, \q \fa (t, \o, y, z) \in [0,T] \times \O \times \hR  \times \hR^d .
 \eeas
  Let $\mu \dfnn  \a \g \vee \beta \vee 1 $.  One can deduce from (H1) that  \dtp
 \bea  \label{eq:c207}
  \wt{f} (t, \o, y, z) \le H(y) \dfnn y \big(\mu+\beta \ln y \big) \b1_{\{y \ge 1\}}  + \mu  \b1_{\{y  < 1\}} , \q  \fa (y, z) \in   \hR \times \hR^d .
 \eea
  Clearly,  $H(\cd)$ is a strictly positive,  increasing, continuous and convex function with $\int_0^\infty \frac{1}{H(y)} dy = \infty$.

 \ss For any  $x \in \hR$ and $\wt{T} \in [0,T]$,  the  ordinary differential equation (ODE)
   \beas   
   \f(t) = e^{\g x} +\int_t^{\wt{T}} H\big(\f(s)\big) ds, \q t \in \big[0,\wt{T}\big]
   \eeas
   can be solved as follows (cf. \cite{BH-06}):

   \ss \no (\,i) When $x  \ge 0$: $ \dis \f^{\wt{T}}_t\neg (x) = \exp\left\{ \mu \vf\big(\wt{T}-t\big) +\g x e^{\beta (\wt{T}-t)}\right\}$, where  $\vf(s) \dfnn \frac{e^{\beta s}-1}{\beta} \b1_{\{\beta > 0\}}+s \b1_{\{\beta = 0\}}$, $\fa s \in  [0,T ]$;

\ss \no (ii) When $x  < 0$:   $   \f^{\wt{T}}_t\neg (x)  =\left\{ \ba{ll}
  e^{ \g x   }+\mu (\wt{T}-t)< 1+ \mu (\wt{T}-t) \le e^{ \mu (\wt{T}-t)} \le e^{ \mu \vf (\wt{T}-t)}  , \q  & \hb{if } e^{\g x}+\mu (\wt{T}-t) < 1 ,  \\
  \exp\left\{ \mu \vf \left(  \wt{T}-t+\frac{e^{\g x} -1}{\mu}  \right)   \right\} \le
  e^{ \mu \vf  (  \wt{T}-t   )  }  & \hb{if } e^{\g x}+\mu (\wt{T}-t) \ge 1.   \ea \right.  $

\ss \no One can    check that

\ss \no ($\f$1) For any $x \in \hR$ and $\wt{T} \in [0,T]$,  $t \to \f^{\wt{T}}_t\neg (x) $ is a  decreasing and continuous function on $\big[0,\wt{T}\big]$;

\ss \no ($\f$2) For any $x \in \hR$ and $t \in [0,T]$,    $\wt{T} \to \f^{\wt{T}}_t\neg (x) $ is an  increasing and continuous function on $\big[t,T]$;

\ss \no ($\f$3) For any $0 \le t \le \wt{T}  \le T $,  $x \to \f^{\wt{T}}_t\neg (x) $ is an  increasing and continuous function on $\hR$;

\ss \no ($\f$4) For any $x \in \hR$ and $0 \le t \le \wt{T}  \le T $,  $\f^{\wt{T}}_t\neg (x) \le  \exp\left\{ \mu \vf(T)  +\g x^+ e^{\beta T}\right\}$.

   \ms  Let $ \wt{\O} \dfnn \{\o \in \O:   ~ L_T(\o) \le \xi(\o) \hb{ and the path }t \to L_t(\o)\hb{ is continuous}   \} \in \cF$, which
  defines a measurable set with probability $1$. Fix  $\o \in \wt{\O}$.
  Theorem 6.2 of \cite{Lep_Xu_2007} shows that the following reflected backward ordinary differential equation
  \beas
     \left\{
       \begin{array}{l}  \dis  e^{\g L_t(\o)} \le \L_t(\o)=  e^{\g \xi(\o)}  + \int_t^T H\big(\L_s(\o)\big) \, ds +   k_T(\o) - k_t(\o) \,, \qq  t \in [0, T] \, ,\\
      \dis \int_0^T   \left(  \L_s(\o)   -  e^{\g L_s(\o)} \right)     d k_s (\o) =0
     \end{array}
     \right.
    \eeas
    admits a unique solution $ \big( \L_\cd(\o), k_\cd(\o) \big)   \in  \hC [0,T] \times \hK [0,T] $, which satisfies
     \bea  \label{eq:c209}
     \L_t(\o)  =   \underset{s \in [t,T]}{\sup}
     \left( \int_t^s H\big(\L_r(\o)\big)dr + e^{\g \xi(\o)}\, \b1_{\{s=T\}}+e^{\g L_s(\o)} \b1_{\{s<T\}}  \right)
     = \underset{s \in [t,T]}{\sup} u^s_t(\o) , \q t \in [0,T] ,
     \eea
     where $\{u^s_r(\o)\}_{r \in [0, s] }$ is the unique solution of  the following ODE
      \beas
       u^s_r  (\o) = e^{\g \xi(\o)}\, \b1_{\{s=T\}}+e^{\g L_s(\o)} \, \b1_{\{s<T\}}+\int_r^s H \big( u^s_a  (\o) \big) d a  , \q r \in [0,s].
      \eeas
      To wit, $ u^s_r  (\o)= \f^s_r  \big(\xi(\o)\, \b1_{\{s=T\}} \neg + \neg  L_s(\o)\, \b1_{\{s<T\}}\big)$.
 Then it follows from \eqref{eq:c209}   and ($\f$4) that
  \bea \label{eq:c217}
   0<  e^{\g L_t(\o)} \le   \L_t(\o) = \underset{s \in [t,T]}{\sup} u^s_t(\o)     \le \exp\Big\{ \mu \vf(T)  +  \g e^{\beta T} \big( \xi^+(\o)  \vee   L^+_*(\o) \big)\Big\}  , \q t \in [0,T]    .
  \eea
For any $0 \le t_1 < t_2 \le T$,   one can deduce from   \eqref{eq:c209}  and ($\f$1)  that
  \bea  \label{eq:c221}
        \L_{t_1}(\o)  = \underset{s \in [t_1,T]}{\sup} u^s_{t_1}(\o) \ge  \underset{s \in [t_2,T]}{\sup} u^s_{t_1}(\o)
   \ge    \underset{s \in [t_2,T]}{\sup} u^s_{t_2}(\o) = \L_{t_2}(\o)  ,
  \eea
Thus   $t \to \L_t(\o) $ is  a decreasing and continuous path.
 Moreover,  for any $t \in [0,T]$ \eqref{eq:c209}  and ($\f$2) imply  that
  \bea  \label{eq:c215}
  \L_t(\o) = \underset{s \in [t,T]}{\sup} u^s_t(\o)= \sup\Big\{  u^s_t(\o): s \in \big( [t,T) \cap \hQ\big) \cup\{T\} \Big\}      .
 \eea

    For any $s \in [0,T]$, since $  \xi   \b1_{\{s=T\}} +  L_s   \b1_{\{s<T\}}$ is an   $\cF_s$-measurable random variable,  the continuity of function
  $\f^s_t (\cd)$ by ($\f$3) implies that    $ \big\{ u^s_t  (\o) \big\}_{\o \in \O}= \f^s_t  \big(\xi   \b1_{\{s=T\}} +  L_s   \b1_{\{s<T\}}\big)$ is also
  an   $\cF_s$-measurable random variable.        Thus   we can deduce  from     \eqref{eq:c215}   that for any  $t \in [0,T]$,
  the random variable $   \L_t $  is $\cF_T$-measurable (however, not necessarily $\cF_t$-measurable).

   \ss Now, let us introduce an $\bF$-adapted process $ \mathfrak{f}_t \dfnn E[  H(\L_t) |\cF_t  ] $, $t \in [0,T]$.
  Since $\L$ is a decreasing  process by \eqref{eq:c221}, and since
  $H(\cd)$ is an increasing  function, it holds for any $0 \le t < s \le T$ that
 \beas
    E[\mathfrak{f}_s  |\cF_t ] =E[  H(\L_s) |\cF_t  ] \le  E[  H(\L_t) |\cF_t  ] = \mathfrak{f}_t, \q \pas
 \eeas
 which implies that $\mathfrak{f}$ is a supermartingale.  As $Y^+ \in \hC^\infty_\bF[0,T]$, it  follows that $(\xi^+, L^+) \in \hL^\infty(\cF_T) \times  \hC^\infty_\bF[0,T]$. Then the continuity of process $H(\L_\cd)$,  \eqref{eq:c217} and the Bounded Convergence Theorem
   imply that
  \beas
   E[ \mathfrak{f}_t ] = E[  H(\L_t)] = \lmt{s \da t}  E[  H(\L_s)] = \lmt{s \da t}   E[ \mathfrak{f}_s ] , \q  t \in [0,T].
  \eeas
 Thanks to  Theorem 1.3.13 of \cite{Kara_Shr_BMSC},   $\mathfrak{f}$ has a   right-continuous modification $\wt{\mathfrak{f}}$.
  Hence, we can regard $\wt{\mathfrak{f}}$
   as a generator that is independent of $( y, z)$.  It follows from Fubini's Theorem, Jensen's inequality  as well as \eqref{eq:c217}  that
    \beas
    \q   E   \int_0^T \neg  \big|\,\wt{\mathfrak{f}}_s \big|^2 ds   =       \int_0^T \neg E \left[ \big|\,\wt{\mathfrak{f}}_s \big|^2 \right] ds
    = \int_0^T \neg E \left[ | \mathfrak{f}_s|^2 \right] ds  \le       \int_0^T  \neg E\Big[E\big[  | H(\L_s)|^2 \big|\cF_s  \big]    \Big] ds
         =       \int_0^T \neg  E\big[    | H(\L_s)|^2    \big] ds < \infty  .
   \eeas
   Since $e^{\g \xi}    \in \hL^\infty(\cF_T)$ and $e^{\g L} \in \hC^\infty_\bF[0,T] $,
      Theorem 5.2 and Proposition 2.3 of \cite{EKPPQ-1997} show that
 the RBSDE$(e^{\g \xi},  \wt{\mathfrak{f}}, e^{\g L} ) $ admits a unique solution $(\cY, \cZ, \cK ) \in \hC^2_\bF[0,T] \times
  \hH^2_\bF\big([0,T]; \hR^d \big) \times \hK^2_\bF[0,T]  $ and that for any $t \in [0,T]$
  \bea  \label{eq:c232}
\cY_t &=& \esup{\t \in \cS_{t,T}}E\left[  \int_t^\t  \wt{\mathfrak{f}}_s \, ds
   + e^{\g \xi }\, \b1_{\{\t=T\}}+e^{\g L_\t }\b1_{\{\t <T\}} \Big|\cF_t\right] , \q \pas
 \eea

   For any $t \in [0,T]$ and $\t \in \cS_{t,T} $,   Fubini's Theorem implies that for any $A \in \cF_t$
\beas
   &&   E\left[ \b1_A \int_t^\t  \wt{\mathfrak{f}}_s \, ds  \right]  = E \int_t^T    \b1_A \b1_{\{s \le \t\}}  \wt{\mathfrak{f}}_s \, ds
  =\int_t^T   E\left[  \b1_A \b1_{\{s \le \t\}}  \wt{\mathfrak{f}}_s  \right]  ds
  =\int_t^T   E\left[  \b1_A \b1_{\{s \le \t\}}   \mathfrak{f}_s  \right]  ds  \\
& &\qq \q     =   \int_t^T   E\left[   \b1_A \b1_{\{s \le \t\}} E[     H(\L_s) |\cF_s  ]   \right]  ds
   =\int_t^T   E\left[   E\big[  \b1_A \b1_{\{s \le \t\}}   H(\L_s) \big|\cF_s  \big]   \right]  ds \\
 & &\qq \q     =  \int_t^T   E\left[     \b1_A \b1_{\{s \le \t\}}   H(\L_s)   \right]  ds
= E    \int_t^T      \b1_A \b1_{\{s \le \t\}}   H(\L_s)    ds      =   E\left[  \b1_A  \int_t^\t          H(\L_s)    ds    \right]    .
 \eeas
 Thus  $ E\left[ \int_t^\t  \wt{\mathfrak{f}}_s \, ds |\cF_t\right]
  =       E\left[   \int_t^\t    H(\L_s) ds  |\cF_t \right ]  $,   \pas ~
  Then   \eqref{eq:c232}, \eqref{eq:c209} and \eqref{eq:c217} imply that for any $t \in [0,T]$
  \bea       \label{eq:c231}
        \cY_t    & \dneg  =& \dneg    \esup{\t \in \cS_{t,T}}E\left[  \int_t^\t   H(\L_s) ds + e^{\g \xi }\, \b1_{\{\t=T\}}
            +e^{\g L_\t  }\b1_{\{\t <T\}} \Big|\cF_t\right]           \nonumber \\
            & \dneg  \le& \dneg   E[\L_t |\cF_t  ] \le   e^{ \mu \vf(T)  }
  E \left[e^{   \g  e^{\beta T }   ( \xi^+  \vee   L^+_* )     } \big|\cF_t \right] \le C_* ,\q \pas,
  \eea
 with $  C_* \dfnn \exp \left\{ \mu \vf(T)  +   \g  e^{\beta T }   \Big( \|\xi^+\|_{\hL^\infty(\cF_T)}  \vee  \| L^+\|_{\hC^\infty_\bF[0,T]} \Big)     \right\}$. By the continuity of process $\cY$, it holds \pas~ that
  \bea \label{eq:c213}
   0<     e^{\g L_t  }   \le \cY_t   \le  e^{ \mu \vf(T)  }
  E \left[e^{   \g  e^{\beta T }   ( \xi^+  \vee   L^+_* )     } \big|\cF_t \right] \le C_*, \q t \in [0,T],
  \eea
   which shows that   $\cY \in \hC^\infty_\bF[0,T]$ with $\|\cY\|_{\hC^\infty_\bF[0,T]} \le C_* $.

\ms To finalize the proof, it suffices to show that $ P \left( \wt{Y}_t \le \cY_t , ~  \fa  t \in [0,T] \right) =1 $.
  To see this, we fix $n \in \hN$ and define the $\bF$-stopping time
  \beas
  \t_n \dfnn \inf\left\{ t \in [0,T]: \int_0^t     |Z_s|^2 ds >n  \right\} \land T .
  \eeas
 Clearly, $\lmtu{n \to \infty} \t_n =T$, \pas ~ Applying  Tanaka's formula to the process $( \wt{Y} - \cY )^+$ yields that
 \bea
         ( \wt{Y}_{\t_n \land t} - \cY_{\t_n \land t} )^+   &\dneg =&\dneg  ( \wt{Y}_{\t_n} - \cY_{\t_n} )^+ + \int_{\t_n \land t}^{\t_n} \neg \b1_{\{\wt{Y}_s > \cY_s\}} \big( \wt{f}(s,\wt{Y}_s,  \wt{Z}_s) - \wt{\mathfrak{f}}_s  \big) ds
    + \int_{\t_n \land t}^{\t_n} \neg \b1_{\{\wt{Y}_s > \cY_s\}} \big( d \wt{K}_s - d \cK_s  \big)    \nonumber  \\
  &\dneg &\dneg  - \int_{\t_n \land t}^{\t_n} \neg \b1_{\{\wt{Y}_s > \cY_s\}} (\wt{Z}_s - \cZ_s  ) \,  dB_s- \frac 12 \int_{\t_n \land t}^{\t_n}  d  \fL_s   , \q     t \in [0,T],
       \label{eq:axa130}
\eea
 where $\fL$ is a real-valued, $\bF$-adapted, increasing and continuous  process known as ``\,local time".

  \ss  Since the function $H(\cd)$ is increasing, continuous and convex,   Jensen's inequality and \eqref{eq:c231} show that
 \bea   \label{eq:axa135}
         H(\wt{Y}_s) - \mathfrak{f}_s      \le          H(\wt{Y}_s) -H \big(  E[\L_s |\cF_s  ]  \big)
     \le            H(\wt{Y}_s)- H ( \cY_s )       \le C_H    \big|\wt{Y}_s - \cY_s\big|           , \q    s \in [0,T] ,
\eea
 where $C_H$ is the  Lipschitz coefficient of function $H (\cd)$ over $\Big\{x \in \hR: |x| \le \|\wt{Y}\|_{\hC^\infty_\bF[0,T]} \vee  \|\cY\|_{\hC^\infty_\bF[0,T]}\Big\}$. Moreover,
 the flat-off condition of $\big(\wt{Y}, \wt{Z}, \wt{K} \big)$ implies that
 \bea    \label{eq:c323}
   \int_0^T  \b1_{\{\wt{Y}_s > \cY_s\}}    d \wt{K}_s = \int_0^T \b1_{\{ e^{\g \neg L_s} = \wt{Y}_s > \cY_s\}}       d \wt{K}_s = 0 , \q   \pas
 \eea

 Taking the expectation in \eqref{eq:axa130}, we can deduce from  \eqref{eq:c207}, Fubini's Theorem, \eqref{eq:axa135}  and \eqref{eq:c323}  that
    \beas
 && \hspace{-1.5cm}     E \left[   ( \wt{Y}_{\t_n \land t} - \cY_{\t_n \land t} )^+      \right] -E \left[     ( \wt{Y}_{\t_n} - \cY_{\t_n} )^+  \right]    \\
     &\le&       E \int_t^T    \b1_{\{s \le \t_n\}}\b1_{\{\wt{Y}_s > \cY_s\}}
     \big(     \wt{f}(s,\wt{Y}_s,  \wt{Z}_s) - \wt{\mathfrak{f}}_s    \big)  ds
        \le     \int_t^T  E \left[ \b1_{\{s \le \t_n\}}\b1_{\{\wt{Y}_s > \cY_s\}}
     \big(   H(\wt{Y}_s)  - \mathfrak{f}_s  \big) \right] ds \\
     &\le&       C_H \int_t^T E \left[ \b1_{\{s \le \t_n\}} \b1_{\{\wt{Y}_s > \cY_s\}}  ( \wt{Y}_s - \cY_s )^+    \right] ds
      \le       C_H \int_t^T E \left[   ( \wt{Y}_{\t_n \land s} - \cY_{\t_n \land s} )^+    \right] ds  , \q  t \in [0,T] .
  \eeas
  Then   Gronwall's inequality shows that  for any $t \in [0,T]$
   \beas
    E \left[   ( \wt{Y}_{\t_n \land t} - \cY_{\t_n \land t} )^+      \right]  \le e^{C_H T} E \left[     ( \wt{Y}_{\t_n} - \cY_{\t_n} )^+  \right]    .
   \eeas
   As $n \to \infty $, the continuity of processes $\wt{Y}$, $ \cY $  and the Bounded Convergence Theorem   imply that
  \beas
   E \left[   ( \wt{Y}_t - \cY_t )^+   \right]=0,  \q \hb{thus} \q    \wt{Y}_t  \le  \cY_t  , \q \pas
  \eeas
  Using the continuity of processes $\wt{Y}$ and $ \cY $ again, we obtain    $ P \left( \wt{Y}_t \le \cY_t , ~  \fa  t \in [0,T] \right) =1 $, which together with \eqref{eq:c213} leads to \eqref{a_priori}. \qed

\ss For a solution $(Y,Z,K) $ of a quadratic  RBSDE$(\xi,f,L)$ such that $L^-_* $ and $Y^+_*$ have exponential moments of certain orders, the next result estimates the norms of $(Z, K)$ in $\hH^{2, 2p}_\bF([0,T];\hR^d)   \times \hK^p_\bF[0,T]$  for some $p \in (1, \infty)$.

 \begin{prop} \label{prop:Z_K_bound}
   Let $(\xi, f, L)$ be a parameter set  such that $f$ satisfies  (H1). If  $(Y,Z,K) 
 $ is a solution of the quadratic  RBSDE$(\xi,   f,L)$ such that $Y \in \hE^{\l \g, \l' \g}_\bF[0,T]$
   for some $\l ,\l'> 1$ with $\frac{1}{\l} + \frac{1}{\l'} <1  $,
  then
    \beas
      \q    E \left[\left( \int_0^T     |Z_s|^2 ds \right)^p + K^p_T\right]
           \le           c_{\l,\l' \neg ,p} \,   E\left[ e^{\l  \g    Y^-_*}  +  e^{   \l'  \g     Y^+_*}    \right]   < \infty , \q  \fa  p \in \left(1, \frac{\l\l'}{\l+\l'}  \right).
       \eeas

\end{prop}

 \ss \no {\bf Proof:}    We set $p_o\dfnn \sqrt{\frac{\l\l'}{p\,(\l+\l')}} \land 2>1$ and define $\bF$-stopping times
     \beas
     \t_n \dfnn \inf\left\{ t \in [0,T]: \int_0^t e^{-p_o \g Y_s}   |Z_s|^2 ds >n  \right\} \land T, \q  \fa n \in \hN.
     \eeas
Since $E \big[ e^{\l \g Y^-_* } \big] < \infty$ and $Z \in \wh{\hH}^{2,loc}_\bF([0,T];\hR^d)$, it holds \pas ~ that
  $     Y^-_*  +  \int_0^T    |Z_s|^2 ds < \infty$.
     Then it follows that
       \beas
        \int_0^T e^{-p_o \g Y_s}   |Z_s|^2 ds \le e^{ p_o \g Y^-_*}   \int_0^T     |Z_s|^2 ds < \infty , \q \pas,
        \eeas
      which implies that   for \pas ~ $\o \in \O$, there exists an $n(\o) \in \hN$ such that $\t_{n(\o)}(\o)=T$.
        For any $n \in \hN$, applying It\^o's formula to the process $   e^{-p_o \g Y  }$ and using the fact that
     \beas
        \a   + \beta x   \le \left(\a \vee \frac{\beta}{(p^2_o-p_o) \g} \right) e^{(p^2_o-p_o) \g x } ,  \q \fa x \ge 0   ,
        \eeas
         we obtain that
    \bea
      && \hspace{-0.6 cm}  e^{-p_o\g  Y_0 }     + \frac12 p_o^2 \g^2  \neg   \int_0^{\t_n} e^{-p_o\g  Y_s}   |Z_s|^2 ds   \nonumber \\
      &&=  e^{-p_o\g  Y_{\t_n}} -p_o\g   \neg \int_0^{\t_n} e^{-p_o\g  Y_s }  f(s, Y_s, Z_s) ds
      - p_o\g   \neg \int_0^{\t_n} e^{-p_o\g  Y_s }  d K_s + p_o\g  \int_0^{\t_n} e^{-p_o\g  Y_s } Z_s d B_s \nonumber   \\
       && \le    e^{p_o \g Y^-_*}  \dneg  +  \neg   p_o\g  \hb{$ \left(\a \neg \vee \neg \frac{\beta}{(p^2_o - p_o) \g}\right)$}
        \dneg \int_0^{\t_n}  \dneg  e^{-p_o\g  Y_s+(p^2_o-p_o) \g|Y_s| }   ds        + \neg   \frac12 p_o \g^2  \dneg \int_0^{\t_n} \neg e^{-p_o\g  Y_s }   |Z_s|^2  ds  \nonumber \\
       && \q + p_o\g  \neg \left|  \int_0^{\t_n} \neg   e^{-p_o\g  Y_s } Z_s d B_s \right|   , \q \pas   \qq \label{eq:c101}
    \eea
     Observe that
    \beas 
     \q  \int_0^{\t_n} e^{-p_o\g  Y_s+(p^2_o-p_o) \g|Y_s| }   ds \le  \int_0^{\t_n} e^{- p^2_o  \g \b1_{\{Y_s < 0 \}} Y_s }   ds
     \le  \int_0^{\t_n} e^{ p^2_o \g \b1_{\{Y_s < 0 \}} Y^-_s }   ds     \le   T e^{ p^2_o  \g  Y^-_*}   ,   ~\;    \pas ,
    \eeas
 which together with     the Burkholder-Davis-Gundy inequality and \eqref{eq:c101} implies that
         \beas
        E \left[\left( \int_0^{\t_n}  e^{-p_o\g  Y_s}   |Z_s|^2 ds \right)^{\l p^{-2}_o}\right]
       &\le &  c_{\l,\l' \neg ,p} \, E \left[     e^{\l  \g  Y^-_*}
       +    \left|  \int_0^{\t_n} e^{-p_o\g  Y_s } Z_s d B_s \right|^{\l p^{-2}_o} \right]     \\
       &\le &  c_{\l,\l' \neg ,p} \, E \left[     e^{\l   \g  Y^-_*}   +
      e^{ \frac{\l}{2 p_{\neg o}}   \g  Y^-_*}     \left( \int_0^{\t_n} e^{-p_o\g  Y_s } |Z_s|^2 d s \right)^{ \frac12 \l p^{-2}_o} \right]  \qq   \\
         &\le &  c_{\l,\l' \neg ,p} \, E\left[      e^{\l   \g  Y^-_*}     \right]
           + \frac12  E \left[ \left( \int_0^{\t_n}  e^{-p_o\g  Y_s}   |Z_s|^2 ds \right)^{\l p^{-2}_o} \right] .
    \eeas
 Since  $E \neg \left[ \left( \int_0^{\t_n}  \neg  e^{-p_o\g  Y_s}   |Z_s|^2 ds \right)^{\l p^{-2}_o} \right] < \infty$, it follows that
 $    E  \neg \left[\left( \int_0^{\t_n}   \neg  e^{-p_o\g  Y_s}   |Z_s|^2 ds \right)^{\l p^{-2}_o}\right]
  \neg \le  \neg   c_{\l,\l' \neg ,p} \, E  \neg \left[      e^{\l  \g  Y^-_*}     \right]$.
      As $n \to \infty$, the Monotone Convergence Theorem gives that
  \beas
     E \left[\left( \int_0^T  e^{-p_o\g  Y_s}   |Z_s|^2 ds \right)^{\l p^{-2}_o}\right]  \le c_{\l,\l' \neg ,p} \,  E\left[      e^{\l  \g  Y^-_*}     \right].
     \eeas
   Observe that $ \frac{\l p_o p}{\l-p^2_o p}  < \frac{\l p^2_o p}{\l-p^2_o p} \le \l'$. Thus,  applying   Young's inequality with $\wt{p}=\frac{\l  }{\l-p^2_o p} $ and $\wt{q}=\frac{\l  }{ p^2_o p} $ yields  that
     \bea  \label{eq:c111}
       &&  E \left[\left( \int_0^T     |Z_s|^2 ds \right)^p \,\right]
         \le     E \left[   e^{p_o p     \g  Y^+_*}  \left( \int_0^T e^{-p_o \g Y_s}   |Z_s|^2 ds \right)^p \,\right]  \nonumber  \\
        && \qq \qq \le    c_{\l,\l' \neg ,p} \,    E \left[   e^{ \frac{\l p_{\neg o} p}{\l-p^2_{\neg o} p}   \g  Y^+_*}
        + \left( \int_0^T e^{-p_o \g Y_s}   |Z_s|^2 ds \right)^{\l p^{-2}_o} \right]
           \le   c_{\l,\l' \neg ,p} \,  E \left[   e^{\l  \g  Y^-_*}  +   e^{  \l'   \g  Y^+_*}         \right]    < \infty  .  \qq \qq
       \eea

     On the other hand,  since   $  Y_*  \le   Y^-_* +   Y^+_*  $,  it holds \pas ~ that
       \beas
       K_T &=&   Y_0-\xi  -     \int_0^T  f(s, Y_s, Z_s)  ds   +   \int_0^T Z_s dB_s   \\
       &\le  &   \a T +(2+\beta T) (Y^-_* +   Y^+_*)+ \frac{\g}{2} \int_0^T     |Z_s|^2 ds + \left|  \int_0^T Z_s dB_s  \right| .
       \eeas
   Then Burkholder-Davis-Gundy inequality and \eqref{eq:c111}  imply that
        \beas
      \hspace{0.8cm}   E \big[ K^p_T \big]  &\le &c_p \, E\left[  1+ \big( Y^-_* \big)^p +   \big( Y^+_* \big)^p +   \left( \int_0^T     |Z_s|^2 ds \right)^p
            +  \left(\int_0^T |Z_s|^2 ds \right)^{\frac{p}{2}} \right] \\
                 &\le & c_{\l,\l' \neg ,p} \,   E\left[  e^{\l  \g    Y^-_*} +  e^{   \l'  \g     Y^+_*}
     +   \left( \int_0^T     |Z_s|^2 ds \right)^p   \,   \right]   \le c_{\l,\l' \neg ,p} \,   E\left[ e^{\l  \g    Y^-_*}  +  e^{   \l'  \g     Y^+_*}    \right]  < \infty .    \hspace{1.5cm} \hb{\qed}
        \eeas

  \section{A Monotone Stability Result}

\begin{thm}    \label{thm:monotone-stable}
 For any $n \in \hN$, let $\big\{(\xi_n, f_n,L^n)\big\}_{n \in \hN} $ be a  parameter set and  let $(Y^n, Z^n, K^n) \in \hC^0_\bF[0,T] \times  \hH^{2,loc}_\bF([0,T];\hR^d) \times \hK_\bF[0,T]$  be a solution of the RBSDE\,$(\xi_n, f_n,L^n) $  such  that

   \ss \no (M1)   All generators $f_n$, $n \in \hN$   satisfy (H1) with the same constants $\a, \beta \ge 0$ and $\g>0$;

   \ss \no (M2)  There exists a function $f: [0,T] \times \O \times \hR \times \hR^d \to \hR$ such that for \dtp ~$(t, \o) \in [0,T] \times \O$,
   the mapping $f(t, \o, \cd, \cd)$ is continuous and   $ f_n(t, \o, y, z)  $
 converges to  $f(t, \o, y, z) $ locally uniformly in $(y, z)$;

  \ss  \no \;\;\,   and  that  for some $L \in \hC^0_\bF[0,T] $ and some real-valued,  $\bF$-adapted process  $Y$,   either of the following two holds:

  \ss \no   (M3a) It holds \pas~ that  for any $ t \in [0,T]$,  $\{L^n_t\}_{n \in \hN} $ and  $\{Y^n_t\}_{n \in \hN} $  are both increasing  sequences in $n$
   with $\lmtu{n \to \infty} L^n_t = L_t$ and $\lmtu{n \to \infty} Y^n_t = Y_t$ respectively;

  \ss \no   (M3b) It holds \pas~ that  for any $ t \in [0,T]$,  $\{L^n_t\}_{n \in \hN} $ and  $\{Y^n_t\}_{n \in \hN} $  are both decreasing  sequences in $n$
   with $\lmtd{n \to \infty} L^n_t = L_t$ and $\lmtd{n \to \infty} Y^n_t = Y_t$ respectively.

   \ms  Denote $  \sL_t   \dfnn (L^1_t)^- \vee L^-_t    $
  and $\sY_t \dfnn  (Y^1_t)^+  \vee  \,  Y^+_t        $, $\fa   t \in [0,T]$.    If  $\, \Xi \dfnn E \left[  e^{ \l   \g    \sL_* }   +e^{ \l'  \g     \sY_*    }    \right]< \infty$
  for some  $\l, \l' >6$ with $\frac{1}{\l}+\frac{1}{\l'} <\frac{1}{6}$,  then $Y \in \hE^{\l \g,\l' \g}_\bF[0,T]$ and there exist
   $(Z,K) \in \underset{p \, \in \big(   1, \frac{\l\l'}{\l+\l'} \neg \big)}{\cap}  \hH^{2, 2 p}_\bF([0,T];\hR^d) \times \hK^{p}_\bF[0,T] $
    such that the triplet $(Y,Z,K)$ is a solution of the   RBSDE\,$(\xi, f, L)$ with $\xi \dfnn Y_T$.

  \end{thm}

\ss \no  {\bf Proof:}     Since it holds \pas  ~  that
    \bea   \label{eq:c151}
           -\sL_t \le L^1_t \land L_t \le L^n_t \le   Y^n_t \le  Y^1_t \vee Y_t \le  \sY_t  , \q t \in [0,T] , \q \fa n \in \hN,
     \eea

  \ms The rest of the proof is divided into several steps.

 \ss \no {\bf 1)}  Let $        \l_o \dfnn 5+\frac12\left(\frac{\l\l'}{\l+\l'}-6\right)<  \frac{\l\l'}{\l+\l'} -1$.
            It follows that $p_o \dfnn  \frac{\l\l'}{\l\l'- \l_o(\l+\l')} \in \Big(1,  \frac{\l\l'}{\l+\l'} \Big) $. For any $n \in \hN $, since
             $    E \neg \left[ e^{\l  \g  (Y^n)^-_*} + e^{\l'  \g   (Y^n)^+_*  } \right]  \le   E \left[  e^{ \l   \g    \sL_* }   +e^{ \l'  \g     \sY_*    }    \right]< \infty $
             by \eqref{eq:c151},
                     applying Proposition \ref{prop:Z_K_bound}  with $p=p_o$ yields that
   \bea \label{eq:c255}
             E \neg \left[ \neg \left( \int_0^T     |Z^n_s|^2 ds \right)^{p_o} \dneg +  \big(K^n_T\big)^{p_o}  \right]
   \neg \le     \neg   c_{\l,\l'}      E \neg \left[ e^{\l  \g  (Y^n)^-_*} + e^{\l'  \g   (Y^n)^+_*  } \right]       \neg     \le    c_{\l, \l'}  \Xi  < \infty   ,
   \eea
    which shows that  $ \{Z^n\}_{n \in \hN}$ is a bounded subset in the reflexive Banach space $\hH^{2, 2 p_o}_\bF([0,T];\hR^d)$.
  Hence Theorem 5.2.1 of \cite{FA_Yosida} implies that $\{Z^n\}_{n \in \hN}$
  has a weakly convergent subsequence (we still denote it by $\{Z^n\}_{n \in \hN}$)
  with   limit $Z \in \hH^{2, 2 p_o}_\bF([0,T];\hR^d)$.

   \ss  Next, we show that this convergence is indeed a strong one  in $\hH^2_\bF([0,T];\hR^d)$. In the second step, we will introduce
  a function that will be useful in establishing this goal and develop several inequalities which will play important roles in the sequel.

   \ms \no {\bf  2)}
     Define a   function   $   \f(x) \dfnn \frac{   1}{\l_o \g } \left( e^{ \l_o \g |x|} -\l_o \g |x| -1 \right)  \ge 0$,  $ \fa    x \in \hR$.
  Fix $n \in \hN$. For  any $m \in \hN$ with   $m \ge n$, since
 $\big|\f'(x) \big|= e^{\l_o \g |x|}-1 $, $ x \in \hR$,   it follows from \eqref{eq:c151} that \pas
   \bea  \label{eq:c401}
       \left| \f ' \left(Y^m_t-Y^n_t\right) \right|  <    e^{\l_o \g   | Y^m_t-Y^n_t |}    \le   e^{ \l_o \g (\sL_t+ \sY_t)  }    ,  \q    t \in [0,T].
     \eea
        Applying It\^o's formula to the process  $\f \big(Y^m_\cd-Y^n_\cd\big)$  yields that
    \bea
     \qq    &&  \hspace{-1.7cm}   \f \left(Y^m_t-Y^n_t\right) + \frac12    \int_t^T  \f '' \left(Y^m_s-Y^n_s\right) |Z^m_s-Z^n_s|^2 ds \nonumber \\
     & &   =    \f \left(\xi_m-\xi_n\right)+     \int_t^T  \f ' \left(Y^m_s-Y^n_s\right)
     \big(f_m (s, Y^m_s, Z^m_s) - f_n (s, Y^n_s, Z^n_s)  \big) ds  \nonumber \\
      &&   \q +      \int_t^T   \f ' \left(Y^m_s-Y^n_s\right) (d K^m_s    -    d K^n_s )
   -\int_t^T  \f ' \left(Y^m_s-Y^n_s\right) \left(Z^m_s-Z^n_s\right) dB_s ,   \q    t \in [0,T].  \qq  \label{eq:c333}
   \eea
  First, we  argue that the stochastic integral term in \eqref{eq:c333} is a martingale. Applying Young's inequality with
    \bea    \label{p1_p3}
 p_1=\frac{\l}{\l_o} ,  \q  p_2=\frac{\l'}{\l_o} \q \hb{and} \q    p_3=\left(1-\frac{1}{p_1}-\frac{1}{p_2} \right)^{-1} =\frac{\l\l'}{\l\l'-\l_o(\l+\l')} = p_o   ,
 \eea
 we can deduce from the Burkholder-Davis-Gundy inequality,  \eqref{eq:c401},  and \eqref{eq:c255} that
   \bea
     \q  &&  \hspace{-2 cm}   E\left[\underset{t \in [0,T]}{\sup}  \left|\int_0^t  \f ' \left(Y^m_s-Y^n_s\right) \left(Z^m_s-Z^n_s\right) dB_s\right|\right]
  \le  c_0 E \left[   \left(\int_0^T  \big|\f ' (Y^m_s-Y^n_s)   \big|^2  | Z^m_s-Z^n_s |^2 d s\right)^{\frac{1}{2}}\right] \nonumber  \\
 &&    \le   c_0 E \left[ \underset{s \in [0,T]}{\sup}  \big| \f '  (Y^m_s-Y^n_s ) \big| \cd \left(1+\int_0^T  | Z^m_s-Z^n_s |^2 d s\right) \right] \nonumber \\
 &&    \le    c_{\l,\l'} \,   E\left[  e^{\l_o p_1  \g    \sL_*  } +e^{\l_o p_2  \g    \sY_*  }
  + \left(1+\int_0^T  | Z^m_s-Z^n_s |^2 d s\right)^{p_o} \right]       \le    c_{\l, \l'}  \big(1+\Xi \big)     < \infty  . \label{eq:c335}
   \eea
 Thus $\int_0^\cd  \f ' \left(Y^m_s-Y^n_s\right) \left(Z^m_s-Z^n_s\right) dB_s $ is a
  uniformly integrable martingale.  Letting $t=0$, taking expectation in \eqref{eq:c333}, and using (H1) we obtain
   \bea
     \q   && \hspace{-1cm}   E \left[ \f \big(Y^m_0-Y^n_0 \big) \right] + \frac12  E \neg \int_0^T \neg \f ''  \big(Y^m_s-Y^n_s\big) |Z^m_s-Z^n_s|^2 ds
        \le   E \left[  \f \big(\xi_m-\xi_n \big) \right] +   E \neg \int_0^T \neg \f ' \big(Y^m_s-Y^n_s\big) (d K^m_s    -    d K^n_s ) ~\;   \nonumber \\
  &&    \hspace{-0.4cm} +   E  \neg \int_0^T \dneg \big|  \f ' (Y^m_s \neg - \neg  Y^n_s ) \big|
     \bigg(2 \a  \neg  + \neg  \beta   |Y^m_s|  \neg + \neg \beta   |Y^n_s|
  \neg  + \neg    \frac12 \g    \Big( \neg   2|Z^m_s  \neg  - \neg  Z^n_s|^2  \neg + \neg  (\l_o  \neg -  \neg  2)|Z_s  \neg - \neg  Z^n_s|^2
  \neg + \neg \big(3 \neg + \neg \hb{$\frac{9}{\l_o \neg - \neg  5}$}\big) | Z_s  |^2 \Big) \neg  \bigg)   ds     ,        \label{eq:c336}
    \eea
    where we used the fact that $ | Z^m_s|^2 +   | Z^n_s|^2  \le  2|Z^m_s-Z^n_s|^2 + 3| Z^n_s|^2 $ and that
    \beas
       | Z^n_s|^2 \le \big(|Z_s - Z^n_s| + | Z_s  |\big)^2 \le  \big(1+\hb{$\frac{\l_o-5}{3}$}\big) |Z_s - Z^n_s|^2 +\big(1+\hb{$\frac{3}{\l_o-5}$}\big) | Z_s  |^2.
    \eeas

    Since it holds \pas~ that
     \beas
     |Y^m_t-Y^n_t|  \le |Y_t-Y^n_t |  \le  |Y_t-Y^1_t|    , \q  t \in [0,T],
      \eeas
    one can deduce from the monotonicity  of functions $\f$ and $|\f'|$ that \pas,     $ \f  (\xi_m- \xi_n  ) \le \f  (\xi - \xi_n  ) $ and
     \bea  \label{eq:c337}
       \big| \f'  (Y^m_t-Y^n_t  ) \big| \le \big| \f '  (Y_t -Y^n_t  )\big| \le   \big| \f '  (Y_t -Y^1_t  )\big|, \q   t \in [0,T].
      \eea
      Similar, it holds \pas~ that
       \bea   \label{eq:c338}
         \big| \f'  (L^m_t - L^n_t  ) \big| \le \big| \f '  (L_t - L^n_t  )\big|  \le \big| \f '  (L_t - L^1_t  )\big|,  \q   t \in [0,T].
      \eea
       We  also see  from \eqref{eq:c335} that
      \bea  \label{eq:c261}
   E \int_0^T   \big| \f'   (Y^m_s-Y^n_s ) \big| \,  |Z^m_s-Z^n_s|^2 ds  \le     E\left[  \underset{s \in [0,T]}{\sup}  \big| \f ' (Y^m_s-Y^n_s ) \big| \,       \int_0^T       | Z^m_s-Z^n_s   |^2          ds    \right]     < \infty    ,
   \eea
       which together with    \eqref{eq:c336}, \eqref{eq:c337} and \eqref{eq:c151} implies that
   \bea
     \q  &&   \hspace{-1.4 cm}
      E \int_0^T  \big(   \f '' - 2 \g | \f' | \big) \left(Y^m_s-Y^n_s\right) |Z^m_s-Z^n_s|^2 ds
    \le   2 E   \left[ \f \left(\xi -\xi_n \right) \right]   +  2  E  \int_0^T \f ' \left(Y^m_s-Y^n_s\right) (d K^m_s    -    d K^n_s ) \qq \q    \nonumber   \\
      & &        +   E \neg  \int_0^{T}   \neg  \big| \f '  (Y_s \neg - \neg Y^n_s ) \big|
      \Big( 4 \a   \neg + \neg  2\beta ( \sL_s + \sY_s    )
           \neg     +    \neg      (\l_o  \neg -  \neg  2)\g |Z_s  \neg - \neg  Z^n_s|^2
  \neg + \neg \big(3 \neg + \neg \hb{$\frac{9}{\l_o \neg - \neg  5}$}\big) \g | Z_s  |^2     \Big)  ds .  \label{eq:c339}
    \eea
\qq  Now we estimate the second term on the right-hand-side of \eqref{eq:c339} by  two cases of assumption (M3).  Assume (M3a) first.
      Since $\f'$ is an increasing and continuous function on $\hR$, the flat-off condition of $(Y^m,Z^m,K^m)$, \eqref{eq:c255}  and \eqref{eq:c338}  imply that
 \bea
      && \hspace{-1.7cm}   E  \int_0^T \f ' \left(Y^m_s-Y^n_s\right) (d K^m_s    -    d K^n_s )
       \le   E  \int_0^{T} \neg \f ' \left(Y^m_s-Y^n_s\right)  d K^m_s
    \le   E  \int_0^{T}  \neg \f ' \left(Y^m_s-L^n_s\right)  d K^m_s   \nonumber  \\
     &=&    E  \int_0^T \b1_{\{Y^m_s = L^m_s\}} \f ' \left(Y^m_s-L^n_s\right) d K^m_s
      =    E  \int_0^T \b1_{\{Y^m_s = L^m_s\}} \f ' \left(L^m_s-L^n_s\right) d K^m_s
        \nonumber \\
   &  \le  &     \left\| K^m_T  \right\|_{\hL^{p_o}(\cF_T)}
   \left\| \f ' \left(L^m-L^n \right) \right\|_{\hC^{\frac{p_o}{p_o -1 }}_\bF[0,T] }
     \le    c_{\l, \l'} \,  \Xi^{\frac{1}{p_o}}     \left\| \f ' \left(L-L^n \right) \right\|_{\hC^{\frac{p_o}{p_o -1 }}_\bF[0,T] }  .  
 \eea
 On the other hand, it holds for the case of (M3b) that
   \bea
       && \hspace{-1.7cm}E  \int_0^T \f ' \left(Y^m_s-Y^n_s\right) (d K^m_s    -    d K^n_s )
        \le  -  E  \int_0^T    \f ' \left(L^m_s-Y^n_s\right) d K^n_s
      =       -  E  \int_0^T  \b1_{\{Y^n_s = L^n_s\}}   \f ' \left(L^m_s-L^n_s\right) d K^n_s     \nonumber  \\
   &&   \le    \left\| K^n_T  \right\|_{\hL^{p_o}(\cF_T)}
   \left\| \f ' \left(L^m_\cd-L^n_\cd \right) \right\|_{\hC^{\frac{p_o}{p_o -1 }}_\bF[0,T] }
      \le   c_{\l, \l'}  \,   \Xi^{\frac{1}{p_o}}   \left\| \f ' \left(L -L^n  \right) \right\|_{\hC^{\frac{p_o}{p_o -1 }}_\bF[0,T] } .   \label{eq:c264}
 \eea
      \bea
 &&\hspace{-1.7cm}\hb{    {\bf 3)} Since  the sequence $\left\{\sqrt{ \big|\f ' (Y^m-Y^n )\big| }   \big(Z^m-Z^n \big) \right\}_{m \ge n  } $
  weakly converges  to}  \hspace{4cm}\nonumber  \\
  && \qq \qq \qq \q ~    \hb{   $\sqrt{ \big|\f ' (Y - Y^n )\big| } \big(Z-Z^n \big) $   in $ \hH^2_\bF ([0,T];\hR^d)$,}  \label{claim1}
\eea
which is proved in  Subsection \ref{appendix_1},  Theorem 5.1.1 ii) of \cite{FA_Yosida} shows that
    \bea  \label{eq:c275}
      E \int_0^T \big|\f '  \big(Y_s-Y^n_s \big) \big| \, |Z_s - Z^n_s |^2    ds  \le     \linf{m \to \infty}
      E \int_0^T \big| \f ' \big(Y^m_s-Y^n_s \big) \big| \,        | Z^m_s-Z^n_s   |^2  ds  .
     \eea
    As  $\hH^{2, 2 p_o}_\bF([0,T];\hR^d) \subset \hH^2_\bF([0,T];\hR^d) $,    the sequence $\{Z^m\}_{m \ge n }$
 also weakly converges to    $Z$ in $\hH^2_\bF([0,T];\hR^d)$. Applying Theorem 5.1.1 ii) of \cite{FA_Yosida} once again,
 we can deduce from  \eqref{eq:c339}-\eqref{eq:c264}  and \eqref{eq:c275} that
   \bea
 \q     && \hspace{-1cm}   \l_o \g     E\int_0^T      |Z_s-Z^n_s|^2 ds
   \le  \l_o \g   \lsup{m \to \infty}   E\int_0^T      |Z^m_s-Z^n_s|^2 ds  \nonumber   \\
    & &  =    \lsup{m \to \infty}   E\int_0^T  \big(   \f '' - \l_o \g  |\f'|  \big)  \left(Y^m_s-Y^n_s\right)   |Z^m_s-Z^n_s|^2 ds
  \q  \big(\because  \f '' (x) - \l_o \g  |\f'(x)| = \l_o \g, ~\fa x \in \hR  \big)    \nonumber  \\
     &  & =     \lsup{m \to \infty}   E\int_0^T  \neg  \big(   \f ''  \neg - \neg  2 \g  |\f'| \big)  \left(Y^m_s  \neg - \neg Y^n_s\right)
       |Z^m_s \neg - \neg Z^n_s|^2 ds          -(\l_o \neg -\neg 2)\g  \linf{m \to \infty}  E \int_0^T   \neg  \big|\f ' (Y^m_s  \neg - \neg Y^n_s ) \big|
           | Z^m_s \neg - \neg Z^n_s   |^2          ds ~\; \;  \nonumber  \\
          &  &  \le  2    E \left[  \f \left(\xi-\xi_n\right) \right]
          +   c_{\l, \l'} \, \Xi^{\frac{1}{p_o}}     \left\| \f ' \left(L-L^n \right) \right\|_{\hC^{\frac{p_o}{p_o -1 }}_\bF[0,T] }  \nonumber  \\
  && \q  +   E \int_0^T \big| \f '  (Y_s-Y^n_s ) \big|
     \Big(   4 \a  + 2\beta  ( \sL_s + \sY_s    )    +  \big(3 \neg + \neg \hb{$\frac{9}{\l_o \neg - \neg  5}$}\big) \g  | Z_s  |^2  \Big) ds . \qq \label{eq:c283}
                 \eea
         Since   $\l_o  <  \frac{\l\l'}{\l+\l'}  $, it follows that $ \l'  >    \frac{\l_o \l}{ \l  - \l_o   } $. Applying Young's inequality with $\wt{p}=\frac{\l}{\l_o} $ and $\wt{q}=\frac{\l}{\l-\l_o}$, we can deduce from   \eqref{eq:c151}       that   \pas
   \bea
      0 &\le&    \f \left(\xi -\xi_n \right)     \le \hb{$\frac{1}{\l_o \g }$}   e^{ \l_o \g | \, \xi -\xi_n| }
   \le  \hb{$\frac{1}{\l_o \g }$}   e^{ \l_o \g ( \sL_* + \sY_*    ) }    \nonumber \\
   &\le&  c_{\l,\l'} \big( e^{ \l \g   \sL_*   }   + e^{ \frac{\l_o \l }{\l-\l_o} \g   \sY_*   } \big)
   \le     c_{\l,\l'} \big( e^{ \l \g   \sL_*   }   + e^{ \l'    \g   \sY_*   } \big) ,   \q  \fa   n   \in   \hN .        \label{eq:c284}
   \eea
  As $E \left[ e^{ \l \g   \sL_*   }   + e^{ \l'    \g   \sY_*   } \right] < \infty $,
  the continuity of function $\f$ and the Dominated Convergence Theorem imply that
  \bea   \label{eq:c285}
  \lmtd{n \to \infty}  E \left[ \f \left(\xi -\xi_n \right) \right] = 0  .
 \eea

   Next, we analyze the convergence of the second term on the right-hand-side of \eqref{eq:c283}.
   In virtue of Dini's Theorem, it holds \pas~ that
   $  \lmt{n \to \infty} \,\underset{t \in [0,T]}{\sup} \big| L_t-L^n_t \big| = 0  $. Then
 the continuity  of function $\f' $ implies that
   \beas
     0&=& \lmt{n \to \infty}   \left| \f ' \left(\underset{t \in [0,T]}{\sup} \big| L_t-L^n_t \big| \right) \right|
      =    \lmt{n \to \infty}  \exp \left\{   \l_o \g \underset{t \in [0,T]}{\sup} \big| L_t-L^n_t \big| \right\} -1 \\
     & =&  \lmt{n \to \infty} \,  \underset{t \in [0,T]}{\sup}  \exp \left\{   \l_o \g  \big| L_t-L^n_t \big| \right\} -1
     =  \lmt{n \to \infty} \, \underset{t \in [0,T]}{\sup} \left |\f '  \big( L_t-L^n_t \big) \right|    , \q \pas
  \eeas
     It follows from \eqref{eq:c338}   that \pas
    \beas
           \underset{t \in [0,T]}{\sup}  \left| \f ' \left(L_t-L^n_t\right) \right|^{\frac{p_o}{p_o -1 }}  \le    \underset{t \in [0,T]}{\sup}  \left| \f ' \left(L_t-L^1_t\right) \right|^{\frac{p_o}{p_o -1 }},     \q  \fa n \in \hN .
     \eeas
   Applying Young's inequality  with $\wt{p}= \frac{\l+ \l'}{\l'}$ and $\wt{q} = \frac{\l+ \l'}{\l}$,   one can deduce from \eqref{eq:c151}    that
   \bea
     E\left[  \underset{t \in [0,T]}{\sup} \left| \f ' \left(L_t-L^1_t\right) \right|^{\frac{p_o}{p_o -1 }} \right]   & \tneg = & \tneg
       E\left[  \underset{t \in [0,T]}{\sup} \left| \f ' \left(L_t-L^1_t\right) \right|^{\frac{\l\l'}{\l_o(\l+\l' \neg ) }} \right]   \le
    E\left[  \underset{t \in [0,T]}{\sup} e^{  \frac{\l\l'}{\l+\l' } \g |L_t-L^1_t|}  \right]                \nonumber \\
    &   \tneg      \le &  \tneg    E \left[    e^{ \frac{\l\l'}{\l+\l' } \g (\sL_*+ \sY_* ) }    \right]        \le     c_{\l,\l'}   E \left[e^{ \l  \g \sL_*}    +  e^{ \l' \g       \sY_*  }     \right]     < \infty .  \qq  \;\;     \label{eq:c265}
   \eea
   The Dominated Convergence Theorem then  implies that
 \bea    \label{eq:c253}
 \lmtd{n \to \infty}  E\left[  \underset{t \in [0,T]}{\sup} \left|\f ' \left(L_t-L^n_t\right) \right|^{\frac{p_o}{p_o -1 }} \right]  =0  .
 \eea
 Similar to \eqref{eq:c265}, one has
  \bea    \label{eq:c265b}
     E\left[  \underset{t \in [0,T]}{\sup} \left| \f ' \left(Y_t-Y^1_t\right) \right|^{\frac{p_o}{p_o -1 }} \right] \le E \left[    e^{ \frac{\l\l'}{\l+\l' } \g (\sL_*+ \sY_* ) }    \right]        \le     c_{\l,\l'}  \Xi      < \infty .
  \eea

 Now we will analyze the convergence of the third term on the right-hand-side of \eqref{eq:c283}.
  We can deduce from   \eqref{eq:c337}, \eqref{eq:c151},   as well as \eqref{eq:c269} that \pas
 \beas
  \qq \qq && \hspace{-2cm}  \big| \f '  (Y_t \neg - \neg Y^n_t )    \big|   \Big(4 \a  \neg  + \neg  2\beta  ( \sL_t + \sY_t    )
               \neg + \neg    \big(3 \neg + \neg \hb{$\frac{9}{\l_o \neg - \neg  5}$}\big) \g   | Z_t  |^2    \Big)
               \le    \big| \f '  (Y_t \neg - \neg Y^1_t )    \big|  \Big(4 \a   \neg + \neg  4 \beta \big(  \sL_t  \neg + \neg  \sY_t  \big)
               \neg + \neg    \big(3 \neg + \neg \hb{$\frac{9}{\l_o \neg - \neg  5}$}\big) \g   | Z_t  |^2    \Big)   \nonumber   \\
       &&   \le c_{\l,\l'}  \big| \f '  (Y_t-Y^1_t )    \big| e^{ \big(\frac{\l \l' }{\l + \l'} - \l_o\big) \g  (  \sL_t + \sY_t   )  }
             +    \big(3 \neg + \neg \hb{$\frac{9}{\l_o \neg - \neg  5}$}\big) \g \big| \f '  (Y_t-Y^1_t )  \big|  | Z_t  |^2  \nonumber   \\
             &&  \le c_{\l,\l'}   e^{ \frac{\l \l' }{\l + \l'} \g (  \sL_t + \sY_t   )  }
             +    \big(3 \neg + \neg \hb{$\frac{9}{\l_o \neg - \neg  5}$}\big) \g \big| \f '  (Y_t-Y^1_t )  \big|  | Z_t  |^2, \q  \fa   t \in [0,T] , ~ \fa n \in \hN.
  \eeas
   Young's inequality, \eqref{eq:c265} and \eqref{eq:c265b} show that
       \beas
 \qq && \hspace{-2cm}  E  \int_0^T e^{ \frac{\l \l' }{\l + \l'} \g (  \sL_t + \sY_t   )  } dt
    + E   \int_0^T \big| \f '  (Y_t-Y^1_t )  \big|  | Z_t  |^2 dt   \\
  &&   \le T  E\left[   e^{ \frac{\l \l' }{\l + \l'} \g (  \sL_* + \sY_*   )  }   \right] +  c_{\l,\l'} E \left[  \underset{t \in [0,T]}{\sup} \left| \f ' \left(Y_t-Y^1_t\right) \right|^{\frac{p_o}{p_o -1 }} + \left( \int_0^T    | Z_t  |^2 dt \right)^{p_o}\right] < \infty.
   \eeas
   Then the continuity  of function $\f' $ and the Dominated Convergence Theorem imply that
     \beas
       \lmt{n \to \infty}  E \int_0^T \big| \f '  (Y_s-Y^n_s ) \big|
     \Big(   4 \a  + 2\beta  ( \sL_s + \sY_s    )    +  \big(3 \neg + \neg \hb{$\frac{9}{\l_o \neg - \neg  5}$}\big) \g  | Z_s  |^2  \Big) ds = 0 ,
     \eeas
   which together with       \eqref{eq:c285} and  \eqref{eq:c253} leads to that
       \bea  \label{eq:c033}
       \lmt{n \to \infty}   E\int_0^T      |Z_s-Z^n_s|^2 ds     = 0  .
      \eea
  Therefore,  the sequence $\left\{Z^n\right\}_{n \in \hN}$
  strongly converges to $Z$ in $\hH^2_\bF([0,T]; \hR^d)$.
  Consequently, Doob's martingale inequality  implies that
      \bea\label{eq:c035}
       \lmt{n \to \infty}   E \left[ \underset{t \in [0,T]}{\sup} \left|\int_0^t     \big(Z_s-Z^n_s \big) d B_s \right|^2 \right] = 0  .
      \eea

  \ss In the next step, we will show   that  $Y \in \hE^{\l \g,\l' \g}_\bF[0,T]$.

  \ms \no {\bf  4)} We  first develop a few auxiliary results.  By \eqref{eq:c033}, we can extract a
  subsequence of $\{ Z^n \}_{n \in \hN}$ (we still denote it by $\{ Z^n \}_{n \in  \hN}$)  such that
  $    \lmt{n \to \infty}  Z^n_t = Z_t $, \dtp ~
   In fact, we can choose this subsequence so that
 $ Z^* \dfnn   \underset{n \in \hN}{\sup}|Z^n|      \in   \hH^2_\bF[0,T] $; see \cite{Lep_San_97} or \cite[Lemma 2.5]{Ko_2000}.
   By (M2), it holds   \dtp ~ that
 \bea  \label{eq:c037}
       f \big(t,\o, y, z \big) = \lmt{n \to \infty}  f_n\big(t, \o, y, z \big), \q \fa (y,z) \in \hR \times \hR^d,
  \eea
   which together with the measurability of $f_n$, $n \in \hN$ implies that $f$
   is also $\sP \times \sB(\hR) \times   \sB(\hR^d)/\sB(\hR)$-measurable.
    Moreover, we see from  \eqref{eq:c037}
    and (M1) that $f$ also satisfies (H1).    
   For \dtp~ $(t, \o) \in [0,T] \times \O$,   the continuity of mapping $ f (t,\o, \cd, \cd)$ shows that
  \bea  \label{eq:c039}
       \lmt{n \to \infty}  \left|  f \big(t, \o, Y^n_t(\o), Z^n_t(\o) \big) -  f \big(t,\o, Y_t(\o), Z_t(\o) \big)  \right| =0 .
  \eea
  On the other hand, (M2) implies that   for \dtp~ $(t, \o) \in [0,T] \times \O$,
  \beas
     0 & \le &  \linf{n \to \infty}  \left|  f_n\big(t, \o, Y^n_t(\o), Z^n_t(\o) \big) -  f \big(t,\o, Y^n_t(\o), Z^n_t(\o) \big)  \right|  \\
        & \le &       \lmt{n \to \infty} \Big( \sup  \Big\{\left|  f_n\big(t, \o, y, z \big) -  f \big(t,\o, y, z \big)  \right| : \,
         |y| \le |Y^1_t(\o)| \vee |Y_t(\o)| <\infty ,\;  |z| \le Z^*_t(\o) <\infty  \Big\}  \Big) = 0  ,
   \eeas
   which together with \eqref{eq:c039} yields that \dtp
    \bea \label{eq:x151}
      \lmt{n \to \infty}  \left|  f_n \big(t, \o, Y^n_t(\o), Z^n_t(\o) \big) -  f \big(t,\o, Y_t(\o), Z_t(\o) \big)  \right| =0   .
    \eea
  Moreover,  (H1) and \eqref{eq:c151}   show that  \dtp
  \bea  \label{eq:x155}
        \left| f_n(t,  Y^n_t, Z^n_t) \neg - \neg   f\big(t, Y_t, Z_t \big) \right|
     & \le &    2  \a  \neg + \neg  \beta | Y^n_t|  \neg + \neg \beta | Y_t|   \neg + \neg   \frac{\g}{2} \big( \big| Z^n_t \big|^2
    \neg + \neg    \big| Z_t \big|^2    \big)  \nonumber  \\
   &         \le &     2  \a  \neg + \neg 2 \beta  (  \sL_* + \sY_*   )
   \neg + \neg  \frac{\g}{2} \big( \big|  Z^*_t \big|^2 \neg + \neg    \big| Z_t \big|^2 \big) , \q \fa n \in \hN  .
  \eea

     Let us assume that except on a $P$-null set $\sN $, \eqref{eq:x151}, \eqref{eq:x155} hold for a.e. $t \in [0,T]$ and
           $ \sL_* + \sY_* + \int_0^T \neg \big( \big|  Z^*_t \big|^2  \neg + \neg  \big|  Z_t \big|^2 \big) dt < \infty $.
     For any $\o \in \sN^c  $,   the Dominated Convergence Theorem implies that
    \bea \label{eq:x159}
    \lmt{n \to \infty}     \int_0^T  \left|  f_n\big(t, \o, Y^n_t(\o), Z^n_t(\o)\big) -  f\big(t, \o, Y_t(\o), Z_t(\o) \big)  \right|   dt =0  .
    \eea
 For any $n \in \hN$,  integrating with respect to $t$ in \eqref{eq:x155} yields that
     \beas
        \int_0^T  \neg  \left|  f_n\big(t, \o, Y^n_t(\o), Z^n_t(\o)\big)  \neg  - \neg   f\big(t, \o, Y_t(\o), Z_t(\o) \big)  \right|   dt \le
         c_{\l,\l'} e^{\frac{\l \l' \g}{(\l+\l') \, p_o}  (  \sL_* (\o)+ \sY_*(\o)   )}  \neg + \neg  \frac{\g}{2} \neg \int_0^T \dneg \Big( \big|  Z^n_t (\o)\big|^2  \neg + \neg  \big|  Z_t (\o)\big|^2 \Big) dt .
     \eeas
     Then it follows from  \eqref{eq:c265} and \eqref{eq:c255} that
     \beas
     E \neg \left[  \neg \left( \int_0^T  \left|  f_n\big(t, ,Y^n_t, Z^n_t\big)  \neg - \neg   f\big(t,  Y_t, Z_t \big)  \right|   dt \right)^{\neg p_o} \, \right]
     & \tneg  \neg \le& \tneg \neg   c_{\l,\l'}   E \left[ e^{\frac{\l \l'}{\l+\l'} \g (  \sL_* + \sY_*   )}
     \neg + \neg  \left( \int_0^T \neg   \big|  Z^n_t  \big|^2 dt \right)^{p_o}
       \neg + \neg  \left( \int_0^T \neg   \big|  Z_t  \big|^2 dt \right)^{p_o} \right]\\
     & \tneg  \neg \le& \tneg \neg  c_{\l,\l'}  \Xi + c_{\l,\l'}  E \left[  \left( \int_0^T \neg  \big|  Z_t  \big|^2 dt \right)^{p_o}   \right]  < \infty ,   \q   \fa n \in \hN,
     \eeas
     which implies that    $\left\{  \left( \int_0^T  \left|  f_n\big(t, ,Y^n_t, Z^n_t\big) \neg - \neg  f\big(t,  Y_t, Z_t \big)  \right|   dt \right)^{\frac{1+p_{ o}}{2}} \right\}_{n \in \hN} $ is uniformly integrable sequence in $\hL^1(\cF_T)$. Hence, one can deduce from \eqref{eq:x159} that
      \bea       \label{eq:c045}
        \lmt{n \to \infty} E \left[ \left(  \int_0^T  \left|  f_n\big(t,  Y^n_t, Z^n_t\big) -  f\big(t,  Y_t, Z_t \big)  \right|   dt
         \right)^{\neg \frac{1+p_{ o}}{2}} \right]= 0 .
      \eea

 Similar to \eqref{eq:c284}, it holds \pas ~  that
  \beas
           \big( \xi -\xi_n \big)^2   \le c_0 e^{\l_o \g | \xi -\xi_n | }
             \le     c_{\l,\l'} \big( e^{ \l \g   \sL_*   }   + e^{ \l'    \g   \sY_*   } \big)  ,   \q  \fa n \in \hN.
   \eeas
 As $E \left[ e^{ \l \g   \sL_*   }   + e^{ \l'    \g   \sY_*   } \right] < \infty $, applying the Dominated Convergence Theorem, we obtain
  \bea   
  \lmtd{n \to \infty}  E \left[ \big( \xi -\xi_n \big)^2   \right] = 0  .
 \eea
 Since $|\f'(x)| = e^{\l_o \g |x| } -1 \ge \l_o \g |x| $, $x \in \hR$, one can deduce from \eqref{eq:c253} that
     \bea    \label{eq:c453}
 \lmtd{n \to \infty}  \| L-L^n \|_{\hC^{\frac{p_o}{p_o -1 }}_\bF[0,T]}      =      0  .
     \eea
    Moreover,  for any $p \in [1, \infty) $,  \eqref{eq:c151} and \eqref{eq:c265} imply that
  \bea  \label{eq:c191}
   \|   Y^n\|^p_{\hC^p_\bF[0,T]}     \le E \Big[ \big( \sL_*+ \sY_* \big)^p \Big]
     \le c_{\l , \l',p}\,  E \left[  e^{ \frac{\l \l'}{\l + \l'} \g ( \sL_*+ \sY_*)  }  \right]  \le c_{\l , \l',p} \, \Xi \,, \q \fa n \in \hN    .
  \eea

   \ss     Now for any $m,  n \in \hN$ with $m \ge  n$,  applying It\^o's formula to the process $ \big(Y^m_\cd - Y^n_\cd\big)^2$  yields that
   \bea
     \qq  && \hspace{-1.5 cm}    ( Y^m_t - Y^n_t)^2 + \int_t^T  \big| Z^m_s- Z^n_s \big|^2 ds      =     (\xi_m-\xi_n)^2  + 2 \int_t^T   \big(Y^m_s - Y^n_s\big)    \big(f_m (s, Y^m_s, Z^m_s) -f_n (s, Y^n_s, Z^n_s)  \big)  ds\nonumber \\
   &&  + 2 \int_t^T  \big(Y^m_s - Y^n_s\big)  (d K^m_s-d K^n_s)
             -  2  \int_t^T  \big(Y^m_s - Y^n_s\big) \left( Z^m_s- Z^n_s   \right)  dB_s  ,  \q   t \in [0,T]. \q  \label{eq:c2958}
   \eea
      Let us estimate the term $\int_t^T \dneg  \big(Y^m_s - Y^n_s\big)  (d K^m_s-d K^n_s) $
      still under two cases of assumption (M3).  Assume  (M3a) first.
      The flat-off condition of $(Y^m,Z^m,K^m)$ implies that     \pas
    \beas
    \q  \int_t^T \dneg   (Y^m_s \neg - \neg  Y^n_s )  (d K^m_s \neg - \neg d K^n_s)  
   \le  \neg \int_t^T \dneg   (Y^m_s  \neg - \neg  L^n_s )  d K^m_s
     =   \neg \int_t^T \neg   (L^m_s  \neg - \neg  L^n_s)  d K^m_s
     \le        K^m_T \underset{s \in [0,T]}{\sup} \big|L^m_s  \neg - \neg  L^n_s\big|    , \q     t \in [0,T] .
 \eeas
 On the other hand, it holds for the case of (M3b) that \pas
 \beas
  \q  \int_t^T \dneg   (Y^m_s \neg - \neg  Y^n_s )  (d K^m_s \neg - \neg d K^n_s)  
   \le  \neg \int_t^T \dneg   (Y^n_s  \neg - \neg  L^m_s )  d K^n_s
    =   \neg \int_t^T \neg   (L^n_s  \neg - \neg  L^m_s)  d K^n_s
     \le        K^n_T \underset{s \in [0,T]}{\sup} \big|L^m_s  \neg - \neg  L^n_s\big|    ,   \q     t \in [0,T] .
 \eeas
             Then        \eqref{eq:c2958},   H\"older's inequality,  \eqref{eq:c255},  the Burkholder-Davis-Gundy inequality and \eqref{eq:c191} imply that
    \beas
           E \neg \left[ \underset{t \in [0,T]}{\sup}|Y^m_t  \neg - \neg  Y^n_t|^2 \right]
     & \dneg  \dneg \le& \dneg  \dneg    E   \big[  (  \xi_m \neg - \neg \xi_n  )^2\big]
     \neg    + \neg  2  \| Y^m  \neg - \neg  Y^n \|_{\hC^{\frac{p_o+1}{p_o - 1}}_\bF[0,T] }
  \big\|  f_m (\cd, Y^m_\cd, Z^m_\cd)  \neg - \neg  f_n (\cd, Y^n_\cd, Z^n_\cd)    \big\|_{\wh{\hH}^{1,\frac{1+p_o}{2}}_\bF([0,T]; \hR)}   \nonumber    \\
    && \dneg  \dneg     +   \,    c_{\l, \l'} \,  \Xi^{\frac{1}{p_o}}  \| L^m \neg  -  \neg L^n \|_{\hC^{\frac{p_o}{p_o - 1}}_\bF[0,T] }
  +   c_0  \,  E \left[    \underset{t \in [0,T]}{\sup}|Y^m_t  \neg - \neg  Y^n_t| \cd \left( \int_0^T    \big| Z^m_s \neg - \neg  Z^n_s   \big|^2  ds \right)^{\frac12} \right]   \nonumber \\
 & \dneg  \dneg \le & \dneg  \dneg   E\left[  (  \xi_m-\xi_n  )^2\right]        +c_{\l,\l'} \Xi^{ \frac{p_o -1}{p_o+1}}
  \left\|  f_m (\cd, Y^m_\cd, Z^m_\cd) -f_n (\cd, Y^n_\cd, Z^n_\cd)    \right\|_{\wh{\hH}^{1,\frac{1+p_o}{2}}_\bF([0,T]; \hR)}     \nonumber   \\
    &&  \dneg  \dneg    +    c_{\l, \l'} \,  \Xi^{\frac{1}{p_o}}  \| L^m \neg - \neg L^n \|_{\hC^{\frac{p_o}{p_o - 1}}_\bF[0,T] }
  +     c_{\l,\l'} \Xi^{\frac12}   \left\|  Z^m- Z^n \right\|_{\hH^2_\bF([0,T]; \hR^d)}   .
  \eeas
       Hence,  we can deduce from  \eqref{eq:c045}-\eqref{eq:c453} and \eqref{eq:c033} that
  $\{Y^n\}_{n \in \hN}$ is a Cauchy sequence in $\hC^2_\bF[0,T]$. Let $\wt{Y}$ be its limit in $\hC^2_\bF[0,T]$. As
   $\lmtd{n \to \infty}E \bigg[ \underset{t \in [0,T]}{\sup}  \big| Y^n_t - \wt{Y}_t   \big|^2 \bigg]=0$,
   there exists a subsequence $\big\{ n_i \big\}_{i \in \hN}$ of $ \hN$ such that
   $   \lmtd{i \to \infty}  \underset{t \in [0,T]}{\sup}  \big| Y^{n_i}_t - \wt{Y}_t   \big| =0 $,    \pas ~
     Then the monotonicity of the sequence $\big\{Y^n\big\}_{n \in \hN}$ by (M3) implies that
     $ \lmtd{n \to \infty}  \underset{t \in [0,T]}{\sup}  \big| Y^n_t - \wt{Y}_t   \big| =0 $, \pas ~ Thus it holds    \pas that
       $         \wt{Y}_t = \lmt{n \to \infty} Y^n_t = Y_t$, $\fa  t \in [0,T]$,
       which shows that processes $  \wt{Y}$ and $Y$ are indistinguishable. To wit, $Y$ is a continuous process that satisfies
   \bea
      \lmtd{n \to \infty}  \underset{t \in [0,T]}{\sup}  \big| Y^n_t - Y_t   \big| =0  ,  \q   \pas               \label{eq:c047}
   \eea
                    Since   $  E \left[  e^{ \l   \g    Y^-_* }  +   e^{ \l'  \g     Y^+_*    }      \right]
      \le  E \left[  e^{ \l   \g    \sL_* }   +e^{ \l'  \g     \sY_*    }    \right] < \infty
      $ by \eqref{eq:c151}, we see that $Y \in \hE^{\l \g,\l' \g}_\bF[0,T]$.

\ms In the next step, we will construct a process $K \in    \hK_\bF[0,T] $ such that $(Y,Z,K)$ is a solution of the quadratic RBSDE$(\xi, f, L)$.

 \ms \no {\bf 5)}      Since $Y$ is a continuous process by  step 4,
  \bea \label{defn_K_process}
  K_t \dfnn Y_0- Y_t   -\int_0^t f(s, Y_s, Z_s) \, ds+ \int_0^t  Z_s dB_s ,  \q  t\in [0,T]
 \eea
  defines an $\bF$-adapted,  continuous process with $K_0=0$. In light of \eqref{eq:c045} and \eqref{eq:c035},
  there exists a subsequence of $\big\{(Y^n, Z^n)\big\}_{n \in \hN}$ \big(we still denote it
  by $\big\{(Y^n, Z^n)\big\}_{n \in \hN}$\big)  such that \pas
  \beas   
       \lmt{n \to \infty}  \left\{       \int_0^T  \left|  f_n(t, Y^n_t, Z^n_t) -  f\big(t, Y_t, Z_t \big)  \right|   dt
   +   \underset{t \in [0,T]}{\sup}  \left|\int_0^t     \big(Z^n_s-Z_s \big) d B_s \right| \,    \right\}  = 0 .
 \eeas
 This  together with \eqref{eq:c047} leads to  that
 \bea   \label{eq:c473b}
  \lmt{n \to \infty} \; \underset{t \in [0,T]}{\sup}    |K^n_t - K_t|   =0  , \q  \pas ,
 \eea
 which implies that  $K$ is also an increasing process. To  wit, $K \in \hK_\bF[0,T]$.
 Letting $n \to \infty $ in  \eqref{eq:c151},    we can deduce from    \eqref{defn_K_process} that \pas
  \beas
   L_t  \le Y_t = \xi  +     \int_t^T f(s, Y_s, Z_s) \, ds + K_T - K_t - \int_t^T  Z_s dB_s\,, \qq t \in [0, T].
  \eeas

 \ss \no {\bf  6)}  It remains to verify that $(Y,Z,K)$ satisfies the flat-off condition \eqref{flat-off}.  For any $p \in [1, \infty) $,    similar to \eqref{eq:c191}, \eqref{eq:c151}       implies  that   \pas
  \beas
     \underset{t \in [0,T]}{\sup}  \big| Y_t - Y^n_t \big|^p \le \big( \sL_*+ \sY_* \big)^p
     \le c_{\l,\l',p}  \,  e^{ \frac{\l \l'}{\l+\l' } \g ( \sL_*+ \sY_*)  } , \q \fa n \in \hN.
  \eeas
   Then one can deduce from \eqref{eq:c047}, \eqref{eq:c265} and the  Dominated Convergence Theorem  that
  \bea  \label{eq:c297}
        \lmtd{n \to \infty}    E \left[\underset{t \in [0,T]}{\sup}  \big| Y_t - Y^n_t \big|^p \right] =0 .
  \eea

   For any $n \in \hN$,  let us show that
 \bea    \label{eq:c485}
   \lmt{n \to \infty} E\int_0^T (Y_t - L_t) dK^n_t =0
 \eea
   by two cases of assumption (M3).   Assume (M3a) first. One can deduce from the flat-off condition of $(Y^n, Z^n, K^n)$  and \eqref{eq:c255}  that
   \beas
  \q  0   &\le&   E\int_0^T \neg (Y_t - L_t) dK^n_t  \le     E\int_0^T (Y_t - L^n_t) dK^n_t  =  E\int_0^T \neg (Y_t - Y^n_t)  dK^n_t \\
  & \le &   \| K^n_T \|_{\hL^{p_o}(\cF_T)} \| Y - Y^n \|_{\hC^{\frac{p_o}{p_o - 1}}_\bF[0,T] }
  \le c_{\l, \l'} \, \Xi^{\frac{1}{p_o}} \| Y - Y^n \|_{\hC^{\frac{p_o}{p_o - 1}}_\bF[0,T] } \,  .
   \eeas
 Thus \eqref{eq:c485} follows from \eqref{eq:c297}.    On the other hand,  it holds for  the case of (M3b) that
     \beas
  \q  0  & \le &  E\int_0^T \neg (Y_t - L_t) dK^n_t  \le     E\int_0^T (Y^n_t - L_t) dK^n_t  =  E\int_0^T \neg (L^n_t - L_t)  dK^n_t \\
  & \le  &  \| K^n_T \|_{\hL^{p_o}(\cF_T)} \| L - L^n  \|_{\hC^{\frac{p_o}{p_o - 1}}_\bF[0,T] }
    \le  c_{\l, \l'} \, \Xi^{\frac{1}{p_o}}  \| L  - L^n   \|_{\hC^{\frac{p_o}{p_o - 1}}_\bF[0,T] }  \,  .
   \eeas
Thus \eqref{eq:c485} follows from \eqref{eq:c453}.

     \ms     Now fix an $\o \in \O$ such that \eqref{eq:c473b} holds and that  $t \to  Y_t(\o) - L_t(\o)$ is a non-negative continuous function on $[0,T] $.
 For any $\e >0$, there exists an $N=N(\o) \in \hN$ such that
  \beas
    0 \le     \int_0^T \big(Y_t(\o) - L_t(\o)\big) dK_t (\o) \le \e+ \sum^N_{j=1}  m_j(\o)  \big(K_{\frac{j}{N}} (\o)- K_{\frac{j-1}{N}} (\o)\big),
 \eeas
   where $m_j(\o) \dfnn  \underset{t \in [\frac{j-1}{N}, \frac{j}{N}] }{\min} \big(Y_t (\o)- L_t (\o)\big) $.  Thus, it follows that
  \beas
    \q && \hspace{-1.2cm} 0 \le  \int_0^T \big(Y_t(\o) - L_t(\o)\big) dK_t (\o)
           \le         \e+ \sum^N_{j=1}  m_j(\o)  \big(K^n_{\frac{j}{N}}(\o)- K^n_{\frac{j-1}{N}}(\o)\big)
     +  2 \underset{t \in [0,T] }{\sup}   |K^n_t (\o)- K_t (\o)| \sum^N_{j=1}  m_j(\o)    \\
     &&   \le    \e+ \int_0^T \big(Y_t (\o) - L_t (\o)\big) dK^n_t (\o)
+  2 \underset{t \in [0,T] }{\sup}   \big|K^n_t (\o)- K_t (\o)\big| \sum^N_{j=1}  m_j(\o)   .
   \eeas
   As $n \to \infty$,   we obtain
  \beas
     0 \le \int_0^T \big(Y_t(\o) - L_t(\o)\big) dK_t (\o) \le  \e+ \linf{n \to \infty}\int_0^T (Y_t  (\o)- L_t (\o)) dK^n_t(\o).
  \eeas
 Then letting $\e \to 0$ yields that
   \beas
    0 \le  \int_0^T   \big(Y_t(\o) - L_t(\o)\big) dK_t (\o) \le    \linf{n \to \infty}\int_0^T   (Y_t  (\o)- L_t (\o)) dK^n_t(\o).
  \eeas
 Eventually,   Fatou's Lemma and \eqref{eq:c485}  imply that
  \beas
      0 \le E \int_0^T \big(Y_t  - L_t \big) dK_t  \le E\left[  \linf{n \to \infty}\int_0^T (Y_t   - L_t ) dK^n_t   \right] \le
  \lmt{n \to \infty}  E  \int_0^T (Y_t   - L_t ) dK^n_t     = 0,
 \eeas
which leads to \eqref{flat-off}.

   \ss \no {\bf  7)}   In the previous steps we constructed a solution of the    quadratic RBSDE$(\xi, f, L)$, namely  $(Y,Z,K)$.
  Since $Y \in \hE^{\l \g,\l' \g}_\bF[0,T]$,     Proposition \ref{prop:Z_K_bound} shows that $(Z,K) \in \hH^{2, 2 p}_\bF([0,T];\hR^d)    \times \hK^p_\bF[0,T] $ for any $p \in \Big(1, \frac{\l\l'}{\l+\l'}\Big)$.
                                    \qed

\section{Existence}


 \begin{thm}   \label{thm:existence}
 Let $(\xi, f, L)$ be a parameter set such that $f$  satisfies  (H1) and that
 \bea  \label{cond:f_conti}
   \hb{For $\dtp ~ (t,\o) \in [0, T] \times \O\,$, the mapping $  f(t, \o, \cd, \cd)$ is continuous.}
 \eea
    If  $\,   E \Big[  e^{ \l   \g    L^-_* }   + e^{ \l'  \g  e^{\beta T  }   (\xi^+ \vee     L^+_*  ) }     \Big]  < \infty$ for some
  $\l , \l' > 6  $ with $\frac{1}{\l}+\frac{1}{\l'} < \frac{1}{6} $,  then the quadratic  RBSDE\,$(\xi,f,L)$ admits a solution
   $(Y,Z,K) \in \underset{p \in \big(1,\frac{\l\l'}{\l+\l'} \neg \big)}{\cap}
     \hE^{\l \g,\l' \g}_\bF[0,T] \times  \hH^{2, 2 p}_\bF([0,T];\hR^d) \times \hK^{p}_\bF[0,T] $  that satisfies \eqref{a_priori}.

       \ss   In addition, if $\xi^+ \vee     L_* \in \hL^e(\cF_T)   $, then  this solution $(Y,Z,K)$ belongs to
     $ \hS^p_\bF[0,T] $ for all $p \in [1, \infty)$. More precisely, for any $ p \in (1, \infty)$ we have
       \bea
     &&    E \neg \left[  e^{ p \g   Y_* }       \right]                    \le   E  \neg \left[  e^{ p \g    L^-_* }      \right]
                    +    c_p \,  E  \neg \left[        e^{ p  \g e^{\beta T} \left ( \xi^+  \vee     L^+_* \right)  }        \right] < \infty,  \nonumber  \\
                          \hb{and}  &&         E \left[\left( \int_0^T     |Z_s|^2 ds \right)^p   + K^p_T \right]
            \le           c_p  \,  E  \neg  \left[          e^{3 p     \g  Y_*}          \right]     < \infty . \label{eq:x225}
       \eea
 \end{thm}

\ss \no  {\bf Proof:}
 Let $i, n  \in \hN $. For any $x \in \hR$,  we define $x^i\dfnn  x \vee (-i)  $ and $x^{i,n}\dfnn \big(x \vee (-i) \big)\land n$. It is plain to check that
    \bea  \label{eq:c169}
   \big(x^i\big)^- \vee \big(x^{i,n}\big)^-   \le x^- \q \hb{and} \q \big(x^i \big)^+ \vee \big(x^{i,n}\big)^+   \le x^+   .
     \eea
   Theorem 1 of \cite{KLQT_RBSDE} shows that the quadratic RBSDE$\big(\xi^{i,n}, f,  L^{i,n} \big)$
  admits a maximal bounded solution $\big(Y^{i,n}, Z^{i,n},   \\    K^{i,n}\big)     \in  \neg
    \hC^\infty_\bF[0,T] \neg \times \neg  \hH^2_\bF([0,T];  \hR^d)  \neg  \times \neg \hK_\bF[0,T]$.
           Then one can deduce from Proposition \ref{prop_estimate} and \eqref{eq:c169}     that \pas
   \bea  \label{eq:c241}
   -L^-_t \le - \big(L^{i,n}_t \big)^- \le     L^{i,n}_t   \le    Y^{i,n}_t   \le    c_0 +  \frac{1}{\g} \ln  E\left[  e^{      \g  e^{\beta T  }
  \big( (\xi^{i,n})^+ \vee    (L^{i,n})^+_*\big)  }\Big|  \cF_t  \right]     \le    c_0 +  \frac{1}{\g} \ln  M_t   ,     ~\;    t \in [0,T] ,
    \eea
     where $M_t \dfnn E\left[  e^{      \g  e^{\beta T  }   (\xi^+ \vee    L^+_*)  }\big|  \cF_t  \right] $.
      Moreover, Proposition \ref{prop_compa_bdd} implies that \pas
   \bea   \label{eq:c251}
     Y^{i+1,n}_t  \le  Y^{i,n}_t \le Y^{i,n+1}_t  , \q t \in [0,T] .
   \eea

   Now fix $i \in \hN$. It is clear that $L^i \in \hC^0_\bF[0,T]$ and that $\{L^{i,n}_t\}_{n \in \hN}$ is an increasing sequence in $n$ with $\lmtu{n \to \infty} L^{i,n}_t =L^i_t$ for any $t \in [0,T]$.
   We see from \eqref{eq:c241} and \eqref{eq:c251} that except on a $P$-null set $\sN_i$, $\{Y^{i,n}_t\}_{n \in \hN}$ is an increasing sequence in $n$
    with an upper bound $ c_0 +  \frac{1}{\g} \ln  M_t $ for any $t \in [0,T]$.
 Thus, one can define a real-valued, $\bF$-adapted  process $Y^i_t (\o) \dfnn \b1_{\{\o \notin \sN_i\}} \lmtu{n \to \infty} Y^{i,n}_t (\o)$, $(t,\o) \in [0,T] \times \O$.
 Note that on $\sN^c_i$
     \bea   \label{eq:x213}
     Y^i_T = \lmtu{n \to \infty} Y^{i,n}_T = \lmtu{n \to \infty} \xi^{i,n} = \xi^i .
     \eea
 Letting $n \to \infty$ in \eqref{eq:c241} yields that   \pas
    \bea \label{eq:x211}
     -L^-_t \le Y^i_t  \le     c_0 +  \frac{1}{\g} \ln  M_t   ,     \q    t \in [0,T].
    \eea
           By \eqref{eq:c169},     $\sL^i_t   \dfnn     \big(L^{i,1}_t\big)^- \vee \big(L^i_t\big)^-  \le  L^-_t $, $\fa t \in [0,T]$.
     Also,     \eqref{eq:c241} and \eqref{eq:x211} imply  that \pas
          \beas
                    \sY^i_t \dfnn \big(Y^{i,1}_t\big)^+ \vee \big(Y^i_t\big)^+  \le    c_0 +  \frac{1}{\g} \ln  M_t ,  \q  t \in [0,T] .
     \eeas
  Then  it follows from    Doob's martingale inequality that
              \bea  \label{eq:c175}
                E \left[  e^{ \l   \g    \sL^i_* }   + e^{ \l'  \g     \sY^i_*    }    \right]  &\le& E \left[  e^{ \l   \g    L^-_* }      \right]+c_{\l'} E \left[     M^{\l'}_*     \right]
               \le E \left[  e^{ \l   \g    L^-_* }      \right] + c_{\l'} E \left[     M^{\l'}_T     \right]  \nonumber \\
               &=&    E \left[  e^{ \l   \g    L^-_* }      \right] +  c_{\l'}  E \left[        e^{ \l' \g  e^{\beta T  }  \left ( \xi^+  \vee     L^+_* \right)  }        \right] < \infty .
                 \eea
              Thus Theorem \ref{thm:monotone-stable} shows that   $Y^i \in \hE^{\l \g,\l' \g}_\bF[0,T]$ and that there exist
   $(Z^i,K^i) \in \underset{p \, \in \big(   1, \frac{\l\l'}{\l+\l'} \neg \big)}{\cap}  \hH^{2, 2 p}_\bF([0,T];\hR^d) \times \hK^{p}_\bF[0,T] $
    such that  $(Y^i,Z^i,K^i)$ is a solution of the quadratic  RBSDE\,$(Y^i_T, f, L^i)$.
  Moreover, letting $n \to \infty$ in    \eqref{eq:c251} yields  that   \pas
    \bea \label{eq:x219}
         Y^{i+1}_t    \le   Y^i_t  ,           \q     t \in [0,T]   .
         \eea

     Clearly,  $\{L^i_t\}_{i \in \hN}$ is a decreasing sequence in $i$ with $\lmtd{i \to \infty} L^i_t =L_t$ for any $t \in [0,T]$.
      We see from \eqref{eq:x211} and \eqref{eq:x219} that
     except on a $P$-null set $\sN$, $\{Y^i_t\}_{i \in \hN}$ is a decreasing sequence in $i$  with a lower bound $ -L^-_t $ for any $t \in [0,T]$.
 Thus, one can define a real-valued, $\bF$-adapted  process $Y_t (\o) \dfnn \b1_{\{\o \notin  \sN\}} \lmtd{i \to \infty} Y^i_t (\o)$, $(t,\o) \in [0,T] \times \O$.
     Letting $i \to \infty$ in \eqref{eq:x213} and \eqref{eq:x211} yields that   \pas
    \bea
         &&        Y_T =   \lmtd{i \to \infty} Y^i_T = \lmtd{i \to \infty} \xi^i = \xi   ,   \label{eq:x229}  \\
    \hb{and}     &&     -L^-_t \le Y_t  \le     c_0 +  \frac{1}{\g} \ln  M_t   ,     \q    t \in [0,T].   \label{eq:x215}
    \eea
           By \eqref{eq:c169},     $\sL_t   \dfnn     \big(L^1_t\big)^- \vee L^-_t  \le  L^-_t $, $\fa t \in [0,T]$.
     Moreover,     \eqref{eq:x211} and \eqref{eq:x215} imply  that \pas
          \beas
                    \sY_t \dfnn \big(Y^1_t\big)^+ \vee Y^+_t  \le    c_0 +  \frac{1}{\g} \ln  M_t ,  \q  t \in [0,T] .
     \eeas
       Similar to \eqref{eq:c175}, one can deduce that
              $                   E \left[  e^{ \l   \g    \sL_* }  +   e^{ \l'  \g     \sY_*    }      \right]
                   \le  E \left[  e^{ \l   \g    L^-_* }      \right]     +    c_{\l'}  E \left[        e^{ \l' \g  e^{\beta T  }  \left ( \xi^+  \vee     L^+_* \right)  }   \right] < \infty   $.
                                 Then         Theorem \ref{thm:monotone-stable}  and \eqref{eq:x229} imply that $Y \in \hE^{\l \g,\l' \g}_\bF[0,T]$ and that there exist
   $(Z,K) \in \underset{p \, \in \big(   1, \frac{\l\l'}{\l+\l'} \neg \big)}{\cap}  \hH^{2, 2 p}_\bF([0,T];\hR^d) \times \hK^{p}_\bF[0,T] $
    such that  $(Y,Z,K)$ is a solution of the quadratic   RBSDE\,$(\xi, f, L)$.

  \ss  Next, let us assume that   $\xi^+ \vee     L_* \in \hL^e(\cF_T)   $.
            For any $p \in (1, \infty)$,     we can deduce from \eqref{eq:a009}, \eqref{eq:x215} and Doob's martingale inequality that
         \beas
           E  \neg \left[  e^{ p \g    Y_* }      \right] &  \dneg \le &  \dneg   E  \neg \left[  e^{ p  \g  Y^-_*    }  +   e^{ p \g   Y^+_*    }  \right]
                    \le     E  \neg \left[  e^{ p  \g   L^-_* }      \right]     +    c_p \,  E \big[        M^{p}_*      \big]
                    \le  E  \neg \left[  e^{ p \g    L^-_* }      \right]     +    c_p \,  E \big[        M^{p }_T      \big]   \\
                 &  \dneg   = &  \dneg     E  \neg \left[  e^{ p  \g   L^-_* }      \right]
                    +     c_p \, E  \neg \left[        e^{ p \g  e^{\beta T  }  \left ( \xi^+  \vee     L^+_* \right)  }        \right]
                    \le c_p \, E  \neg \left[        e^{ p \g  e^{\beta T  }  \left ( \xi^+  \vee     L_* \right)  }        \right] < \infty  ,
           \eeas
               which shows that       $Y \in \hE^{p \g}_\bF[0,T]$. Finally,  an application of   Proposition \ref{prop:Z_K_bound} with  $\l =\l'=3 p  $ leads to \eqref{eq:x225}.    \qed

 \section{Comparison}

A function   $f:  [0,T] \times \O \times \hR \times \hR^d \rightarrow  \hR $ is said to be {\it convex} (resp. {\it concave}) in $z$ if it holds \dtp~ that
 \bea
   &&\hspace{-1.5cm}     f \big(t,\o, y, \th z_1 \neg + \neg (1-\th)z_2 \big)   \nonumber  \\
  && \le \hb{(resp. $\ge$) } \th  f  (t,\o, y,   z_1) \neg + \neg  (1-\th) f  (t,\o, y,   z_2), ~
 \fa (\th, y) \in (0,1) \times  \hR, \; \fa z_1, z_2 \in \hR^d.  \label{eq:axa011}
 \eea

 \ss  In the rest of the paper, we impose two more hypotheses on generator $f$ which together imply \eqref{cond:f_conti}.

\ms \no  {\bf (H2)} $f$ is {\it Lipschitz} in $y$: For some  $\k \ge 0$,   it holds \dtp~ that
  \bea       \label{eq:axa015}
   |f(t , \o, y_1,z)-f(t , \o, y_2,z) |  \le \k |y_1-y_2| , \q  \fa  y_1, y_2 \in \hR  , ~ \fa  z \in \hR^d.
  \eea

   \ss  \no  {\bf (H3)} $f$ is either convex or concave in $z$.

   \ms    From now on, for any $\l \ge 0$ the generic constant $c_\l$ also depends on $\k$ implicitly.
   Inspired by the ``$\th$-difference" method introduced in \cite{BH-07},
  we obtain  two  comparison theorems for quadratic RBSDEs with unbounded obstacles.

   \if{0}
    \ms For \dtp~$(t,\o) \in [0,T] \times \O$, let $\big\{(y_n,z_n)\big\}_{n \in \hN}\subset \hR \times \hR^d$ be a convergent sequence
   with  limit $(y,z) \in \hR \times \hR^d$. It holds for any $n \in \hN$ that
    \bea   \label{eqn-x01}
     |f(t , \o, y_n ,z_n )  -f(t , \o, y ,z)  | &\le & |f(t , \o, y_n ,z_n )  -f(t , \o, y ,z_n)  |+  |f(t , \o, y ,z_n)  -f(t , \o, y ,z)  |    \nonumber   \\
     &\le & \k | y_n - y | + |f(t , \o, y ,z_n)  -f(t , \o, y ,z)  |
    \eea
 As a  real-valued convex (resp. concave) function on $\hR^d$, $ f(t , \o, y , \cd)$ is definitely a continuous one.  Letting $n \to \infty$ in \eqref{eqn-x01}
 yields that $\lmt{n \to \infty} f(t , \o, y_n ,z_n )  =f(t , \o, y ,z) $. Therefore, (H2) and (H3) imply \eqref{cond:f_conti}.
  \fi

  \begin{thm}   \label{thm:comparison}
  Let $(\xi,f,L)$, $(\hat{\xi}, \hat{f},\wh{L})$ be two parameter sets  and let $(Y,Z,K) \big(\hb{resp. }  (\wh{Y}, \wh{Z},\wh{K})\big) 
 $ be a solution of  RBSDE$(\xi,f,L)$ \big(resp.\,RBSDE$(\hat{\xi}, \hat{f},\wh{L})$\big)  such that

 \ss \no (C1)   It holds \pas ~that  $    \xi  \le  \hat{\xi}   $ and that $  L_t  \le  \wh{L}_t    $ for any $  t  \in [0,T]$;

 \ss \no  (C2)    $   E \left[  e^{\l  Y^+_*} + e^{\l \wh{Y}^-_*}        \right] < \infty$ for all $  \l   \in (1, \infty) $  and $K \in   \hK^p_\bF[0,T]$ for some $p \in (1, \infty)$.

 \ss \no  For $\a,    \beta,  \k    \ge   0$, $\g   >    0$,   if  either of the following two  holds:

  \ss \no  (\,i)    $f$   satisfies (H1), (H2),  $f$   is convex in $z$,
 and $\D f(t) \dfnn  f\big(t, \wh{Y}_t, \wh{Z}_t  \big) - \hat{f} \big(t, \wh{Y}_t, \wh{Z}_t \big)   \le   0  $, \dtp;

  \ss \no  (ii)   $\hat{f}$ satisfies (H1), (H2),   $\hat{f}$ is convex in $z$, and
  $\D f(t) \dfnn  f(t, Y_t, Z_t  ) - \hat{f}(t, Y_t, Z_t ) \le 0  $, \dtp;

  \ss \no   then it holds \pas ~ that $Y_t \le \wh{Y}_t  $ for any $ t  \in [0,T]$.

    \end{thm}

\ss  \no  {\bf Proof:}
    Fix $\th \in (0,1) $.  We set $U \dfnn Y \neg - \th \wh{Y}$, $V \dfnn Z - \th \wh{Z}$      and define an $\bF$-progressively measurable process
   \bea
   \qq   a_t &\dfnn& \b1_{\{Y_t \ge 0\}} \left(  \b1_{\{Y_t \ne \wh{Y}_t \}}\frac{\fF \big( t,Y_t,\wh{Z}_t \big)
   -\fF \big( t,  \wh{Y}_t,\wh{Z}_t \big)}{Y_t - \wh{Y}_t      }
      - \k \b1_{\{Y_t =\wh{Y}_t   \}} \right)   - \k \b1_{\{   Y_t < 0 \le \wh{Y}_t \}}  \nonumber   \\
   && + \b1_{\{ Y_t \vee \wh{Y}_t <  0 \}} \left( \b1_{\{U_t \ne 0 \}}\frac{\fF(t,Y_t,Z_t) -\fF \big(t,\th \wh{Y}_t,Z_t \big)}{U_t
      } - \k \b1_{\{U_t = 0 \}} \right)    ,    \q   t \in [0,T] ,   \label{eq:axf011}
      \eea
     where  $ \fF$ stands for $f$ if (i) holds, and for $\hat{f}$ otherwise.
     It follows that $A_t \dfnn \int_0^t a_s \, ds$, $   t \in [0,T]$ is an $\bF$-adapted process.
  By  (H2), it holds \dtp~ that  $|a_t| \le \k$. Thus   $A_* \dfnn \underset{t \in [0,T]}{\sup} |A_t| \le \int_0^T \neg   | a_s |   ds \le  \k T$, \pas ~
  In light of (H1) and the convexity of $\fF$ in $z$, it holds \dtp~ that
  \bea  \label{eq:d131}
  \fF(t, y, Z_t)     \le  \th \fF \big(t, y,   \wh{Z}_t \big)  \neg +  \neg   (1 \neg - \neg \th) \fF  \neg  \left(t, y,   \frac{V_t}{1-\th}\right)
    \neg    \le    \th \fF \big(t,  y,   \wh{Z}_t \big)   \neg  + \neg    (1  \neg  -  \neg   \th)  \big(\a  \neg  + \neg   \beta |y|\big)
      \neg  + \neg   \hb{$\frac{\g}{2(1-\th)}$} |V_t|^2   , \q   \fa y \in \hR.
     \eea

\ss   Given $n \in \hN$, we define the $\bF$-stopping time
  \bea  \label{eq:axd071}
    \t_n \dfnn \inf\bigg\{ t \in [0,T]: \int_0^t  \big(       |Z_s|^2 + \big| \wh{Z}_s \big|^2  \big) ds >n  \bigg\} \land T.
    \eea
    Clearly, $\lmtu{n \to \infty} \t_n = T$, \pas ~  Let   $ \z_\th  \dfnn \frac{ \g e^{\k T}}{1-\th}$.    Applying It\^o's formula to the process
    $ \G_t \dfnn \exp\big\{\z_\th  e^{A_t }U_t \big\}  $, $t \in [0,T]$   yields that
  \beas
  \G_{\t_n \land t} = \G_{\t_n} + \int_{\t_n \land t}^{\t_n}  G_s ds + \z_\th   \int_{\t_n \land t}^{\t_n} \G_s e^{A_s}     (dK_s-\th d \wh{K}_s) - \z_\th   \int_{\t_n \land t}^{\t_n} \G_s e^{A_s}  V_s dB_s, \q t \in [0,T] ,
  \eeas
  where   $   G_t  =   \z_\th \,  \G_t e^{A_t}  \left(  f(t,Y_t,Z_t)-\th \hat{f}(t, \wh{Y}_t, \wh{Z}_t) -a_t U_t- \hb{$\frac12$} \z_\th  e^{A_t}|V_t|^2 \right) $.
   Clearly, it holds \dtp~ that
         \bea   \label{eq:x251}
    G_t  \le   \z_\th \,  \G_t e^{A_t}  \left(   \fF(t,Y_t,Z_t)-\th \fF(t, \wh{Y}_t, \wh{Z}_t) -a_t U_t- \hb{$\frac{\g}{2(1-\th)}$} |V_t|^2 \right)
  \eea
  whether  (i)  or (ii) holds. Moreover, let us show by 3 cases that
   \bea    \label{eq:d141}
   G_t \le   \g e^{2 \k T} \, \G_t \Big(\a + (\beta + \k) \big(Y^+_t + \wh{Y}^-_t \big)\Big) , \q \dtp
   \eea

   \ss \no 1) For \dtp ~ $(t,\o) \in \{ Y_t (\o) \ge 0 \}  $,   applying \eqref{eq:d131} with $y=Y_t$, we can deduce from    \eqref{eq:x251} and (H2)      that
    \beas
         G_t
       &    \tneg    \le   &    \tneg   \z_\th  \,  \G_t e^{A_t}   \neg \left(       \th \fF \big(t, Y_t,   \wh{Z}_t \big) - \th \fF \big(t, \wh{Y}_t, \wh{Z}_t \big)
        -a_t U_t+  (1-\th)  \big(\a+ \beta |Y_t|  \big)   \right)    \\
          &    \tneg = &    \tneg     \z_\th  \,  \G_t e^{A_t}    \Big(     \neg         (\th  \neg - \neg 1) \, a_t  Y_t
          \neg +  \neg  (1  \neg -  \neg \th)     \big(\a  \neg +  \neg  \beta |Y_t| \big)   \neg   \Big)
         \neg    \le  \neg    \g e^{2 \k T}  \G_t      \big(\a  \neg +  \neg  (\beta  \neg +  \neg \k) | Y_t | \big)
          \neg  =   \neg    \g e^{2 \k T}  \G_t      \big(\a  \neg +  \neg  (\beta  \neg +  \neg \k) Y^+_t \big)  .
   \eeas

     \ss \no 2) For \dtp ~ $(t,\o) \in \{ Y_t (\o) < 0 \le \wh{Y}_t (\o)   \}  $,   applying \eqref{eq:d131} with $y=0$,
     we see from    \eqref{eq:x251} and (H2)      that
                 \beas
        G_t           & \neg  \tneg \le& \neg  \tneg \z_\th \, \G_t e^{A_t}   \neg \left(  \big| \fF(t, Y_t,Z_t)  \neg - \neg \fF(t, 0,Z_t) \big| \neg + \neg \fF(t, 0,Z_t)  \neg  -  \neg \th \fF \big(t, \wh{Y}_t, \wh{Z}_t \big)
        \neg + \neg \k \big( Y_t  \neg - \neg  \th \wh{Y}_t \big)   \neg  -  \neg \hb{$\frac{\g}{2(1-\th)}$}  |V_t|^2   \right)    \\
       & \neg \tneg \le& \neg \tneg \z_\th \,  \G_t e^{A_t}   \neg \left(    \th \big| \fF \big(t, 0,   \wh{Z}_t \big)  \neg - \neg  \fF \big(t, \wh{Y}_t, \wh{Z}_t \big) \big|  \neg + \neg   (1 \neg - \neg \th)  \a
        \neg - \neg  \k    \th \wh{Y}_t           \right)  \le    \a   \g e^{2 \k T}  \G_t         .
   \eeas

       \ss \no 3) For \dtp ~ $(t,\o) \in \{Y_t  (\o) \vee \wh{Y}_t  (\o)  <  0   \}  $,   applying \eqref{eq:d131} with $y=\wh{Y}_t$,
        we see  from    \eqref{eq:x251} and (H2)      that
          \beas
   G_t        &  \dneg \dneg \le &\dneg  \dneg
       \z_\th \,  \G_t e^{A_t}  \left(  \fF(t, \th \wh{Y}_t,Z_t)-\th \fF(t, \wh{Y}_t, \wh{Z}_t) - \hb{$\frac{\g}{2(1-\th)}$} |V_t|^2 \right)       \\
       & \dneg  \dneg   \le& \dneg  \dneg     \z_\th  \,   \G_t e^{A_t} \neg \left(   \big| \fF(t,\th \wh{Y}_t,Z_t) \neg  - \neg  \fF(t,  \wh{Y}_t,Z_t)  \big|
        \neg  +  \neg \fF(t,  \wh{Y}_t,Z_t) \neg - \neg  \th \fF \big(t, \wh{Y}_t, \wh{Z}_t \big)
       \neg  - \neg \hb{$\frac{\g}{2(1-\th)}$} |V_t|^2 \right) \\
          & \dneg  \dneg \le   & \dneg  \dneg  \g e^{2 \k T}  \G_t      \big(\a \neg + \neg (\beta \neg + \neg \k) \big|\wh{Y}_t \big| \big)
          = \g e^{2 \k T}  \G_t      \big(\a \neg + \neg (\beta \neg + \neg \k) \wh{Y}^-_t \big)  .
   \eeas

  Now, we   define   a    process
        \bea  \label{def_D_process}
           D_t \dfnn \exp \left\{   \g e^{2 \k T}   \int_0^t        \Big(\a+  (\beta+\k) \big( Y^+_s + \wh{Y}^-_s \big) \Big)       ds  \right\} ,  \q  t \in [0,T] .
           \eea
     \if{0}
    Clearly, it holds \pas~  that
 \bea \label{bound_D}
  \frac{D_s}{D_t}
    =    \exp \left\{  \g e^{2\k T}   \int_t^s       \big(\a+ (\beta+\k) |\wh{Y}_r| \big)          dr  \right\}
       \le c_0 \,  \exp \left\{  (\beta+\k)   \g  T   e^{2\k T}       \wh{Y}_*        \right\}  ,  \q  \fa 0 \le t \le s \le T.
  \eea
  \fi
 Fix $t \in [0,T]$. Integration by parts and    \eqref{eq:d141} imply  that \pas
 \bea
    \G_{\t_n \land t}  \le    D_{\t_n \land t} \G_{\t_n \land t}   \le     D_{\t_n} \G_{\t_n}   + \z_\th  \int_{\t_n \land t}^{\t_n} D_s \G_s e^{A_s}       d K_s   - \z_\th  \int_{\t_n \land t}^{\t_n} D_s \G_s e^{A_s}  V_s dB_s.   
  \label{eq:d151}
    \eea
     Since it holds \pas~ that   $ L_s  \le  \wh{L}_s \le \wh{Y}_s$ for any $ s \in [0,T]$,
    the flat-off condition of $(Y,Z,K)$ implies that \pas
  \bea         \label{eq:axd014}
     \int_0^T \neg   D_s    \G_s    \,     d K_s    =   \neg     \int_0^T    \neg  \b1_{\{ Y_s = L_s \}}    D_s    \G_s  \,  d K_s
            \le   \neg   \int_0^T \neg \b1_{\{Y_s \le \wh{Y}_s\}}    D_s    \G_s \,  d K_s
              \le \neg   \int_0^T \neg       D_s     e^{  \g e^{2 \k T}      Y^+_s  }   d K_s               \le   \eta     K_T
 \eea
 with $\eta \dfnn  \exp \Big\{  \neg  (\beta \neg + \neg \k)   \g  T   e^{2\k T}     \wh{Y}^-_*    \neg  + \neg  \big( 1 \neg + \neg (\beta \neg + \neg \k)   T  \big) \g e^{2 \k T}        Y^+_*  \neg    \Big\}    $.      H\"older's inequality and (C2)  imply that
   \bea    \label{eq:axc015}
        &&        E\big[  \eta  (1+K_T)  \big]
        \le        \left\|  \eta \right\|_{\hL^{\frac{p}{ p-1}}(\cF_T) }  \neg  \left(1 \neg + \neg   \big\|  K_T \big\|_{\hL^p(\cF_T)} \right) \neg < \neg  \infty  ,   \\
          \hb{and}  &&     E \left[ D_T  \G_*  \right]  \le              E \left[ \exp \Big\{   \big( (\beta+\k)   \g T    e^{2\k T}    + \z_\th  e^{\k T} \big)
      \big(  Y^+_*   + \wh{Y}^-_* \big)  \Big\}   \right]   < \infty      .          \label{eq:axd131}
         \eea
                Then  the Burkholder-Davis-Gundy inequality    shows  that
  \bea
     E\left[\underset{t \in [0,T]}{\sup}  \left|   \int_0^{\t_n \land t}  D_s \G_s e^{A_s}  V_s dB_s \right|\right]
        \le    c_0 E \left[ D_T  \G_* e^{ A_*}   \left(\int_0^{\t_n}         | V_s |^2 d s\right)^{\frac{1}{2}}\right]
      \le  c_0     \sqrt{n} \,  E \left[ D_T  \G_*  \right]   
        <    \infty   .    \qq   \label{eq:d165}
 \eea
    Thus $\int_0^{\t_n \land \cd}   D_s \G_s e^{A_s}  V_s  dB_s $ is a       uniformly integrable martingale.

\ms   Taking  $E\big[\cd\big|\cF_{\t_n \land t}\big]$ in  \eqref{eq:d151}, we can deduce from   \eqref{eq:axd014}-\eqref{eq:axd131}  that \pas
  \beas
 ~\;\;  \G_{\t_n \land t}   \neg    \le   \neg      E\big[   D_{\t_n} \G_{\t_n}  \neg   \big|\cF_{\t_n \land t}\big]
   \neg +  \neg        c_0   \z_\th    E\big[ \eta K_T    \big|\cF_{\t_n \land t}\big]
       \neg   =   \neg     \b1_{\{\t_n < t\} }     D_{\t_n}  \neg  \G_{\t_n}
     \dneg +  \neg  \b1_{\{\t_n \ge t\} } E\big[   D_{\t_n}  \neg  \G_{\t_n} \neg \big|\cF_t \big]   \neg
     +  \neg    c_0   \z_\th      E\big[ \eta K_T    \big|\cF_{\t_n \land t}\big] .
  \eeas
  As $n \to \infty$, the Dominated Convergence Theorem, \eqref{eq:axd131} and \eqref{eq:axc015} imply that \pas
    \beas
  \q  \G_t   \neg   \le    \neg        E\big[   D_T \G_T  \big|\cF_t \big]    \neg  +   \neg     c_0   \z_\th     E\big[ \eta K_T    \big|\cF_t\big]
      \neg   \le   \neg    E\Big[   D_T \,   e^{     \g  e^{2 \k T}  \xi^+    }  \Big|\cF_t \Big]
   \neg  +    \neg    c_0   \z_\th       E\big[ \eta K_T    \big|\cF_t\big]
   \neg   \le   \neg    c_0  (1   \neg  \vee   \neg   \z_\th)      E\big[ \eta (1  \neg  +  \neg   K_T )   \big|\cF_t\big]  ,
 \eeas
      which leads to that
      \bea
     Y_t   -   \th \wh{Y}_t      \le  \neg \hb{$\frac{1-\th}{\g}$}  \ln \neg  \left( \neg 1 \vee  \hb{$\frac{ \g e^{\k T}}{1-\th}$}   \right) e^{-\k T-A_t}
      \neg  +\neg   \hb{$\frac{1-\th}{\g}$}     \big(c_0 \neg+   \ln   E  [      \eta (1+K_T)    |\cF_t ]   \big) e^{-\k T-A_t} , \q \pas     \label{eq:x255} 
    \eea
  Letting $\th \to 1$ yields that  $  Y_t -   \wh{Y}_t \le 0$, \pas~
Then the continuity of processes $Y$ and $\wh{Y}$ proves the theorem. \qed

   \begin{thm}   \label{thm:comparison2}
   Let $(\xi,f,L)$, $(\hat{\xi}, \hat{f},\wh{L})$ be two parameter sets  and let $(Y,Z,K) \big(\hb{resp. }  (\wh{Y}, \wh{Z},\wh{K})\,\big) 
 $ be a solution of  RBSDE$(\xi,f,L)$ \big(resp.\,RBSDE$(\hat{\xi}, \hat{f},\wh{L})$\,\big)  such that  \($C1$\), \($C2$\) hold.
For $\a,    \beta,  \k    \ge   0$, $\g   >    0$,    if  either of the following two  holds:

\ss \no (\,i)    $f$   satisfies (H1), (H2), $f$  is concave in $z$,  and $ \D f(t) \dfnn  f(t, \wh{Y}_t, \wh{Z}_t  ) - \hat{f}(t, \wh{Y}_t, \wh{Z}_t ) \le 0 $, \dtp;

\ss \no (ii)   $\hat{f}$ satisfies (H1), (H2), $\hat{f}$ is concave in $z$,   and $ \D f(t) \dfnn  f(t, Y_t, Z_t  ) - \hat{f}(t, Y_t, Z_t ) \le 0 $, \dtp;

 \ss \no then it holds \pas ~ that $Y_t \le \wh{Y}_t  $ for any $ t  \in [0,T]$.

   \end{thm}

\ss  \no  {\bf Proof:}    Fix $\th \in (0,1) $.  We set $\wt{U} \dfnn \th Y \neg - \wh{Y}$, $\wt{V} \dfnn \th Z-  \wh{Z}$
    and define  an $\bF$-progressively measurable process
      \beas
   \qq   \wt{a}_t &\dfnn& \b1_{\{Y_t \ge 0\}} \left(  \b1_{\{\wt{U}_t  \ne 0 \}}\frac{\fF \big(t, \th Y_t, \wh{Z}_t\big)
     -   \fF \big(t, \wh{Y}_t, \wh{Z}_t\big)}{\wt{U}_t        }         - \k \b1_{\{\wt{U}_t  =0   \}} \right)   - \k \b1_{\{   Y_t < 0 \le \wh{Y}_t \}}     \\
   && + \b1_{\{ Y_t \vee \wh{Y}_t <  0 \}} \left( \b1_{\{ Y_t \ne \wh{Y}_t \}}\frac{ \fF (t, Y_t,Z_t)
   -  \fF   \big(t, \wh{Y}_t,   Z_t \big) }{Y_t - \wh{Y}_t      }         - \k \b1_{\{ Y_t = \wh{Y}_t \}} \right)   ,   \q   t \in [0,T] ,
      \eeas
     where  $ \fF$ stands for $f$ if (i) holds, and for $\hat{f}$ otherwise.   It follows that $\wt{A}_t \dfnn \int_0^t \wt{a}_s \, ds$, $   t \in [0,T]$ is an $\bF$-adapted process
    with $\wt{A}_* \dfnn \underset{t \in [0,T]}{\sup} \big| \wt{A}_t \big| \le \int_0^T \neg   \big| \wt{a}_s \big|   ds \le  \k T$, \pas ~
       In light of (H1) and the concavity of $\fF$ in $z$, it holds \dtp~ that
          \bea  \label{eq:d131b}
   \fF \big(t, y, \wh{Z}_t \big)         \ge   \th  \fF  (t, y,   Z_t )  \neg + \neg   (1 \neg - \neg \th)   \fF    \bigg(t, y,   \frac{-\wt{V}_t}{1 \neg - \neg \th} \bigg) \neg
        \ge    \th  \fF  (t, y,   Z_t  )  \neg - \neg   (1 \neg - \neg \th)  \big(\a \neg + \neg  \beta |y|\big)  \neg - \neg  \hb{$\frac{\g}{2(1-\th)}$} \big| \wt{V}_t \big|^2   , \q \fa y \in \hR .
  \eea

    \ss   Given $n \in \hN$, we still define    the $\bF$-stopping time $\t_n$ as in \eqref{eq:axd071}.  Let   $ \z_\th  \dfnn \frac{ \g e^{\k T}}{1-\th}$.
 Applying It\^o's formula to the process $\wt{\G}_t \dfnn \exp\big\{\z_\th  e^{\wt{A}_t }\wt{U}_t \big\}$, $t \in [0,T]$ yields that
  \beas
  \wt{\G}_{\t_n \land t} = \wt{\G}_{\t_n} + \int_{\t_n \land t}^{\t_n}  \wt{G}_s \, ds + \z_\th  \int_{\t_n \land t}^{\t_n} \wt{\G}_s e^{\wt{A}_s}     (\th dK_s- d \wh{K}_s) - \z_\th  \int_{\t_n \land t}^{\t_n} \wt{\G}_s e^{\wt{A}_s}  \wt{V}_s dB_s, \q t \in [0,T] ,
  \eeas
  where   $   \wt{G}_t  =   \z_\th \,  \wt{\G}_t e^{\wt{A}_t}  \left( \th f(t,Y_t,Z_t)-  \hat{f}(t, \wh{Y}_t, \wh{Z}_t) -\wt{a}_t \wt{U}_t- \hb{$\frac12$} \z_\th  e^{\wt{A}_t}\big| \wt{V}_t \big|^2 \right) $. Clearly, it holds \dtp~that
  \bea   \label{eq:x251b}
    \wt{G}_t  \le   \z_\th \,  \wt{\G}_t e^{\wt{A}_t}  \left(   \th \fF(t,Y_t,Z_t)-  \fF(t, \wh{Y}_t, \wh{Z}_t) -\wt{a}_t \wt{U}_t- \hb{$\frac{\g}{2(1-\th)}$}  \big| \wt{V}_t \big|^2   \right)
  \eea
  whether (i)  or (ii) holds. Moreover, let us show by 3 cases that
      \bea    \label{eq:d141b}
   \wt{G}_t \le  \g e^{2 \k T} \, \wt{\G}_t \Big(\a + (\beta + \k) \big(Y^+_t + \wh{Y}^-_t \big)\Big) , \q \dtp
   \eea

   \ss \no 1) For \dtp ~ $(t,\o) \in \{ Y_t (\o) \ge 0 \}  $,   applying \eqref{eq:d131b} with $y=Y_t$, we can deduce from    \eqref{eq:x251b} and (H2)      that
       \beas
   \wt{G}_t       &  \dneg \dneg\le &\dneg  \dneg          \z_\th \, \wt{\G}_t e^{\wt{A}_t} \neg  \left(    \th \fF(t,Y_t,Z_t) \neg - \neg  \fF \big( t, \th Y_t, \wh{Z}_t \big)  \neg  - \neg  \hb{$\frac{\g}{2(1-\th)}$} \big|\wt{V}_t\big|^2 \right)  \\
       & \dneg \dneg  \le& \dneg \dneg    \z_\th \,  \wt{\G}_t e^{\wt{A}_t} \neg \left(   \big|\fF (t, Y_t, \wh{Z}_t)     - \neg  \fF (t, \th Y_t, \wh{Z}_t)   \big|
          + \neg (1-\th)  \big(\a+ \beta |Y_t|\big) \right)         \\
        & \dneg \dneg        \le  & \dneg \dneg   \g e^{2 \k T}  \wt{\G}_t      \big(\a  +  (\beta + \k) |Y_t| \big)
         =   \g e^{2 \k T}  \wt{\G}_t      \big(\a  +  (\beta + \k) Y^+_t \big) .
   \eeas

     \ss \no 2) For \dtp ~ $(t,\o) \in \{ Y_t (\o) < 0 \le \wh{Y}_t (\o)   \}  $,   applying \eqref{eq:d131b} with $y=0$,
     we see from    \eqref{eq:x251b} and (H2)      that
                    \beas
        \wt{G}_t     &\dneg \tneg \le& \dneg \tneg \z_\th \,  \wt{\G}_t e^{\wt{A}_t}
          \neg \left( \th \big| \fF (t, Y_t,Z_t) \neg- \neg \fF (t, 0,Z_t) \big|
        \neg+\neg \th \fF (t, 0,Z_t)  \neg - \neg \fF  \big(t, \wh{Y}_t, \wh{Z}_t \big)
       \neg+\neg \k \big( \th Y_t \neg-\neg \wh{Y}_t \big)  \neg -\neg \hb{$\frac{\g}{2(1 \neg - \neg \th)}$}  \big| \wt{V}_t \big|^2   \right)    \\
       & \tneg \le& \tneg \z_\th \,  \wt{\G}_t e^{\wt{A}_t}   \neg \left(     \big| \fF  \big(t, 0,   \wh{Z}_t \big) - \fF  \big(t, \wh{Y}_t, \wh{Z}_t \big) \big| +  (1-\th)  \a
       -\k     \wh{Y}_t              \right)  \le    \a   \g e^{2 \k T}  \wt{\G}_t         .
   \eeas

       \ss \no 3) For \dtp ~ $(t,\o) \in \{Y_t  (\o) \vee \wh{Y}_t  (\o)  <  0   \}  $,   applying \eqref{eq:d131b} with $y=\wh{Y}_t$,
        we see  from    \eqref{eq:x251b} and (H2)      that
     \beas
       \wt{G}_t                &\tneg   \le   &\tneg     \z_\th \,  \wt{\G}_t e^{\wt{A}_t} \neg \left(   \th \fF(t, Y_t,Z_t) - \th \fF   \big(t, \wh{Y}_t,   Z_t \big)
       -\wt{a}_t \wt{U}_t   +       (1 \neg - \neg \th)  \big(\a \neg + \neg  \beta \big| \wh{Y}_t \big|\big)          \right) \\
       &\tneg   =   &\tneg         \z_\th  \,  \wt{\G}_t e^{\wt{A}_t} \neg \left(   (1-\th)  \wt{a}_t   \wh{Y}_t
           +       (1 \neg - \neg \th)  \big(\a \neg + \neg  \beta \big| \wh{Y}_t \big| \big)          \right)
           \le     \g e^{2 \k T}  \wt{\G}_t      \big(\a \neg + \neg (\beta \neg + \neg \k) \big| \wh{Y}_t \big| \big)
          = \g e^{2 \k T}  \wt{\G}_t      \big(\a \neg + \neg (\beta \neg + \neg \k) \wh{Y}^-_t \big)  .
   \eeas

   Let  $D$ be the $\bF$-adapted process defined in \eqref{def_D_process}. Fix $t \in [0,T]$. Similar to \eqref{eq:d151},
             integration by parts and \eqref{eq:d141b} imply that \pas
 \bea
      \wt{\G}_{\t_n \land t}
     \le  D_{\t_n } \wt{\G}_{\t_n }     +\z_\th         \int_{\t_n \land t}^{\t_n }    D_s  \wt{\G}_s e^{\wt{A}_s}  d K_s
       - \z_\th  \int_{\t_n \land t}^{\t_n } D_s \wt{\G}_s e^{\wt{A}_s}  \wt{V}_s dB_s    .  \label{eq:axd121}
 \eea
 Similar to \eqref{eq:axd014},      the flat-off condition of $(Y,Z,K)$ implies that
  $       \int_0^T \neg    D_s  \wt{\G}_s     d K_s        \neg         \le    \neg     \wt{\eta}     K_T$,  \pas~with
       $\wt{\eta} \dfnn  \exp \Big\{  \neg  (\beta \neg + \neg \k)   \g  T   e^{2\k T}       Y^+_*   \neg  + \neg  \big( 1 \neg + \neg (\beta \neg + \neg \k)   T  \big) \g e^{2 \k T}    \wh{Y}^-_*   \neg    \Big\}$.
                H\"older's inequality   and (C2)   imply   that
    \bea      \label{eq:axc015b}
          E\big[  \wt{\eta} (1+K_T)  \big]
       \le       \left\|   \wt{\eta} \right\|_{\hL^{\frac{p}{ p-1}}(\cF_T) }  \neg  \left(1 \neg + \neg   \big\|  K_T \big\|_{\hL^p(\cF_T)} \right) \neg < \neg  \infty
         \eea
     and     $        E \left[ D_T  \wt{\G}_*  \right] \neg  \le       \neg         E \left[ \exp \Big\{   \big( (\beta \neg  + \neg  \k)   \g T    e^{2\k T}     \neg  + \neg   \z_\th  e^{\k T} \big)     \big(  Y^+_*   \neg   +  \neg  \wh{Y}^-_* \big)  \Big\}   \right]   < \infty$. Similar to  \eqref{eq:d165},
              the Burkholder-Davis-Gundy inequality  shows that
        $\int_0^{\t_n \land \cd} \neg   D_s \wt{\G}_s e^{\wt{A}_s}  \wt{V}_s  dB_s $ is a      uniformly integrable martingale.
     Taking  $E\big[\cd\big|\cF_{\t_n \land t}\big]$ in  \eqref{eq:axd121} and using the similar arguments to those that lead to  \eqref{eq:x255},
we can deduce from \eqref{eq:axc015b}   that
       \beas
      \th Y_t   -    \wh{Y}_t
    \le \neg \hb{$\frac{1-\th}{\g}$}  \ln \neg  \left( \neg 1 \vee  \hb{$\frac{   \g  e^{\k T}}{1-\th}$}   \right) e^{-\k T-\wt{A}_t} \neg  +\neg   \hb{$\frac{1-\th}{\g}$}     \big(c_0 \neg+   \ln   E  [      \wt{\eta} (1+K_T)   |\cF_t ]   \big) e^{-\k T-\wt{A}_t} , \q \pas
    \eeas
   Letting $\th \to 1$ yields that  $  Y_t -   \wh{Y}_t \le 0$, \pas~
Then the continuity of processes $Y$ and $\wh{Y}$ proves the theorem. \qed

 \ms   Using  Theorem \ref{thm:existence}, Theorem \ref{thm:comparison} and Theorem \ref{thm:comparison2}, we obtain  the following uniqueness result for quadratic RBSDEs.

\begin{cor} \label{cor:unique}
    Let $(\xi, f, L)$ be a parameter set such that $f$  satisfies  (H1)-(H3).
         If $\xi^+ \vee     L_* \in \hL^e(\cF_T)   $,   then
    the quadratic  RBSDE\,$(\xi, f, L)$ admits a unique solution    $(Y,Z,K)$ in $ \underset{p \in [1, \infty)}{\cap} \hS^p_\bF[0,T]$
    that satisfies \eqref{a_priori}.
\end{cor}

 \ss \no {\bf Proof:} The existence results from Theorem \ref{thm:existence}. Let $(\wh{Y},\wh{Z},\wh{K})$ be another solution of
  the quadratic  RBSDE\,$(\xi, f, L)$ such that $(\wh{Y},\wh{Z},\wh{K})  \in   \hS^p_\bF[0,T] $ for all $p \in [1, \infty)$. One can deduce from  Theorem \ref{thm:comparison} or Theorem \ref{thm:comparison2} that $Y $ and $\wh{Y}$ are
   indistinguishable, which implies that
        \bea
   \q   0 &\dneg=&\dneg Y_0 - Y_t  - ( \wh{Y}_0    - \wh{Y}_t  )  =  \int_0^t \big( f(s, Y_s, Z_s) -f(s, \wh{Y}_s, \wh{Z}_s) \big) \, ds
     +   K_t - \wh{K}_t  - \int_0^t  (Z_s - \wh{Z}_s) dB_s  \nonumber \\
         &\dneg=&\dneg  \int_0^t \big( f(s, Y_s, Z_s) -f(s, Y_s, \wh{Z}_s) \big) \, ds
     +   K_t - \wh{K}_t  - \int_0^t  (Z_s - \wh{Z}_s) dB_s\,, \q     t \in [0, T].   \label{eq:d175}
  \eea
   Since the set of continuous martingales and that of finite variation processes only intersect at constants, one can deduce that $Z_t = \wh{Z}_t$, \dtp~
   Putting it back into \eqref{eq:d175} shows that  $K$ and $\wh{K}$ are indistinguishable.  \qed

\section{A Uniqueness Result}

When the generator $f$ is concave in $z$, we have    the following more general uniqueness result than Corollary \ref{cor:unique}
by a Legendre-Fenchel transformation    argument, which was   used in   \cite{DHR_2010},
   \cite[Section 7]{EKPPQ-1997} and \cite[Section 4]{OS_CRM}.

\begin{thm}   \label{thm:uniqueness}
      Let $(\xi, f, L)$ be a parameter set  such that  $f$ satisfies \(H2\) and is concave in $z$. Assume that
        for three constants $\a, \beta \ge 0$ and $\g>0$, it holds \dtp~ that
\bea    \label{eq:axa201}
    f(t , \o, y,z)  \ge - \a  - \beta |y| - \frac{\g}{2} |z|^2, \q \;\;  \fa  (y,z) \in \hR \times \hR^d ,
\eea
   Then the RBSDE$(\xi, f, L)$ has at most one solution $(Y,Z,K)   \in \hE^{\l,\l'}_\bF [0,T] \times \wh{\hH}^{2,loc}_\bF([0,T];\hR^d) \times \hK_\bF[0,T]$  with $\l \in (\g,\infty)$ and   $\l' \in (0,\infty)$.
      \end{thm}

   \ss \no {\bf Proof:}
     Suppose that the RBSDE$(\xi, f, L)$ has two solutions $\big\{(Y^i, Z^i, K^i ) \big\}_{i=1,2}   \subset  \hE^{\l_i,\l'_i}_\bF [0,T] \times \wh{\hH}^{2,loc}_\bF([0,T];\hR^d) \times \hK_\bF[0,T]$  with   $\l_i \in (\g,\infty)$ and   $\l'_i \in (0,\infty)$.
     We set  $\l \dfnn \l_1 \land \l_2$    and      $\l' \dfnn  \l'_1 \land \l'_2$.

     \ms  Clearly, $ -f$ is convex in $z$.  For any $(t,\o,y) \in [0,T] \times \O \times \hR  $, it is well-known that the Legendre-Fenchel transformation of $ f(t,\o,y,\cd)$:
 \beas
 \wh{f}(t, \o, y, \fq) \dfnn \underset{z \in  \hR^d}{\sup} \big( \lan  \fq, z   \ran +  f(t,\o,y,z)\big),   \q    \fa   \fq  \in  \hR^d
 \eeas
 is   an  $\hR \cup \{\infty\}$-valued, convex and lower semicontinuous function.
   Let $\fN $ be the $dt \otimes dP$-null set   except on which
     \eqref{eq:axa011}, \eqref{eq:axa015} and \eqref{eq:axa201} hold.  Given  $(t,\o ) \in \fN^c  $, $\wh{f} $ has the following properties:
\bea \label{eq:axb015}
     (1)~ \hb{By \eqref{eq:axa201}, }   \wh{f}(t,\o,y,\fq) \ge -\a - \beta |y| +\frac{1}{2 \g} |\fq|^2, ~ \fa (y,z) \in   \hR \times \hR^d. \hspace{6.2cm}
\eea
    (2) For any $ \fq  \in  \hR^d $, if  $\wh{f}(t,\o,y, \fq) < \infty$ for some $y \in \hR$,  then (H2) implies that  for any $     y'  \in \hR$,
 \bea   \label{eq:axa019}
      \wh{f}(t,\o,y', \fq) < \infty \q \hb{and} \q
 \big| \wh{f}(t,\o,y, \fq)-\wh{f}(t,\o,y', \fq)\big| \le \k |y-y'| .
 \eea

 \no (3) For any $ y \in   \hR $,  since  $- f(t,\o,y,\cd)$ is convex on  $\hR^d$,   the conjugacy relation shows that
 \bea  \label{eq:axa021}
  - f(t, \o, y, z) =\underset{\fq \in \hR^d}{\sup} \big(  \lan z, \fq  \ran- \wh{f}(t, \o, y, \fq) \big) ,    \q    \fa   z  \in   \hR^d .
  \eea
 Moreover, the convexity of $- f(t,\o,y,\cd)$ on  $\hR^d$  implies its continuity on  $\hR^d$, thus
   \beas
  \wh{f}(t, \o, y, \fq) = \underset{z \in  \hQ^d}{\sup} \big( \lan  \fq, z   \ran +   f(t,\o,y,z)\big),   \q    \fa   \fq  \in  \hR^d ,
 \eeas
 which implies that $\wh{f}$ is $\sP \times \sB(\hR) \times   \sB(\hR^d)/\sB(\hR)$-measurable.

 \ss \no (4)  For any $(y,z) \in   \hR \times \hR^d $, let $\pa  (- f)  (t,\o,y, z)$ denote the subdifferential of the function $ -  f(t,\o,y,\cd)$
  at $z $ (see e.g. \cite{rock}).  It  is a non-empty  convex compact subset of $\fq \in \hR^d$ such that $-  f(t,\o,y, z' ) +   f(t,\o,y, z) \ge \lan \fq, z'- z \ran$ for any $z' \in \hR^d$, to wit,
   \bea  \label{eq:axa025}
     \lan \fq,   z \ran  +     f(t,\o,y, z)   = \underset{z' \in \hR^d}{\sup} \big(  \lan \fq,   z' \ran +   f(t,\o,y, z')  \big) =   \wh{f}(t, \o, y, \fq).
        \eea

    Let $i=1,2$.  For any $(t,\o ) \in \fN^c  $, we choose a $ \fq^i(t,\o) \in \pa (  - f) \big(t,\o,Y^i_t(\o) , Z^i_t(\o)\big)  $. By     \eqref{eq:axa025},
    \bea  \label{eq:axb019}
     \wh{f}\big(t, \o, Y^i_t(\o), \fq^i(t,\o)\big)  = \big\lan     Z^i_t(\o), \fq^i(t,\o)  \big\ran +    f \big(t,\o,Y^i_t(\o), Z^i_t(\o)\big)      < \infty     .
         \eea
              Thanks to the Measurable Selection Theorem (see e.g. Lemma 1 of  \cite{Benes_1970} or  Lemma 16.34 of
\cite{Elliott_1982}\big),  there exists an $\bF$-progressively measurable process $\wt{\fq}^{i}$ such that
    \bea   \label{eq:axa033}
  \qq       f \big(t,\o,Y^i_t(\o), Z^i_t(\o)\big)  =      \wh{f}\big(t, \o, Y^i_t(\o),  \wt{\fq}^{i}_t(\o) \big) - \big\lan     Z^i_t(\o), \wt{\fq}^{i}_t(\o)  \big\ran  ,
  \q  \fa  (t,\o ) \in \fN^c ,
        \eea
    which together with   \eqref{eq:axb015}  leads to that
 \bea  \label{eq:axb011}
          f \big(t,\o,Y^i_t(\o), Z^i_t(\o)\big)
    \ge -\a - \beta |Y^i_t(\o)| +  \frac{1}{2 \g} \big|\wt{\fq}^{i}_t(\o)\big|^2 - \frac12 \Big(     2 \g \big|Z^i_t(\o)\big|^2 +  \frac{1}{2 \g}  \big|\wt{\fq}^{i}_t(\o)\big|^2\Big) , ~\;  \fa  (t,\o ) \in \fN^c.
 \eea
 Since $ (Y^i, Z^i, K^i ) \in \hE^{\l_i,\l'_i}_\bF [0,T] \times \wh{\hH}^{2,loc}_\bF([0,T];\hR^d) \times \hK_\bF[0,T]$   solves  the RBSDE$(\xi, f, L)$,
 it holds \pas ~ that
 \beas
   Y^i_*  +  \int_0^T\big|Z^i_t \big|^2 dt +        \bigg| \int_0^T f \big(t, Y^i_t , Z^i_t \big) dt \bigg|<\infty .
  \eeas
   Then it follows from \eqref{eq:axb011} that
  \bea  \label{eq:axb023}
  \q  \frac{1}{4 \g} \int_0^T\big|\wt{\fq}^{i}_t \big|^2 dt \le   \neg \int_0^T      f \big(t, Y^i_t , Z^i_t\big) dt
   +\big(\a \neg + \neg   \beta  Y^i_* \big)\,T   +\g \int_0^T\big|Z^i_t\big|^2 dt < \infty     , \q \pas \q  \;\;\;
  \eea

Next, let us pick up an $N \in  \hN $ such that
\bea  \label{eq:axb030}
  \frac{T}{N} \le  \frac{ \l \l'}{2 \beta (\l+\l' )} \Big(\frac{1}{\g}-\frac{1}{\l}\Big)    .
  \eea
Let $t_0 \dfnn 0$. For $j \in \{ 1, \cds ,N\}$, we set $t_j \dfnn \frac{j T}{N}$ and  define the process
 \beas
 M^{i,j}_t \dfnn \exp \left( - \int^t_0 \b1_{\{s \ge t_{j-1}\}} \wt{\fq}^{i}_s dB_s - \frac12 \int^t_0 \b1_{\{s \ge t_{j-1}\}} \big|\wt{\fq}^{i}_s\big|^2 ds  \right)  ,
 \q \fa t \in [0, t_j] .
 \eeas

  Given $n \in \hN$, we define the $\bF$-stopping time
  \beas
    \t^j_n \dfnn \inf\bigg\{ t \in [t_{j-1},t_j]: \int_{t_{j-1}}^t  \big(   |Z^1_s|^2 + |Z^2_s|^2+ |\wt{\fq}^{1}_s|^2 + |\wt{\fq}^{2}_s|^2 \big) ds >n  \bigg\} \land t_j.
    \eeas
    Clearly, $\lmtu{n \to \infty} \t^j_n = t_j$, \pas ~ by \eqref{eq:axb023}, and  $ \Big\{  M^{i,j}_{\t^j_n \land t} \Big\}_{t \in [0, t_j]} $ is a uniformly integrable martingale thanks to Novikov's Criterion.
    Hence, $ \dis \frac{d  Q^{i,j}_n}{dP} \dfnn  M^{i,j}_{\neg \t^j_n}$ induces a probability $Q^{i,j}_n$ 
    that is equivalent to $P$.
       Girsanov Theorem shows that $ \big\{B^{i,j,n}_t\dfnn B_t + \int_0^t \b1_{\{ t_{j-1} \le s \le \t^j_n \}} \wt{\fq}^{i}_s ds
    \big\}_{t \in [0,t_j]}$ is a Brownian Motion under  $Q^{i,j}_n$ and
   \bea
 E  \left[    M^{i,j}_{\t^j_n} \,  \ln M^{i,j}_{\t^j_n} \right]  &=& E_{Q^{i,j}_n} \neg  \left[       \ln M^{i,j}_{\t^j_n} \right]
   =  E_{Q^{i,j}_n}  \neg  \left[     - \int^{\t^j_n}_{t_{j-1}} \wt{\fq}^{i}_s dB_s -
\frac12 \int^{\t^j_n}_{t_{j-1}} \big|\wt{\fq}^{i}_s\big|^2 ds    \right]  \nonumber \\
  &=&      E_{Q^{i,j}_n}  \neg   \left[      - \int^{\t^j_n}_{t_{j-1}} \wt{\fq}^{i}_s dB^{i,j,n}_s
  +  \frac12  \int^{\t^j_n}_{t_{j-1}}     |\wt{\fq}^{i}_s|^2 ds       \right]
     =    \frac12    E_{Q^{i,j}_n}   \neg  \left[       \int^{\t^j_n}_{t_{j-1}}     |\wt{\fq}^{i}_s|^2 ds       \right] .   \label{eq:axb027}
  \eea

 It is well-known that  for any $(x,\mu) \in \hR \times (0,\infty)$
 \bea  \label{eq:axb042}
       x \mu \le e^x + \mu (\ln \mu - 1) \le e^x + \mu \ln \mu ,    \q  \hb{thus} \q
         x\mu = \l x \frac{\mu}{\l} \le e^{\l x} + \frac{\mu}{\l} \big(\ln \mu - \ln \l \big) ,
 \eea
 which together with \eqref{eq:axb027} implies that for $k=1,2$
  \bea
   E_{Q^{i,j}_n} \left[ \, \underset{t \in [0, t_j]}{\sup} \big( Y^k_t\big)^-\right] &=& E  \left[ \, \underset{t \in [0, t_j]}{\sup} \big( Y^k_t\big)^- M^{i,j}_{\t^j_n}\right]
  \le  E  \left[  e^{\l (Y^k)^-_*} \right]+ \frac{1}{\l}  E_{Q^{i,j}_n}  \Big[      \ln M^{i,j}_{\t^j_n} - \ln \l  \,    \Big] \nonumber \\
  &\le &  \wt{c}^{\,k}_{\l}  +  \frac{1}{2\l}  E_{Q^{i,j}_n} \neg  \left[       \int^{\t^j_n}_{t_{j-1}}     |\wt{\fq}^{i}_s|^2 ds       \right] . \label{eq:axb036}
  \eea
 where $\wt{c}^{\,k}_\l \dfnn    E  \left[  e^{\l (Y^k)^-_*} \right] +\frac{( \ln \l )^-}{\l}$. Similarly,
  \bea    \label{eq:axb039}
  && E_{Q^{i,j}_n} \left[ \, \underset{t \in [0, t_j]}{\sup} \big( Y^k_t\big)^+\right]
   \le     \wt{c}^{\,k}_{\l'}  +  \frac{1}{2\l'}  E_{Q^{i,j}_n} \neg  \left[       \int^{\t^j_n}_{t_{j-1}}     |\wt{\fq}^{i}_s|^2 ds       \right].
  \eea

 We can deduce from  \eqref{eq:axa033}, \eqref{eq:axb015} and Girsanov Theorem that
 \bea
 Y^i_{t_{j-1}} \neg - \neg  Y^i_{\t^j_n}  & \dneg \tneg =&\tneg \dneg      \int_{t_{j-1}}^{\t^j_n}  \neg f(s, Y^i_s, Z^i_s) ds    + \neg    K^i_{\t^j_n} \neg- \neg K^i_{t_{j-1}}  \dneg -  \neg  \int_{t_{j-1}}^{\t^j_n}  \neg  Z^i_s  d B_s          \ge        \int_{t_{j-1}}^{\t^j_n}   \dneg \Big(\wh{f}\big(s,   Y^i_s ,  \wt{\fq}^{i}_s  \big)  \neg - \neg   \big\lan     Z^i_s , \wt{\fq}^{i}_s   \big\ran  \Big)     ds
         \neg -  \neg  \int_{t_{j-1}}^{\t^j_n}  \neg  Z^i_s  d B_s  ~   \nonumber \\
   &\dneg \tneg \ge &\tneg \dneg     \int_{t_{j-1}}^{\t^j_n}  \neg   \Big( - \neg  \a  \neg  -  \neg  \beta \big|Y^i_s  \big|
    \neg + \neg  \frac{1}{2 \g}   \big|  \wt{\fq}^{i}_s  \big|^2    \Big)     ds      - \neg   \int_{t_{j-1}}^{\t^j_n}  \neg  Z^i_s  d B^{i,j,n}_s , \q \pas \label{eq:axb033}
 \eea
 By Bayes' rule (see e.g., \cite[Lemma 3.5.3]{Kara_Shr_BMSC}),  $E_{Q^{i,j}_n}  [Y^i_{t_{j-1}}] =  E  [ Y^i_{t_{j-1}} M^{i,j}_{t_{j-1}}] = E  [ Y^i_{t_{j-1}}]   $.
 Then taking $E_{Q^{i,j}_n}$ in \eqref{eq:axb033}, one can deduce from \eqref{eq:axb036} and \eqref{eq:axb039}  that
  \beas
          \frac{1}{2 \g}  E_{Q^{i,j}_n} \bigg[ \int_{t_{j-1}}^{\t^j_n}  \big|  \wt{\fq}^{i}_s  \big|^2          ds \bigg]
      & \tneg \le& \tneg     E  \big[ Y^i_{t_{j-1}} \big] + E_{Q^{i,j}_n}  \Big[ \Big( Y^i_{\t^j_n} \Big)^{\neg -} \,  \Big]+\frac{\a T}{N}   +\frac{\beta T}{N}    E_{Q^{i,j}_n}  \bigg[ \, \underset{t \in [0, t_j]}{\sup} \big| Y^i_t\big|  \,  \bigg]     \\
    & \tneg \le& \tneg  E   \Big[  \big(Y^i \big)^+_* \Big] +  \frac{\a T}{N}   + \Big(1+\frac{\beta T}{N} \Big)   E_{Q^{i,j}_n}  \bigg[\, \underset{t \in [0, t_j]}{\sup} \big( Y^i_t\big)^-\bigg] +
     \frac{\beta T}{N}   E_{Q^{i,j}_n}  \bigg[  \,  \underset{t \in [0, t_j]}{\sup} \big( Y^i_t\big)^+ \bigg] \\
    & \tneg \le& \tneg    \Xi +    \bigg(\frac{1}{2\l} +\frac{\beta T}{N}\Big(\frac{1}{2\l}+\frac{1}{2\l'}  \Big)\bigg)    E_{Q^{i,j}_n}  \neg \left[       \int^{\t^j_n}_{t_{j-1}}     |\wt{\fq}^{i}_s|^2 ds       \right]  ,
  \eeas
  where $\Xi \dfnn \frac{1}{\l'} E  \Big[ e^{\l'  (Y^i  )^+_*} \Big]  +  \frac{\a T}{N}    +  (1+ \beta T ) \wt{c}^{\,i}_\l +    \beta T  \wt{c}^{\,i}_{\l'}$\,.
  It follows from   \eqref{eq:axb027} and \eqref{eq:axb030}    that
 $$
 ~   \frac{1}{2} \Big(  \frac{1}{ \g} - \frac{1}{\l} \Big)  E  \left[    M^{i,j}_{\t^j_n}  \,  \ln M^{i,j}_{\t^j_n} \right] = \frac{1}{4} \Big(  \frac{1}{ \g} - \frac{1}{\l}\Big)    E_{Q^{i,j}_n} \neg \left[\int_{t_{j-1}}^{\t^j_n}  \big|  \wt{\fq}^{i}_s  \big|^2          ds   \right]
 \le  \Xi   .
  $$
 In light of    de la Vall\'ee-Poussin's lemma,
    $\big\{ M^{i,j}_{\t^j_n}\big\}_{n \in \hN}$ is uniformly integrable. Hence, $ E\Big[M^{i,j}_{t_j} \Big] \neg = \neg \lmt{n \to \infty} E\Big[M^{i,j}_{\t^j_n} \Big] \neg =1$, which shows that $M^{i,j}$ is a martingale. Thus $ \dis \frac{d  Q^{i,j} }{dP} \dfnn  M^{i,j}_{t_j}$ induces a probability $Q^{i,j} $
    that is equivalent to $P$,  and $ \big\{B^{i,j}_t\dfnn B_t + \int_0^t \b1_{\{s \ge  t_{j-1}    \}} \wt{\fq}^{i}_s ds    \big\}_{t \in [0,t_j]}$
    is a Brownian Motion under   $Q^{i,j}$.
    Then   Fatou's lemma implies that
      \beas
         ~     E_{Q^{i,j}}    \neg    \left[       \int^{t_j}_{t_{j-1}}  \neg    |\wt{\fq}^{i}_s|^2 ds       \right]  \neg =   \neg       E  \neg    \left[    M^{i,j}_{t_j}    \int^{t_j}_{t_{j-1}}  \neg   |\wt{\fq}^{i}_s|^2 ds       \right]
        \neg   \le   \neg       \linf{n \to \infty}     E   \neg   \left[    M^{i,j}_{\t^j_n}    \int^{\t^j_n}_{t_{j-1}}  \neg      |\wt{\fq}^{i}_s|^2 ds       \right]
           \neg   =    \neg        \linf{n \to \infty}     E_{Q^{i,j}_n}  \neg  \left[         \int^{\t^j_n}_{t_{j-1}}    \neg    |\wt{\fq}^{i}_s|^2 ds       \right]     \neg  \le  \neg
      \frac{4 \l \g \Xi}{\l -\g }.
      \eeas
    And  an analogy to   \eqref{eq:axb027} shows that
    \bea   \label{eq:axb045}
    ~    E_{Q^{i,j}}  \neg  \left[     \ln M^{i,j}_{t_j} \right]        =         \frac12        E_{Q^{i,j}}    \neg    \left[       \int^{t_j}_{t_{j-1}}     |\wt{\fq}^{i}_s|^2 ds       \right]
           \le      \frac{2 \l \g \Xi}{\l -\g }   \,   .
      \eea
      \if{0}
    Then similar to \eqref{eq:axb036} and \eqref{eq:axb039}, it holds for $k=1,2$ that
         \beas
   E_{Q^{i,j}} \left[ \underset{t \in [t_{j-1}, t_j]}{\sup} \big( Y^k_t\big)^-\right]
   \le    \wt{c}^k_{\l}  +  \frac{1}{2\l}  E_{Q^{i,j}}  \left[       \int^{t_j}_{t_{j-1}}     |\wt{\fq}^{i}_s|^2 ds       \right] , \q
       E_{Q^{i,j}} \left[ \underset{t \in [t_{j-1}, t_j]}{\sup} \big( Y^k_t\big)^+\right]
    \le      \wt{c}^k_{\l'}  +  \frac{1}{2\l'}  E_{Q^{i,j}}  \left[       \int^{t_j}_{t_{j-1}}     |\wt{\fq}^{i}_s|^2 ds       \right]
  \eeas
   \fi

  Now for any $n \in \hN$, applying  Tanaka's formula to  the process $\big( Y^1 -  Y^2  \big)^+$, we can deduce from \eqref{eq:axa021},   \eqref{eq:axa033}, the flat-off condition of $(Y^1,Z^1,K^1)$, \eqref{eq:axb019}, \eqref{eq:axa019}   as well as Girsanov Theorem   that
    \bea
         \Big( Y^1_{\t^n_j \land t}  \neg  -  \neg  Y^2_{\t^n_j \land t}  \Big)^{\neg +} \dneg  &  \tneg\neg = &   \tneg\dneg  \Big( Y^1_{\t^n_j}  \neg - \neg   Y^2_{\t^n_j}  \Big)^{\neg +}   + \int_{\t^n_j \land t}^{\t^n_j}  \neg \b1_{\{ Y^1_s   > Y^2_s   \}} \big(  f (s, Y^1_s, Z^1_s )  - f (s, Y^2_s, Z^2_s ) \big) ds     \label{eq:axb051} \\
    &  \tneg \neg  &  \tneg\dneg  + \int_{\t^n_j \land t}^{\t^n_j} \neg \b1_{\{ Y^1_s   > Y^2_s \}}      (  d K^1_s- d K^2_s ) - \int_{\t^n_j \land t}^{\t^n_j} \neg \b1_{\{ Y^1_s   > Y^2_s \}}    ( Z^1_s-Z^2_s)      \,  dB_s- \frac 12 \int_{\t^n_j \land t}^{\t^n_j} d  \fL_s   \nonumber \\
     &  \neg\tneg \le &  \dneg\tneg     \Big( Y^1_{\t^n_j}  \neg -  \neg  Y^2_{\t^n_j}  \Big)^{\neg +} +  \int_{\t^n_j \land t}^{\t^n_j} \neg \b1_{\{ Y^1_s   > Y^2_s   \}} \Big( \wh{f}\big(s, Y^1_s,  \wt{\fq}^2_s \big) - \big\lan     Z^1_s, \wt{\fq}^2_s  \big\ran
     - \wh{f}\big(s, Y^2_s,  \wt{\fq}^2_s \big) +  \big\lan     Z^2_s , \wt{\fq}^2_s  \big\ran   \Big)   ds \nonumber \\
    &  \neg\tneg &  \dneg\tneg   + \int_{\t^n_j \land t}^{\t^n_j} \neg \b1_{\{   L_s =Y^1_s > Y^2_s \}}   d     K^1_s   - \int_{\t^n_j \land t}^{\t^n_j} \neg \b1_{\{ Y^1_s   > Y^2_s \}}   ( Z^1_s-Z^2_s)      \,  dB_s     \nonumber   \\
       &  \neg\tneg  \le & \dneg\tneg     \Big( Y^1_{\t^n_j}  \neg - \neg   Y^2_{\t^n_j}  \Big)^{\neg +}  \dneg + \neg   \k  \neg  \int_{\t^n_j \land t}^{\t^n_j} \neg \b1_{\{ Y^1_s   > Y^2_s   \}}    \big( Y^1_s  \neg -  \neg  Y^2_s  \big)^{\neg +}  \neg \,    ds         - \neg  \int_{\t^n_j \land t}^{\t^n_j} \neg \b1_{\{ Y^1_s   > Y^2_s \}}   ( Z^1_s  \neg - \neg  Z^2_s)      \,  dB^{2,j}_s
        , ~     t \in [t_{j-1},t_j],  \nonumber
\eea
 where $\fL$ is a real-valued, $\bF$-adapted, increasing and continuous  process known as ``\,local time".
  Taking the expectation $E_{Q^{2,j}}$ and using Fubini Theorem, we obtain
  \beas
     \q     E_{Q^{2,j}}\bigg[ \Big( Y^1_{\t^n_j \land t} -  Y^2_{\t^n_j \land t}  \Big)^+   \bigg] & \le&    E_{Q^{2,j}} \bigg[ \Big( Y^1_{\t^n_j} -  Y^2_{\t^n_j}  \Big)^+\bigg]  +\k E_{Q^{2,j}}\bigg[ \int_t^{t_j} \neg \b1_{\{s \le \t^n_j\}}  \b1_{\{ Y^1_s   > Y^2_s \}}  \Big( Y^1_s -  Y^2_s  \Big)^+   ds  \bigg]  \qq \\
     & \le&    E_{Q^{2,j}} \bigg[ \Big( Y^1_{\t^n_j} -  Y^2_{\t^n_j}  \Big)^+\bigg]  +\k \int_t^{t_j} \neg   E_{Q^{2,j}}\bigg[ \Big( Y^1_s -  Y^2_s  \Big)^+   \bigg]     ds   , \q     t \in [t_{j-1},t_j].
  \eeas
 Then an application of Gronwall's inequality yields that
  \bea   \label{eq:axb048}
                   E_{Q^{2,j}}\bigg[ \Big( Y^1_{\t^n_j \land t} -  Y^2_{\t^n_j \land t}  \Big)^+   \bigg]
       \le      e^{\k T }  E_{Q^{2,j}} \bigg[ \Big( Y^1_{\t^n_j} -  Y^2_{\t^n_j}  \Big)^+\bigg] , \q     t \in [t_{j-1},t_j].
       \eea
Similar to   \eqref{eq:axb036} and \eqref{eq:axb039}, one can deduce from \eqref{eq:axb042} and \eqref{eq:axb045} that
      $$
          E_{Q^{2,j}} \neg \bigg[ \, \underset{t \in [0,t_j]}{\sup}   \Big( Y^1_t -  Y^2_t  \Big)^{\neg +}   \bigg]   \neg  \le  \neg
          E_{Q^{2,j}}  \neg \bigg[ \, \underset{t \in [0,t_j]}{\sup}  \big( Y^1_t\big)^{\neg +}   + \underset{t \in [0,t_j]}{\sup} \big( Y^2_t\big)^{\neg - }   \bigg]
             \neg  \le      \wt{c}^{\,1}_{\l'} + \wt{c}^{\,2}_{\l}
              +  \Big(\frac{1}{2\l'}+   \frac{1}{2\l} \Big)  E_{Q^{2,j}}   \neg  \left[       \int^{t_j}_{t_{j-1}}     |\wt{\fq}^2_s|^2 ds       \right]  \neg  <  \neg   \infty.
      $$

  If $Y^1_{t_j} \le Y^2_{t_j}$, \pas,  as $n \to \infty$ in \eqref{eq:axb048}, Dominated convergence theorem implies that for any
   $t \in [t_{j-1},t_j]$
    \bea   \label{eq:axb055}
     E_{Q^{2,j}}\Big[ \big( Y^1_t -  Y^2_t  \big)^+   \Big]  =0,   \q \hb{thus} \q     Y^1_t        \le Y^2_t,     \q \pas
  \eea
  In particular, $Y^1_{t_{j-1}}        \le Y^2_{t_{j-1}} $,      \pas  ~ On the other hand, if $Y^2_{t_j} \le Y^1_{t_j}$, \pas,
   interchanging   $(Y^1,Z^1,Z^1)$ with $(Y^2,Z^2,Z^2)$ and estimating under $Q^{1,j} $  in the above arguments \big(from \eqref{eq:axb051} to \eqref{eq:axb055}\big) give that  for any    $t \in [t_{j-1},t_j]$, $Y^2_t \le Y^1_t$, \pas ~
   Therefore, starting from $ Y^1_T = Y^2_T = \xi$, \pas,    we can use backward induction to conclude that for any    $t \in [0,T]$,
   $Y^1_t = Y^2_t$, \pas  ~
  Then the continuity of processes $Y^1$ and $Y^2$   shows that  $Y^1$ and $Y^2$ are indistinguishable.
   Similar to the proof of Corollary \ref{cor:unique}, it follows that $Z^1_t = Z^2_t$, \dtp~
   as well as that  $K^1$ and $K^2$ are indistinguishable.  \qed

\begin{rem}
This uniqueness result via Legendre-Fenchel transformation may not work for  the convex-generator  case:
     In fact,    if the generator $f$ is convex in the z-variable \big(while requiring that it holds \dtp~ that $ f(t,y,z) \le \a +\beta |y| +\frac{\g}{2} |z|^2$,
     $\fa (y,z) \in \hR \times \hR^d$\big),     we have to alternatively define
      \beas
       \wh{f}(t, \o, y, \fq) \neg \dfnn  \neg  \underset{z \in  \hR^d}{\sup} \big(  \lan  \fq, z   \ran   - \neg   f(t,\o,y,z)\big),
       \q \fa (t, \o, y, \fq) \in [0,T] \times \O \times \hR  \times \hR^d.
       \eeas
        Correspondingly, equation \eqref{eq:axb033} becomes
  \beas
 Y^i_{t_{j-1}} \neg - \neg  Y^i_{\t^j_n}  & \dneg \tneg =&\tneg \dneg
      \int_{t_{j-1}}^{\t^j_n}  \neg f(s, Y^i_s, Z^i_s) ds    + \neg    K^i_{\t^j_n} \neg- \neg K^i_{t_{j-1}}  \neg -  \neg  \int_{t_{j-1}}^{\t^j_n}  \neg  Z^i_s  d B_s  \\
        & \dneg \tneg =& \dneg \tneg  \int_{t_{j-1}}^{\t^j_n}  \neg \Big( \big\lan     Z^i_s , \wt{\fq}^{i}_s   \big\ran \neg -\neg \wh{f}\big(s,   Y^i_s ,  \wt{\fq}^{i}_s  \big)       \Big) ds    + \neg    K^i_{\t^j_n} \neg- \neg K^i_{t_{j-1}}  \neg -  \neg  \int_{t_{j-1}}^{\t^j_n}  \neg  Z^i_s  d B_s          \\
   &\dneg \tneg \le &\tneg \dneg     \int_{t_{j-1}}^{\t^j_n}  \neg   \Big(     \a  \neg  +  \neg  \beta \big|Y^i_s  \big|
    \neg - \neg  \frac{1}{2 \g}   \big|  \wt{\fq}^{i}_s  \big|^2    \Big)     ds    + \neg    K^i_{\t^j_n} \neg- \neg K^i_{t_{j-1}}   - \neg   \int_{t_{j-1}}^{\t^j_n}  \neg  Z^i_s  d \wt{B}^{i,j,n}_s , \q \pas ,
 \eeas
 where $ \big\{\wt{B}^{i,j,n}_t  \neg \dfnn  \neg   B_t    -    \int_0^t    \b1_{\{ t_{j-1} \le s \le \t^j_n \}} \wt{\fq}^{i}_s ds
    \big\}_{t \in [0,t_j]}$ is a Brownian Motion under the probability $\wt{Q}^{i,j}_n$,  which is induced by
      $  \frac{d  \wt{Q}^{i,j}_n}{dP}   \dfnn         \exp  \neg  \Big(  \neg \int^{\t^j_n}_0 \b1_{\{s \ge t_{j-1}\}} \wt{\fq}^{i}_s dB_s
      - \frac12 \int^{\t^j_n}_0 \b1_{\{s \ge t_{j-1}\}} \big|\wt{\fq}^{i}_s\big|^2 ds  \Big)  $. Hence,
      $E_{\wt{Q}^{i,j}_n} \Big[ \int_{t_{j-1}}^{\t^j_n}  \big|  \wt{\fq}^{i}_s  \big|^2          ds \Big]$ also depends on $E_{\wt{Q}^{i,j}_n} \Big[ K^i_{\t^j_n} \neg- \neg K^i_{t_{j-1}} \Big]$, which in turn depends on  $E_{\wt{Q}^{i,j}_n}  \Big[\, \underset{s \in [t_{j-1}, t_j]}{\sup} \big| Y^i_s\big| \Big]$ and $E_{\wt{Q}^{i,j}_n} \Big[ \int_{t_{j-1}}^{\t^j_n} |Z^i_s|^2 ds\Big]$ due to the the structure of the quadratic RBSDE.
       If we  estimate $E_{\wt{Q}^{i,j}_n} \Big[ \int_{t_{j-1}}^{\t^j_n} |Z^i_s|^2 ds\Big]$ similar to \eqref{eq:axb036},
        then  $      E\Big[\exp\big\{\wt{\l}\int_{t_{j-1}}^{\t^j_n} |Z^i_s|^2 ds \big\} \Big]$      is supposed to be uniformly bounded in $n \in \hN$    for some $\wt{\l} > 0$. To wit,    $E\Big[\exp\big\{\wt{\l}\int_{t_{j-1}}^{t_j} |Z^i_s|^2 ds \big\} \Big] < \infty$.
   However,   $\int_{t_{j-1}}^{t_j} |Z^i_s|^2 ds$  may not even in  $\hL^p (\cF_{t_j})$ for   $p \ge \frac{\l \l'}{\l + \l'}$  according to   Proposition \ref{prop:Z_K_bound}.
\end{rem}

\section{An Optimal Stopping Problem for Quadratic $g$-Evaluations} \label{sec:g-evaluation}

\ms In this section, we will solve an optimal stopping problem in which the objective of the stopper is to determine an
optimal stopping time  $\t_*$ that satisfies
  \bea    \label{eq:op-gE}
  \underset{ \t  \in  \cS_{0,T}}{\sup} \cE^g_{0, \t}\big[\cR_\t\big]  =\cE^g_{0, \t_*}\big[\cR_{\t_* }\big] ,
  \eea
 where  $\cE^g$ is a ``quadratic $g$-evaluation" (a type of non-linear expectation to be defined below), and $\cR$ is a reward process that we will specify shortly.

 \ms Let  $g:  [0,T] \times \O \times \hR \times \hR^d \rightarrow  \hR $ be a $\sP \times \sB(\hR) \times   \sB(\hR^d)/\sB(\hR)$-measurable function that satisfies (H1)-(H3).  For any  $\t \in \cS_{0,T}$,
It is clear that $ g_\t(t,\o, y,z) \dfnn  \b1_{\{t < \t\}} g(t,\o,y,z)$, $(t,\o, y,z) \in  [0,T] \times \O \times \hR \times \hR^d$
 is also a $\sP \times \sB(\hR) \times   \sB(\hR^d)/\sB(\hR)$-measurable function that satisfies (H1)-(H3). Thus, we know from Corollary 6 of \cite{BH-07}
    that  for any $  \xi \in \hL^{\neg e}(\cF_T)$, the following quadratic BSDE
    \bea \label{BSDEtau}
Y_t=\xi+\int_t^T \b1_{\{s < \t
\}}g(s,Y_s,Z_s)ds-\int_t^T Z_s  dB_s,\qq  t \in [0,T]
 \eea
 admits a unique solution $(Y^{\t, \xi}, Z^{\t, \xi}) $\, in $ \underset{p \in (1, \infty)}{\cap}
 \hE^p_\bF[0,T] \times  \hH^{2,2p}_\bF([0,T];\hR^d) $. Moreover, if $ \xi \in \hL^{\neg e}(\cF_\t) $,
       one can deduce that
     \bea \label{eq:xx25}
         P \Big( Y^{\t, \xi}_t=Y^{\t, \xi}_{\t \land t}, \,   \fa t \in [0,T] \Big) =1 \q \hb{and} \q  Z^{\t, \xi}_t=\b1_{\{t < \t \}}Z^{\t, \xi}_t , \q      \dtp~
    \eea

\begin{deff}    \label{def:qg_evalu}
A ``quadratic $g$-evaluation" with domain $\hL^e(\cF_T) $ is a family of operators $\big\{ \cE^g_{\nu,\t}:
\hL^e(\cF_\t) \mapsto \hL^e(\cF_\nu) \big\}_{\nu \in \cS_{0,T},    \t \in \cS_{\nu,T}}$ such that  $ \cE^g_{\nu,\t}[\xi]
\dfnn Y^{\t, \xi}_\nu $, $    \fa \xi \in \hL^e(\cF_\t)   $.
            In particular, for any $\xi \in \hL^e(\cF_T)   $,    we can
define the ``quadratic $g$-expectation" of $\xi$ at a stopping time   $\nu \in
\cS_{0,T}$ by $\cE^g[\xi|\cF_\nu] \dfnn \cE^g_{\nu,T}[\xi]$.
 \end{deff}

  \ss    The $g$-evaluation  was introduced by \cite{Pln} for Lipschitz generators over $\hL^2(\cF_T)$.
  Then \cite{MaYao_2010}  extended the notion for
 quadratic generators, however, on $\hL^\infty(\cF_T) $.
  Thanks to Theorem 5 of \cite{BH-07} and the uniqueness of the solution $(Y^{\t, \xi}, Z^{\t, \xi}) $,
 the quadratic $g$-evaluation $\cE^g_{\nu,\t}$ introduced in
  Definition \ref{def:qg_evalu} has the following properties:   

\ss \no (1) {\it ``Monotonicity":}~~For any $\xi,\eta \in \hL^e(\cF_\t)$ with $\xi \geq \eta$,
$\pas$, we have $\cE^g_{\nu,\t}[\xi] \geq \cE^g_{\nu,\t}[\eta]$, $\pas$;

\ms \no (2) {\it ``Time-Consistency":}~~For any $\nu_1, \nu_2, \t \in \cS_{0,T}$ with $\nu_1 \le \nu_2 \le \t$, \pas,
   and for any $\xi \in \hL^e(\cF_\t)$, we have
 $    \cE^g_{\nu_1,\nu_2}\big[\cE^g_{\nu_2,\t}[\xi]\big]=\cE^g_{\nu_1,\t}[\xi]$,   \pas;

\ms \no (3) {\it ``Constant-Preserving":} ~~If it holds \dtp~ that $g(t,y,0)=0$,   $\fa  y \in \hR$,
   then  for any $\xi \in      \hL^e(\cF_\nu)$, we have      $\cE^g_{\nu,\t}[\xi]=\xi$, \pas;

\ms \no (4) {\it ``Zero-one Law":}~~For any $\xi \in \hL^e(\cF_\t)$
and $A \in \cF_\nu$, we have $\b1_A\cE^g_{\nu,\t}[\b1_A \xi]=\b1_A
\cE^g_{\nu,\t}[\xi]$, $\pas$. Moreover, if $g(t,0,0)=0$, $\dtp$,
then $\cE^g_{\nu,\t}[\b1_A \xi]=\b1_A \cE^g_{\nu,\t}[\xi]$, $\pas$;

\ms \no (5) {\it ``Translation Invariance":}~~If $g$ is independent of $y$, then for any $ \xi \in \hL^e(\cF_\t)$ and  $\eta \in \hL^e(\cF_\nu)$,
 we have $   \cE^g_{\nu,\t}[\xi+\eta]=\cE^g_{\nu,\t}[\xi]+\eta$,  \pas

 \ms   Now, we assume that the reward process $\cR$ is in the form of
  \bea   \label{def:reward}
  \cR_t \dfnn \b1_{\{ t < T\}} \cL_t + \b1_{\{ t = T\}} \xi , \q   t \in [0,T],
  \eea
   for some $ \cL   \in   \hC^0_\bF[0,T]$ and $\xi \in \hL^0(\cF_T) $ with
  $\cL_T \le \xi$, \pas~ One can regard $\cL$ as the running reward and $\xi$ as the final reward with a possible bonus.

  \ms   When      $\xi^+ \vee     \cL_* \in \hL^e(\cF_T)   $,
    the quadratic  RBSDE\,$(\xi, g, \cL)$ admits a unique solution    $(\cY,\cZ,\cK)$ in $ \underset{p \in [1, \infty)}{\cap} \hS^p_\bF[0,T]$ thanks to Corollary \ref{cor:unique}.  In fact, the continuous process $\cY$ is  the snell envelope of the reward process $\cR$ under the quadratic $g$-evaluation, and the first time process  $\cY$ meets process $\cR$     after time $t = 0$ is an   optimal
    stopping time for \eqref{eq:op-gE}. More precisely, we have the following result.

\begin{thm} \label{thm:op-gE}
Let  $g:  [0,T] \times \O \times \hR \times \hR^d \rightarrow  \hR $ be a $\sP \times \sB(\hR) \times   \sB(\hR^d)/\sB(\hR)$-measurable function that satisfies (H1)-(H3), and let $\cR$ be a reward process in the form of \eqref{def:reward}. If $\xi^+ \vee     \cL_* \in \hL^e(\cF_T)   $, then for any $\nu \in \cS_{0,T}$,
\beas
     \cY_\nu =   \esup{  \t \in  \cS_{\nu,T}}    \cE^g_{\nu, \t}\big[\cR_\t\big]
  =\cE^g_{\nu, \t_*(\nu)}\big[\cR_{\t_*(\nu) }\big]        , \q \pas,
     \eeas
  where $\cY$ is the first component  of the unique solution to the quadratic  RBSDE\,$(\xi, g, \cL)$ and $\t_*(\nu) \dfnn     \inf\big\{t \in [\nu, T]
:\, \cY_t=\cR_t\big\} \in \cS_{\nu,T}   $.
  \end{thm}

 \ss This theorem     extends Section 3 of \cite{Morlais_2008}, it also extends Theorem 5.3 of \cite{OSNE2} except that the continuity condition on the reward process $\cR$ is  strengthened.
 The proof of  Theorem \ref{thm:op-gE} depends on the following two comparison theorems for quadratic BSDEs,
 which  generalize Theorem 5 of \cite{BH-07}.

 \begin{prop}  \label{prop:comparion_ext}
Let $f, \hat{f}:    [0,T] \times \O \times \hR \times \hR^d \rightarrow  \hR $ be two
 $\sP \times \sB(\hR) \times   \sB(\hR^d)/\sB(\hR)$-measurable functions, and let $(Y,Z,V) , (\wh{Y},\wh{Z},\wh{V})\in
 \hC^0_\bF[0,T] \times  \wh{\hH}^{2,loc}_\bF([0,T];\hR^d)    \times \hV_\bF[0,T]$ solve  the following two BSDEs
  \bea
    &&  Y_t= Y_T  + \int_t^T f(s, Y_s, Z_s) \, ds +   V_T - V_t - \int_t^T  Z_s dB_s\,, \q \;\;\;    t \in [0, T]     \label{VBSDE} \\
 \hb{and}  &&     \wh{Y}_t= \wh{Y}_T  + \int_t^T \hat{f} (s, \wh{Y}_s, \wh{Z}_s) \, ds +   \wh{V}_T - \wh{V}_t - \int_t^T  \wh{Z}_s dB_s\,, \q \;\;\;    t \in [0, T]
 \label{VBSDE2}
  \eea
  respectively such that    $ Y_T \le \wh{Y}_T   $, \pas, that
  \bea   \label{eq:axf123}
     E \left[  e^{\l  Y^+_*} + e^{\l \wh{Y}^-_*}        \right] < \infty, \q \fa    \l   \in (1, \infty) ,
     \eea
   and\;that\;for\;some\;$\th_0 \neg \in  \neg (0, 1)$,\;$V  \neg - \neg  \th \wh{V}$\;is\;a\;decreasing\;process\;for\;any\;$\th  \neg \in \neg  (\th_0,1)$.\;If either\;of\;the\;following\;two\;holds:

 \ss \no \(\,i\)  $f$ satisfies (H1'),   (H2);  $f$ is convex in $z$; and  $\D f(t) \dfnn    f (t, \wh{Y}_t, \wh{Z}_t )-\hat{f}  (t, \wh{Y}_t, \wh{Z}_t )  \le 0  $, \dtp;

\ss \no \(ii\)  $\hat{f}$ satisfies (H1'),   (H2);   $\hat{f}$ is convex in $z$;  and  $ \D f(t) \dfnn   f(t, Y_t, Z_t )-\hat{f}  (t, Y_t, Z_t )   \le 0 $, \dtp;

  \ss \no \big(where (H1') is an extension of (H1) in that   the constant $\a$  is   replaced by an $\bF$-progressively measurable,
  non-negative process $\{\a_t\}_{t \in [0,T]}$
  such that $E \big[\exp \big\{  p  \int_0^T \a_r dr \big\} \big] < \infty$ for some $p >\g   e^{2\k T} $\big)
      then it holds \pas ~ that $Y_t \le \wh{Y}_t     $ for any $ t  \in [0,T]$.

      \ms In addition, if $ Y_\t  =    \wh{Y}_\t $, \pas ~ for some $\t \in \cS_{0,T}$,    then
  \bea  \label{eq:f333}
  P \left( Y_T = \wh{Y}_T,  \int_\t^T  \D f(s) ds = 0   \right)  >  0.
  \eea
 \end{prop}

 \ss \no {\bf Proof:}   For any $\th \in (\th_0,1)$,  we set $U \dfnn Y \neg - \th \wh{Y}$, $\U \dfnn Z - \th \wh{Z}$      and still define
the $\bF$-progressively measurable process $\{a_t\}_{t \in [0,T]}$ as in $ \eqref{eq:axf011}$.
      It follows that $A_t \dfnn \int_0^t a_s \, ds$, $   t \in [0,T]$ is an $\bF$-adapted process
     with     $A_* \dfnn \underset{t \in [0,T]}{\sup} |A_t| \le \int_0^T \neg   | a_s |   ds \le  \k T$, \pas ~
    Similar to \eqref{eq:d131},   (H1) and the convexity of $\fF$ in $z$ show that \dtp
  \bea   \label{eq:axf111}
  \fF(t, y, Z_t)   
         \le    \th \fF \big(t,  y,   \wh{Z}_t \big)   \neg  + \neg    (1  \neg  -  \neg   \th)  \big(\a_t  \neg  + \neg   \beta |y|\big)
      \neg  + \neg   \hb{$\frac{\g}{2(1-\th)}$} |\U_t|^2   , \q   \fa y \in \hR.
     \eea

    Given $n \in \hN$, we still define the $\bF$-stopping time \eqref{eq:axd071}.  Let $ \z_\th  \dfnn \frac{ \g e^{\k T}}{1-\th}$.
 Applying It\^o's formula to the process $\G_t \dfnn \exp\big\{\z_\th  e^{A_t }U_t \big\}$, $t \in [0,T]$, yields that
  \beas
  \G_{\t_n \land t} = \G_{\t_n} + \int_{\t_n \land t}^{\t_n}  G_s ds + \z_\th  \int_{\t_n \land t}^{\t_n} \G_s e^{A_s}     \big( d V_s -\th d \wh{V}_s \big)
  - \z_\th  \int_{\t_n \land t}^{\t_n} \G_s e^{A_s}  \U_s dB_s, \q t \in [0,T] ,
  \eeas
  where   $   G_t  =   \z_\th  \G_t e^{A_t}  \Big(  f(t,Y_t,Z_t)-\th \hat{f} (t, \wh{Y}_t, \wh{Z}_t) -a_t U_t- \hb{$\frac12$} \z_\th  e^{A_t}|\U_t|^2 \Big) $.
  Let  $\D_\th f (t) \dfnn \th \D f (t)$ if  (i) holds and $\D_\th f (t) \dfnn \D f (t)$ if  (ii) holds for $t \in [0,T]$.
      Then
  \beas
    G_t  =   \z_\th \,  \G_t e^{A_t}  \left(\D_\th f(t) +  \fF(t,Y_t,Z_t)-\th \fF(t, \wh{Y}_t, \wh{Z}_t) -a_t U_t- \hb{$\frac{\g}{2(1-\th)}$} |\U_t|^2 \right), \q t \in [0,T]
  \eeas
  whether  (i)  or (ii) holds. Similar to \eqref{eq:d141}, using \eqref{eq:axf111} and (H2), one can show  by 3 cases that
   \bea    \label{eq:d141x}
   G_t \le \z_\th \,  \G_t e^{A_t}   \D_\th f(t) + \g e^{2 \k T} \, \G_t \Big(\a_t + (\beta + \k) \big(Y^+_t + \wh{Y}^-_t \big)\Big) , \q \dtp
   \eea

 Now, we   define   a    process
        \beas
           D_t \dfnn \exp \left\{ \z_\th \,  \int_0^t e^{A_s}   \D_\th f(s) ds+ \g e^{2 \k T}   \int_0^t        \Big(\a_t+  (\beta+\k) \big( Y^+_s + \wh{Y}^-_s \big) \Big)       ds  \right\} ,  \q  t \in [0,T] .
           \eeas
           Then   it holds \pas~  that
 \beas
  \frac{D_s}{D_t}       & \dneg =& \dneg    \exp \left\{ \z_\th \,  \int_t^s e^{A_r}   \D_\th f(r) dr+ \g e^{2 \k T}   \int_t^s        \Big(\a_r+  (\beta+\k) \big( Y^+_r + \wh{Y}^-_r \big) \Big)       dr  \right\}   \nonumber \\
      & \dneg \le & \dneg        \exp \left\{  \g      e^{2\k T} \bigg( \int_0^T \a_rdr \neg  + \neg     (\beta \neg  + \neg  \k)     T \big(       Y^+_* \neg  + \neg   \wh{Y}^-_*  \big)   \bigg)       \right\}
       \dfnn \eta,  \q  \fa \,  0 \le s \le t \le T.
  \eeas
 As $D_0=1$, we in  particular have $ D_* \le \eta $, \pas~
 Let $q= \frac{p}{\g   }e^{-2\k T}$ and $\wt{\eta} \dfnn \eta \exp  \big\{     \g  e^{2 \k T}  Y^+_*    \big\} $. H\"older's inequality and  \eqref{eq:axf123}  imply that
   \bea
   E\big[D_* \G_* \big] &\tneg \le&\tneg E\big[\eta \G_* \big]   \le   E \bigg[\exp \bigg\{  \g   e^{2\k T}   \int_0^T \a_r dr
    \neg  +  \neg  \big(  (\beta  \neg +  \neg \k)   \g  T   e^{2\k T}    \neg +  \neg  \z_\th  e^{\k T}  \big)   \big( Y^+_*   \neg +  \neg  \wh{Y}^-_* \big)            \bigg\} \bigg]  \nonumber \\
   &\tneg \le&\tneg   \left\| \exp \bigg\{  \g   e^{2\k T}  \neg  \int_0^T  \neg  \a_r dr            \bigg\}\right\|_{\hL^q(\cF_T)}  \neg  \left\| \exp \bigg\{    \big(  (\beta  \neg + \neg \k)   \g  T   e^{2\k T}   \neg + \neg  \z_\th  e^{\k T}  \big)   \big( Y^+_*  \neg + \neg  \wh{Y}^-_* \big)            \bigg\}\right\|_{\hL^{\frac{q}{q-1}}(\cF_T)}
     \neg < \neg \infty.  \label{eq:axf115}     \qq \;\;\;
   \eea
   and  that
   \bea
    E \big[   \wt{\eta}    \big]  \neg   \le   \neg
     \left\| \exp \bigg\{  \g   e^{2\k T}  \neg  \int_0^T  \neg  \a_r dr            \bigg\}\right\|_{\hL^q(\cF_T)}   \neg
      \left\| \exp \bigg\{       (\beta  \neg + \neg \k)   \g  T   e^{2\k T}    \wh{Y}^-_*
        \neg +\neg \big(1 \neg + \neg  (\beta  \neg + \neg \k) T \big) \g      e^{2\k T}   Y^+_*               \bigg\}\right\|_{\hL^{\frac{q}{q-1}}(\cF_T)}
                   \neg < \neg \infty .  \q   \label{eq:axf126}
  \eea
  Fix $t \in [0,T]$.    Integration by parts and \eqref{eq:d141x} imply that    \pas
 \bea
         D_{\t_n \land t} \G_{\t_n \land t}   &\le&     D_{\t_n} \G_{\t_n}   + \z_\th  \int_{\t_n \land t}^{\t_n} D_s \G_s e^{A_s}    \big( d V_s -\th d \wh{V}_s \big)   - \z_\th  \int_{\t_n \land t}^{\t_n} D_s \G_s e^{A_s}  \U_s dB_s \nonumber  \\
     &\le & D_{\t_n} \G_{\t_n}     - \z_\th  \int_{\t_n \land t}^{\t_n} D_s \G_s e^{A_s}  \U_s dB_s       .
   \label{eq:axf119}
 \eea
 Similar to \eqref{eq:d165},      the Burkholder-Davis-Gundy inequality   and \eqref{eq:axf115} show  that
    $\int_0^{\t_n \land \cd}   D_s \G_s e^{A_s}  \U_s  dB_s $ is a       uniformly integrable martingale.

\ms   Taking  $E\big[\cd\big|\cF_{\t_n \land t}\big]$ in  \eqref{eq:axf119}, we see from \eqref{eq:axf115}  that
  \beas
       D_{\t_n \land t} \,  \G_{\t_n \land t}    & \dneg \le  & \dneg   E\big[   D_{\t_n} \G_{\t_n}  \big|\cF_{\t_n \land t}\big]
       =     \b1_{\{\t_n < t\} }     D_{\t_n} \G_{\t_n}
     \dneg +  \neg  \b1_{\{\t_n \ge t\} } E\big[   D_{\t_n} \G_{\t_n}  \big|\cF_t \big]   \neg, \q \pas
  \eeas
   Let $\D \, Y_T \dfnn Y_T -  \wh{Y}_T $. As $n \to \infty$,       Dominated Convergence Theorem,  \eqref{eq:axf115} and \eqref{eq:axf126}  imply that
  \bea
 D_t \G_t      \le        E\big[   D_T \G_T  \big|\cF_t \big]
    \le      E\big[   D_T \exp  \big\{   \th \z_\th\, e^{A_T} \D Y_T +   \g  e^{2 \k T}  Y_T^+    \big\}  \big|\cF_t \big]
  , \q \pas   \label{eq:axd011}
 \eea
Then it follows from   \eqref{eq:axf126} that \pas
 \beas
    \G_t          \le       E\bigg[ \frac{D_T}{ D_t} \exp  \big\{ \th  \z_\th\, e^{A_T} \D Y_T  \neg +  \neg   \g  e^{2 \k T}  Y^+_*    \big\}  \Big|\cF_t \bigg]
     \neg \le \neg     E \Big[   \eta \exp  \big\{     \g  e^{2 \k T}  Y^+_*    \big\}     \big| \cF_t  \Big]   =
      E \big[   \wt{\eta}       \big| \cF_t  \big]   ,
 \eeas
     which leads to that
      \beas
     Y_t   -   \th \wh{Y}_t      \le      \hb{$\frac{1-\th}{\g}$}         \ln   E  \big[         \wt{\eta}     \big| \cF_t \big]
        e^{-\k T-A_t} , \q \pas      
    \eeas
  Letting $\th \to 1$ yields that  $  Y_t -   \wh{Y}_t \le 0$, \pas~
Then the continuity of processes $Y$ and $\wh{Y}$ shows that $P\big( Y_t \le \wh{Y}_t  , \; \fa    t  \in [0,T]  \big) =1$.

\ms In addition, suppose that $ Y_\t  =   \wh{Y}_\t $, \pas ~ for some $\t \in \cS_{0,T}$,
however that   \eqref{eq:f333} does not hold   \big(i.e., $ \D Y_T + \int_\t^T      \D  f(s) ds < 0 $, \pas\big).
 Letting $\th \in (\th_0, 1)$ and taking  $t=\t$  in \eqref{eq:axd011} yield that
 \beas
 0 & \dneg   <&\dneg    \exp  \big\{ \g e^{\k T + A_\t}  Y_\t \big\} = \G_\t      \le        E\bigg[  \frac{ D_T}{D_\t} \exp  \big\{ \th  \z_\th\, e^{A_T} \D Y_T +   \g  e^{2 \k T}  Y^+_*    \big\}  \Big|\,\cF_\t \bigg]  \nonumber \\
 &\dneg    \le &  \dneg     E\bigg[  \wt{\eta}  \exp  \Big\{ \z_\th \, e^{-\k T} \Big( \th \D Y_T  \neg +\neg  \int_\t^T \neg      \D_\th f(s) ds  \Big)         \Big\}   \Big|\,\cF_\t \bigg]     , \q \pas
 \eeas
 As   $\th \to 1$, Dominated Convergence Theorem and \eqref{eq:axf126} imply that
  \beas
 0    <     \exp  \big\{ \g e^{\k T + A_\t}  Y_\t \big\}      \le  \lmtd{\th \to 1}     E\bigg[  \wt{\eta}  \exp  \Big\{ \z_\th \, e^{-\k T} \Big( \th \D Y_T  \neg +\neg  \int_\t^T \neg      \D_\th f(s) ds  \Big)         \Big\}   \Big|\,\cF_\t \bigg]  =0    , \q \pas
 \eeas
  A contradiction appears. \qed

 \begin{prop}  \label{prop:comparion_ext2}
Let $f, \hat{f}:    [0,T] \times \O \times \hR \times \hR^d \rightarrow  \hR $ be two
 $\sP \times \sB(\hR) \times   \sB(\hR^d)/\sB(\hR)$-measurable functions, and let $(Y,Z,V) , (\wh{Y},\wh{Z},\wh{V})\in
 \hC^0_\bF[0,T] \times  \wh{\hH}^{2,loc}_\bF([0,T];\hR^d)    \times \hV_\bF[0,T]$ solve  the   BSDEs \eqref{VBSDE} and BSDEs \eqref{VBSDE2}
 respectively   such that    $ Y_T \le \wh{Y}_T   $, \pas, that \eqref{eq:axf123} holds, and that for some   $\th_0 \in (0, 1)$, $ \th V - \wh{V}$ is a decreasing process for any $\th \in (\th_0,1)$.    If either of the following two holds:

 \ss \no \(\,i\)  $f$ satisfies (H1'),   (H2);  $f$ is concave in $z$; and  $\D f(t) \dfnn   f  (t, \wh{Y}_t, \wh{Z}_t )-\hat{f} (t, \wh{Y}_t, \wh{Z}_t )  \le 0  $, \dtp;

\ss \no \(ii\)  $\hat{f}$ satisfies (H1'),   (H2);  $\hat{f}$ is concave in $z$;  and  $ \D f(t) \dfnn  f  (t, Y_t, Z_t )-\hat{f} (t, Y_t, Z_t )   \le 0 $, \dtp;

\ss \no      then it holds \pas ~ that $Y_t \le \wh{Y}_t     $ for any $ t  \in [0,T]$.   In addition, if $ Y_\t  =    \wh{Y}_\t $, \pas ~ for some $\t \in \cS_{0,T}$,    then \eqref{eq:f333} holds.

 \end{prop}

  \ss \no {\bf Proof:} The triplet $(\cY, \cZ, \cV) \dfnn \big(-\wh{Y}, -\wh{Z}, -\wh{V}\big)$ solves the BSDE\eqref{VBSDE}
  with generator $\mathfrak{f}(t,\o,y,z)\dfnn - \hat{f} (t,\o,-y,-z)$, $\fa (t,\o,y,z) \in [0,T] \times \O \times \hR \times \hR^d$, and
  triplet $\big(\wh{\cY}, \wh{\cZ}, \wh{\cV}\big) \dfnn  (-Y, -Z, -V )$ solves the BSDE\eqref{VBSDE2}
  with generator $\hat{\mathfrak{f}}(t,\o,y,z)\dfnn - f (t,\o,-y,-z)$, $\fa (t,\o,y,z) \in [0,T] \times \O \times \hR \times \hR^d$.
    One can  check that all  conditions in Proposition   \ref{prop:comparion_ext} are satisfied by the new settings.
    Therefore, an application of  Proposition  \ref{prop:comparion_ext}   gives rise to the conclusion.    \qed

 \ss \no {\bf Proof of Theorem \ref{thm:op-gE} :}
 Fix $\nu \in \cS_{0, T}$. For any $\t  \in \cS_{\nu,T}$,  it holds \pas ~ that
  \bea
      \cY_{\t \land t}
     &=&  \cY_\t       + \int_{\t \land t}^\t  g \big(s,  \cY_s , \cZ_s     \big)   ds  +   \cK_\t  - \cK_{\t \land t} - \int_{\t \land t}^\t  \cZ_s dB_s  \nonumber  \\
     &=&   \cY_\t    + \int_t^T  \b1_{\{s < \t \}} g \big(s,  \cY_{\t \land s} ,  \b1_{\{s < \t \}} \cZ_s     \big)   ds  +   \cK_\t  - \cK_{\t \land t} - \int_t^T  \b1_{\{s < \t \}} \cZ_s dB_s, \q  t   \in    [0, T] .   \q
     \label{eq:xx21}
  \eea
     Since  $  \cY_\t  \ge     \b1_{\{\t < T\}}  \cL_\t  +\b1_{\{\t = T\}} \xi   =   \cR_\t      $, \pas,  applying Proposition \ref{prop:comparion_ext}
     and Proposition \ref{prop:comparion_ext2}     with $(Y, Z, V)=\big(Y^{\t, \cR_\t}, Z^{\t, \cR_\t} ,0\big)$
     and $\big(\wh{Y}, \wh{Z}, \wh{V}\big)=\big\{\big(\cY_{\t \land t},  \b1_{\{t < \t \}} \cZ_t,  \cK_{\t \land t} \big)\big\}_{t \in [0,T]} $ yields that \pas
     \beas
      \cY_{\t \land t} \ge Y^{\t, \cR_\t}_t                , \q    t \in [0,T] .
      \eeas
    In particular, we have
    $     \cY_{\nu} \ge Y^{\t, \cR_\t}_\nu = \cE^g_{\nu, \t}[\cR_\t]  $  ,   \pas   ~

    \ms So it remains to show that    $\cY_\nu    =\cE^g_{\nu, \t_*(\nu)}\big[\cR_{\t_*(\nu) }\big]   $, \pas~ To see this, we
     define
      \beas
        \wt{\cY}_t \dfnn \b1_{\{t < \nu\}} Y^{\nu, \cY_\nu}_t+\b1_{\{t \ge \nu    \}} \cY_{\t_*(\nu) \land t}
    \;\;   \hb{ and } \;\;     \wt{\cZ}_t \dfnn \b1_{\{t < \nu\}} Z^{\nu, \cY_\nu}_t+\b1_{\{  \nu \le  t < \t_*(\nu)  \}} \cZ_t   , \q      \fa t \in [0,T]      .
    \eeas
      Clearly,       $ \big(\wt{\cY}, \wt{\cZ} \big) \in \underset{p \in (1, \infty)}{\cap}
 \hE^p_\bF[0,T] \times  \hH^{2,2p}_\bF([0,T];\hR^d)$.
     The flat-off condition of   $(\cY,\cZ,\cK)$ and the continuity of $\cK$ imply that \pas
      \beas
        0= \neg  \int_{ [\nu, \t_*(\nu))}  \neg   \b1_{\{ \cY_s   >  \cL_s  \}}   d\cK_s = \int_{ [\nu, \t_*(\nu))}  \neg   \b1_{\{ \cY_s   >  \cR_s  \}}   d\cK_s
      =\neg  \int_{  [\nu, \t_*(\nu) )} \neg   d\cK_s      =   \underset{s \nearrow \t_*(\nu)}{\lim}\cK_s   -  \cK_{\nu}
      =\cK_{\t_*(\nu)}  -  \cK_{\nu}  .
      \eeas
      Hence, taking $\t = \t_*(\nu)$ and $t =  \nu \vee t       $ in \eqref{eq:xx21}, we can deduce that    \pas
       \bea
   \q  \; \;      \cY_{(\nu \vee t) \land \t_*(\nu) }
     &=& \cY_{\t_*(\nu)}          + \int_{ \nu \vee t  }^T \b1_{\{s < \t_*(\nu) \}}   g \big(s,  \cY_{ \t_*(\nu) \land s} ,  \b1_{\{s <  \t_*(\nu) \}} \cZ_s     \big)   ds
                   -  \int_{ \nu \vee t  }^T \b1_{\{s < \t_*(\nu) \}} \cZ_s dB_s  \nonumber  \\
    &=& \cR_{\t_*(\nu)}          + \int_{ \nu \vee t  }^T \b1_{\{s < \t_*(\nu) \}}   g \big(s,  \wt{\cY}_s,  \wt{\cZ}_s     \big)   ds
                -  \int_{ \nu \vee t  }^T \wt{\cZ}_s dB_s  , \q t \in [0,T] . \label{eq:xx27}
       \eea
     In particular, we have
     \bea  \label{eq:xx29}
      \cY_\nu
     = \cR_{\t_*(\nu)}          + \int_\nu^T \b1_{\{s < \t_*(\nu) \}}   g \big(s,  \wt{\cY}_s,  \wt{\cZ}_s     \big)   ds
                -  \int_\nu^T \wt{\cZ}_s dB_s  , \q  \pas
     \eea

    Fix $t \in [0,T]$. One can deduce from \eqref{eq:xx25} and \eqref{eq:xx29} that
        \beas
          \b1_{\{t < \nu\}} Y^{\nu, \cY_\nu}_t &=& \b1_{\{t < \nu\}}  \cY_\nu
          + \b1_{\{t < \nu\}}  \int_t^\nu  g \Big(s,Y^{\nu, \cY_\nu}_s,  Z^{\nu, \cY_\nu}_s \Big)ds-\b1_{\{t < \nu\}}  \int_t^\nu    Z^{\nu, \cY_\nu}_s  dB_s \\
  &=& \b1_{\{t < \nu\}}  \cY_\nu + \b1_{\{t < \nu\}} \int_t^\nu  g\big(s,\wt{\cY}_s,\wt{\cZ}_s\big)ds
  - \b1_{\{t < \nu\}} \int_t^\nu    \wt{\cZ}_s  dB_s\\
   &=&  \b1_{\{t < \nu\}}  \cR_{\t_*(\nu)}  + \b1_{\{t < \nu\}} \int_t^T \b1_{\{s < \t_*(\nu) \}}  g \big(s,\wt{\cY}_s,\wt{\cZ}_s\big)ds
  - \b1_{\{t < \nu\}} \int_t^T    \wt{\cZ}_s  dB_s   ,
        \eeas
      which together with \eqref{eq:xx27} implies that \pas
       \bea \label{eq:xx35}
        \wt{\cY}_t = \cR_{\t_*(\nu)}  +   \int_t^T \b1_{\{s < \t_*(\nu) \}}  g \big(s,\wt{\cY}_s,\wt{\cZ}_s\big) ds
  -   \int_t^T    \wt{\cZ}_s  dB_s   .
       \eea
      The continuity of process   $\wt{\cY}_t$ further shows that \pas, \eqref{eq:xx35} holds for any $t \in [0,T]$.  To wit,
       $(\wt{\cY}, \wt{\cZ})  \in \underset{p \in (1, \infty)}{\cap}
 \hE^p_\bF[0,T] \times  \hH^{2,2p}_\bF([0,T];\hR^d) $ is the unique solution of the BSDE \eqref{BSDEtau} with $(\t, \xi)= \left(\t_*(\nu), \cR_{\t_*(\nu)}\right) $. Therefore, it follows that
        $           \cY_\nu = \wt{\cY}_\nu  = \cE^g_{\nu,  \t_*(\nu)} \left[ \cR_{\t_*(\nu)} \right]$.
              \qed

\section{Stability}

   \begin{thm}    \label{thm:stable}
Let   $\left\{(\xi_m, f_m,L^m) \right\}_{m \in \hN_0}$ be a sequence of parameter sets such that

  \ss   \no   (S1)  With the same constants $\a,    \beta,  \k  \ge 0$ and $\g > 0$, $f_0$ satisfies \(H1\) and $\{f_n\}_{n \in \hN}$    satisfy   \(H1\)-\(H3\);

  \ss   \no   (S2)  It holds \pas~ that $\xi_n$ converges to $\xi_0$ and that $L^n_t$ converges to $L^0_t$ uniformly in $t \in [0,T]$;

 \ss   \no   (S3)  $   \Xi(p)   \dfnn   \underset{m \in  \hN_0}{\sup}  E \left[ e^{ p \,   (\xi^+_m \vee L^m_*  )   }    \right] < \infty  $ for all $p \in (1, \infty)$.

\ss We let $(Y^0, Z^0, K^0) \in  \underset{p \in [1, \infty)}{\cap} \hS^p_\bF[0,T] $  be a solution of the quadratic RBSDE$(\xi_0, f_0,L^0) $
 , and for any $n \in  \hN$ we let $(Y^n, Z^n, K^n)  $  be the unique solution of the quadratic RBSDE$(\xi_n, f_n,L^n) $
in $ \underset{p \in [1, \infty)}{\cap} \hS^p_\bF[0,T]$.
 If   $  f_n \big(t, Y^0_t,   Z^0_t \big) $ converges \dtp~ to $  f_0\big(t, Y^0_t, Z^0_t \big) $,  then   for any $p \in [1, \infty) $
 \beas
   \lmt{n \to \infty} E \neg \left[ \exp\bigg\{ p \cd \underset{t \in [0,T]}{\sup} |Y^n_t -Y^0_t|\bigg\} \right] =1 \;\;\; \hb{and} \;\;\;
          \lmt{n \to \infty} E \neg \left[ \left( \int_0^T |Z^n_s - Z^0_s |^2 ds\right)^p \, \right] = 0     .
 \eeas
 Moreover, if it holds \dtp~ that  $ f_n(t, \o, y, z)  $
 converges to  $f_0(t, \o, y, z) $ locally uniformly in $(y, z)$, 
 then up to a subsequence, we further have
    \bea  \label{eq:e141}
   \lmt{n \to \infty} E \Big[ \big|K^n_T - K^0_T \big|^p \Big] = 0  , \q  \fa  p \in [1, \infty)  .
 \eea

 \end{thm}

  \ss   In (S1) of Theorem \ref{thm:stable}, the convexity/concavity does not need to be the same for all generators $f_n's$,
   for example, it can be alternate as in the following example.

  \begin{eg}
  Let $d =1$. For any $m \in \hN_0$, the function
  \beas
  f_m(t, \o, y,z ) \dfnn  (-1)^m    z^2     , \q  \fa (t, \o, y,z ) \in [0,T] \times \O \times \hR \times \hR
  \eeas
  is   $\sP \times \sB(\hR) \times   \sB(\hR)/\sB(\hR)$-measurable and   satisfies (H1), (H2) with $(\a,\beta, \g, \k)= (0,0,2,0)$. Moreover, $f_m$ is convex (resp.\;concave) in $z$  when $m$ is even (resp.\;odd).   Clearly,
     $ (0,0,0)$ is the unique solution of RBSDE$(0,f_0, 0)$ in $\underset{p \in [1, \infty)}{\cap} \hS^p_\bF[0,T]$. For any $n \in  \hN$ we set $L^n_t \dfnn
      \frac{T-t}{n} $, $t \in [0,T]$, and let $(Y^n, Z^n, K^n)  $  be the unique solution of the quadratic RBSDE$(0, f_n,  L^n ) $
in $ \underset{p \in [1, \infty)}{\cap} \hS^p_\bF[0,T]$. As $f_m(\cd,\cd,0,0) \equiv 0$ for all $m \in \hN_0$,
 the first part of Theorem \ref{thm:stable} yields that
 \beas
   \lmt{n \to \infty} E \left[ e^{ p  Y^n_*         }         \right] = 1    \q \hb{and} \q \lmt{n \to \infty} E \left[ \left( \int_0^T |Z^n_s   |^2 ds\right)^p \,  \right] = 0  , \q  \fa  p \in [1, \infty)  .
 \eeas
\end{eg}

\ss  \no  {\bf Proof of Theorem \ref{thm:stable}:}   {\bf 1)} Fix $n \in \hN$, $\th \in (0,1)$ and $\e >0$.  We first show that \pas
\bea \label{eq:e191}
     | Y^n_t - Y^0_t    |  \le  (1-\th) \big(|Y^0_t| + |Y^n_t|\big) +  \hb{$\frac{1-\th  }{\g}$} \ln \left( \sum^4_{i=1} I^{n,i}_t \right) , \q t \in [0,T]   ,
   \eea
 where  $I^{n, i}_t \dfnn E \Big[  I^{n, i}_T \big| \cF_t \Big]$ for $i=1,2,3,4$ such that
 \beas
     I^{n,1}_T &\dfnn&     D_{T} \eta_n  \; \hb{ with } \, D_t \dfnn \exp \left\{ \g e^{2 \k T}   \int_0^t        \big(\a+  (\beta+\k) |Y^0_s| \big)       ds  \right\}  ,  ~   t \in [0,T]  \, \hb{ and }  \\
    & &    \q      \eta_n \dfnn \exp\Big\{ \z_\th e^{\k T} \big( | \xi_n-\th \xi_0 | \vee | \xi_0-\th \xi_n | \big)  \Big\};   \\
      I^{n,2}_T  &\dfnn&   \z_\th e^{\k T}     D_T \U_n \neg \int_0^T \neg       \big|\D_n f(s)   \big| ds   \,
      \hb{ with } \, \z_\th  \dfnn \frac{ \g e^{\k T}}{1-\th}, ~  \U_n \dfnn   \exp\Big\{ \z_\th e^{\k T} \big( Y^n_* + Y^0_* \big)   \Big\} \, \hb{ and }  \\
       &&    \q     \D_n f(t) \dfnn  f_n \big(t, Y^0_t,   Z^0_t \big) \neg - \neg  f_0\big(t, Y^0_t, Z^0_t \big)  , ~ t \in [0,T] ; \\
       I^{n,3}_T &\dfnn& \Big( 1   + \z_\th \exp \big\{\k T+ \e \z_\th e^{\k T}\big\} \Big)  \Big( 1+      D_T  \exp \left\{\g e^{2\k T}  \big( Y^0_*  + Y^n_* \big)  \right\}
      \big(  K^0_T + K^n_T \big) \Big) ;   \\
        I^{n,4}_T &\dfnn&   \frac{\z_\th}{\e}   e^{\k T}       D_T \U_n      \Big(   \underset{t \in [0,T]}{\sup} |L^n_t - L^0_t |  \Big) \big(  K^0_T + K^n_T \big) \, .
    \eeas
 {\bf Case 1:}   $f_n$ is convex in $z$.   We set $U^n \dfnn Y^n \neg - \th Y^0$, $V^n \dfnn Z^n  \neg -\th Z^0$
     and define two processes
   \beas
   \qq   a^n_t \dfnn       \b1_{\{U^n_t \ne 0 \}}\frac{f_n(t,Y^n_t,Z^n_t) -f_n(t,\th Y^0_t,Z^n_t)}{U^n_t
      } - \k \b1_{\{U^n_t = 0 \}} ,
          \q    A^n_t \dfnn \int_0^t a^n_s ds, \q   t \in [0,T].
      \eeas
          Applying It\^o's formula to the process $\G^{1,n}_t \dfnn \exp\big\{\z_\th  e^{A^n_t }U^n_t \big\}$, $t \in [0,T]$, yields that
  \beas
  \G^{1,n}_t = \G^{1,n}_T + \int_t^T  G^{1,n}_s ds + \z_\th  \int_t^T \G^{1,n}_s e^{A^n_s}    (d K^n_s-\th  d K^0_s)
   - \z_\th  \int_t^T \G^{1,n}_s e^{A^n_s}  V^n_s dB_s, \q t \in [0,T] ,
  \eeas
  where   $   G^{1,n}_t  =   \z_\th \,  \G^{1,n}_t e^{A^n_t}  \left(  f_n(t,Y^n_t,Z^n_t)-\th f_0(t, Y^0_t, Z^0_t) -a^n_t U^n_t- \hb{$\frac12$} \z_\th  e^{A^n_t}|V^n_t|^2 \right) $.
    Similar to \eqref{eq:d131},  (H1) and the convexity of $f_n$ in $z$  show    that \dtp
  \beas
  f_n(t, Y^0_t,Z^n_t)
           \le    \th f_n \big(t, Y^0_t,   Z^0_t \big) +  (1-\th)  \left(\a+ \beta |Y^0_t|\right) + \hb{$\frac{\g}{2(1-\th)}$}  |V^n_t|^2 ,
     \eeas
    which together with (H2)  implies that  \dtp
    \beas
        G^{1,n}_t      &=& \z_\th  \G^{1,n}_t e^{A^n_t}  \left(   f_n(t,\th Y^0_t,Z^n_t) -\th f_0(t, Y^0_t, Z^0_t)  - \hb{$\frac12$} \z_\th  e^{A^n_t}|V^n_t|^2 \right)  \\
       &\le&\z_\th  \G^{1,n}_t e^{A^n_t} \left(   \big|f_n(t,\th Y^0_t,Z^n_t) -f_n(t,  Y^0_t,Z^n_t)  \big| + f_n(t,  Y^0_t,Z^n_t) - \th f_0\big(t, Y^0_t, Z^0_t \big)
        -\hb{$\frac{\g}{2(1-\th)}$} |V^n_t|^2 \right)       \nonumber  \\
         &\le &   \g e^{2\k T}  \G^{1,n}_t       \big(\a+ (\beta+\k) |Y^0_t| \big)         +   \z_\th   e^{\k T}  \G^{1,n}_t      \big|\D_n f(t)   \big|   .
   \eeas
            Integration by parts gives that
 \bea
   \G^{1,n}_t  & \tneg \le &  \tneg    D_t \G^{1,n}_t    \le     D_{T} \G^{1,n}_{T}
    + \z_\th e^{\k T} \neg \int_t^T \neg   D_s    \G^{1,n}_s     \big|\D_n f(s)   \big| ds
    +   \z_\th  \int_t^{T} \neg  D_s \G^{1,n}_s   e^{A^n_s}  d K^n_s     - \z_\th  \neg \int_t^{T} \dneg  D_s \G^{1,n}_s e^{A^n_s}  V^n_s dB_s  \nonumber \\
    &   \tneg    \le &  \tneg   I^{n,1}_T +   I^{n,2}_T
    +   \z_\th e^{\k T}   D_T \neg \int_0^T \neg \G^{1,n}_s  dK^n_s     - \z_\th   \neg \int_t^{T} \dneg  D_s \G^{1,n}_s e^{A^n_s}  V^n_s dB_s  , \q   t \in [0,T] . \qq
 \label{eq:e151}
 \eea
   The flat-off condition of $(Y^n, Z^n, K^n)$ implies that
 \bea
 \int_0^T \neg \G^{1,n}_s  dK^n_s &=& \int_0^T \neg \b1_{\{Y^n_s = L^n_s\}} \G^{1,n}_s  dK^n_s= \int_0^T \neg \b1_{\{Y^n_s = L^n_s \le L^0_s + \e \}}\G^{1,n}_s  dK^n_s + \int_0^T \neg \b1_{\{Y^n_s = L^n_s > L^0_s + \e \}}\G^{1,n}_s  dK^n_s  \nonumber \\
  &\le &   \int_0^T \neg \b1_{\{Y^n_s   \le Y^0_s + \e \}}  \exp \left\{\g e^{2\k T}   |Y^0_s| +\e \z_\th e^{\k T}    \right\}   dK^n_s
   + \U_n \int_0^T \neg \b1_{\{  |L^n_s - L^0_s |>  \e \}}       dK^n_s \nonumber \\
    &\le &          \exp \left\{\g e^{2\k T}   Y^0_*   +   \e \z_\th e^{\k T}   \right\}      K^n_T
     + \frac{1}{\e}   \U_n   \Big(   \underset{t \in [0,T]}{\sup} |L^n_t - L^0_t |  \Big)  K^n_T, \q \pas   \label{eq:e153}
 \eea

      For each $p \in (1, \infty)$,  Theorem \ref{thm:existence}   and (S3) imply that
  \beas
                 \underset{n' \in \hN}{\sup}      E \neg \left[e^{ p \g   Y^{n'}_* } \dneg + \dneg  \left( \int_0^T \neg    |Z^{n'}_s|^2 ds \right)^p
                  \neg + \neg \big( K^{n'}_T\big)^p \right]     \neg       \le           c_p  \underset{n' \in \hN}{\sup}
                  E \neg \left[        e^{ 3 p  \g e^{\beta T} \left ( \xi^+_{n'}  \vee     L^{n'}_* \right)  }        \right]
                  \neg  \le \neg c_p\, \Xi \Big(3 p  \g e^{\beta T}\Big) .
       \eeas
       Thus, it follows that
        \bea   \label{eq:e155}
                 \underset{m \in \hN_0}{\sup}      E \neg \left[e^{ p \g   Y^m_* } \dneg + \dneg  \left( \int_0^T \neg    |Z^m_s|^2 ds \right)^p
                  \neg + \neg \big( K^m_T\big)^p \right]     \neg       \le  \neg         c_p   \, \Xi \Big(3 p  \g e^{\beta T}\Big) \neg
                  +\neg  E \neg \left[e^{ p \g   Y^0_* } \dneg + \dneg  \left( \int_0^T \neg    |Z^0_s|^2 ds \right)^p  \neg + \neg \big( K^0_T\big)^p \right]
                  \neg \dfnn  \wt{\Xi}(p)  ,
       \eea
       which together with   (S1) implies that
          \bea
               E [ \eta^p_n]    & \tneg \neg \le & \tneg  \neg   E  \Big[ e^{ p \, \z_\th e^{\k T}  ( | \xi_n \neg |+| \xi_0 |   )   }  \Big]     \le   \frac12    E \left[  e^{ 2 p \, \z_\th e^{\k T}   (\xi^+_n\vee L^n_* )   } \neg + \neg  e^{ 2 p \, \z_\th e^{\k T}  (\xi^+_0\vee L^0_* )   } \right] \neg \le  \Xi \Big( 2 p \, \z_\th e^{\k T}  \Big) ,  \label{eq:e157a} \qq \\
                  E [ \U^p_n]         &  \tneg  \neg  \le & \tneg  \neg    \frac12    E \left[  e^{ 2 p \, \z_\th e^{\k T}   Y^n_*   } + e^{ 2 p \, \z_\th e^{\k T}   Y^0_*   } \right]
                \le \wt{\Xi} \Big(\hb{$\frac{2 p   }{1-\th}$}  e^{2\k  T} \Big) ,  \label{eq:e157b} \\
         E \neg \left[ \neg \left( \int_0^T \neg       \big|\D_n f(s)   \big| ds \right)^p \, \right]
       &   \tneg  \neg  \le & \tneg   \neg     E \neg \left[ \neg \left(  2T(  \a  + \beta Y^0_* )   \neg + \neg \g \int_0^T \neg   | Z^0_s |^2    ds  \right)^p \, \right] \le    c_p  E \neg  \left[  e^{ p \g Y^0_* }       \neg + \neg    \left(  \int_0^T \neg           |Z^0_s|^2    ds  \right)^p \, \right]  , \label{eq:e157c}  \\
         E \neg \left[ \underset{t \in [0,T]}{\sup} |L^n_t - L^0_t |^p \, \right]
       &   \tneg  \neg  \le & \tneg  \neg   c_p  E  \Big[ \big(L^n_*\big)^p +  \big(L^0_*\big)^p  \Big] \le c_p  E  \Big[ e^{p L^n_* } +  e^{p L^0_* }  \Big]
       \le   c_p \,  \Xi (p) .  \label{eq:e157d}
         \eea
         Since $       D_T        \le   c_0    \exp\Big\{  \g       (\beta \neg + \neg \k) T e^{2 \k T}  \, Y^0_*        \Big\} $,  \pas,
         we also see that $D_T \in \hL^p(\cF_T)$.  Thus, one can deduce from  Young's inequality and  \eqref{eq:e155}-\eqref{eq:e157d} that
         random variables $I^{n,i}_T $, $i=1,2,3,4$           are all integrable.  
                 Moreover,   the Burkholder-Davis-Gundy inequality and H\"older's inequality imply that
  \bea
     E \neg \left[\underset{t \in [0,T]}{\sup}  \left|   \int_0^t  \neg D_s \G^{1,n}_s e^{A^n_s}  V^n_s dB_s \right|\right]
      & \tneg \neg \le & \tneg \neg  c_0 E  \neg \left[   \left(\int_0^T    \dneg  \big( D_s \G^{1,n}_s \big)^2  e^{2 A^n_s}   | V^n_s |^2 d s\right)^{\frac{1}{2}}\right]
   \le c_0 E \neg \left[ D_T \U_n  \neg \left( \int_0^T \neg  | V^n_s |^2 d s\right)^{\frac{1}{2}}  \right]  \nonumber  \\
 & \tneg \neg  \le & \tneg \neg   c_0     \left\| D_T \right\|_{\hL^4(\cF_T)}          \left\| \U_n \right\|_{\hL^4(\cF_T)}
   \left\|  V^n \right\|_{ \hH^2_\bF([0,T]; \hR^d)} < \infty ,  \label{eq:e159}
 \eea
    thus $\int_0^\cd   D_s \G^{1,n}_s e^{A^n_s}  V^n_s  dB_s $ is a  uniformly integrable martingale.

   \ms  For any $t \in [0,T]$,       taking  $E[\cd|\cF_t]$ in \eqref{eq:e153} and \eqref{eq:e151} yields that
   $ \dis  \G^{1,n}_t     \le      \sum^4_{i=1} I^{n,i}_t$, \pas  ~ It then follows that
     \beas
            Y^n_t-\th Y^0_t \le  \hb{$\frac{1-\th }{\g}$}  e^{-\k T -A^n_t}  \ln  \neg \left(  \sum^4_{i=1} I^{n,i}_t \right)
           \le \hb{$\frac{1-\th  }{\g}$}   \ln \neg \left(   \sum^4_{i=1} I^{n,i}_t  \right) , \q \pas,
     \eeas
 which implies that
         \bea     \label{eq:e163}
           Y^n_t-  Y^0_t \le  (1-\th) |Y^0_t| +  \hb{$\frac{1-\th  }{\g}$}      \ln  \left(  \sum^4_{i=1} I^{n,i}_t\right) ,   \q \pas
         \eea

               To show the other half of \eqref{eq:e191},  we set $\wt{U}^n \dfnn Y^0-\th Y^n$,
           $\wt{V}^n \dfnn Z^0-\th Z^n$  and define two processes
   \beas
   \qq   \tilde{a}^n_t \dfnn       \b1_{\{Y^0_t \ne Y^n_t \}}\frac{f_n(t,Y^0_t,Z^n_t) -f_n(t, Y^n_t,Z^n_t)}{Y^0_t - Y^n_t
      } - \k  \b1_{\{Y^0_t = Y^n_t \}}  ,
          \q    \wt{A}^n_t \dfnn \int_0^t \tilde{a}^n_s ds, \q   t \in [0,T].
      \eeas
   Applying It\^o's formula to the process $\wt{\G}^{1,n}_t \dfnn \exp\big\{\z_\th  e^{\wt{A}^n_t }\wt{U}^n_t \big\}$, $t \in [0,T]$, yields that
  \beas
  \wt{\G}^{1,n}_t = \wt{\G}^{1,n}_T + \int_t^T  \wt{G}^{1,n}_s ds + \z_\th  \int_t^T \wt{\G}^{1,n}_s e^{\wt{A}^n_s}     ( d K^0_s-\th d K^n_s) - \z_\th  \int_t^T \wt{\G}^{1,n}_s e^{\wt{A}^n_s}  \wt{V}^n_s dB_s, \q t \in [0,T] ,
  \eeas
  where   $   \wt{G}^{1,n}_t  =   \z_\th \,  \wt{\G}^{1,n}_t e^{\wt{A}^n_t}  \left(  f_0(t,Y^0_t,Z^0_t)-\th f_n(t, Y^n_t, Z^n_t) -\tilde{a}^n_t \wt{U}^n_t- \hb{$\frac12$} \z_\th  e^{\wt{A}^n_t}|\wt{V}^n_t|^2 \right) $.
   Similar to \eqref{eq:d131},  (H1) and the convexity of $f_n$ in $z$   show that \dtp
  \beas
  f_n(t, Y^0_t,Z^0_t)         \le    \th f_n \big(t, Y^0_t,   Z^n_t \big)  +  (1-\th)  \left(\a+ \beta |Y^0_t|\right) + \hb{$\frac{\g}{2(1-\th)}$} \big|\wt{V}^n_t \big|^2 ,
     \eeas
    which together with    (H2) implies that \dtp
    \beas
        \wt{G}^{1,n}_t      & \tneg \le & \tneg \z_\th  \wt{\G}^{1,n}_t e^{\wt{A}^n_t} \neg \left( - \D_n f(t)    +f_n(t, Y^0_t,Z^0_t)   -\th f_n \big(t, Y^n_t, Z^n_t \big)        -\tilde{a}^n_t \wt{U}^n_t         - \hb{$\frac{\g}{2(1-\th)}$}  |\wt{V}^n_t|^2 \right)  \\
       & \tneg \le& \tneg \z_\th  \wt{\G}^{1,n}_t e^{\wt{A}^n_t}   \neg \left(   \big| \D_n f(t) \big| + \th f_n \big(t, Y^0_t,   Z^n_t \big) - \th f_n \big(t, Y^n_t, Z^n_t \big)
        -\tilde{a}^n_t \wt{U}^n_t+  (1-\th)  \left(\a+ \beta |Y^0_t|\right)   \right)    \\
         & \tneg= & \tneg  \z_\th  \wt{\G}^{1,n}_t e^{\wt{A}^n_t}  \neg \left(   \big| \D_n f(t) \big|  \neg
         +  \neg (\th  \neg - \neg 1) \tilde{a}^n_t  Y^0_t
          \neg +  \neg  (1  \neg -  \neg \th)   \neg \left(\a  \neg +  \neg  \beta |Y^0_t|\right)   \right)
          \le   \g e^{2 \k T}  \wt{\G}^{1,n}_t      \big(\a  \neg +  \neg  (\beta  \neg +  \neg \k) |Y^0_t| \big)
          \neg +  \neg  \z_\th e^{\k T}    \wt{\G}^{1,n}_t     \big|\D_n f(t)   \big|   .
   \eeas
    Similarly to  \eqref{eq:e151},   integration by parts  gives that
 \bea
    \wt{\G}^{1,n}_t     
          \le     I^{n,1}_T +   I^{n,2}_T
    \neg +    \z_\th e^{\k T}   D_T \neg \int_0^T  \neg \wt{\G}^{1,n}_s d K^0_s \neg
    - \neg \z_\th   \int_t^{T} \neg  D_s \wt{\G}^{1,n}_s e^{\wt{A}^n_s}  \wt{V}^n_s dB_s   ,    \q      t \in [0,T] ,
 \label{eq:e181}
 \eea
 where $\int_0^\cd   D_s \wt{\G}^{1,n}_s e^{\wt{A}^n_s}  \wt{V}^n_s  dB_s $ is a
  uniformly integrable martingale, which can be shown by using similar arguments to those that lead to  \eqref{eq:e159}.
  And similar to \eqref{eq:e153}, the flat-off condition of $(Y^0, Z^0, K^0)$ implies that
 \bea  \label{eq:e185}
    \int_0^T \neg \wt{\G}^{1,n}_s  dK^0_s     \le        \exp \left\{\g e^{2\k T}   Y^n_*  + \e \z_\th e^{\k T}   \right\}   K^0_T
     + \frac{1}{\e} \U_n      \Big(   \underset{t \in [0,T]}{\sup} |L^n_t - L^0_t |  \Big) K^0_T   ,  \q  \pas
 \eea

 For any $t \in [0,T]$,  taking  $E[\cd|\cF_t]$ in \eqref{eq:e185} and \eqref{eq:e181} yields that
         \beas
           Y^0_t-  Y^n_t \le  (1-\th) |Y^n_t| + \hb{$\frac{1-\th  }{\g}$}    \ln  \left(  \sum^4_{i=1} I^{n,i}_t \right)     , \q \pas,
     \eeas
  which together with  \eqref{eq:e163} as well as  the continuity of processes  $Y^n$, $Y^0$ and $\dis \sum^4_{i=1} I^{n,i}$  implies \eqref{eq:e191}.

   \ss \no  {\bf Case 2:}   $f_n$ is concave in $z$.    Applying It\^o's formula to the process
 $\G^{2,n}_t   \dfnn  \big( \G^{1,n}_t\big)^{-1} =   \exp \neg \big\{\neg - \neg \z_\th  e^{A^n_t }U^n_t \big\}$, $t \in [0,T]$ yields that
  \beas
  \G^{2,n}_t = \G^{2,n}_T + \int_t^T  G^{2,n}_s ds + \z_\th  \int_t^T \G^{2,n}_s e^{A^n_s}    (\th  d K^0_s - d K^n_s )
   + \z_\th  \int_t^T \G^{2,n}_s e^{A^n_s}  V^n_s dB_s, \q t \in [0,T] ,
  \eeas
  where   $   G^{2,n}_t  =   \z_\th \,  \G^{2,n}_t e^{A^n_t}  \left( \th f_0(t, Y^0_t, Z^0_t) - f_n(t,Y^n_t,Z^n_t)  + a^n_t U^n_t- \hb{$\frac12$} \z_\th  e^{A^n_t}|V^n_t|^2 \right) $.
   Similar to \eqref{eq:d131b},  (H1) and the concavity of $f_n$ in $z$ show    that \dtp
  \bea  \label{eq:d133}
     f_n(t, Y^0_t,Z^n_t)
                \ge     \th f_n \big(t, Y^0_t,   Z^0_t \big)   -   (1- \th)  \neg \left(\a + \beta |Y^0_t|\right) \neg - \neg  \hb{$\frac{\g}{2(1-\th)}$}  |V^n_t|^2,    \q
     \eea
    which together with (H2)  implies that  \dtp
    \beas
        G^{2,n}_t      & \tneg =& \tneg  \z_\th  \G^{2,n}_t e^{A^n_t}  \left(  - \th \D_n f(t) + \th f_n(t, Y^0_t, Z^0_t) - f_n(t, \th Y^0_t, Z^n_t)   - \hb{$\frac12$} \z_\th  e^{A^n_t}|V^n_t|^2 \right)  \\
       & \tneg  \le& \tneg  \z_\th  \G^{2,n}_t e^{A^n_t} \left( |\D_n f(t)| + \th f_n(t, Y^0_t, Z^0_t)  - f_n(t, Y^0_t,Z^n_t) +  |f_n(t,  Y^0_t,Z^n_t)  -  f_n(t,\th Y^0_t,Z^n_t)   |
        -\hb{$\frac{\g}{2(1-\th)}$} |V^n_t|^2 \right)       \nonumber  \\
         & \tneg  \le & \tneg     \g e^{2\k T}  \G^{2,n}_t       \big(\a+ (\beta+\k) |Y^0_t| \big)         +   \z_\th   e^{\k T}  \G^{2,n}_t      \big|\D_n f(t)   \big|   .
   \eeas
  Similar to \eqref{eq:e151},  integration by parts gives that
 \bea
  \qq  \G^{2,n}_t  
       \le       I^{n,1}_T +   I^{n,2}_T
    +     \z_\th e^{\k T}   D_T \neg \int_0^T \neg \G^{2,n}_s  dK^0_s
     + \neg  \z_\th   \neg \int_t^{T} \dneg  D_s \G^{2,n}_s e^{A^n_s}  V^n_s dB_s ,     \q      t \in [0,T] , \qq
 \label{eq:e151b}
 \eea
  where $\int_0^\cd   D_s \G^{2,n}_s e^{A^n_s}  V^n_s  dB_s $ is a  uniformly integrable martingale, which can be shown
   by using similar arguments to those lead to \eqref{eq:e159}.  And  similar to \eqref{eq:e153}, the flat-off condition of $(Y^0, Z^0, K^0)$ implies that
 \bea
 \int_0^T \neg \G^{2,n}_s  dK^0_s  
     \le            \exp \left\{\g e^{2\k T}   Y^n_* + \e \z_\th e^{\k T}   \right\}      K^0_T
     + \frac{1}{\e}   \U_n   \Big(   \underset{t \in [0,T]}{\sup} |L^n_t - L^0_t |  \Big)  K^0_T, \q \pas   \label{eq:e153b}
 \eea
    For any $t \in [0,T]$,     taking  $E[\cd|\cF_t]$ in \eqref{eq:e153b} and \eqref{eq:e151b} yields that
   $ \dis  \G^{2,n}_t     \le      \sum^4_{i=1} I^{n,i}_t$, \pas  ~ It then follows that
      \bea  \label{eq:e163b}
       Y^0_t-  Y^n_t \le  (1-\th) |Y^0_t| + \th Y^0_t -   Y^n_t
           \le (1-\th) |Y^0_t| +\hb{$\frac{1-\th  }{\g}$}   \ln \neg \left(   \sum^4_{i=1} I^{n,i}_t  \right) , \q \pas
     \eea

               It remains to show the other half of \eqref{eq:e191} for Case 2.
   Applying It\^o's formula to the process  $\wt{\G}^{2,n}_t \dfnn  \big(\wt{\G}^{1,n}_t\big)^{-1} = \exp\big\{-\z_\th  e^{\wt{A}^n_t }\wt{U}^n_t \big\}$,
   $t \in [0,T]$, yields that
  \beas
  \wt{\G}^{2,n}_t = \wt{\G}^{2,n}_T + \int_t^T  \wt{G}^{2,n}_s ds + \z_\th  \int_t^T \wt{\G}^{2,n}_s e^{\wt{A}^n_s}     ( \th d K^n_s- d K^0_s) + \z_\th  \int_t^T \wt{\G}^{2,n}_s e^{\wt{A}^n_s}  \wt{V}^n_s dB_s, \q t \in [0,T] ,
  \eeas
  where   $   \wt{G}^{2,n}_t  =   \z_\th \,  \wt{\G}^{2,n}_t e^{\wt{A}^n_t}  \left(  \th f_n(t,Y^n_t,Z^n_t)-  f_0(t, Y^0_t, Z^0_t) + \tilde{a}^n_t \wt{U}^n_t- \hb{$\frac12$} \z_\th  e^{\wt{A}^n_t} \big| \wt{V}^n_t \big|^2 \right) $.
   Similar to \eqref{eq:d133},   (H1) and the concavity of $f_n$ in $z$  show   that \dtp
  \beas
  f_n(t, Y^0_t,Z^0_t)         \ge    \th f_n \big(t, Y^0_t,   Z^n_t \big)  -  (1-\th) \neg  \left(\a+ \beta |Y^0_t|\right) - \hb{$\frac{\g}{2(1-\th)}$} \big| \wt{V}^n_t \big|^2 ,
     \eeas
    which together with    (H2) implies that \dtp
    \beas
        \wt{G}^{2,n}_t      & \tneg \le & \tneg \z_\th  \wt{\G}^{2,n}_t e^{\wt{A}^n_t} \neg \left( \th f_n \big(t, Y^n_t, Z^n_t \big) - f_n(t, Y^0_t,Z^0_t)  + \D_n f(t)               + \tilde{a}^n_t \wt{U}^n_t         - \hb{$\frac{\g}{2(1-\th)}$}  |\wt{V}^n_t|^2 \right)  \\
       & \tneg \le& \tneg \z_\th  \wt{\G}^{2,n}_t e^{\wt{A}^n_t}   \neg \left(  \th f_n \big(t, Y^n_t, Z^n_t \big) - \th f_n \big(t, Y^0_t,   Z^n_t \big)+   \big| \D_n f(t) \big|         + \tilde{a}^n_t \wt{U}^n_t+  (1-\th)  \left(\a+ \beta |Y^0_t|\right)   \right)    \\
         & \tneg= & \tneg  \z_\th  \wt{\G}^{2,n}_t e^{\wt{A}^n_t}  \neg \left(       (1 \neg - \neg \th    ) \tilde{a}^n_t  Y^0_t   \neg
         +  \neg  \big| \D_n f(t) \big|          \neg +  \neg  (1  \neg -  \neg \th)   \neg \left(\a  \neg +  \neg  \beta |Y^0_t|\right)   \right)
          \le   \g e^{2 \k T}  \wt{\G}^{2,n}_t      \big(\a  \neg +  \neg  (\beta  \neg +  \neg \k) |Y^0_t| \big)
          \neg +  \neg  \z_\th e^{\k T}    \wt{\G}^{2,n}_t     \big|\D_n f(t)   \big|   .
   \eeas
    Similarly to  \eqref{eq:e151b},   integration by parts  gives that
 \bea
    \wt{\G}^{2,n}_t     
          \le     I^{n,1}_T +   I^{n,2}_T      +     \z_\th e^{\k T}   D_T \neg \int_0^T  \neg \wt{\G}^{2,n}_s d K^n_s \neg
    + \neg \z_\th   \int_t^{T} \neg  D_s \wt{\G}^{2,n}_s e^{\wt{A}^n_s}  \wt{V}^n_s dB_s   ,    \q      t \in [0,T] ,
 \label{eq:e181b}
 \eea
 where $\int_0^\cd   D_s \wt{\G}^{2,n}_s e^{\wt{A}^n_s}  \wt{V}^n_s  dB_s $ is a
  uniformly integrable martingale, which can be shown by using similar arguments to those lead to \eqref{eq:e159}.
  And similar to \eqref{eq:e153}, the flat-off condition of $(Y^n, Z^n, K^n)$ implies that
 \bea  \label{eq:e185b}
    \int_0^T \neg \wt{\G}^{2,n}_s  dK^n_s     \le        \exp \left\{\g e^{2\k T}   Y^0_*  + \e \z_\th e^{\k T}   \right\}   K^n_T
     + \frac{1}{\e} \U_n      \Big(   \underset{t \in [0,T]}{\sup} |L^n_t - L^0_t |  \Big) K^n_T   ,  \q  \pas
 \eea

 For any $t \in [0,T]$,  taking  $E[\cd|\cF_t]$ in \eqref{eq:e185b} and \eqref{eq:e181b} yields that
 $ \dis   \wt{\G}^{2,n}_t     \le          \sum^4_{i=1} I^{n,i}_t$, \pas ~It then follows that
         \beas
           Y^n_t-  Y^0_t \le  (1-\th) |Y^n_t| + \th Y^n_t -   Y^0_t \le  (1-\th) |Y^0_t| + \hb{$\frac{1-\th  }{\g}$}    \ln  \left(  \sum^4_{i=1} I^{n,i}_t \right)     , \q \pas,
     \eeas
  which together with  \eqref{eq:e163b} as well as the continuity of processes $Y^n$, $Y^0$ and $\dis \sum^4_{i=1} I^{n,i}$  implies \eqref{eq:e191}.

  \ss  \no {\bf 2)}     For any $\d > 0$,  \eqref{eq:e191}, \eqref{eq:e155}, \eqref{eq:e157b},  Doob's martingale inequality and H\"older's inequality imply that
 \bea
   \qq   && \hspace{-1.5cm}  P \bigg(  \underset{t \in [0,T]}{\sup} |Y^n_t  \neg -\neg Y^0_t |   \ge  \d  \bigg)   \le
      P \Big( (1-\th) \big(Y^0_* + Y^n_* \big) \ge \d /2 \Big)   +     P  \neg \left( \hb{$\frac{1-\th  }{\g}$}
      \ln  \bigg(  \sum^4_{i=1} I^{n,i}_*   \bigg)  \ge \d /2 \right)  \label{eq:e193}  \\
       & &    \le         2   \hb{$\frac{1-\th}{\d}$} E\big[ Y^0_* + Y^n_* \big]
       +     \sum^4_{i=1} P  \neg \left(  I^{n, i}_* \ge  \hb{$\frac14$} e^{\frac{\d \g}{2(1-\th)}}  \right)
      \le    \hb{$\frac{1-\th}{\d \g}$} E\Big[   e^{2 \g    Y^0_*    } \neg + \neg e^{2 \g     Y^n_*  } \Big]
         \neg  +  \neg   4 e^{ \frac{ - \d \g}{2(1-\th)}} \sum^4_{i=1} E \Big[  I^{n, i}_T \Big]   \qq   \nonumber  \\
        &  &  \le    2  \hb{$\frac{1-\th}{\d \g}$}   \wt{\Xi} (2)    \neg + \neg  4 e^{\k T} e^{ \frac{-\d \g}{2(1-\th)}}
       C \bigg(  \big\| \eta_n  \big\|_{\hL^2(\cF_T)}  \neg+ \neg  \z_\th      \left\{\wt{\Xi} \Big(\hb{$\frac{8   }{1-\th}$}  e^{2\k  T} \Big)\right\}^{\frac14}
       \big\| \hb{$\int_0^T \neg    |\D_n f(s)   | ds $} \big\|_{\hL^4(\cF_T)}
         + 1  +    \z_\th e^{  \e \z_\th e^{\k T}}    \nonumber  \\
         &&     \q      +  \frac{\z_\th}{\e}  \left\{\wt{\Xi} \Big(\hb{$\frac{8   }{1-\th}$}  e^{2\k  T} \Big)\right\}^{\frac14}
         \big\|L^n \neg-\neg L^0\big\|_{\hC^4_\bF[0,T]}    \bigg) ,    \nonumber
        \eea
        with        $ C = 1+ \big\|  D_T     \big\|_{\hL^2(\cF_T)} + \underset{n \in \hN}{\sup}
           \bigg(         E \neg \left[       D_T  e^{\g e^{2\k T}  ( Y^0_*  \neg + \neg Y^n_*  )  }      \big(  K^0_T \neg +\neg  K^n_T \big) \right]
             \neg  +    \big\|    D_T         \big(  K^0_T + K^n_T \big) \big\|_{\hL^2(\cF_T)} \neg  \bigg)$.
        H\"older's inequality  and \eqref{eq:e155}  show that  $C $ is a finite constant.

     \ms   The convergence of $ \D_n f $ to $0$ and (S1) imply  that  \dtp
            \bea  \label{eq:e201}
              \lmt{n \to \infty} \D_n f(t,\o) =0 \q \hb{and} \q |\D_n f(t,\o)|   \le 2\a + 2 \beta Y^0_*(\o)+ \g \big|Z^0_t (\o) \big|^2, \q \fa n \in \hN.
            \eea
            Hence, for \pas ~ $\o \in \O$ we may assume that \eqref{eq:e201} holds for a.e. $t \in [0,T]$, and that $ Y^0_*(\o)+ \int_0^T  \big|Z^0_s (\o) \big|^2 ds < \infty$.   The Dominated convergence theorem then yields that $ \lmt{n \to \infty}  \int_0^T \big|\D_n f(s,\o)\big| ds =0 $.
            By (S2), it also holds \pas ~ that
             \beas
              \lmt{n \to \infty} \eta_n = e^{\g e^{2 \k T}|\xi_0| } \q \hb{and} \q
              \lmt{n \to \infty} \Big(\,\underset{t \in [0,T]}{\sup} |L^n_t - L^0_t | \, \Big)= 0    .
              \eeas
              Using
            \eqref{eq:e157a},   \eqref{eq:e157c} and \eqref{eq:e157d} with any $p>4$ shows that $\big\{ \eta_n^2 \big\}_{n \in \hN}$,
            $\left\{    \left( \int_0^T \neg       \big|\D_n f(s)   \big| ds \right)^4 \right\}_{n \in \hN}$
            and $\bigg\{    \underset{t \in [0,T]}{\sup} |L^n_t - L^0_t |^4 \bigg\}_{n \in \hN}$ are all uniformly integrable sequences in $\hL^1(\cF_T)$,
            which leads    to  that
            \beas
               \lmt{n \to \infty} E\big[ \eta_n^2 \big] = E \Big[  e^{2 \g e^{2 \k T}|\xi_0| }  \Big] \q \hb{and}   \q  \lmt{n \to \infty} E \left[  \left( \int_0^T \neg       \big|\D_n f(s)   \big| ds \right)^4 + \underset{t \in [0,T]}{\sup} |L^n_t - L^0_t |^4  \right] = 0  .
            \eeas

        Hence,  letting $n \to \infty$ in \eqref{eq:e193} and then letting $\e \to 0$ yield   that
           \beas
           \lsup{n \to \infty}  P \bigg(  \underset{t \in [0,T]}{\sup} |Y^n_t  \neg -\neg Y^0_t |   \ge  \d  \bigg)
           \le    2  \hb{$\frac{1-\th}{\d \g}$}   \wt{\Xi} (2)     \neg + \neg  4 e^{\k T} e^{\frac{- \d \g}{2(1-\th)}}
       C \left(    1 +   \big\| e^{ \g e^{2 \k T}|\xi_0| }  \big\|_{\hL^2(\cF_T)}        +    \hb{$\frac{\g e^{\k T}}{1-\th}$}   \right) .
           \eeas
   As $\th \to 1$, we obtain   $ \lmt{n \to \infty}  P \bigg(  \underset{t \in [0,T]}{\sup} |Y^n_t  \neg -\neg Y^0_t |   \ge  \d  \bigg) =0$,
    which implies that for any $p \in [1, \infty)$,  $    \exp\Big\{ p \g \cd \underset{t \in [0,T]}{\sup} |Y^n_t -Y^0_t|\Big\} $
           converges to $1$ in probability.

  \ss  \no {\bf 3)}       Fix $p \in [1, \infty)$.  Since
    $      E \bigg[ \exp\Big\{2  p \g \cd \neg \underset{t \in [0,T]}{\sup} |Y^n_t -Y^0_t|\Big\} \bigg] \le \frac12
                E \left[ e^{4  p \g  Y^n_*} +  e^{4  p \g Y^0_*} \right]  \le   \wt{\Xi} (4p)   $
              holds  for any $  n \in \hN$ thanks to \eqref{eq:e155},
                we see that $\bigg\{ \exp\Big\{ p \g \cd \underset{t \in [0,T]}{\sup} |Y^n_t -Y^0_t|\Big\}   \bigg\}_{n \in \hN} $
                is a uniformly integrable sequence in $\hL^1(\cF_T)$. Then it follows that
               $      \lmt{n \to \infty } E \left[ \exp\bigg\{ p \g \cd \underset{t \in [0,T]}{\sup} \big|Y^n_t -Y^0_t\big|\bigg\}\right] =1$,
               which in particular implies that
                     \bea  \label{eq:e211}
                     \lmt{n \to \infty } E \left[        \underset{t \in [0,T]}{\sup} |Y^n_t -Y^0_t|^q \right] =0, \q \fa q \in [1, \infty).
                     \eea
 For any $n \in \hN$, applying It\^o's formula to the process $ | Y^n - Y^0 |^2 $, we can deduce from (S1)  that
  \beas
      \int_0^T \neg  |Z^n_s-Z^0_s|^2 ds &\tneg =& \tneg  |\xi_n - \xi_0|^2 - |Y^n_0-Y^0_0|^2 +2 \int_0^T ( Y^n_s-Y^0_s )
   \big(f_n(s, Y^n_s, Z^n_s) - f_0(s, Y^0_s, Z^0_s)\big) \, ds \\
  & \tneg& \tneg     +2 \int_0^T ( Y^n_s-Y^0_s )     ( d K^n_s - d K^0_s)  - 2 \int_0^T ( Y^n_s-Y^0_s )   ( Z^n_s  -  Z^0_s)  \, d B_s   \\
   & \tneg \le& \tneg   2  \underset{t \in [0,T]}{\sup}
   |Y^n_t \neg -\neg Y^0_t|  \neg \left( 2 \a T \neg+\neg \beta T  \big( Y^n_* \neg+\neg     Y^0_* \big)
    \neg+\neg \frac{\g}{2} \int_0^T \neg \big( |Z^n_s|^2 \neg+\neg   |Z^0_s|^2  \big)  ds \neg + \neg  K^n_T \neg + \neg  K^0_T \right)    \\
   &\tneg &  \tneg      +  \underset{t \in [0,T]}{\sup} |Y^n_t \neg -\neg Y^0_t|^2  + 2 \left| \int_0^T ( Y^n_s-Y^0_s )
   ( Z^n_s  -  Z^0_s)  \, d B_s  \right|  , ~ \;  \pas
  \eeas
 Then  the Burkholder-Davis-Gundy inequality, H\"older's inequality, and  \eqref{eq:e155} imply that
   \beas
 && \hspace{-1cm}  E\left[ \left( \int_0^T |Z^n_s-Z^0_s|^2 ds \right)^p \, \right] \le
          c_p E\left[ \underset{t \in [0,T]}{\sup} |Y^n_t -Y^0_t|^{2p}\right] + c_p  E \left[ \underset{t \in [0,T]}{\sup} |Y^n_t -Y^0_t|^{p} \cd
          \left( \int_0^T   | Z^n_s  -  Z^0_s|^2  \, d s  \right)^{\frac{p}{2}} \, \right] \\
  && \q + c_p   \left\{ E\left[\underset{t \in [0,T]}{\sup} |Y^n_t \neg - \neg  Y^0_t|^{2p}\right]   \right\}^{\frac12} \left\{
    \underset{m \in \hN_0}{\sup}  E\left[    e^{2p \g  Y^m_*}      \neg  + \neg    \left(   \int_0^T  \neg   |Z^m_s|^2  ds \right)^{2p}
    \neg  + \neg     \big(K^m_T\big)^{2p}   \right]   \right\}^{\frac12}  \\
      && \le         c_p E\left[ \underset{t \in [0,T]}{\sup} |Y^n_t \neg -\neg Y^0_t|^{2p}\right]
      \neg +\neg  \frac12 E\left[ \left( \int_0^T \neg |Z^n_s \neg -\neg Z^0_s|^2 ds \right)^{\neg p} \, \right]
     \neg +\neg  c_p \sqrt{     \wt{\Xi}(2p)    }   \left\{ E\left[\underset{t \in [0,T]}{\sup} |Y^n_t \neg - \neg  Y^0_t|^{2p}\right]   \right\}^{\frac12}    .
  \eeas
    It is clear that $  E\left[ \left( \int_0^T |Z^n_s-Z^0_s|^2 ds \right)^p \, \right]  < \infty $ as $Z^n, Z^0 \in   \hH^{2, 2p}_\bF([0,T];\hR^d) $.
    Hence, it follows that
   \beas
    E\left[ \left( \int_0^T |Z^n_s-Z^0_s|^2 ds \right)^p \, \right]  \le   c_p E\left[ \underset{t \in [0,T]}{\sup} |Y^n_t \neg -\neg Y^0_t|^{2p}\right]
     \neg +\neg  c_p \sqrt{     \wt{\Xi}(2p)    }   \left\{ E\left[\underset{t \in [0,T]}{\sup} |Y^n_t \neg - \neg  Y^0_t|^{2p}\right]   \right\}^{\frac12} .
   \eeas
As $n \to \infty$, \eqref{eq:e211} implies that
 \bea  \label{eq:e215}
 \lmt{n \to \infty} E\left[ \left( \int_0^T |Z^n_s-Z^0_s|^2 ds \right)^p \, \right] =0 .
 \eea
 In particular, we have
 \bea  \label{eq:e215b}
 \lmt{n \to \infty} E  \int_0^T |Z^n_s-Z^0_s|^2 ds   =0 .
 \eea

 \ss  \no {\bf 4)} Let us further assume that \dtp,  $ f_n(t, \o, y, z)  $
 converges to  $f_0(t, \o, y, z) $ locally uniformly in $(y, z)$. By \eqref{eq:e211} and \eqref{eq:e215b}, $\big\{(Y^n,Z^n)\big\}_{n \in \hN}$ has a subsequence (we still denote it by $\big\{(Y^n,Z^n)\big\}_{n \in \hN}$) such that
  \bea  \label{eq:e219}
    \lmt{n \to \infty}   \underset{t \in [0,T]}{\sup} |Y^n_t -Y^0_t|    = 0, \q \pas  \q \hb{and} \q    \lmt{n \to \infty} Z^n_t = Z^0_t , \q \dtp
  \eea
  In fact, we can choose this subsequence so that
 $ Z^* \dfnn   \underset{n \in \hN}{\sup}|Z^n|      \in   \hH^2_\bF[0,T] $; see \cite{Lep_San_97} or \cite[Lemma 2.5]{Ko_2000}.
 Hence, except on a $dt \otimes dP$-null set of $[0,T] \times \O $, one may suppose the following statements hold:

 \ss \no  (\,\,i) $ \lmt{n \to \infty}  Y^n_t(\o) = Y^0_t(\o)$ and  $\lmt{n \to \infty} Z^n_t (\o)= Z^0_t(\o)$,

 \ss \no  (\,ii)  The mapping $  f_0(t, \o, \cd, \cd)$ is continuous,

 \ss \no  (iii)  For any compact subset $\sK$ of $\hR \times \hR^d$, $\lmt{n \to \infty} \Big( \underset{(y,z) \in \sK}{ \sup}   \left|  f_n\big(t, \o, y, z \big) -  f_0 \big(t,\o, y, z \big)  \right|    \Big) = 0$.

   \ms \no   Let $ \sK(t,\o)  \dfnn \Big\{(y,z) \in \hR \times \hR^d : | y| \le \underset{n \in \hN}{\sup} | Y^n_t(\o) | < \infty
  \hb{ and }  | z| \le  Z^*_t(\o)    < \infty   \Big\}$, which is clearly  a compact subset
  of $\hR \times \hR^d$.  Since
    \beas
 \qq \qq   &&  \hspace{-2cm} \big|f_n(t, \o, Y^n_t, Z^n_t) \neg - \neg  f_0(t, \o, Y^0_t, Z^0_t) \big|     \le  \big|f_n(t, \o, Y^n_t, Z^n_t)
   \neg  -  \neg  f_0(t, \o, Y^n_t, Z^n_t) \big| \neg + \neg  \big|f_0(t, \o, Y^n_t, Z^n_t)  \neg - \neg  f_0(t, \o, Y^0_t, Z^0_t) \big| \\
  && \le  \underset{(y,z) \in \sK(t,\o)}{\sup}\big|f_n(t, \o, y, z) - f_0(t, \o, y, z) \big| +
  \big|f_0(t, \o, Y^n_t, Z^n_t) - f_0(t, \o, Y^0_t, Z^0_t) \big|, \q \fa  n \in \hN,
 \eeas
    letting $n \to \infty$ yields that
    \bea  \label{eq:e229}
     \lmt{n \to \infty} f_n(t, \o, Y^n_t, Z^n_t) = f_0(t, \o, Y^0_t, Z^0_t) .
     \eea
       By (S1),  it also holds \dtp~ that
  \bea  \label{eq:e233}
   \q  \big|f_n(t,  Y^n_t, Z^n_t)- f_0(t,  Y^0_t, Z^0_t) \big|   \le  2\a + 2  \beta    \underset{m \in \hN_0}{\sup} Y^m_*
        + \frac{\g}{2} \Big(  \big|Z^*_t  \big|^2 +\big|Z^0_t  \big|^2    \Big),  \q  \fa n \in   \hN   ,
    \eea
    where  $ \underset{m \in \hN_0}{\sup} Y^m_* < \infty $, \pas ~thanks to \eqref{eq:e219}.
    Thus, for \pas ~ $\o \in \O$ we may assume that \eqref{eq:e229} and \eqref{eq:e233} hold for a.e. $t \in [0,T]$, as well as that
       $         \underset{m \in \hN_0}{\sup} Y^m_*(\o)     + \int_0^T  \Big(  \big|Z^*_s (\o) \big|^2 + \big|Z^0_s (\o) \big|^2  \Big)     ds < \infty $.
       The Dominated convergence theorem then yields that $ \lmt{n \to \infty}  \int_0^T \big|f_n(s, \o, Y^n_s, Z^n_s)- f_0(s, \o, Y^0_s, Z^0_s) \big|  ds =0 $.

  \ss Fix $p \in [1, \infty)$.    For any $n \in \hN$, (S1) and \eqref{eq:e155} shows that
          \beas
            E \left[ \left( \int_0^T \big|f_n(s, Y^n_s, Z^n_s) - f_0(s, Y^0_s, Z^0_s) \big|  ds \right)^{2p} \right] &  \tneg   \le &  \tneg
            c_p E \left[ \left( 2 \a T \neg+\neg \beta T \big( Y^n_* \neg+\neg    Y^0_* \big)
    \neg+\neg \frac{\g}{2} \int_0^T \neg \big( |Z^n_s|^2 \neg+\neg   |Z^0_s|^2  \big)  ds \right)^{2p} \right] \\
  & \tneg   \le  &  \tneg      c_p  \underset{m \in \hN_0}{\sup}  E\left[     e^{2p \g  Y^m_*}      \neg  + \neg    \left(   \int_0^T  \neg   |Z^m_s|^2  ds \right)^{2p}
        \right]   \le c_p  \,  \wt{\Xi}(2p) ,
            \eeas
            which implies that $\left\{ \left( \int_0^T \big|f_n(s, Y^n_s, Z^n_s) - f_0(s, Y^0_s, Z^0_s) \big|  ds \right)^p\right\}_{n \in \hN}$ is a uniformly integrable sequence in $\hL^1(\cF_T)$. Hence, it follows that
            \bea  \label{eq:e241}
           \lmt{n \to \infty}   E \left[ \left( \int_0^T \big|f_n(s, Y^n_s, Z^n_s) - f_0(s, Y^0_s, Z^0_s) \big|  ds \right)^p \,  \right] =0.
            \eea
            For any $n \in \hN$, it holds \pas ~ that
            \beas
  K^n_T-K^0_T = Y^n_0 -Y^0_0 - (\xi_n -\xi_0) - \int_0^T \big(f_n(s, Y^n_s, Z^n_s) - f_0(s, Y^0_s, Z^0_s)\big) \, ds
  +   \int_0^T   ( Z^n_s  -  Z^0_s)  \, d B_s.
  \eeas
 The Burkholder-Davis-Gundy inequality  then implies  that
   \beas
  E \Big[ \big| K^n_T-K^0_T \big|^p \Big]
 &\dneg   \le & \dneg    c_p  E \neg \left[ \underset{t \in [0,T]}{\sup} |Y^n_t \neg -\neg Y^0_t|^p \right]
   \neg+c_p  E \neg \left[ \left( \int_0^T \neg \big|f_n(s, Y^n_s, Z^n_s) \neg -\neg  f_0(s, Y^0_s, Z^0_s) \big|  ds \right)^p \, \right]  \\
   &\dneg    & \dneg   + \,    c_p E \left[ \left( \int_0^T |Z^n_s-Z^0_s|^2 ds \right)^{\frac{p}{2}}  \right],
   \eeas
 where $ E \left[ \left( \int_0^T |Z^n_s-Z^0_s|^2 ds \right)^{\frac{p}{2}}  \right] \le \left\{ E \neg \left[ \left( \int_0^T \neg  |Z^n_s-Z^0_s|^2 ds \right)^p  \right] \right\}^{\frac12} $ due to  H\"older's inequality. As $n \to \infty$,   \eqref{eq:e211}, \eqref{eq:e241} and \eqref{eq:e215} lead to \eqref{eq:e141}.  \qed

\section{An Obstacle Problem for PDEs.}

  \ms  In this section, we show that in the Markovian case,  quadratic RBSDEs with unbounded obstacles provide a probabilistic interpretation of solutions of some obstacle  problem for semi-linear parabolic PDEs, in which the non-linearity appears as the square of the gradient.

   \ms  For any $ t \in  [0, \infty)$,   $ B^t = \{B^t_s \dfnn B_{t+s}-B_t \}_{s \in [0,\infty)} $ is also a $d$-dimensional standard Brownian Motion on the probability space $(\O, \cF, P)$. Let $\bF^t $ be the augmented filtration generated by $B^t$, i.e.,
   \beas
 \bF^t \neg = \Big\{\cF^t_s \dfnn \si \Big(\si\big(B^t_r; r\in [0,s]\big) \cup \cN \Big) \Big\}_{s \in [0, \infty) }.
 \eeas

     Let $k \in \hN$,  $\k \ge 0$ and $\varpi \in [1,2)$. We consider the following functions:

    \ss  \no 1)     $b: [0,T] \times  \hR^k \to  \hR^k$ and $\si: [0,T] \times  \hR^k \to  \hR^{k \times d} $
     are two  continuous functions
     such  that $\si_*  \dfnn \underset{(t,x) \in [0,T] \times \hR^k}{\sup}|\si(t,x)|  \\ < \infty$, and that
   \bea  \label{b_si_Lip}
      |b(t,x)-b(t,x')|+|\si(t,x)-\si(t,x') | \le  \k    |x-x'|, \q \fa t \in [0,T],  ~ \fa x, x' \in \hR^k.
   \eea

      \ss  \no 2)  $h: \hR^k \to \hR $ and $l: [0,T] \times   \hR^k \to \hR$ are two  continuous functions such that
        \bea  \label{hl_varipi}
 l(T,x) \le h(x), \q  \fa  x \in \hR^k \q \hb{and} \q    |h(x)|  \vee   |l(t,x)|       \le \k \big( 1+|x|^\varpi \big)  , \q     \fa (t,x) \in [0,T] \times \hR^k.
   \eea

  \ss \no 3)    $f: [0,T] \times \hR^k \times \hR \times \hR^d \to \hR$ is a jointly continuous function   that satisfies

  \ss  \no \q    \,i)  There exist   $\a, \beta  \ge 0$ and $\g > 0$ such that
   for any $(t,x, z) \in [0,T] \times \hR^k    \times \hR^d  $  and $y, y' \in \hR$
   \bea   \label{cond:f_basic}
     |f(t,x,y,z)|   \le \a+\beta |y|+\frac{\g}{2} |z|^2    \q \hb{and} \q        |f(t,x,y,z)-f(t,x, y',z)| \le \k |y-y'| \, ;
   \eea
        \q   ii)  The mapping $z \to f(t,x,y,z)$ is
 \bea
         \bullet  &\tneg & \tneg \dneg  \hb{either {\it convex\,}  for   all $(t, x, y) \in   [0,T] \times \hR^k    \times \hR $,} \hspace{7.2cm}  \label{cond:f_convex}\\
      \bullet  &\tneg&\tneg \dneg \hb{ or {\it concave\,}    for all $(t, x, y) \in   [0,T] \times   \hR^k    \times \hR $.}  \label{cond:f_concave}
  \eea

  \ss  \no    For any $\l \ge 0$,    we let  $\wt{c}_\l$ denote a generic constant,  depending on $\l,  \a, \beta, \g, \k, \varpi, T,   \si_* $ and  on $b_0 \dfnn \underset{t \in [0,T]}{\sup} |b  (t, 0) | < \infty$, whose form may vary from line to line.

    \ms     Given $(t,x) \in  [0,T] \times \hR^k$, it is well-known that the SDE
  \bea  \label{FSDE}
     X_s= x + \int_t^s b(r, X_r) dr + \int_t^s \si(r, X_r) dB_r,    \q s \in [t, T]
  \eea
  admits a  unique  solution $\left\{ X^{t, x}_s \right\}_{s \in [t,T]}$, an $\hR^k$-valued 
  continuous process,   such that $X^{t, x}_s \in \cF^t_{s-t} \subset \cF_s$ for any $ s \in [t,T]$.
  In addition, we set $X^{t, x}_s \dfnn x $, $\fa s \in [0,t]$.

  \bs
  The following lemma gives an estimate for the exponential moments of process  $\big\{ |X^{t,x}_s|^\varpi \big\}_{s \in [t,T]}$\,.

\begin{lemm}   \label{lemm:esti_X}
Let $p \in [1, \infty)$. For any $(t,x) \in  [0,T] \times \hR^k$, we have
 \beas   
 E\left[ \exp\left\{ p \underset{s \in [t,T]}{\sup} \left| X^{t,x}_s \right|^\varpi \right\}\right]
 \le \wt{c}_p  \exp\left\{p \, 3^{\varpi-1}  e^{\k \varpi T}   |x|^\varpi \right\} .
 \eeas
\end{lemm}

 \ss \no {\bf Proof:}        One can deduce from \eqref{FSDE} and \eqref{b_si_Lip}  that \pas
   \beas
  \underset{s \in [t,t'\,]}{\sup}   \left| X^{t,x}_s \right| \le |x|+b_0 T+\k \int_t^{t'}  \underset{s \in [t, r]}{\sup}    \left| X^{t,x}_s \right|  dr
  +  \underset{s \in [t,T]}{\sup} \left| \int_t^s \si\left(r, X^{t,x}_r \right) dB_r \right|, \q  t' \in [t,T] .
   \eeas
     Then Gronwall's inequality implies that  \pas
    \beas
       \underset{s \in [t,t'\,]}{\sup}   \left| X^{t,x}_s \right|  \le e^{\k T}  \left( |x|+b_0 T+  \underset{s \in [t,T]}{\sup} \left| \int_t^s \si\left(r, X^{t,x}_r \right) dB_r \right| \right)  ,  \q   t' \in [t,T].
    \eeas
    Letting $t'=T$ and taking power of $\varpi$ yield that
     \bea  \label{eq:f121}
         \underset{s \in [t,T]}{\sup}   \left| X^{t,x}_s \right|^\varpi  \le  3^{\varpi-1}  e^{\k \varpi T} \left(  |x|^\varpi + ( b_0 T )^\varpi
         + \underset{s \in [t,T]}{\sup} \left| \int_t^s \si\left(r, X^{t,x}_r \right) dB_r \right|^\varpi \right)  , \q \pas
     \eea

   Let $\b1_{k \times d}$ denote the $k \times d$ matrix whose entries are all 1's.  We define an $\hR^k$-valued process
    \beas
        M_s \dfnn \int_0^s  \Big(\b1_{\{ t\le r \le T \}} \si\left( r, X^{t,x}_{r} \right) + \b1_{\{r > T \}} \b1_{k \times d}  \Big)\, dB_r  ,  \q \fa s \in [0, \infty).
        \eeas
     Given $i \in \{1,\cds,k\}$, it is clear that $M^i $ is an $\hR$-valued continuous  martingale such that $\lmt{s \to \infty}\lan M^i \ran_s =\infty$.  For any  $s \in [0, \infty)$, we define an $\bF$-stopping time
     $     \t^i_s \dfnn \inf\big\{ r \in [0, \infty) : \lan M^i \ran_r > s \big\}$.
     In light of the Dambis-Dubins-Schwarz Theorem (see, e.g., Theorem 3.4.6 of \cite{Kara_Shr_BMSC}),
   $W^i_s\dfnn M^i_{\t^i_s}$, $s \in [0, \infty)$ defines an $1$-dimensional standard Brownian Motion
   on the probability space $(\O, \cF, P)$  with respect to  the filtration $ \left\{  \cF_{\t^i_s } \right\}_{s \in [0, \infty)}$,
   and it holds \pas~that $\dis M^i_s=W^i_{\lan M^i\ran_s}$ for any $s \in [0, \infty)$.

     The convexity of function $y \to  e^{|y|^\varpi }$ on $\hR$ and Jensen's inequality imply that
  $ \Big\{ \neg \exp\big\{     \big| W^i_s \big|^\varpi  \big\} \neg \Big\}_{s \in [0, \infty)}$ is a continuous positive submartingale with respect to
  the filtration $\left\{  \cF_{\t^i_s } \right\}_{s \in [0, \infty)} $.  Applying Doob's Martingale Inequality, we obtain
     \bea
       E \left[ \underset{s \in [0,  \si^2_* T]}{\sup}  \bigg( \neg \exp\left\{     \left| W^i_s  \right|^\varpi  \right\} \neg \bigg)^{pk} \right]
   & \tneg \le &  \tneg \hb{$\big(\frac{pk}{pk-1}\big)^{pk}$}
   E \left[   \left(\exp\left\{     \left| W^i_{\si^2_* T}  \right|^\varpi  \right\} \right)^{pk} \right]      \nonumber \\
  &  \tneg = &  \tneg   \hb{$\big(\frac{pk}{pk-1}\big)^{pk}$}
  E \left[    \exp\left\{  pk  \big(\si_* \sqrt{T}\big)^\varpi  \left| W^i_1  \right|^\varpi  \right\}  \right]
 \le  \wt{c}_p \, E \left[    \exp\left\{   \frac14  \left| W^i_1  \right|^2  \right\}  \right].   \qq  \label{eq:f137}
   \eea
      As $W^i_1$ is a standard normal random variable under $P$, we have
   \bea  \label{eq:f139}
    E \left[    \exp\left\{   \frac14 \left| W^i_1  \right|^2  \right\}  \right]= \frac{1}{\sqrt{2\pi} }\int_{-\infty}^\infty e^{-\frac14  y^2} dy = \sqrt{2} .
    \eea

     For any $p \in (1, \infty)$,  since $ \dis   \lan M^i\ran_s = \int_t^s  \sum^d_{j=1}\big(\si_{ij}\big(  r, X^{t,x}_r \big) \big)^2 dr \le \si^2_*   T $ for any  $s \in [t, T]$,  one can deduce from \eqref{eq:f137}, \eqref{eq:f139} and H\"older's inequality that
   \bea
   \qq &&\hspace{-1.5cm}  E \left[   \exp\left\{ p  \underset{s \in [t,T]}{\sup} \left| \int_t^s \si\left(r, X^{t,x}_{r}\right) dB_r \right|^\varpi  \right\}  \right]
      =E \left[  \underset{s \in [t,T]}{\sup}  \exp\left\{ p  \left| M_s \right|^\varpi  \right\}  \right]
      \le E \left[  \underset{s \in [t,T]}{\sup}  \exp\left\{ p \sum^k_{i=1} \left| M^i_s \right|^\varpi  \right\}  \right]   \nonumber  \\
    &&= 
       E \left[ \underset{s \in [t, T]}{\sup}  \prod^k_{i=1}  \exp\left\{ p   \left| W^i_{\lan M^i\ran_s}  \right|^\varpi  \right\}  \right]
    \le   E \left[  \prod^k_{i=1} \left( \underset{s \in [t, T]}{\sup}   \exp\left\{ p   \left| W^i_{\lan M^i\ran_s}  \right|^\varpi  \right\} \right)  \right]   \nonumber \\
     && \le   E \left[ \prod^k_{i=1} \left(\underset{s \in [0,  \si^2_* T]}{\sup}  \exp\left\{ p   \left| W^i_s  \right|^\varpi  \right\}  \right) \right]
     \le \left\{ \prod^k_{i=1} E \left[   \underset{s \in [0,  \si^2_* T]}{\sup}  \exp\left\{ p k  \left| W^i_s  \right|^\varpi  \right\}    \right] \right\}^{\frac{1}{k}}
     \le \wt{c}_p   ,      \label{eq:f131}
     \eea
     where we used in the first inequality the fact that $\dis   |x|^\varpi =\bigg(\sum^k_{i=1} |x^i|^2\bigg)^{\frac{\varpi}{2}}
     \le \sum^k_{i=1} \big(|x^i|^2\big)^{\frac{\varpi}{2}}     =\sum^k_{i=1} |x^i|^\varpi$ for any $x \in \hR^k$.
  Plugging it back into   \eqref{eq:f121} proves the lemma.      \qed

 \ms Our objective in this section is to find a unique viscosity solution of the following obstacle problem for semi-linear parabolic PDEs:
 \bea\label{eq:PDE}
  \left\{\ba{l}
  \dis  \dneg \min\left\{ \neg (u \neg- \neg l)(t,x),  -\frac{\pa u}{\pa t}   (t,x)- \neg \cL u(t,x)    - \neg f \neg
   \left(t,x, u(t,x),    (\si^T \neg \cd \neg  \nabla_x u  )     (t,x)\right) \neg \right\} =0, ~ \fa (t,x) \in (0,T) \times \hR^k, \vspace{2mm} \\
   u(T,x)=h(x) , ~  \fa x \in \hR^k,
   \ea
   \right.
 \eea
    where  $\si^T$ denotes the transpose of $\si$ and
     $  \cL u(t,x)  \dfnn \frac12  \hb{trace} \big(  ( \si  \si^T  D^2_x u  ) (t,x) \big) + \lan b(t,x), \nabla_x u(t,x) \ran$.

  \ms  Now let us consider the  obstacle problem for   PDEs in a  more general form:
   \bea\label{eq:PDE_gen}
  \left\{\ba{l}
  \dis  \dneg \min\left\{ \neg (u \neg- \neg \mathfrak{l})(t,x),  -\frac{\pa u}{\pa t}   (t,x)-F\left(t,x,   u(t,x),  \nabla_x u   (t,x), D^2_x u(t,x) \right)  \neg \right\} =0, ~\; \fa (t,x) \in (0,T) \times \hR^k, \vspace{2mm} \\
   u(T,x)=\mathfrak{h}(x) , ~  \fa x \in \hR^k,
   \ea
   \right.
 \eea
 where $\mathfrak{h}: \hR^k \to \hR $,  $\mathfrak{l}: [0,T] \times   \hR^k \to \hR$, and $F: [0,T] \times   \hR^k \times \hR  \times   \hR^k
  \times \hS^k \to \hR$ are all (jointly)  continuous functions with  $\hS^k$ denoting the set of all real symmetric $k \times k$ matrices.

    \begin{deff}
  A function $u \in C([0, T] \times \hR^k )$    is called  a viscosity subsolution \(resp. viscosity supersolution\) of \eqref{eq:PDE_gen} if   $u(T, x) \le $ \(resp. $\ge$\) $ \mathfrak{h}(x)$,   $\fa x \in   \hR^k$,
  and if    for any $(t_0,x_0, \vf ) \in (0,T) \times \hR^k \times C^{1,2}\big([0,T] \times \hR^k\big)$   such that $u(t_0,x_0) = \vf (t_0,x_0)$ and that
  $u-\vf$ attains a local  maximum (resp. local  minimum) at $(t_0,x_0)$, we have
 \beas  
  ~   \min\left\{ \neg (u \neg - \neg \mathfrak{l}) (t_0,x_0   ), - \frac{\pa \vf }{\pa t}(t_0,x_0   )-\neg
   F\Big(t_0,x_0,   u(t_0,x_0),  \nabla_x  \vf   (t_0,x_0), D^2_x  \vf(t_0,x_0) \Big)
   \neg \right\} \le (\hb{resp.}\, \ge) ~0. \q
 \eeas
   A function $u \in C([0, T] \times \hR^k )$  is called a viscosity solution of \eqref{eq:PDE_gen} if it is both a viscosity subsolution and a viscosity supersolution of \eqref{eq:PDE_gen}.
  \end{deff}

   One can alternatively define  viscosity subsolutions/supersolutions of \eqref{eq:PDE_gen} in term of second-order superjets/subjets (see \cite{User_1992}).
  \begin{deff} 1\) For a function $u :   [0, T] \times \hR^k  \to \hR$, its second-order superjet \(resp. subjet\) at some $(t_0,x_0) \in (0, T) \times \hR^k$, denoted by $\cP^{2,+}u(t_0,x_0)$ \big(resp. $\cP^{2,-} u(t_0,x_0) $\big),
  is a collection of all triplets $ (p,q,W) \in \hR \times \hR^k \times   \hS^k$ such that as $(t,x) \to (t_0,x_0)$ in $(0, T) \times \hR^k$,
     \beas
     u(t,x) \le (\hb{resp.\,} \ge) ~ u(t_0,x_0) + p(t-t_0) + \lan q, x - x_0 \ran + \frac12 \big\lan W(x-x_0), x-x_0 \big\ran + o \big(|t-t_0|+|x-x_0|^2\big) .
   \eeas

  \ms \no 2\) For a function $u :   [0, T] \times \hR^k  \to \hR$ and some $(t_0,x_0) \in (0, T) \times \hR^k$, we define
   $\ol{\cP}^{2,+}u(t_0,x_0)$ \(resp. $\ol{\cP}^{2,-} u(t_0,x_0) $\) as the
   collection of all triplets $ (p,q,W) \in \hR \times \hR^k \times   \hS^k$ such that for some sequence $\big\{(t_n,x_n,p_n,q_n,  \\  W_n)\big\}_{n \in \hN}
    \subset  (0, T) \times \hR^k \times \hR \times \hR^k \times   \hS^k$,
     \beas
     &&  ( p_n,q_n,W_n) \in \cP^{2,+}u(t_n,x_n)\; \big(\hb{resp.\,} \cP^{2,-} u(t_n,x_n)\big), ~\fa n \in \hN \\
       \hb{and}  &&      \big(t_0,x_0, u(t_0,x_0), p,q,W \big)= \lmt{n \to \infty}  \big(t_n,x_n, u(t_n,x_n), p_n,q_n,W_n \big) .
     \eeas
  \end{deff}

   \begin{deff} \label{def:viscosity_alternative}
      A function $u \in C([0, T] \times \hR^k )$ is called  a viscosity subsolution \(resp. viscosity supersolution\) of \eqref{eq:PDE_gen} if   $u(T, x) \le $ \(resp. $\ge$\) $ \mathfrak{h}(x)$,   $\fa x \in   \hR^k$,
      and if    for any $(t_0,x_0 ) \in (0,T) \times \hR^k  $ and  $ (p,q,W) \in \ol{\cP}^{2,+}u(t_0,x_0)$ \big(resp. $\ol{\cP}^{2,-} u(t_0,x_0) $\big), we have
 \beas
    ~\;  \min\Big\{ \neg (u \neg - \neg \mathfrak{l}) (t_0,x_0   ), - p   -  F\big(t_0,x_0,   u(t_0,x_0),  q , W \big)
    \neg \Big\} \neg \le      (\hb{resp.}\; \ge)    ~  0.
 \eeas
  A function $u \in C([0, T] \times \hR^k )$   is called a viscosity solution of \eqref{eq:PDE_gen} if it is both a viscosity subsolution and a viscosity supersolution of \eqref{eq:PDE_gen}.
   \end{deff}

        \ss    For any $(t,x) \in  [0,T] \times \hR^k$, let $ \sP^{\,t}$ denote the $\bF^t$-progressively measurable $\si$-field on $ [0,T-t] \times \O$.
        Since $\wt{X}^{t,x}_s \dfnn    X^{t,x}_{t+s} $, $s \in [0,T-t]$ is an $\bF^t$-adapted continuous process,
          the joint continuity of $f$ implies that
         \beas
         \tilde{f}^{t,x}(s,\o,y,z )\dfnn   f \Big(t+s,  \wt{X}^{t,x}_s  (\o), y,z \Big)  ,  \q    \fa (s,\o,y,z ) \in   [0,T-t] \times \O \times \hR \times \hR^d
          \eeas
          is a  $\sP^{\,t} \times \sB(\hR) \times   \sB(\hR^d)/\sB(\hR)$-measurable function, namely,  it is  a  generator with respect to $\bF^t$ over the period $[0,T-t]$.
          By \eqref{cond:f_basic}-\eqref{cond:f_concave},  $\tilde{f}^{t,x}$ also satisfies (H1)-(H3).
          On the other hand,   \eqref{hl_varipi}
           shows that $\big\{\wt{L}^{t, x}_s  \dfnn  l \big(t+s, \wt{X}^{t, x}_s   \big) \big\}_{s \in [0,T-t]}$  is also
       an $\bF^t$-adapted continuous process such that $\wt{L}^{t, x}_{T-t} = l \big( T, X^{t, x}_{T}   \big)  \le h\big( X^{t, x}_{T} \big) \in \cF^t_{T-t}$.
       For any $p \in [1, \infty)$,  \eqref{hl_varipi} and Lemma \ref{lemm:esti_X} imply that
        \bea
         && \hspace{-1.5cm}  E \left[ \exp \left\{ p \Big( \big|h\big( X^{t, x}_T \big) \big|   \vee   \wt{L}^{t, x}_* \Big) \right\} \right]
           \le     E \left[ \exp \left\{ p\k \bigg(1 + \underset{s \in [t,T]}{\sup}  \big| X^{t, x}_s  \big|^\varpi  \bigg) \right\} \right] \nonumber   \\
        &&  \le   e^{p \k} E \left[ \exp \left\{ (1 \vee p\k) \underset{s \in [t,T]}{\sup}  \big| X^{t, x}_s  \big|^\varpi   \right\} \right]
        \le    \wt{c}_p \exp\left\{(1 \vee p\k) \, 3^{\varpi-1}  e^{\k \varpi T}   |x|^\varpi \right\}    .   \label{eq:f147}
        \eea
        Hence, Corollary \ref{cor:unique}  shows that the quadratic RBSDE$\Big(h ( X^{t, x}_{T}  ) , \tilde{f}^{t,x}, \wt{L}^{t, x} \Big)$
        with respect to $B^t$ over the period $[0,T-t]$ admits a  unique solution $ \Big(\wt{Y}^{t,x}, \wt{Z}^{t,x}, \wt{K}^{t,x} \Big)$ in $ \underset{p \in [1, \infty)}{\cap} \hS^p_{\bF^t}[0,T-t] $.

      \ms   The continuity of process $\big\{X^{t,x}_s\big\}_{s \in [0,T]}$ and \eqref{cond:f_basic}-\eqref{cond:f_concave} imply  that
      $$      ~ \; f^{t,x}(s, \o, y,z) \dfnn  \b1_{\{s \ge t\}}  \tilde{f}^{t,x}(s-t, \o, y,z)
      = \b1_{\{s \ge t\}}  f \big(s, X^{t,x}_s(\o),y,z \big),  ~    \fa (s,\o,y,z ) \in   [0,T] \times \O \times \hR \times \hR^d $$
       is a  $\sP \times \sB(\hR) \times   \sB(\hR^d)/\sB(\hR)$-measurable function that   
       satisfies (H1)-(H3) with the same   constants $\a, \beta,   \k  \ge 0$ and $\g > 0$ as $f$.
       Let   $    L^{t,x}_s  \dfnn    \wt{L}^{t,x}_{(s-t)^{\neg +}}        =  l \big(s\vee t , X^{t,x}_{s\vee t}\big)$ , $s \in [0,T] $,
       which is clearly an $\bF$-adapted continuous process  with $L^{t,x}_T  =
     \wt{L}^{t,x}_{T-t} \le h\big( X^{t, x}_{T} \big)$.  Then  one can show that
        \beas
         \big( Y^{t,x}_s , Z^{t,x}_s , K^{t,x}_s  \big)
  \dfnn   \left(  \wt{Y}^{t,x}_{(s-t)^{\neg +}}\,, \b1_{\{s \ge t\}} \wt{Z}^{t,x}_{s-t}\,, \b1_{\{s \ge t\}} \wt{K}^{t,x}_{s-t} \right) , \q   s \in [0,T]
        \eeas
        satisfies the quadratic RBSDE$\Big(h ( X^{t, x}_{T}  ) ,  f^{t,x},  L^{t, x} \Big)$ over the period $[0,T]$, and that $\big( Y^{t,x} , Z^{t,x} , K^{t,x}  \big)
        \in  \underset{p \in [1, \infty)}{\cap} \hS^p_\bF[0,T]$. Since $ E \left[ \exp \left\{ p \Big( \big|h\big( X^{t, x}_{T} \big) \big|   \vee   L^{t, x}_* \Big) \right\} \right]
         < \infty$ by  \eqref{eq:f147},  Corollary \ref{cor:unique} again shows that $ \big( Y^{t,x} , Z^{t,x} ,     K^{t,x}  \big) $
         is the unique solution of  the quadratic RBSDE$\Big(h ( X^{t, x}_{T}  ) ,  f^{t,x},  L^{t, x} \Big)$ in $ \underset{p \in [1, \infty)}{\cap} \hS^p_\bF[0,T]$.

       \ss  The main objective of this section is to demonstrate that
        \bea  \label{def_utx}
                 u(t,x) \dfnn   \wt{Y}^{t,x}_0   =   Y^{t,x}_t   ,   \q  \fa (t,x) \in [0,T] \times \hR^k
                 \eea
                 is a  viscosity solution of   \eqref{eq:PDE}. First, we recall a well-known moment estimate of diffusion process $X^{t,x}$, without proof, in order
                 to show that  $u$ is a continuous function.

  \begin{lemm}
   For any $(t,x), (t',x') \in [0,T] \times \hR^k$, we have
   \bea
     E\left[  \underset{r \in [t,s]}{\sup}  \big| X^{t,x}_{r} - x  \big|^2 \right]   &\le&  \wt{c}_0 \big(1+|x|^2\big) (s \neg - \neg t) ,
     \q \fa s \in [t, T],     \label{eq:esti_X_1} \\
     E \neg \left[  \underset{s \in [\tilde{t} \vee t,T]}{\sup}  \big| X^{\tilde{t},\tilde{x}}_s \neg - \neg X^{t,x}_s  \big|^2 \right]
          & \le & \wt{c}_0     E \neg \left[     \big| X^{\tilde{t},\tilde{x}}_{\tilde{t} \vee t} \neg - \neg X^{t,x}_{\tilde{t} \vee t}  \big|^2 \right] .   \label{eq:esti_X_2}
   \eea
  \end{lemm}

\begin{prop}  \label{prop:u_continuous}
   The function $u $ defined in \eqref{def_utx} is  a continuous one such that $  |u(t, x)|  \le \wt{c}_0 \big(1 + |x|^\varpi\big) $
   for any $ (t,x) \in [0,T] \times \hR^k$.
\end{prop}

\ss \no {\bf Proof:} Given $(t,x) \in [0,T] \times \hR^k $,  we let
 $\left\{(t^0_n, x^0_n)\right\}_{n \in \hN} \subset [0,T] \times \hR^k$ be an arbitrary sequence that converges to $(t,x)$.
   Without loss of generality, we assume that     $ \{x^0_n \}_{n \in \hN} \subset \ol{B}_1(x)   \dfnn \{\tilde{x} \in \hR^k: |\tilde{x}-x| \le 1 \}$.
   To see  $    \lmt{n \to \infty} u\big( t^0_n,x^0_n\big)  = u(t,x)$, we only need to show that any subsequence
         $\left\{u(t_n, x_n)\right\}_{n \in \hN} $  of $\left\{u(t^0_n, x^0_n)\right\}_{n \in \hN} $ has in turn a subsequence that converges  to $u(t,x)$.
         For any $n \in \hN$,   \eqref{eq:esti_X_2} shows that
   \bea  \label{eq:f161}
   E \neg \left[  \underset{s \in [ t_n \vee t, T]}{\sup}  \big|  X^{t_n, x_n}_s  \neg - \neg   X^{t,x}_s  \big|^2 \right]
    \le \wt{c}_0 E \neg \left[     \big| X^{t_n, x_n}_{t_n \vee t} \neg - \neg X^{t,x}_{t_n \vee t}  \big|^2 \right]
               \le \wt{c}_0   E \neg \left[  \underset{s \in [t_n \land t, t_n \vee t]}{\sup}  \big|  X^{t_n, x_n}_s  \neg - \neg   X^{t,x}_s   \big|^2 \right]     .
  \eea
  When  $t_n \le t  $, \eqref{eq:esti_X_1} implies that
 \bea
    E \neg \left[  \underset{s \in [t_n \land t, t_n \vee t]}{\sup}  \big|  X^{t_n, x_n}_{s}  \neg - \neg   X^{t,x}_{s}   \big|^2 \right]
     & \tneg = & \tneg   E \neg \left[  \underset{s \in [t_n, t]}{\sup}  \big|  X^{t_n,x_n}_{s}  -x  \big|^2 \right]
     \le   2 |x_n-x|^2 +  2  E \neg \left[  \underset{s \in [t_n, t]}{\sup}  \big| X^{t_n,x_n}_{s} \neg - x_n      \big|^2 \right]  \nonumber     \\
     & \tneg \le & \tneg    2 |x_n-x|^2 +  \wt{c}_0   \big(1+|x_n |^2 \big) (t \neg - \neg t_n) .      \label{eq:f163}
   \eea
   Similarly, when $t_n > t  $,
 \bea
    E \neg \left[  \underset{s \in [t_n \land t, t_n \vee t]}{\sup}  \big|  X^{t_n, x_n}_{s}  \neg - \neg   X^{t,x}_{s}   \big|^2 \right]
      \le    2 |x_n-x|^2 +    \wt{c}_0   \big(1+|x_n |^2 \big) (t_n \neg - \neg t)   .     \label{eq:f165}
 \eea
 Since   $X^{t_n, x_n}_s  \neg - \neg   X^{t,x}_s= x_n-x$ for any $  s \in [0, \  t_n \land t]$,
 \eqref{eq:f161} and \eqref{eq:f163} (or \eqref{eq:f165})   imply that
    \beas
    E \neg \left[  \underset{s \in [0, T]}{\sup}  \big|  X^{t_n, x_n}_s  \neg - \neg   X^{t,x}_s   \big|^2 \right] \le \wt{c}_0 |x_n-x|^2
    + \wt{c}_0   \big(1+|x_n |^2 \big) |t_n \neg - \neg t| \to 0,  \;  \hb{ as } n \to \infty.
    \eeas
 Hence, we can extract a subsequence of $\{(t_n,x_n)\}_{n \in \hN}$  \big(we still denote it by $\{(t_n,x_n)\}_{n \in \hN}$\big)
 such that except on a $P$-null set $\sN$,
     \bea \label{eq:f169}
  \lmt{n \to \infty} \Big( \underset{s \in [0, T]}{\sup}  \big|  X^{t_n, x_n}_s    \neg - \neg   X^{t,x}_s    \big| \Big)=0 \q \hb{and} \q \hb{the path $s \to  X^{t,x}_s$ is continuous.}
      \eea

   \ss  To apply Theorem \ref{thm:stable} to the sequence $\big\{\big( Y^{t_n,x_n} , Z^{t_n,x_n} , K^{t_n,x_n}  \big)\big\}_{n \in \hN}$,
 let us check the assumptions of this theorem first. We have seen that $f^{t,x}$ together with   $\big\{f^{t_n,x_n}\big\}_{n \in \hN} $ satisfy (S1).

    \ms Fix $\o \in \sN^c$.  For any $\e > 0$, the continuity of functions $h$ and $l$   assures that there exists a $\d (\o)  \in (0,1) $ such that
  \beas
     |h(\tilde{x}) \neg - \neg h(x') | \vee  |l(\tilde{s}, \tilde{x}) \neg -\neg l(s', x') | < \e, ~  \fa (\tilde{s},\tilde{x}), (s',x') \in [0,T] \times \sD (\o)
     \hb{\;with\;} |\tilde{s}-s'|^2 \neg +    |\tilde{x}-x'|^2      < \d^2 (\o) ,
  \eeas
    where $ \sD (\o) \dfnn \Big\{\tilde{x} \in \hR^k: |\tilde{x}| \le 1 + \underset{s \in [0, T]}{\sup}  \big|     X^{t,x}_s (\o)   \big|< \infty \Big\}$.
      Moreover, in light of \eqref{eq:f169}, there exists  an $ N(\o ) \in \hN$ such that for any $n \ge N(\o )$,
   \beas
     |t_n-t|  \vee \underset{s \in [0,T]}{\sup} \big|  X^{t_n, x_n}_s (\o)-    X^{t,x}_s (\o)   \big| < \frac{\d(\o)}{2}    .
     \eeas
     Then for any $n \ge N(\o ) $,  one can deduce the following statements:

     \ss   \no $\bullet$  $\big| h\big(X^{t_n,x_n}_T(\o)\big) -h\big(X^{t,x}_T(\o)\big)  \big| < \e$,

      \ss \no $\bullet$  for any $s \in [t_n \vee t,T]$,
     $    \big| L^{t_n,x_n}_s(\o)-L^{t,x}_s(\o) \big| =  \big|    l \big(s, X^{t_n,x_n}_s(\o)\big)-l \big(s, X^{t,x}_s(\o)\big) \big| < \e$;

    \ss \no $\bullet$ for any $s \in [0, t_n \vee t]$, if $t_n \le t$,
      $ \big| L^{t_n,x_n}_s(\o)-L^{t,x}_s(\o) \big|  =   \big| l \big(s \vee t_n, X^{t_n,x_n}_{s \vee t_n}(\o) \big)-l \big(t, X^{t,x}_t(\o) \big) \big|
     =  \big| l \big(s \vee t_n, X^{t_n,x_n}_{s \vee t_n}(\o) \big)-l \big(t, X^{t,x}_{s \vee t_n}(\o) \big) \big| < \e  $ since $s \vee t_n \in [t_n,t] $;
      on the other hand, if $t_n > t$, one can similarly deduce that $ \big| L^{t_n,x_n}_s(\o)-L^{t,x}_s(\o) \big| < \e  $.

   \ss \no     Thus   (S2) is satisfied.

   \ms For any $p \in [1, \infty)$,  we have seen from \eqref{eq:f147} that
    \beas
      \qq     E \left[ \exp \left\{ p \,  \Big( \big|h\big( X^{t, x}_{T} \big) \big|   \vee   L^{t, x}_* \Big) \right\} \right]
          =   E \left[ \exp \left\{ p \, \Big( \big|h\big( X^{t, x}_{T} \big) \big|   \vee   \wt{L}^{t, x}_* \Big) \right\} \right]
           \le       \wt{c}_p \exp\left\{   (1 \vee p\k) \, 3^{\varpi-1}  e^{\k \varpi T}   |x|^\varpi \right\}       .
        \eeas
           Similarly, it holds for any $n \in \hN$ that
         \beas
          E \left[ \exp \left\{ p \, \Big( \big|h\big( X^{t_n, x_n}_{T} \big) \big|   \vee   L^{t_n, x_n}_* \Big) \right\} \right]
         &  \le &   \wt{c}_p \exp\left\{(1 \vee p\k) \, 3^{\varpi-1}  e^{\k \varpi T}   |x_n|^\varpi \right\}  \\
         &   \le  &  \wt{c}_p \exp\left\{(1 \vee p\k) \, 3^{\varpi-1}  e^{\k \varpi T}  \big(1+ |x|\big)^\varpi \right\}      .
        \eeas
        Thus (S3) also holds.

        \ms  Given $(s, \o, y,z) \in [0,T] \times \sN^c \times \hR \times \hR^d$, it holds  for any $n \in \hN$ that
            \beas
        \big|f^{t_n,x_n}(s, \o, y, z) - f^{t,x}(s, \o, y, z)\big|
     & \dneg\le &  \dneg  \big| f \big(s, X^{t_n,x_n}_s(\o), y, z\big) -  f \big(s, X^{t,x}_s(\o), y, z\big)    \big|  \q \\
    &&     + \big|\b1_{\{s \ge t_n\}} - \b1_{\{s \ge t\}} \big| \cd  \big| f \big(s, X^{t,x}_s(\o), y, z\big)    \big|   .
  \eeas
  As $n \to \infty$,  the continuity of $f$ and \eqref{eq:f169} imply that    $        \lmt{n \to \infty}    f^{t_n,x_n}(s, \o, y, z) =  f^{t,x}(s, \o, y, z)$,
    which in particular shows that for \dtp ~ $(s, \o) \in [0,T] \times \O$,   $ f^{t_n,x_n}\big(s, \o, Y^{t,x}_s(\o), Z^{t,x}_s(\o)\big)  $
    converges to  $f^{t,x} \big(s, \o, Y^{t,x}_s(\o), Z^{t,x}_s(\o)\big) $.

        \ms   Now, applying Theorem \ref{thm:stable} yields that $ \lmt{n \to \infty} E \left[ \exp\bigg\{   \underset{s \in [0,T]}{\sup} \big|Y^{t_n,x_n}_s -Y^{t,x}_s \big|\bigg\}  \right] = 1 $, which allows us to extract   a  subsequence
        $\left\{(t_{n_i}, x_{n_i})\right\}_{i \in \hN} $ from $\left\{(t_n, x_n)\right\}_{n \in \hN} $ such that
        $ \lmt{i \to \infty} \,  \underset{s \in [0,T]}{\sup} \big|Y^{t_{n_i},x_{n_i}}_s -Y^{t,x}_s \big|  = 0 $, \pas~ In particular,
        one has
         \beas
              \lmt{i \to \infty} u\big( t_{n_i},x_{n_i}\big) =    \lmt{i \to \infty}  \wt{Y}^{t_{n_i},x_{n_i}}_0  =     \lmt{i \to \infty}  Y^{t_{n_i},x_{n_i}}_0 = Y^{t,x}_0  = \wt{Y}^{t,x}_0 = u(t,x),
         \eeas
          which shows that $u$ is a continuous function. Moreover, Corollary \ref{cor:unique}, \eqref{hl_varipi} and Lemma \ref{lemm:esti_X} imply that
             \beas
     \hspace{1cm} - \k \big(1+|x|^\varpi \big)  & \le&  l(t, x)    =   \wt{L}^{t,x}_0    \le     u(t,x)  =   \wt{Y}^{t,x}_0  \le \wt{c}_0 + \frac{1}{\g} \ln E \bigg[ \exp \left\{\g e^{\beta T} \big( \big|h(X^{t,x}_T)\big|
           \vee  \wt{L}^{t,x}_*\big)\right\} \bigg]  \\
           &\le& \wt{c}_0 + \frac{1}{\g} \ln E \left[ \exp \left\{\g \k e^{\beta T}    \underset{s \in [t,T]}{\sup} |X^{t,x}_s|^\varpi \right\} \right]
           \le \wt{c}_0 \big(1+|x|^\varpi\big).  \hspace{3.4cm} \hb{\qed}
           \eeas

  \ms For any  $ \xi  \in \cO^{t,x} \dfnn \big\{ \xi  \in \hL^0(\cF_T):    \xi \ge L^{t,x}_T , \pas
  \hb{ and }    E \big[   e^{ p  \,   \xi^+     }   \big] < \infty, ~\fa  p \in (1, \infty) \big\}$, Corollary  \ref{cor:unique} guarantees  a unique solution
  $ \big( Y^{t,x,\xi} , Z^{t,x, \xi} , K^{t,x, \xi}  \big) $ of  the quadratic RBSDE$\Big(\xi ,  f^{t,x},  L^{t, x} \Big)$
  in $ \underset{p \in [1, \infty)}{\cap} \hS^p_\bF[0,T]$. For each $s \in [0,T]$, we can regard $  \cE^{t,x}[\xi|\cF_s] \dfnn Y^{t,x,\xi}_s$, $\xi \in\cO^{t,x}$
  as a nonlinear conditional expectation on $\cO^{t,x}$ with respect to $\cF_s$ \big(cf. $g$-expectations in the case of BSDEs, see e.g. \cite{Peng-97}, \cite{MaYao_2010}, Subsection 5.4 of \cite{OSNE2} and Section \ref{sec:g-evaluation} of the current paper\big).
    Then the diffusion $X^{t,x}$ has the following Markov property under $\cE^{t,x}$:

  \begin{prop} \label{prop:u_markov}
  Let $u$ be the function defined  in \eqref{def_utx}. For any $(t,x) \in [0,T] \times \hR^k $ it holds \pas ~ that
    \bea  \label{eq:f201}
        u \big( s   , X^{t,x}_s  \big)      = Y^{t,x}_s   =  \wt{Y}^{t,x}_{s-t} ,  \q s \in [t,T] .
    \eea
  \end{prop}

 \ss \no {\bf Proof:}    {\bf 1)}  We fix $s \in [\,t,T]$ and denote  $ \Th^0_{t'} \dfnn \Th^{t,x}_{t'}$, $  t'  \in [s, T]$ for $\Th=X, Y,Z,K$.
  Given $n \in \hN$, there exist a finite subset $\big\{x^n_i \big\}_{ i=1}^{ j_n} $ of $B_{2^n}(0) \dfnn \{x \in \hR^k: |x| < 2^n \}$
    and a disjoint partition $ \big\{\cI^n_i  \big\}_{ i=1}^{ j_n} $ of $B_{2^n}(0) $ such that $x^n_i \in \cI^n_i \in \sB(\hR^k)$
    and $\cI^n_i \subset \ol{B}_{2^{-n}} (x^n_i)$ 
          for $i=1, \cds ,j_n$. Let
 \beas
      \cA^n_i \dfnn \left\{   X^0_{s}  \in \cI^n_i \right\} \in \cF_s , \q \fa i =1, \cds ,j_n \q \hb{and} \q \cA^n_0 \dfnn\big\{  X^0_{s} \in B^c_{2^n}(0)
         \big\} \in \cF_s .
 \eeas
        For any $t' \in [s, T]$ and $\Th=X, Y,Z,K$,  we define
   $\dis  \Th^n_{t'} \dfnn  \sum_{i=0}^{j_n}  \b1_{ \cA^n_i}  \Th^{s, x^n_i}_{t'} \in \cF_{t'}$ with $x^n_0 \dfnn 0$.         Then   for any $i=0, \cds ,j_n$,
    \beas
   \b1_{ \cA^n_i}  X^{s,   x^n_i}_{t'} &=&   x^n_i  \b1_{ \cA^n_i} +  \int_{s}^{t'}   \b1_{ \cA^n_i}  b\Big( r,   X^{s, x^n_i}_r\Big) dr
  + \int_{s}^{t'} \b1_{ \cA^n_i}  \si \Big(  r,    X^{s, x^n_i}_r  \Big)  dB_r     \\
    &=&  x^n_i \b1_{ \cA^n_i} + \int_{s}^{t'}  \b1_{ \cA^n_i}  b \big(r,    X^n_r\big) dr
  +\int_{s}^{t'}  \b1_{ \cA^n_i}  \si \big(r,    X^n_r\big) \,  dB_r , \q \pas ;
  \eeas
  and that
  \beas
     & \tneg &  \hspace{-0.5cm}    \b1_{ \cA^n_i}  l \big(t',  X^n_{t'}   \big) = \b1_{ \cA^n_i}  l \big(t',  X^{s,   x^n_i}_{t'}   \big)     =   \b1_{ \cA^n_i}  L^{s,   x^n_i}_{t'}  \le   \b1_{ \cA^n_i}  Y^{s,   x^n_i}_{t'}  \\
      & \tneg &    =\b1_{ \cA^n_i}   h \big( X^{s,   x^n_i}_T\big) \neg + \neg  \int _{t'}^{T} \neg \b1_{ \cA^n_i}   f \big(r, X^{s,   x^n_i}_r\neg , Y^{s,   x^n_i}_r \neg , Z^{s,   x^n_i}_r \neg \big) dr
        \neg + \neg \b1_{ \cA^n_i}  K^{s,   x^n_i}_T \dneg - \neg \b1_{ \cA^n_i}  K^{s,   x^n_i}_{t'}  \dneg - \neg \int_{t'}^{T} \neg  \b1_{ \cA^n_i}  Z^{s,   x^n_i}_r dB_r \\
        & \tneg &    =\b1_{ \cA^n_i}   h \big( X^n_T\big)  + \neg  \int _{t'}^{T} \neg \b1_{ \cA^n_i}   f \big(r, X^n_r  , Y^n_r   , Z^n_r \big) dr
         + \neg \b1_{ \cA^n_i}  K^n_T  - \neg \b1_{ \cA^n_i}  K^n_{t'}   - \neg \int_{t'}^{T} \neg  \b1_{ \cA^n_i}  Z^n_r dB_r , \q  \pas
   \eeas
    Summing up both expressions over $i=0, \cds ,j_n $, one can deduce from  the continuity of function $l$
    as well as the continuity of processes  $ \{X^n_{t'}\}_{t' \in [s,T]}$,
    $ \{Y^n_{t'}\}_{t' \in [s,T]}$ and $ \{K^n_{t'}\}_{t' \in [s,T]}$ that \pas
    \bea
 X^n_{t'} & \tneg =& \tneg  X^n_s    + \int_{s}^{t'}     b(r,   X^n_r) dr  +\int_{s}^{t'}    \si(r,    X^n_r) dB_r ,    \q    t' \in  [s, T] ;  \label{eq:f221}    \\
  l \big(t',   X^n_{t'}   \big)    \le        Y^n_{t'}
       & \tneg =& \tneg  h( X^n_T) + \int _{t'}^T     f (r,    X^n_r, Y^n_r, Z^n_r) dr
          +     K^n_T -    K^n_{t'}      - \int_{t'}^T       Z^n_r dB_r, \q    t' \in  [s, T] .  \label{eq:f222a}
    \eea
     Moreover, we also have
  \bea
    \int_s^T \big(Y^n_r - l  (r,   X^n_r ) \big) d K^n_r =  \sum_{i=0}^{j_n} \b1_{ \cA^n_i}    \int_s^T \Big(Y^{s,   x^n_i}_r \neg - L^{s,   x^n_i}_r\Big) \, d K^{s,   x^n_i}_r =0, \q \pas   \label{eq:f222b}
  \eea

    By \eqref{FSDE},  it holds \pas ~ that
     \beas
 X^0_{t'} =  X^0_s    + \int_{s}^{t'}     b(r,   X^0_r) dr  +\int_{s}^{t'}    \si(r,    X^0_r) dB_r ,    \q    t' \in  [s, T] .
     \eeas
  Subtracting it  from \eqref{eq:f221}, we can  deduce from \eqref{b_si_Lip} that \pas
      \bea \label{eq:f269}
          \underset{s' \in [s, t']}{\sup}\big| X^n_{s'} \neg  - \neg X^0_{s'}  \big|  \le \neg  \big|X^n_s  \neg  - \neg  X^0_s\big|
   \neg + \neg   \k \neg \int_{s}^{t'} \neg    \big| X^n_r  \neg  -  \neg  X^0_r  \big| dr
   \neg + \dneg \underset {s' \in [s, t']}{\sup} \neg \left| \, \int_{s}^{s'} \dneg \big( \si(r,    X^n_r)  \neg  -  \neg  \si(r,    X^0_r)\big) dB_r \right|, ~   t' \in  [s, T].
    \eea
    By similar arguments to those that lead to \eqref{eq:f121},
     we can deduce from \eqref{eq:f269}    that \pas
      \beas
         \underset{t' \in [s,T]}{\sup}   \left| X^n_{t'}   -X^0_{t'} \right|^\varpi  \le  2^{\varpi-1}  e^{\k \varpi T}
         \left(  \big|X^n_s-X^0_s\big|^\varpi  + \underset{t' \in [s,T]}{\sup} \left| \int_s^{t'}  \big( \si(r,    X^n_r)-\si(r,    X^0_r)\big) dB_r \right|^\varpi \right).
     \eeas
   And using similar arguments to those that lead to \eqref{eq:f131}, we can deduce that for any $p \in (1, \infty)$
    \beas
      E \left[   \exp\left\{ p  \underset{t' \in [s,T]}{\sup} \left| \int_s^{t'}  \big( \si(r,    X^n_r)-\si(r,    X^0_r)\big) dB_r \right|^\varpi  \right\}  \right]
      \le \wt{c}_p   .
     \eeas
  Then H\"older's inequality implies that
   \beas
    && \hspace{-1cm}  E  \left[   \exp\left\{ p    \underset{t' \in [s,T]}{\sup}   \left|  X^n_{t'}   -X^0_{t'}  \right|^\varpi   \right\}  \right]     \le
    \left\{  E   \left[  \exp\left\{ p \, 2^{\varpi}  e^{\k \varpi T}    \underset{t' \in [s,T]}{\sup} \left| \int_s^{t'}  \big( \si(r,    X^n_r)
   \neg - \neg \si(r,    X^n_r)\big) dB_r \right|^\varpi     \right\}  \right]\right\}^{\frac12}         \\
   &     &       \times  \Big\{  E \Big[   \exp\Big\{ p \, 2^{\varpi}  e^{\k \varpi T}
     \big|X^n_s-X^0_s\big|^\varpi   \Big\}   \Big]\Big\}^{\frac12}    \le        \wt{c}_p  
       \Big\{  E   \Big[   \exp    \Big\{ p \, 2^{2\varpi-1}  e^{\k \varpi T}     \big| X^0_s\big|^\varpi  \Big\}   \Big] \Big\}^{\frac12} ,
   \eeas
    where we used the fact that
        \bea \label{eq:f273}
                       \big|X^n_s \neg - \neg  X^0_s\big|    =    \b1_{\{|X^0_s| < 2^n \}} \big|X^n_s \neg  -  \neg   X^0_s\big| \neg  +  \neg
         \b1_{\{|X^0_s| \ge 2^n \}} \big|  X^0_s\big| \le 2^{-n}  \neg  +  \neg    \b1_{\{|X^0_s| \ge 2^n \}} \big|  X^0_s\big|  .
          \eea
          Thus it follows that for any   $p \in [1, \infty)$
          \bea
          E  \left[   \exp\left\{ p    \underset{t' \in [s,T]}{\sup}   \left|  X^n_{t'}     \right|^\varpi   \right\}  \right] & \tneg \le & \tneg
        \frac12  E  \left[   \exp\left\{ p  2^\varpi  \underset{t' \in [s,T]}{\sup}   \left|  X^n_{t'}   -X^0_{t'}  \right|^\varpi   \right\}  \right]
        +\frac12  E  \left[   \exp\left\{ p  2^\varpi  \underset{t' \in [s,T]}{\sup}   \left|   X^0_{t'}  \right|^\varpi   \right\}  \right] \nonumber \\
        & \tneg  \le & \tneg    \wt{c}_p +     E \neg \left[   \exp \neg  \left\{ p \, 2^{3\varpi-1}  e^{\k \varpi T}
           \underset{t' \in [s,T]}{\sup}   \left|   X^0_{t'}  \right|^\varpi  \right\}   \right]   .    \label{eq:f275}
          \eea
     As $   \Big\{\big( Y^{s, x^n_i} , Z^{s, x^n_i} , K^{s, x^n_i}  \big) \Big\}_{i=0, \cds ,j_n}
        \subset   \underset{p \in [1, \infty)}{\cap}     \hS^p_\bF[0,T]$,  one can deduce that for any $p \in [1, \infty)$
          \bea
   &&  \hspace{-1.5cm}    E \neg \left[ \exp \neg \bigg\{ p \underset{t' \in [s,T]}{\sup} |Y^n_{t'}|   \bigg\}
    \neg+ \neg \left(\int_s^T |Z^n_r|^2 dr \right)^p +\big(K^n_T\big)^p \, \right]  \nonumber    \\
     &&   = \neg  E\left[ \,  \sum_{i=0}^{j_n} \neg  \b1_{ \cA^n_i} \left( \exp \neg \bigg\{ p \underset{t' \in [s,T]}{\sup} \big|Y^{s, x^n_i}_{t'}\big|   \bigg\}
      \neg+ \neg \bigg( \neg \int_s^T \big|Z^{s, x^n_i}_r \neg \big|^2 dr \bigg)^p   \neg +\neg \Big(K^{s, x^n_i}_T\Big)^p \right) \right]  \nonumber\\
  && \le \sum_{i=0}^{j_n} \neg  E \neg \left[   \exp \neg \bigg\{ p \underset{t' \in [s,T]}{\sup} \big|Y^{s, x^n_i}_{t'}\big|\bigg\}
  \neg+ \neg \bigg( \neg \int_s^T \big|Z^{s, x^n_i}_r \neg \big|^2 dr \bigg)^p   \neg +\neg \Big(K^{s, x^n_i}_T\Big)^p \right] < \infty .  \label{eq:f225}
  \eea

     \ss \no {\bf 2)}    Fix $m \in \hN_0$.  Since $ \cX^m_{t'}  \dfnn  \b1_{\{t' <s\}}  E[X^m_s|\cF_{t'} ] \neg+ \neg \b1_{\{t' \ge s\}} X^m_{t'} $,
 $t' \in [0,T]$  is  an $\bF$-adapted continuous process, the continuity of function $l$ and $f$ shows that
  $\cL^m_{t'} \neg \dfnn    l \big(t', \cX^m_{t'} \big)$, $t' \in [0,T]$ is also an $\bF$-adapted continuous process, and that
  \beas
   f_m(t', \o, y,z)   \dfnn   f\big(t',  \cX^m_{t'} (\o) ,y,z\big), \q \fa (t', \o, y,z) \in [0,T] \times \O \times \hR^k \times \hR^d
  \eeas
  is a    $\sP \times \sB(\hR) \times   \sB(\hR^d)/\sB(\hR)$-measurable function. Moreover,  \eqref{cond:f_basic}-\eqref{cond:f_concave} show that $f_m$ satisfies (H1)-(H3)  with the same   constants $\a, \beta,  \k  \ge 0$ and $\g > 0$ as $f$.       For any $p \in (1, \infty)$,   the convexity of function $y \to  e^{|y|^\varpi }$ on $\hR$ and Jensen's inequality imply that
  $ \bigg\{ \neg \exp\Big\{     \Big(  E\big[ |X^n_s| \big|\cF_{t'} \big] \Big)^\varpi  \Big\} \neg \bigg\}_{t' \in [0, \infty)}$ is a continuous positive submartingale.  Doob's Martingale Inequality then shows that
     \beas
       E \left[ \underset{t' \in [0,  s]}{\sup}  \bigg( \exp\Big\{     \Big(  E\big[ |X^m_s| \big|\cF_{t'} \big] \Big)^\varpi  \Big\} \bigg)^p \right]
   & \tneg \le &  \tneg \hb{$\big(\frac{p}{p-1}\big)^p$}
   E \bigg[   \Big( \exp\big\{      |X^m_s|^\varpi  \big\} \Big)^p \bigg]   ,
   \eeas
 which together with \eqref{eq:f275} and Lemma \ref{lemm:esti_X}  leads to that
   \beas
        E \left[ \exp\Big\{   p  \big(  \cX^m_*   \big)^\varpi  \Big\}   \right]
       & \tneg \le  & \tneg   E \left[ \underset{t' \in [0,  s]}{\sup}  \exp\Big\{   p \, \Big(  E\big[ |X^m_s| \big|\cF_{t'} \big]   \Big)^\varpi  \Big\}   \right] +
      E \bigg[  \, \underset{t' \in [s,T]}{\sup}   \exp\Big\{    p \,  \big|X^m_{t'}\big|^\varpi  \Big\}  \bigg] \\
        & \tneg \le  & \tneg   \wt{c}_p  
          E \Bigg[     \exp\bigg\{    p \underset{t' \in [s,T]}{\sup}  |X^m_{t'}|^\varpi  \bigg\}  \Bigg]
          \le \wt{c}_p +  \wt{c}_p   E \neg \left[   \exp \neg  \left\{ p \, 2^{3\varpi-1}  e^{\k \varpi T}
           \underset{t' \in [s,T]}{\sup}   \left|   X^0_{t'}  \right|^\varpi  \right\}   \right]  \\
         & \tneg \le  & \tneg    \wt{c}_p +  \wt{c}_p     \exp \neg  \left\{ p \, 2^{3\varpi-1} 3^{\varpi-1}    e^{2\k \varpi T}      | x |^\varpi  \right\}    .
   \eeas
 Hence it follows from \eqref{hl_varipi} that
       \bea
             E \neg \left[ \exp \neg \left\{ p \Big( \big|h\big( \cX^m_T \big) \big|  \neg \vee \neg  \cL^m_* \Big) \neg \right\} \right]
           \neg \le   e^{p \k} E \neg \left[ \exp \neg \left\{ \neg (1\neg  \vee \neg p\k)  \big( \cX^m_*  \big)^\varpi   \right\} \right]
           \neg  \le   \wt{c}_p \neg + \neg  \wt{c}_p     \exp \neg  \left\{   (1 \neg \vee \neg  p\k) \,  2^{3\varpi-1} 3^{\varpi-1}    e^{2\k \varpi T}      | x |^\varpi  \right\} .  \q \label{eq:f229}
        \eea
  As $Y^{t,x} \in \hE^p_\bF[0,T]$, we  also see from \eqref{eq:f225} that $ E  \Big[ e^{ p \,  |Y^m_s|   } \,  \Big]   < \infty $.
  Since $Y^m_s \ge  l \big(s, X^m_s\big) =  l \big(s, \cX^m_s\big) = \cL^m_s$, \pas, Corollary \ref{cor:unique} implies that
 the quadratic RBSDE$( Y^m_s, f_m, \cL^m )$ over time interval $[0,s]$ admits a unique solution $\big\{ (\cY^m_r, \cZ^m_r, \cK^m_r ) \big\}_{r \in [0,s]}   $
  in $ \underset{p \in [1, \infty)}{\cap} \hS^p_\bF[0,s] $.

   \ms   We extend the processes $ (\cY^m, \cZ^m, \cK^m) $ to the period $(s,T]$ by setting:  $\fa t' \in (s,T]$,
 \beas
     (\cY^m_{t'}, \cZ^m_{t'} )  \dfnn   ( Y^m_{t'},  Z^m_{t'} ) \q \hb{and} \q
      \cK^m_{t'}   \dfnn   \left\{\ba{ll}
      \cK^0_s+ K^0_{t'} -K^0_s    ,\q & \hb{if } m=0;  \\
         \cK^m_s+ K^m_{t'}   ,\q & \hb{if } m \in \hN.
      \ea \right.
 \eeas
   Then    \eqref{eq:f222a} and \eqref{eq:f222b} imply that $\{(\cY^m_{t'}, \cZ^m_{t'}, \cK^m_{t'})\}_{t' \in [0,T]}  $ is a solution of
     the quadratic RBSDE$\big(h(\cX^m_T),    f_m, \cL^m \big)$. As    $(Y^{t,x}, Z^{t,x}, K^{t,x})  \in  \underset{p \in [1, \infty)}{\cap}    \hS^p_\bF[0,T] $, we see from
   \eqref{eq:f225} that $(\cY^m, \cZ^m, \cK^m)  \in  \underset{p \in [1, \infty)}{\cap}    \hS^p_\bF[0,T] $.
   Moreover, Corollary \ref{cor:unique} and \eqref{eq:f229}   show that $(\cY^m, \cZ^m, \cK^m)$
   is the unique solution of the quadratic RBSDE$\big(h(\cX^m_T), f_m, \cL^m \big)$ in $\underset{p \in [1, \infty)}{\cap}    \hS^p_\bF[0,T] $.

  \ss \no {\bf 3)} Squaring both sides of \eqref{eq:f269}, one can deduce from  H\"older's inequality,   Doob's martingale inequality, Fubini's Theorem and  \eqref{b_si_Lip} that
    \beas
      E \neg \left[    \underset{s' \in [s, t']}{\sup}\big| X^n_{s'} \neg  - \neg X^0_{s'}  \big|^2 \right]
        &\tneg \le& \tneg  3 E \neg \left[ \big|X^n_s \neg - \neg X^0_s\big|^2 \right] \neg + \neg   3 \k^2  T E\neg \int_s^{t'} \neg
             \big| X^n_r \neg - \neg  X^0_r  \big|^2   dr  \neg
         + \neg 12 \,  E  \neg   \int_s^{t'} \neg  \big|  \si(r,    X^n_r) \neg - \neg \si(r,    X^0_r) \big|^2  dr   \\
    & \tneg  \le  &  \tneg   3 E \neg \left[ \big|X^n_s \neg - \neg X^0_s\big|^2 \right]  \neg + \neg    3 \k^2  (T \neg + \neg 4) \neg \int_s^{t'} \neg
        E \neg \left[ \underset{s' \in [s, r]}{\sup} \big| X^n_{s'} \neg - \neg X^0_{s'}  \big|^2 \right] dr  ,   ~\;   t' \in [s,T].
       \eeas
 Then Gronwall's inequality and \eqref{eq:f273} imply that
        \beas
                   E \neg \left[    \underset{t' \in [s, T]}{\sup}\big| X^n_{t'} \neg  - \neg X^0_{t'}  \big|^2 \right]
         \le  3 E \neg \left[ \big|X^n_s \neg - \neg X^0_s\big|^2 \right] \,  e^{ 3 \k^2 (T^2+4T)   }
         \le  \wt{c}_0   \bigg( 2^{-2 n}  \neg + \neg E \neg \left[     \b1_{\{|X^0_s| \ge 2^n \}} \big|  X^0_s\big|^2   \right] \bigg) .
         \eeas
   As $E \Big[\big|  X^{t,x}_s\big|^2\Big]< \infty$ by \eqref{eq:esti_X_1},  letting $n \to \infty$
     yields that
  $   \lmt{n \to \infty}  E \neg \left[    \underset{t' \in [s, T]}{\sup}\big| X^n_{t'} \neg  - \neg X^0_{t'}  \big|^2 \right]  = 0$. Since Doob's martingale inequality
  implies that
   \beas
     E \neg \left[    \underset{t' \in [0, T]}{\sup}\big| \cX^n_{t'} \neg  - \neg \cX^0_{t'}  \big|^2 \right] &\dneg \le & \dneg
      E \neg \left[    \underset{t' \in [0, s]}{\sup}\Big| E[ X^n_s - X^0_s | \cF_{t'}]   \Big|^2 \right]+
       E \neg \left[    \underset{t' \in [s, T]}{\sup}\big| X^n_{t'} \neg  - \neg X^0_{t'}  \big|^2 \right] \\
       &\dneg \le & \dneg 4  E \neg \left[     \big|   X^n_s  \neg  - \neg   X^0_s     \big|^2 \right] +   E \neg \left[    \underset{t' \in [s, T]}{\sup}\big| X^n_{t'} \neg  - \neg X^0_{t'}  \big|^2 \right] \le 5 E \neg \left[    \underset{t' \in [s, T]}{\sup}\big| X^n_{t'} \neg  - \neg X^0_{t'}  \big|^2 \right]   ,
   \eeas
   it follows that $\lmt{n \to \infty}  E \neg \left[    \underset{t' \in [0, T]}{\sup}\big| \cX^n_{t'} \neg  - \neg \cX^0_{t'}  \big|^2 \right]  = 0$.
  Hence, we can pick up a subsequence of $\{\cX^n\}_{n \in \hN}$ (we still denote it by $\{\cX^n\}_{n \in \hN}$) such that
  except on a $P$-null set $\sN$,
      \bea  \label{eq:f301}
  \lmt{n \to \infty}  \Big( \underset{t' \in [0, T]}{\sup}\big| \cX^n_{t'} \neg  - \neg \cX^0_{t'}  \big|  \Big)  =0
  \q \hb{and} \q \hb{the path $t'   \to   \cX^0_{t'}$ is continuous}.
    \eea

   \ss  To apply Theorem \ref{thm:stable} to the sequence $\big\{\big( \cY^n , \cZ^n , \cK^n  \big)\big\}_{n \in \hN}$,
 let us check the assumptions of this theorem first. We have seen that the sequence $ \{  f_m  \}_{m \in \hN_0 } $ satisfies (S1),
 and that \eqref{eq:f229} justifies (S3).

    \ms
    Fix $\o \in \sN^c$.  For any $\e > 0$, the continuity of $h$ assures that there exists  a $\d (\o)  \in (0,1) $ such that
    \beas
     |h(\tilde{x}) \neg - \neg h(x') | \vee  |l(\tilde{s}, \tilde{x}) \neg -\neg l(s', x') | < \e, ~  \fa (\tilde{s},\tilde{x}), (s',x') \in [0,T] \times \wt{\sD} (\o)
     \hb{\;with\;} |\tilde{s}-s'|^2 \neg +    |\tilde{x}-x'|^2      < \d^2 (\o) ,
  \eeas
    where $ \wt{\sD} (\o) \dfnn \Big\{\tilde{x} \in \hR^k: |\tilde{x}| \le 1 + \underset{t' \in [0, T]}{\sup}  \big| \cX^0_{t'}        (\o)   \big|< \infty \Big\}$.
    Moreover, in light of \eqref{eq:f301},   there exists an $ N(\o ) \in \hN$ such that for any $n \ge N(\o)$,
     $  \underset{t' \in [0, T]}{\sup}\big| \cX^n_{t'} (\o)  -   \cX^0_{t'} (\o)  \big| <  \d (\o)  $.
 Then it holds  for any $n \ge N(\o)$ that
      \beas
      \big| h\big( \cX^n_T (\o) \big) \neg - \neg h\big( \cX^0_T (\o) \big) \big| < \e \q \hb{and} \q
            \big| \,  \cL^n_{t'} (\o) \neg - \neg \cL^0_{t'} (\o)  \big|= \big| \,  l \big(t' \cX^n_{t'} (\o) \big) \neg - \neg l \big(t' \cX^0_{t'} (\o) \big)  \big| < \e, \q  t' \in [0,T].
      \eeas
 Thus    (S2) is satisfied.

 \ms      \ms  Given $(t', \o) \in [0,T] \times \sN^c$,
   the continuity of $f$ and \eqref{eq:f301} imply that
   \beas
           \lmt{n \to \infty}    f_n\big(t', \o, \cY^0_{t'}(\o),  \cZ^0_{t'}(\o) \big)
          &=&   \lmt{n \to \infty}    f \big(t', \cX^n_{t'}(\o), \cY^0_{t'}(\o),  \cZ^0_{t'}(\o) \big)
           =  f \big(t', \cX^0_{t'}(\o), \cY^0_{t'}(\o),  \cZ^0_{t'}(\o) \big) \\
           &=&   f_0\big(t', \o, \cY^0_{t'}(\o),  \cZ^0_{t'}(\o) \big) .
           \eeas

        \ms   Now, applying Theorem \ref{thm:stable} yields that $ \lmt{n \to \infty} E \left[ \exp\bigg\{   \underset{t' \in [0,T]}{\sup}
        \big|\cY^n_{t'} -\cY^0_{t'} \big|\bigg\}  \right] = 1 $, which allows us to extract   a  subsequence of $\{ \cY^n \}_{n \in \hN}$
        \big(we still denote it by     $\{ \cY^n \}_{n \in \hN}$\big) such that
        $ \lmt{n \to \infty} \,  \underset{t' \in [0,T]}{\sup}  \big|\cY^n_{t'} -\cY^0_{t'} \big|  = 0 $, \pas~ In particular,
       it holds \pas ~ that
         \bea    \label{eq:f331}
           \lmt{n \to \infty}   Y^n_s =  \lmt{n \to \infty}  \cY^n_s =   \cY^0_s = Y^0_s = Y^{t,x}_s   ,
         \eea
           where
           \beas
             Y^n_s =     \sum_{i=0}^{j_n}  \b1_{ \cA^n_i}  Y^{s, x^n_i}_s = \sum_{i=0}^{j_n}  \b1_{ \cA^n_i} u\big(s,  x^n_i\big)
                 = \sum_{i=0}^{j_n}  \b1_{ \cA^n_i}  u\big(s,  X^n_s \big)  = u\big(s,  X^n_s \big)  ,  \q  \fa n \in \hN .
           \eeas
          Since $\lmt{n \to \infty} X^n_s = X^0_s =X^{t,x}_s$, \pas ~  by \eqref{eq:f273}, Proposition \ref{prop:u_continuous} and \eqref{eq:f331}
         then imply that
           \beas
            Y^{t,x}_s =  \lmt{n \to \infty} u\big(s,  X^n_s \big) =  u\big(s,  X^{t,x}_s \big) ,\q \pas
           \eeas
           Eventually, the continuity of processes $X^{t,x}$, $Y^{t,x}$ and Proposition \ref{prop:u_continuous} leads to  \eqref{eq:f201}.   \qed

  \begin{thm}   \label{thm:viscos_exist}
   The function $u$ defined in \eqref{def_utx} is a viscosity solution of \eqref{eq:PDE}.
  \end{thm}

 \ss \no {\bf Proof:} {\bf 1) } For any $x \in \hR^k$, it is clear that $u(T,x)= \wt{Y}^{T,x}_0= h\big( X^{T,x}_T\big) = h(x)$. We first show that $u$  is a viscosity subsolution of \eqref{eq:PDE}. Let $(t_0,x_0, \vf) \in (0,T) \times \hR^k \times C^{1,2}\big([0,T] \times \hR^k\big)$ be  such that $u(t_0,x_0) = \vf (t_0,x_0)$ and that
  $u-\vf$ attains a local  maximum  at $(t_0,x_0)$. We prove by contradiction. Suppose that
   \beas
    \q  \e \dfnn \frac12 \min \left\{    (u-l) (t_0, x_0),  -   \frac{\pa \vf }{\pa t}(t_0,x_0)
     - \neg \cL \vf(t_0,x_0) - \neg f \neg \left(t_0,x_0, \vf(t_0,x_0),
   (\si^T \nabla_x \vf  ) (t_0,x_0)\right)  \right\} > 0.
   \eeas
 Since $\vf \in C^{1,2}\big([0,T] \times \hR^k\big)$,  the continuity of functions $u, l, f$ and $\si$ as well as the  assumption on local  maximum of $u-\vf$ assure that there exists   a $\d \in \big(0,  T - t_0 \big] $ such that for any $ t \in [t_0 , t_0+\d ]  $  and any $x \in \hR^k$ with $|x-x_0| \le \d $
   \bea
  &&   |u  (t,x) \neg - \neg u(t_0,x_0) | \le \frac13 \e, \q   (u - l) (t,x) \ge    \e ,     \q  (u - \vf) (t,x) \le  0 ,  \hspace{2cm} \label{eq:f351} \\
    \hb{ and } &&   -   \frac{\pa \vf }{\pa t}(t,x)       - \neg \cL \vf(t,x) - \neg f \neg \left(t,x,   \vf(t,x),     (\si^T \nabla_x \vf  ) (t,x)\right) \ge \e . \label{eq:f353}
   \eea
 Since $\big\{\wt{X}^{t_0,x_0}_s\big\}_{s \in [0, T-t_0]} $ and $\wt{Y}^{t_0,x_0}$ are both $\bF^{t_0}$-adapted continuous processes,
   \bea   \label{eq:f355}
   \qq  \nu  \dfnn   \inf\Big\{ s \in [0,  \d] : \big| \wt{X}^{t_0,x_0}_s - x_0 \big|  > \d  \Big\} \land
    \inf\Big\{ s \in [0,   \d] :   \big| \wt{Y}^{t_0,x_0}_s - \wt{Y}^{t_0,x_0}_0 \big|> \frac13 \e    \Big\}  \land   \d
   \eea
   defines  an $\bF^{t_0}$-stopping time such that $\nu > 0 $, \pas ~
    For any $\o \in \O$ and $s \in [0, \nu(\o)]$, one can deduce from \eqref{eq:f351} that
  \beas
    \wt{Y}^{t_0,x_0}_s (\o) & \dneg \ge&  \dneg   \wt{Y}^{t_0,x_0}_0  - \frac13 \e =u (t_0,x_0) - \frac13 \e
    \ge u  \big(t_0+s, \wt{X}^{t_0,x_0}_s (\o) \big) - \frac23 \e \\
    &  \dneg  \ge&  \dneg    l \big(t_0+s,  \wt{X}^{t_0,x_0}_s (\o) \big)  + \frac13 \e = \wt{L}^{t_0, x_0}_s(\o) + \frac13 \e .
  \eeas
  Because $ \big(\wt{Y}^{t_0,x_0}, \wt{Z}^{t_0,x_0}, \wt{K}^{t_0,x_0} \big) \in  \underset{p \in [1, \infty)}{\cap} \hS^p_{\bF^{t_0}}[0,T-t_0] $
   solves the quadratic RBSDE$\Big(h \big( X^{t_0,x_0}_{T}  \big) , \tilde{f}^{t_0,x_0}, \wt{L}^{t_0,x_0} \Big)$ with respect to $B^{t_0}$ over the period $[0,T-t_0]$,
   its flat-off condition  implies that    \pas, $\wt{K}^{t_0,x_0}_s  =0$ for any $  s \in [0, \nu ]$.
         Hence, it holds  \pas ~ that
     \beas
       \wt{Y}^{t_0,x_0}_{\nu \land s} =  \wt{Y}^{t_0,x_0}_\nu + \int_{\nu \land s}^\nu  \tilde{f}^{t_0,x_0}\Big(r, \wt{Y}^{t_0,x_0}_r,
         \wt{Z}^{t_0,x_0}_r \Big) dr -  \int_{\nu \land s}^\nu   \wt{Z}^{t_0,x_0}_r   d B^{t_0}_r, \q  s \in [0,\d].
 \eeas
 In other words, the  processes $\big(\cY, \cZ\big) \dfnn \left\{ \left(\wt{Y}^{t_0,x_0}_{\nu   \land    s}, \b1_{\{s < \nu\}} \wt{Z}^{t_0,x_0}_s\right)
  \right\}_{s \in [0,\d]} \in \hC^\infty_{\bF^{t_0}}[0,\d] \times  \underset{p \in [1, \infty)}{\cap} \hH^{2, 2p}_{\bF^{t_0}}([0,\d];\hR^d)$ solves the BSDE
 \bea
       \cY_s  & \tneg =& \tneg \wt{Y}^{t_0,x_0}_\nu        + \int_s^\d \mathfrak{f}  \left(r, \cY_r,
       \cZ_r \right) dr -  \int_s^\d  \cZ_r  d B^{t_0}_r, \q  s \in [0,\d],  \nonumber \\
   \hb{with} \q  \mathfrak{f} (s, \o, y, z )  & \tneg \dfnn & \tneg  \b1_{\{ s \, < \nu(\o) \}}\tilde{f}^{t_0,x_0} (s, \o, y,z) ,
   \q   \fa (s, \o, y, z ) \in [0,\d] \times \O \times \hR \times \hR^d .  \hspace{2cm}   \label{eq:f363}
  \eea
  Like $\tilde{f}^{t_0,x_0}$,   $\mathfrak{f}$ is a  generator with respect to $\bF^{t_0}$ over the period $[0,\d]$ that    satisfies (H1)-(H3).

  \ms  On the other hand, since
     \beas
     \wt{X}^{t_0,x_0}_s= x + \int_0^s  b(r+t_0, \wt{X}^{t_0,x_0}_r) dr
      + \int_0^s  \si(r+t_0, \wt{X}^{t_0,x_0}_r) dB^{t_0}_r,    \q s \in [0, T-t_0] ,
    \eeas
  applying It\^o's formula to the process $\vf(t_0+\cd, \wt{X}^{t_0,x_0}_\cd)$  yields that
    \beas
       \vf \big( t_0+\nu  \land s, \wt{X}^{t_0,x_0}_{\nu  \land s} \big) & \dneg =& \dneg  \vf \big(t_0+\nu , \wt{X}^{t_0,x_0}_{\nu}\big)
       - \int_{ \nu  \land s}^{ \nu}  \big(\hb{$\frac{\pa \vf }{\pa t}$} + \cL \vf\big) \neg \big(t_0+r, \wt{X}^{t_0,x_0}_r\big)dr   \nonumber \\
       && -  \int_{ \nu  \land s}^{ \nu}   (\si^T \nabla_x \vf  ) \big(t_0+r, \wt{X}^{t_0,x_0}_r\big)dB^{t_0}_r , \q   s \in [0, \d] .  
    \eeas
   Namely,  $   \big(\cY', \cZ'\big) \dfnn \Big\{ \Big(   \vf \big( t_0  \neg +  \neg  \nu  \land s, \wt{X}^{t_0,x_0}_{    \nu  \land s}\big),
  \b1_{\{s < \nu\}}    (\si^T \nabla_x \vf  ) \big(t_0  \neg +  \neg  s, \wt{X}^{t_0,x_0}_s  \big)\Big)  \Big\}_{s \in [0,\d]}$
  solves the BSDE
     \beas  
       \cY'_s  =   \vf \big(t_0+\nu , \wt{X}^{t_0,x_0}_{\nu}\big)   + \int_s^\d   \mathfrak{f}'_r  dr   -  \int_s^\d    \cZ'_r    dB^{t_0}_r , \q   s \in [0, \d] ,
    \eeas
  where $\mathfrak{f}'_s    \dfnn   - \b1_{\{ s \, < \nu  \}}   \big(\hb{$\frac{\pa \vf }{\pa t}$} + \cL \vf\big) \big(t_0+s,
       \wt{X}^{t_0,x_0}_s \big) $, $\fa  s  \in [0,\d]$.   Since $\wt{X}^{t_0,x_0} $ is an $\bF^{t_0}$-adapted continuous process, and since
      $\vf \in   C^{1,2}\big([0,T] \times \hR^k\big)$,  the continuity of function $\si$ implies that
      $\cY'$ is an  $\bF^{t_0}$-adapted continuous process as well as that $\cZ'$ and $\mathfrak{f}'$
      are both $\bF^{t_0}$-progressively measurable processes.  Moreover,
      since $   \big| \wt{X}^{t_0,x_0}_s - x_0   \big|  \le    \d  $ holds for \pas~ $\o \in \O$ and $s \in \big[0, \nu(\o)\big]$, and since
      $\vf \in   C^{1,2}\big([0,T] \times \hR^k\big)$, we further see from the continuity of function $b$ and the boundedness of function $\si$ that
       $\cY'$, $\cZ'$ and $\mathfrak{f}' $ are all bounded processes.

   \ms One can deduce from  Proposition \ref{prop:u_markov} and  \eqref{eq:f351} -\eqref{eq:f355}  that
       \beas
        \wt{Y}^{t_0,x_0}_\nu =  u \big( t_0+\nu   , \wt{X}^{t_0,x_0}_{ \nu   }\big)       \le  \vf \big( t_0+\nu   , \wt{X}^{t_0,x_0}_{ \nu   }\big)  , \q  \pas ,
    \eeas
 and that on $\O$
 \bea
  \q && \hspace{-1 cm} \mathfrak{f}'_s -\mathfrak{f} \big(s,    \cY'_s,  \cZ'_s\big)
     =   -   \b1_{\{ s \, < \nu  \}}   \big(\hb{$\frac{\pa \vf }{\pa t}$} + \cL \vf\big) \big(t_0+s,
       \wt{X}^{t_0,x_0}_s  \big)   \nonumber \\
     &    &         -   \b1_{\{ s < \nu  \}}f \Big(t_0 \neg + \neg  s, \wt{X}^{t_0,x_0}_s  ,    \vf \big( t_0  \neg +  \neg  s, \wt{X}^{t_0,x_0}_s  \big),
       (\si^T \nabla_x \vf  ) \big(t_0  \neg +  \neg  s, \wt{X}^{t_0,x_0}_s   \big) \Big)  \ge \e \b1_{\{ s \, < \nu  \}} ,  ~\; \fa  s \in [0,\d]  .  \qq   \label{eq:f359}
 \eea
 The first part of Proposition \ref{prop:comparion_ext} or that of Proposition \ref{prop:comparion_ext2} implies that \pas,   $\cY'_s \ge \cY_s$ for any $   s \in [0,\d]$. Since  $      \cY'_0  = \vf(t_0, x_0)  = u(t_0,x_0) =    \wt{Y}^{t_0,x_0}_0 = \cY_0$, the second part of   Proposition \ref{prop:comparion_ext}
 or that of Proposition \ref{prop:comparion_ext2} further shows that
   $        P\Big(  \int_0^\d \big( \mathfrak{f}'_s -\mathfrak{f} \big(s,    \cY'_s,  \cZ'_s\big)  \big)   ds =0\Big)    >0$.
             However,  \eqref{eq:f359} and \eqref{eq:f355} show that  \pas,   $
       \int_0^\d \Big( \mathfrak{f}'_s -\mathfrak{f} \big(s,    \cY'_s,  \cZ'_s\big)  \Big)   ds 
       \ge  \e \nu >0$,   which leads to a contradiction.

 \ms \no {\bf 2)}  Next, we show that $u$  is a viscosity supersolution of \eqref{eq:PDE}. Let $(t_0,x_0, \vf) \in (0,T) \times \hR^k \times C^{1,2}\big([0,T] \times \hR^k\big)$ be  such that $u(t_0,x_0) = \vf (t_0,x_0)$ and that
  $u-\vf$ attains a local  minimum   at $(t_0,x_0)$.  Since
 \beas
 u(t_0,x_0) = Y^{t_0,x_0}_{t_0} \ge L^{t_0,x_0}_{t_0} = l  \big(t_0, X^{t_0,x_0}_{t_0} \big)= l(t_0, x_0) ,
 \eeas
 it suffices to show that
  \beas
 - \frac{\pa \vf }{\pa t}(t_0,x_0)
     - \neg \cL \vf(t_0,x_0) - \neg f \neg \left(t_0,x_0, \vf(t_0,x_0),
   (\si^T \nabla_x \vf  ) (t_0,x_0)\right) \ge 0.
  \eeas
 To make a  contradiction, we assume that
   \beas
    \q  \e \dfnn \frac12   \left(          \frac{\pa \vf }{\pa t}(t_0,x_0)
     + \neg \cL \vf(t_0,x_0) + \neg f \neg \left(t_0,x_0, \vf(t_0,x_0),
   (\si^T \nabla_x \vf  ) (t_0,x_0)\right)  \right) > 0.
   \eeas
 Since $\vf \in C^{1,2}\big([0,T] \times \hR^k\big)$,  the continuity of functions $f$ and $\si$ as well as the  assumption on local  minimum of $u-\vf$ assures that there exists   a $\d  \in \big(0,  T - t_0 \big] $ such that for any $ t \in [t_0 , t_0+\d ]  $  and any $x \in \hR^k$ with $|x-x_0| \le \d $, one has
   \bea
            \frac{\pa \vf }{\pa t}(t,x)       + \neg \cL \vf(t,x)  +  \neg f \neg \left(t,x,   \vf(t,x) ,     (\si^T \nabla_x \vf  ) (t,x)\right) \ge \e
        \q \hb{and} \q (u - \vf) (t,x) \ge  0. \label{eq:f453}
   \eea
 We still define   the $\bF^{t_0}$-stopping time $\nu$ as in    \eqref{eq:f355}. It is easy to see  that
  the  processes
   \beas
 \q   \big(\cY, \cZ, \cV \big) \dfnn \left\{ \left(\wt{Y}^{t_0,x_0}_{\nu   \land    s}, \b1_{\{s < \nu\}} \wt{Z}^{t_0,x_0}_s,
  \wt{K}^{t_0,x_0}_{\nu   \land    s} \right)  \right\}_{s \in [0,\d]} \in \hC^\infty_{\bF^{t_0}}[0,\d] \times  \underset{p \in [1, \infty)}{\cap} \hH^{2, 2p}_{\bF^{t_0}}([0,\d];\hR^d) \times  \underset{p \in [1, \infty)}{\cap} \hK^p_{\bF^{t_0}}[0,\d]
   \eeas
    solves the BSDE \eqref{VBSDE2}  with generator  $\hat{f}= \mathfrak{f}$ as defined in \eqref{eq:f363} over the period $[0, \d]$.
       Let $(\cY',\cZ')$ be the pair of processes considered in part 1.
     Proposition \ref{prop:u_markov}, \eqref{eq:f453} and the definition of $\nu$ imply    that
       \beas
        \wt{Y}^{t_0,x_0}_\nu =  u \big( t_0+\nu   , \wt{X}^{t_0,x_0}_{ \nu   }\big)
       \ge  \vf \big( t_0+\nu   , \wt{X}^{t_0,x_0}_{ \nu   }\big)   , \q  \pas ,
    \eeas
 and that on $\O$
 \beas
    \mathfrak{f} \big(s,   \cY' ,  \cZ' \big)-  \mathfrak{f}'_s
    & \tneg =& \tneg   \b1_{\{ s < \nu  \}}   f \Big(t_0 \neg + \neg  s, \wt{X}^{t_0,x_0}_s  ,
        \vf \big( t_0  \neg +  \neg  s, \wt{X}^{t_0,x_0}_s  \big),
       (\si^T \nabla_x \vf  ) \big(t_0  \neg +  \neg  s, \wt{X}^{t_0,x_0}_s   \big) \Big) \\
     & \tneg  & \tneg + \b1_{\{ s < \nu  \}}  \big(\hb{$\frac{\pa \vf }{\pa t}$} + \cL \vf\big) \big(t_0+s,
       \wt{X}^{t_0,x_0}_s  \big)    \ge        \e  \b1_{\{ s < \nu  \}} ,  \q \fa s \in [0,\d]   .
 \eeas
 Using similar arguments to those that follow \eqref{eq:f359}, we  reach a  contradiction.  \qed

For the uniqueness of the viscosity solution of \eqref{eq:PDE}, we first establish a comparison principle between its viscosity subsolution and
 viscosity supersolution:

\begin{lemm} \label{lem:tilde_u}
 Let  $\mathfrak{a} > 0 $ and $\z \in \hR$.  If $u \in C\big([0,T] \times \hR^k\big)$ is a viscosity subsolution \(resp. viscosity supersolution\)  of \eqref{eq:PDE},
  then
  \beas       
   \tilde{u} (t,x)   \dfnn     \mathfrak{a} e^{\z t} u (t,x)  , \q  \fa  (t,x) \in [0,T] \times \hR^k
  \eeas
  becomes a viscosity subsolution \(resp. viscosity supersolution\) of the following   obstacle problem of  semi-linear parabolic PDE
 \bea  \label{eq:PDE2}
  \left\{\ba{l}
  \dis  \dneg \min\bigg\{ \tilde{u}     (t,x)-\mathfrak{a} e^{\z t}l(t,x) ,    -\frac{\pa  \tilde{u} }{\pa t}   (t,x)- \neg \cL \tilde{u} (t,x)
    - \neg \tilde{f}_\mathfrak{a}   \big(t,x,  \tilde{u} (t,x) ,         \nabla_x    \tilde{u} (t,x)
         \big)        \bigg\}
    =0,     ~      \fa (t,x) \in (0,T) \times \hR^k, \vspace{2mm} \\
   \tilde{u}(T,x)=\mathfrak{a} e^{\z T} h(x)    , ~  \fa x \in \hR^k,
   \ea
   \right.
 \eea
  where          $  \tilde{f}_\mathfrak{a}    (t,x,y,z) \dfnn   - \z y +  \mathfrak{a} e^{\z t} f \Big(t,x,\hb{$\frac{1}{\mathfrak{a}}$}e^{-\z t}   y ,
  \hb{$\frac{1}{\mathfrak{a}}$}e^{-\z t}  \,   \si^T (t,x) \cd     z          \Big)    ,     ~    \fa (y,z) \in \hR \times \hR^k$.
\end{lemm}

\ss \no {\bf Proof:}   We first assume that $u$ is a viscosity subsolution of \eqref{eq:PDE}.
   Clearly,    $\tilde{u}  \in C([0, T] \times \hR^k )$ and   $\tilde{u}(T,x) = \mathfrak{a} e^{\z T}  u (T,x)  \le \mathfrak{a} e^{\z T}  h(x)  $, $\fa x \in \hR^k$.
   Let $(t_0,x_0,  \wt{ \vf}) \in (0,T) \times \hR^k \times C^{1,2}\big([0,T] \times \hR^k\big)$ be  such that $\tilde{u}(t_0,x_0) =  \wt{ \vf} (t_0,x_0)$ and that
  $\tilde{u}-  \wt{ \vf}$ attains a local  maximum  at $(t_0,x_0)$. Then
   \beas
      \vf (t,x) \dfnn \frac{1}{\mathfrak{a}}  e^{-\z t}   \wt{\vf} (t,x)     , \q     \fa (t,x)   \in [0,T] \times \hR^k
  \eeas
   is   a $C^{1,2}\big([0,T] \times \hR^k\big)$ function such that $u(t_0,x_0) = \vf (t_0,x_0)$ and that
  $u-\vf$ attains a local  maximum  at $(t_0,x_0)$.
 Thus,
   \bea
     \min\left\{ \neg (u \neg - \neg l) (t_0,x_0  ), - \frac{\pa \vf }{\pa t}(t_0,x_0  )-\neg \cL \vf(t_0,x_0  )
     - \neg f \neg \left(t_0,x_0, u(t_0,x_0  ),       (\si^T \neg \cd \neg \nabla_x \vf  )        (t_0,x_0  )\right) \neg \right\} \le    0.  \label{eq:f511}
 \eea
  Suppose that  $\tilde{u}     (t_0,x_0)-\mathfrak{a} e^{\z t_0}  l(t_0,x_0)  > 0$, or equivalently, $u(t_0,x_0)-l(t_0,x_0)>0$.  By \eqref{eq:f511},
   \bea   \label{eq:f513}
   -   \frac{\pa \vf }{\pa t}(t_0,x_0)   -   \cL \vf (t_0,x_0)
     -  f \Big(t_0,x_0,    u(t_0,x_0) ,       (\si^T \neg \cd \neg \nabla_x \vf  ) (t_0,x_0)        \Big) \le 0 .
   \eea
   For any $(t,x) \in (0,T) \times \hR^k$, one can compute that
    \beas
  \qq     \frac{\pa \vf  }{\pa t}   (t,x)   =    \frac{1}{\mathfrak{a}}   e^{-\z t}    \Big( \frac{\pa \wt{\vf} }{\pa t}(t,x)
    - \z   \wt{\vf}  (t,x)  \Big) , ~ \;\;
   \nabla_x   \vf  (t,x)    =  \frac{1}{\mathfrak{a}}   e^{-\z t}       \nabla_x \wt{\vf}  (t,x),
           ~ \;\;    \hb{and}   ~ \;\;        \cL  \vf  (t,x)     =   \frac{1}{\mathfrak{a}}  e^{-\z t}    \cL \wt{\vf} (t,x)     .
 \eeas
 Plugging them into \eqref{eq:f513} yields  that
    \beas
      -\frac{\pa  \wt{\vf}  }{\pa t}   (t_0,x_0)+   \z    \tilde{u}   (t_0,x_0) - \neg \cL \wt{\vf}  (t_0,x_0)
    - \neg \mathfrak{a}    e^{ \z t_0} f   \Big(t_0,x_0,  \frac{1}{\mathfrak{a}}  e^{-\z t_0} \tilde{u} (t_0,x_0)  ,     \frac{1}{\mathfrak{a}}  e^{-\z t_0}    (\si^T \neg \cd \neg \nabla_x    \wt{\vf} )  (t_0,x_0)      \Big)   \le 0   .
    \eeas
    Hence, we have
    \beas
    \q   \min\bigg\{   \tilde{u}     (t_0,x_0)-\mathfrak{a} e^{\z t_0}  l(t_0,x_0),      -\frac{\pa  \wt{\vf}  }{\pa t}   (t_0,x_0) - \neg \cL \wt{\vf}  (t_0,x_0)
     - \neg \tilde{f}   \Big(t_0,x_0,  \tilde{u} (t_0,x_0)  ,         \nabla_x    \wt{\vf}  (t_0,x_0)         \Big)
        \bigg\}    \le  0   ,
    \eeas
    which means that $ \tilde{u}$  is a  viscosity subsolution of \eqref{eq:PDE2}.
    For the case of  viscosity supersolution, we can argue similarly. \qed

  Inspired by   Theorem 3.1 of \cite{Da_Lio_Ley_2010}, we have the following  comparison theorem, which together with
  Theorem \ref{thm:viscos_exist}  shows that   \eqref{eq:PDE} admits a unique viscosity solution.

\begin{thm} \label{thm:viscos_comparison}
    Suppose that there exists an increasing function $\fM: (0,\infty) \to (0,\infty)$ such that   for any $R > 0$,
    \bea  \label{cond:f_lip_x}
     \big|f(t,x,y,z) -f(t,x',y,z) \big| \le \fM(R) \big(1+|z|\big) |x-x'|
    \eea
     holds for any $(t,x,x',y,z) \in [0,T] \times \hR^k \times \hR^k \times \hR \times \hR^d$ with $|x| \vee |x'| \vee |y| \le R$.
    Let $u \in  C\big([0,T] \times \hR^k\big)$ \Big(resp. $v \in  C\big([0,T] \times \hR^k\big)$\Big) be a viscosity subsolution \(resp. viscosity supersolution\) of \eqref{eq:PDE}   such that for  some $\wt{\k} > 0 $,
  \bea  \label{eq:uv_varpi_growth}
   |u(t,x)| \vee |v(t,x)| \le \wt{\k} (1+|x|^\varpi), \q \fa (t,x) \in [0,T] \times \hR^k  .
  \eea
       Then $u(t,x) \le v(t,x)$ for all $(t,x) \in [0, T] \times \hR^k$.
 \end{thm}

 \ss \no {\bf Proof:}  For any $\th \in (0,1]$, we define
        \beas
        \tilde{u}_\th (t,x)   \dfnn     \th e^{\k t} u (t,x)  \q \hb{and} \q \tilde{v}_\th (t,x)   \dfnn     \th e^{\k t} v (t,x), \q  \fa  (t,x) \in [0,T] \times \hR^k.
        \eeas
        Lemma \ref{lem:tilde_u} shows that $ \tilde{u}_\th$ (resp. $\tilde{v}_\th$) is a viscosity subsolution \(resp. viscosity supersolution\) of \eqref{eq:PDE2} with $(\mathfrak{a},\z)=(\th, \k)$.

    \ms   Let $\l \dfnn     8  (b_0+\k) +4  (1+ 4 \g  e      )  \si^2_* + 2     (\a+    4 \k  \wt{\k} )  e^{\k T} $. Suppose that we have proven the following statement:
  \bea  \label{eq:f505}
     \ba{ll}\hb{For any $[T_1,T_2] \subset [0,T]$ with $T_2-T_1 \le \frac{1}{\l}$,  if $u(T_2,x) \le v(T_2,x)$,  $\fa  x \in \hR^k$, then} \vspace{2mm}\\
      \hb{\hspace{2.5cm}   $u(t,x) \le v(t,x)$, $\fa  (t,x) \in [T_1, T_2] \times \hR^k$.}
      \ea
  \eea
 Set $N \dfnn  \lceil \l T \rceil  $  and  $t_i \dfnn i\frac{T}{N}$, for $i=0,1,\cds, N$. Since  $u (T,x) \le h(x) \le v (T,x)$, $\fa x \in \hR^k$,   \eqref{eq:f505} shows that $u(t,x) \le v(t,x)$, $\fa  (t,x) \in [t_{N-1}, t_N] \times \hR^k$, in particular, $u(t_{N-1},x) \le v(t_{N-1},x)$, $\fa x \in \hR^k$. Again by \eqref{eq:f505}, we  have $u(t,x) \le v(t,x)$, $\fa  (t,x) \in [t_{N-2}, t_{N-1}] \times \hR^k$, in particular, $u(t_{N-2},x) \le v(t_{N-2},x)$, $\fa x \in \hR^k$. Iteratively, one can show that $u(t,x) \le v(t,x)$ for all $(t,x) \in [0, T] \times \hR^k$. Therefore, it suffices to show \eqref{eq:f505}.

   \ms  Assume that  \eqref{eq:f505} does not hold, i.e., there exists a time interval $[T_1,T_2] \subset [0,T]$ with $T_2-T_1 \le \frac{1}{\l}$ such that  $u(T_2,x) \le v(T_2,x)$,  $\fa  x \in \hR^k$ and that
       $            u(\hat{t},\hat{x})-  v(\hat{t},\hat{x}) > \d $
              for some $(\hat{t},\hat{x}) \in [T_1,T_2) \times \hR^k$ and some $\d>0$. By the continuity of functions $u$ and $v$, we may assume that
        $\hat{t} > T_1$.

         \bs We divide the proof into two cases.

    \ss  \no {\bf Case 1:}  The mapping $z \to f(t,x,y,z)$ is
   convex for   all $(t, x, y) \in   [0,T] \times \hR^k    \times \hR $.

  \ss \no  We fix a  $\, \th \in (0,1)$  such that
        \bea   \label{eq:f521}
      \big|     e^{\k \hat{t}}  u(\hat{t},\hat{x}) \big|   \vee  \big|  e^{\k \hat{t}}   v(\hat{t},\hat{x}) \big|   \vee   e^{\l(T_2-\hat{t})}       (1+2| \hat{x} |^2 )   < \frac{\d}{4(1-\th)} ,
        \eea
 and fix a $\varrho \in \big(0, \frac{\d}{4}   (\hat{t}-T_1)\big) $.   For any $\e >0  $, we define
  \beas
           \F_\e(t, x, x') & \dneg \dneg\dfnn & \dneg \dneg \frac{\varrho}{t-T_1} + e^{\l(T_2-t)} \left( \frac{  |x-x'|^2}{\e} +  (1-\th) (1+|x|^2+|x'|^2) \right)   \q \fa  t  \in (T_1,T_2], ~ \fa x,x' \in \hR^k   ,  \\
 \hb{and}   \q    M_\e & \dneg \dneg \dfnn & \dneg \dneg \underset{(t,x,x')  \in (T_1,T_2] \times  \hR^k \times \hR^k   }{\sup} \big\{  \tilde{u}_1 (t,x)  -\tilde{v}_\th      (t,x')    -   \F_\e(t,x,x')\big\} .
  \eeas

  \ss
    Since $\dis r^2 \ge \frac{4 \wt{\k} e^{\k T} }{1-\th}  (1+r^\varpi)$ holds for any $\dis r \ge R_\th \dfnn 1 \vee \bigg(\frac{8 \wt{\k} e^{\k T} }{1-\th}\bigg)^{\frac{1}{2-\varpi}} $,   \eqref{eq:uv_varpi_growth} shows that   for any $(t,x,x')  \in [T_1,T_2] \times \hR^k \times \hR^k $ with $   |x|\vee  |x'|  \ge R_\th$
    \bea
        \tilde{u}_1 (t,x)  -  \tilde{v}_\th (t,x')  &\le& e^{\k T}  \big(|u (t,x) |+|v(t,x') | \big)  \le 2\wt{\k} e^{\k T} \big(1+ (|x|\vee  |x'|)^\varpi \big)  \le \frac12 (1-\th)   \big(|x|\vee  |x'|\big)^2  \nonumber \\
     &\le&  \frac12  e^{\l(T_2-t)}   (1-\th) \big(1+|x|^2+|x'|^2\big), \label{eq:f533}
    \eea
    which implies that
     \beas    
        \lmt{ \frac{1}{t-T_1}\vee |x|\vee  |x'|  \to \infty  }      \big( \tilde{u}_1 (t,x)  -  \tilde{v}_\th (t,x')    -   \F_\e(t,x,x') \big) = -\infty .
     \eeas
    Hence,  one can deduce that    the supremum  $M_\e$    is finite and attainable at some $(t_\e,x_\e, x'_\e) \in   (T_1,T_2] \times  \hR^k \times \hR^k$.
      \if{0}
       For any $n \in \hN$, there exists some $(t_n,x_n, x'_n) =\big(t_n(\e),x_n(\e), x'_n(\e)\big) \in   (T_1,T_2] \times  \hR^k \times \hR^k$ such that
     \beas
     \tilde{u}_1 (t_n,x_n)  -  \tilde{v}_\th (t_n,x'_n)    -   \F_\e(t_n,x_n,x'_n) > \left\{ \ba{ll}  M_\e -\frac{1}{n}  & \q \hb{if }  M_\e \in (-\infty, \infty),  \\
     n   & \q \hb{if }  M_\e =\infty.  \ea \right.
     \eeas
     Since $ \lmt{n \to \infty}  \big( \tilde{u}_1 (t_n,x_n)  -  \tilde{v}_\th (t_n,x'_n)    -   \F_\e(t_n,x_n,x'_n)\big) =M_\e > - \infty$,     \eqref{eq:f527} implies that $\big\{(t_n, x_n,x'_n )\big\}_{n \in \hN}$ does not have a subsequence $\big\{(t_{n_i}, x_{n_i},x'_{n_i} )\big\}_{j \in \hN}$ such that   $\lmt{j \to \infty} \Big( \frac{1}{t_{n_i}-T_1}\vee |x_{n_i}|\vee  |x'_{n_i}|  \Big) = \infty$,
       which further implies that $\{t_n\}_{n \in \hN} \subset [\,\wt{T}_1 ,T_2]$ for some $\wt{T}_1 = \wt{T}_1 (\e) \in (T_1, T_2)$ and that $\{x_n\}_{n \in \hN}$, $\{x'_n\}_{n \in \hN}$
       are both bounded sequences in $\hR^k$. Hence, $\{(t_n, x_n,x'_n)\}_{n \in \hN}   $ has a convergent subsequence $\big\{(t_{n_i}, x_{n_i},x'_{n_i} )\big\}_{j \in \hN}$ with the limit $(t_\e,x_\e, x'_\e) \in [\,\wt{T}_1 ,T_2] \times \hR^k \times \hR^k$.

       \ms By the continuity of functions $u$ and $v$,
      \bea
      M_\e &=&  \lmt{j \to \infty}  \big(\tilde{u}_1 (t_{n_i}, x_{n_i})  -  \tilde{v}_\th (t_{n_i},x'_{n_i})    -   \F_\e(t_{n_i},x_{n_i},x'_{n_i})\big)   \nonumber\\
       &=&   \tilde{u}_1 (t_\e,x_\e)  -    \tilde{v}_\th(t_\e,x'_\e)    -   \F_\e(t_\e,x_\e,x'_\e) < \infty   .  \label{eq:f515}
     \eea
     \fi
 Then   it follows from \eqref{eq:f521} that
  \bea  \label{eq:f529}
     \qq   \qq   && \hspace{-2cm}   \tilde{u}_1 (t_\e,x_\e)  -    \tilde{v}_\th(t_\e,x'_\e)    \neg   -   \neg   \F_\e(t_\e,x_\e,x'_\e) =M_\e   \ge  \tilde{u}_1 (\hat{t} ,\hat{x})    \neg -    \neg   \tilde{v}_\th(\hat{t} ,\hat{x})     \neg  -    \neg  \frac{\varrho}{\hat{t}   \neg -  \neg T_1}  \neg  -  \neg  e^{\l(T_2    -    \hat{t})}    (1   \neg -   \neg  \th) (1  \neg +  \neg 2|\hat{x}|^2 )   \nonumber    \\
       && \ge        u (\hat{t} ,\hat{x})  \neg -    \neg     v(\hat{t} ,\hat{x})   \neg+ \neg (1\neg - \neg \th) e^{\k \hat{t}}\,  v(\hat{t} ,\hat{x})        \neg  -    \neg  \frac{\varrho}{\hat{t}   \neg -  \neg T_1}  \neg  -  \neg  e^{\l(T_2    -    \hat{t})}    (1   \neg -   \neg  \th) (1  \neg +  \neg 2|\hat{x}|^2 )   > \frac{\d}{4} ,
  \eea
  which implies that
  \bea   \label{eq:f535}
    \frac{\d}{4}  +  e^{\l(T_2-t)} \left( \frac{  |x_\e-x'_\e|^2}{\e} +  (1-\th) \big(1+|x_\e|^2+|x'_\e|^2\big)  \right)  <  \tilde{u}_1 (t_\e,x_\e)  -    \tilde{v}_\th(t_\e,x'_\e)      .
  \eea
  Hence, we see  from \eqref{eq:f533}  that
  \bea
  |x_\e|\vee  |x'_\e|  < R_\th  . \label{eq:f539}
  \eea
  As      $ \Big\{(t_\e,x_\e,x'_\e): \e > 0   \Big\}
   \subset (T_1,T_2] \times B_{R_\th}(0)  \times B_{R_\th}(0) $,
   we can pick up a sequence $\{ \e_n  \}_{n \in \hN} \subset  (0, \infty   )$ with $\lmtd{n \to \infty} \e_n =0 $ such  that the sequence $\big\{(t_{ \e_n},x_{ \e_n},x'_{ \e_n}) \big\}_{n \in \hN}$ converges to some $(t_*,x_*,x'_*) \in [T_1,T_2] \times \ol{B}_{R_\th}(0)  \times \ol{B}_{R_\th}(0)$.  Then  \eqref{eq:f529} and the continuity of function $u$ and $v$ imply that
    \beas
    \lmt{n \to \infty} \,\frac{\varrho}{t_{ \e_n}\dneg- \neg T_1}  \le \lsup{n \to \infty}  \F_{ \e_n}\neg(t_{ \e_n}\neg,x_{ \e_n}\neg,x'_{ \e_n}\neg)  \le
      \tilde{u}_1 (t_*,x_*)  \neg -   \neg  \tilde{v}_\th(t_*,x'_*) \neg -   \neg \frac{\d}{4} < \infty,
    \eeas
    which implies that $t_* =\lmt{n \to \infty}t_{ \e_n} > T_1$, i.e., $t_* \in (T_1,T_2]$.

   \ms One can  also deduce from \eqref{eq:f535} that
   \beas
     \lsup{n \to \infty}  \frac{  |x_{\e_n}  \neg-  \neg x'_{\e_n}\neg |^2}{\e_n}   \le  \tilde{u}_1 (t_*,x_*)     -    \tilde{v}_\th(t_*,x'_*) < \infty   ,
   \eeas
     which   leads to that     $  \lmt{n \to \infty}   |x_{\e_n}  \neg-  \neg x'_{\e_n}\neg |  =0$, namely,     $ x_* = x'_*$.         For any $n \in \hN$,
    \beas
      && \hspace{-2cm}     \tilde{u}_1 (t_{\e_n}\neg,x_{\e_n}\neg)  \neg -   \neg    \tilde{v}_\th(t_{\e_n}\neg,x'_{\e_n}\neg)    \neg   -   \neg   \F_{\e_n}          \neg   (t_{\e_n}\neg,x_{\e_n}\neg,x'_{\e_n}\neg) =M_{\e_n}  \\
     &&        \ge   \tilde{u}_1 (t_* ,x_*\neg)    \neg -    \neg   \tilde{v}_\th(t_* ,x_*\neg)     \neg  -    \neg  \frac{\varrho}{t_*   \neg -  \neg T_1}  \neg  -  \neg  e^{\l(T_2    -    t_*)}    (1   \neg -   \neg  \th) (1  \neg +  \neg 2|x_*|^2 )    .
  \eeas
 As $n \to \infty$,  the continuity of functions $u$ and $v$ implies that
  \bea  \label{eq:f551b}
     \lmt{n \to \infty}  \frac{  |x_{\e_n}  \neg-  \neg x'_{\e_n}|^2}{\e_n}   = 0.
  \eea

    \ss  Now we claim that
        \bea
          \hb{$\{\e_n\}_{n \in \hN} $ has a subsequence $\{\wt{\e}_n\}_{n \in \hN} $ such that for any $n \in \hN$,  either }   t_{\wt{\e}_n}=T_2
           \hb{ or }   u(t_{\wt{\e}_n}\neg,x_{\wt{\e}_n}\neg) \le l(t_{\wt{\e}_n}\neg,x_{\wt{\e}_n}\neg).  \q   \label{claim3}
        \eea
 Assume the contrary. Then   there exists an $n^o \in \hN$ such that for any $n \ge n^o$,  $t_{\e_n} \in (T_1,T_2)    $ and  $u(t_{\e_n}\neg,x_{\e_n}\neg) > l(t_{\e_n}\neg,x_{\e_n}\neg)$.

  \ms  Fix $n \ge n^o$. The continuity of functions  $u$ and $l$ shows that $(t_{\e_n}\neg,x_{\e_n}\neg)$ has an open neighborhood
   $ \cO_n \dfnn  (t_{\e_n}\neg-   r_n, t_{\e_n}\neg+   r_n) \times B_{r_n}(x_{\e_n}\neg) \subset (T_1,T_2) \times \hR^k$ for some $r_n > 0$ such that  $ u(t,x)  >   l(t,x)  $ for any $  (t,x) \in  \cO_n $. Then $\tilde{u}(t,x)  $   becomes a viscosity subsolution of \eqref{eq:PDE2}  without obstacle and terminal condition over $ \cO_n$, i.e.
    \bea
      -\frac{\pa \tilde{u} }{\pa t}   (t,x) - \neg \cL \tilde{u} (t,x) + \k \tilde{u}(t,x)  -   \neg e^{\k t}  f   \Big(t,x,    e^{-\k t}   \tilde{u} (t,x)   ,    e^{-\k t}     (\si^T\neg \cd \neg \nabla_x \tilde{u}) (t,x)      \Big)      \neg     =0 ,        ~    \fa (t,x) \in    \cO_n .  \q  \label{eq:f545a}
    \eea
    As $\tilde{v}_\th  $ is  a viscosity supersolution of \eqref{eq:PDE2}, it  is clearly a  viscosity supersolution of \eqref{eq:PDE2} without obstacle and terminal condition over $(0,T) \times \hR^k$ \big(thus over $ \cO'_n \dfnn  (t_{\e_n}\neg-r_n, t_{\e_n}\neg+r_n) \times B_{r_n}(x'_{\e_n}\neg)$\big), i.e.
  \bea
        -  \frac{\pa  \tilde{v} }{\pa t}   (t,x)  - \neg \cL \tilde{v} (t,x)   + \k \tilde{v} (t,x)   -  \th  e^{\k t}   f   \Big(t,x,  \hb{$\frac{1}{\th}$}   e^{-\k t}    \tilde{v} (t,x) ,
           \hb{$\frac{1}{\th}$}  e^{-\k t}   (\si^T\neg \cd \neg \nabla_x \tilde{v} ) (t,x)       \Big) \neg      =    0 ,    ~    \fa (t,x) \in  \cO'_n.    \q   \label{eq:f545b}
    \eea
    Since the mapping $(t,x,x')  \to \tilde{u}_1 (t,x) - \tilde{v}_\th (t,x')   -\F_{\e_n}(t,x,x') $ is maximized at $(t_{\e_n}\neg,x_{\e_n}\neg,x'_{\e_n})$ over
  $(T_1,T_2] \times \hR^k \times \hR^k$ \big(thus over $(t_{\e_n}\neg-   r_n, t_{\e_n}\neg+   r_n) \times B_{r_n}(x_{\e_n} ) \times B_{r_n}(x'_{\e_n} )$\big),      Theorem 8.3 of \cite{User_1992} shows that there exist $p_n, p'_n \in \hR$ and $W_n, W'_n \in \hS^k $   such that
    \bea
      &\tneg \neg &  \hspace{-2cm}    \;\, \q \left( p_n,  \nabla_x\F_{\e_n} (t_{\e_n}\neg,x_{\e_n}\neg, x'_{\e_n}\neg), W_n \right) \in \ol{\cP}^{2,+} \tilde{u}_1 (t_{\e_n}\neg,x_{\e_n}\neg),     \qq \qq \label{eq:f531a}\\
      &\tneg \neg &  \hspace{-2cm}    \q     \left(  p'_n,  - \nabla_{x'}\F_{\e_n} (t_{\e_n}\neg,x_{\e_n}\neg,x'_{\e_n}\neg), W'_n  \right) \in \ol{\cP}^{2,-} \tilde{v}_\th (t_{\e_n}\neg,x'_{\e_n}\neg),  \qq \qq \label{eq:f531b}\\
   p_n - p'_n  &\tneg \neg =&     \tneg \neg    \frac{\pa\F_{\e_n}}{\pa t}(t_{\e_n}\neg,x_{\e_n}\neg,x'_{\e_n}\neg)
     =-\frac{\varrho}{(t_{\e_n}\dneg- \neg T_1)^2} - \l  \F_{\e_n} (t_{\e_n}\neg, x_{\e_n}\neg) ,  \label{eq:f531c} \\
     \hb{and}   \q  \left( \neg  \ba{cc} W_n &0 \\ 0& - W'_n  \ea \dneg \right)  &\tneg \neg \le & \tneg \neg   D^2_{x,x'}\F_{\e_n} (t_{\e_n}\neg,x_{\e_n}\neg,x'_{\e_n}\neg) + \e^3_n
       \Big( D^2_{x,x'}\F_{\e_n} (t_{\e_n}\neg,x_{\e_n}\neg,x'_{\e_n}\neg)\Big)^2 . \qq \qq \qq
                  \label{eq:f531d}
    \eea
          As  $\tilde{u}_1$  is a viscosity subsolution of \eqref{eq:f545a},  one can deduce from \eqref{eq:f531a}    that
     \bea
     &\dneg \dneg& \dneg \dneg          -     p_n
      \neg -\neg \frac12   \hb{trace}\big(W_n \neg \cd \neg (\si \si^T) (t_{\e_n}\neg, x_{\e_n}\neg)\big)
         \neg   -  \neg  2 e^{\l(T_2  -  t_{\e_n}\neg)}  \bigg\lan \neg b(t_{\e_n}\neg, x_{\e_n}\neg),     \frac{x_{\e_n} \dneg- \neg x'_{\e_n}}{\e_n}  \neg
         + \neg (1 \neg - \neg \th) x_{\e_n} \neg \bigg\ran    \neg +\neg \k  e^{\k t_{\e_n}} u (t_{\e_n}\neg,x_{\e_n} \neg )           \nonumber  \\
        &\dneg \dneg&   \dneg \dneg   \qq     -      e^{\k t_{\e_n}}        f   \left(t_{\e_n}\neg, x_{\e_n}\neg,     u (t_{\e_n}\neg,x_{\e_n} \neg )  ,     2 e^{-\k t_{\e_n} +\l(T_2  -  t_{\e_n}\neg)}        \si^T(t_{\e_n}\neg,x_{\e_n} \neg) \cd \left(     \frac{x_{\e_n} \dneg- \neg x'_{\e_n}}{\e_n}  \neg
         + \neg (1\neg - \neg \th) x_{\e_n} \neg     \right)  \right)  \le 0  .  \hspace{1cm}  \label{eq:f537}
         \eea
        Since $\tilde{v}_\th$ is a viscosity supersolution of \eqref{eq:f545b},  it follows from  \eqref{eq:f531b}   that
       \bea
       &\dneg \dneg& \dneg \dneg         -  p'_n     \neg   -  \neg  \frac12 \hb{trace}\big(W'_n  \neg \cd  \neg  (\si \si^T) (t_{\e_n}\neg, x'_{\e_n}\neg)\big)  \neg   -  \neg   2 e^{\l(T_2  -  t_{\e_n}\neg)} \bigg\lan \neg b(t_{\e_n}\neg, x'_{\e_n}\neg),
           \frac{x_{\e_n} \dneg- \neg x'_{\e_n}}{\e_n}  \neg         - \neg (1\neg - \neg \th) x'_{\e_n}       \neg  \bigg\ran \neg +\neg \th \k  e^{\k t_{\e_n}} \, v(t_{\e_n}\neg,x'_{\e_n} \neg )
                \nonumber  \\
        &\dneg \dneg& \dneg     \qq       -      \th      e^{\k t_{\e_n}}   f   \left(t_{\e_n}\neg,x'_{\e_n}\neg,  v (t_{\e_n}\neg,x'_{\e_n}\neg )  ,      \frac{2}{\th}  e^{-\k t_{\e_n} + \l(T_2  -  t_{\e_n}\neg)}   \si^T(t_{\e_n}\neg,x'_{\e_n}\neg ) \cd\left(     \frac{x_{\e_n} \dneg- \neg x'_{\e_n}}{\e_n}  \neg
         - \neg (1\neg - \neg \th) x'_{\e_n} \neg          \right)   \right)  \ge 0   .   \qq   \qq    \label{eq:f541}
  \eea
 Subtracting \eqref{eq:f541} from  \eqref{eq:f537}, we see from \eqref{eq:f531c}    that
     \bea
   \frac{\varrho}{(t_{\e_n}\dneg- \neg T_1)^2} + \l  \F_{\e_n} (t_{\e_n}\neg, x_{\e_n}\neg)   \le    I^1_n + 2 e^{\l(T_2  -  t_{\e_n}\neg)}  I^2_n+   e^{\k t_{\e_n}}    \sum^6_{j=3} I^j_n    , \q   \label{eq:f571}
          \eea
          where
       \beas
         I^1_n & \dneg \dneg \dfnn&  \dneg \dneg    \frac12 \hb{trace}\big(W_n \cd (\si \si^T) (t_{\e_n}\neg, x_{\e_n}\neg)\big)  -  \frac12 \hb{trace}\big(W'_n \cd (\si \si^T) (t_{\e_n}\neg, x'_{\e_n}\neg)\big)    ,   \\
        I^2_n & \dneg \dneg \dfnn&  \dneg \dneg     \bigg\lan \neg b(t_{\e_n}\neg, x_{\e_n}\neg)- b(t_{\e_n}\neg, x'_{\e_n}\neg),     \frac{x_{\e_n} \dneg- \neg x'_{\e_n}}{\e_n}    \neg \bigg\ran
             + (1 \neg - \neg \th)  \Big( \big\lan   b(t_{\e_n}\neg, x_{\e_n}\neg),           x_{\e_n}\neg  \big\ran
           +     \big\lan   b(t_{\e_n}\neg, x'_{\e_n}\neg),           x'_{\e_n}\neg  \big\ran \Big) ,   \\
             I^3_n & \dneg \dneg \dfnn&  \dneg \dneg        -  \k   u (t_{\e_n}\neg,x_{\e_n} \neg )   + \th   \k    v(t_{\e_n}\neg,x'_{\e_n} \neg )   ,   \\
         I^4_n &  \dneg \dneg  \dfnn&    \dneg \dneg            \left[  f   \Big(t_{\e_n}\neg,x_{\e_n}\neg,     u (t_{\e_n}\neg,x_{\e_n}  \neg)  ,        J_n    \Big)
         -          f   \Big(t_{\e_n}\neg,x_{\e_n}\neg,     v (t_{\e_n}\neg,x'_{\e_n}\neg)  ,        J_n    \Big) \right],  \hb{\; with }     \\
       & \dneg \dneg &  \dneg \dneg    \q     J_n \dfnn     2 e^{-\k t_{\e_n} +\l(T_2  -  t_{\e_n}\neg)}        \si^T(t_{\e_n}\neg,x_{\e_n} \neg) \cd \left(     \frac{x_{\e_n} \dneg- \neg x'_{\e_n}}{\e_n}  \neg         + \neg (1\neg - \neg \th) x_{\e_n} \neg     \right),    \\
         I^5_n & \dneg \dneg  \dfnn&  \dneg \dneg        \left[                   f   \Big(t_{\e_n}\neg,x_{\e_n}\neg,     v (t_{\e_n}\neg,x'_{\e_n}\neg)  ,
                        J_n    \Big)   -  f   \Big(t_{\e_n}\neg, x'_{\e_n},     v (t_{\e_n}\neg,x'_{\e_n}\neg)  ,        J_n    \Big) \right] , \\
        I^6_n & \dneg \dneg \dfnn&    \dneg \dneg          f   \Big(t_{\e_n}\neg, x'_{\e_n}\neg,     v (t_{\e_n}\neg,x'_{\e_n}\neg)  ,        J_n    \Big)
               \neg   -  \neg      \th      f   \neg   \left(t_{\e_n}\neg,x'_{\e_n}\neg,  v (t_{\e_n}\neg,x'_{\e_n}\neg )  ,     \frac{1}{\th}   J'_n   \right)  , \hb{\; with }     \\
       & \dneg \dneg &  \dneg \dneg    \q     J'_n   \dfnn  2   e^{-\k t_{\e_n}  \neg  + \l(T_2  -  t_{\e_n}\neg)}   \si^T(t_{\e_n}\neg,x'_{\e_n}\neg ) \neg  \cd \neg   \left(     \frac{x_{\e_n} \dneg- \neg x'_{\e_n}}{\e_n}  \neg
         - \neg (1\neg - \neg \th) x'_{\e_n} \neg          \right)    .
       \eeas
     $\bullet$    One can deduce from \eqref{eq:f531d}  and \eqref{b_si_Lip} that
  \bea
     I^1_n  & \dneg =&  \dneg \frac12  \left( \dneg \ba{c}  \si (t_{\e_n}\neg, x_{\e_n}\neg)   \\   \si (t_{\e_n}\neg, x'_{\e_n}\neg)     \ea \dneg \right)^T   \left( \ba{cc} W_n &0 \\ 0& - W'_n  \ea \dneg \right) \left( \dneg \ba{c}  \si (t_{\e_n}\neg, x_{\e_n}\neg)   \\   \si (t_{\e_n}\neg, x'_{\e_n}\neg)     \ea \dneg \right) \nonumber \\
     &  \dneg \le &  \dneg      \Big(  \frac{1}{\e_n}  e^{\l(T_2-t_{\e_n})} \neg  +  \neg   4  \e_n   e^{2 \l(T_2-t_{\e_n})} \neg  +  \neg  4 \e^2_n (1 \neg  - \neg  \th) e^{2 \l(T_2-t_{\e_n})} \Big)     \big|\si (t_{\e_n}\neg,x_{\e_n}\neg)-\si (t_{\e_n}\neg,x'_{\e_n}\neg) \big|^2  \nonumber \\
       &\tneg \neg   &     +  \,    \Big(      (1  \neg  - \neg   \th)  e^{\l(T_2-t_{\e_n})}  \neg  + \neg   2 \e^3_n  (1  \neg  -  \neg   \th)^2  e^{2 \l(T_2-t_{\e_n})} \Big) \left( \big|\si (t_{\e_n}\neg,x_{\e_n}\neg)\big|^2+ \big| \si (t_{\e_n}\neg,x'_{\e_n}\neg) \big|^2 \right)
            \nonumber \\
    &  \dneg \le &  \dneg   e    \k^2  \frac{       | x_{\e_n} \dneg - \neg x'_{\e_n}|^2 }{\e_n} + 2  (1  \neg  - \neg   \th) e^{\l(T_2-t_{\e_n})}        \si^2_*  +  c_{\si_*} \big( \e_n + \e^2_n+\e^3_n \big) .  \label{eq:f573}
  \eea

   \ms \no  $\bullet$ It follows from \eqref{b_si_Lip} that
   \bea
     I^2_n  \le      \k   \frac{       | x_{\e_n} \dneg - \neg x'_{\e_n}|^2 }{\e_n} + (1 \neg - \neg  \th)   \big(b_0 |x_{\e_n}\neg| +b_0 |x'_{\e_n}\neg| +\k |x_{\e_n}\neg|^2+\k|x'_{\e_n}\neg|^2\big).
   \eea

   \ms \no  $\bullet$  We see from \eqref{eq:f535} that $ u (t_{\e_n}\neg,x_{\e_n}\neg)  -     \th v(t_{\e_n}\neg,x'_{\e_n}\neg) >0$. Then    \eqref{cond:f_basic} shows that
      \beas
       I^4_n       \le      \k      \big|u (t_{\e_n}\neg,x_{\e_n}  \neg)- v (t_{\e_n}\neg,x'_{\e_n}\neg )\big|
      \le  \k      \big(u (t_{\e_n}\neg,x_{\e_n}  \neg)-\th v (t_{\e_n}\neg,x'_{\e_n}\neg )\big) +  \k  (1-\th)    \big|  v (t_{\e_n}\neg,x'_{\e_n}\neg )  \big|     .
      \eeas
      Thus,
      \bea
       I^3_n + I^4_n       \le \k  (1-\th)    \big|  v (t_{\e_n}\neg,x'_{\e_n}\neg )  \big|  .
      \eea

  \ss \no  $\bullet$  \eqref{eq:f539} and \eqref{eq:uv_varpi_growth} imply  that
    $\underset{i \in \hN}{\sup} \big\{| x_{\e_i}| \vee |x'_{\e_i}| \vee \big|v (t_{\e_i}\neg,x'_{\e_i}\neg)\big| \big\}
   \neg \le \neg  \wt{R}_\th \dfnn  (1 \vee   \wt{\k}) \big(1 \neg + \neg  R^\varpi_\th \big)$. Then \eqref{cond:f_lip_x} shows that
   \bea
    I^5_n                          \le     \fM(\wt{R}_\th) \big(1+ |J_n|\big)|x_{\e_n}\neg -x'_{\e_n}  |
    \le  \fM(\wt{R}_\th)  \Big(1  \neg +2 e   \si_*        \frac{ |x_{\e_n} \dneg- \neg x'_{\e_n}|}{\e_n}   \neg +2 e   \si_*  R_\th
       \Big) |x_{\e_n}\neg -x'_{\e_n}  |   .
   \eea

   \ms \no  $\bullet$ The convexity of the mapping
   $z \to   f   \Big(t_{\e_n}\neg,x'_{\e_n}\neg,     v (t_{\e_n}\neg,x'_{\e_n}  \neg)  , z \Big)$,
     \eqref{cond:f_basic}    implies    that
  \beas
     I^6_n  & \tneg \le & \tneg    (1- \th)     f   \bigg(t_{\e_n}\neg,x'_{\e_n}\neg,     v (t_{\e_n}\neg,x'_{\e_n}  \neg)  ,    \frac{J_n-J'_n}{1\neg- \neg \th}     \bigg)
            \le   (1- \th)     f   \bigg(t_{\e_n}\neg,x'_{\e_n}\neg,     0  ,      \frac{J_n-J'_n}{1\neg- \neg \th}   \bigg)
          +  \k (1- \th)     | v (t_{\e_n}\neg,x'_{\e_n}  \neg) |       ,
  \eeas
  where
    \beas
   \q    \frac{J_n-J'_n}{1\neg- \neg \th}     =  \frac{2     e^{-\k t_{\e_n}  \neg  +\neg  \l(T_2  -  t_{\e_n}\neg)} }{1 - \th} \neg    \bigg(   \neg  \Big\lan \si(t_{\e_n}\neg,x_{\e_n} \neg) \neg- \neg \si(t_{\e_n}\neg,x'_{\e_n}\neg ) ,       \frac{x_{\e_n} \dneg- \neg x'_{\e_n}}{\e_n}    \neg   \Big\ran \neg + \neg   (1 \neg - \neg \th)  \Big(     \si^T \neg (t_{\e_n}\neg, x_{\e_n}\neg) \neg \cd \neg           x_{\e_n}\dneg
           +  \neg \si^T \neg (t_{\e_n}\neg, x'_{\e_n}\neg) \neg \cd \neg           x'_{\e_n}\neg  \Big)   \neg    \bigg).
        \eeas
   Since $ \dis \left|\Big\lan \si(t_{\e_n}\neg,x_{\e_n} \neg)-\si(t_{\e_n}\neg,x'_{\e_n}\neg ) ,       \frac{x_{\e_n} \dneg- \neg x'_{\e_n}}{\e_n}       \Big\ran  \right| \le   \k  \frac{       | x_{\e_n} \dneg - \neg x'_{\e_n}|^2 }{\e_n} $, \eqref{eq:f551b} and the continuity of function $\si$ that
  \bea  \label{eq:f579}
   \lmt{n \to \infty}   \frac{J_n-J'_n}{1\neg- \neg \th}   =        4     e^{-\k t_*    +\neg  \l(T_2  -  t_*)}  \si^T \neg (t_*, x_*)   \cd   x_*  .
  \eea
 Letting $n \to \infty$ in \eqref{eq:f571} and using the continuity of all functions involved, we can deduce  from \eqref{eq:f551b}, \eqref{eq:f573} through \eqref{eq:f579}  that
    \bea
      \l  (1 \neg-\neg \th) e^{\l(T_2-t_*)}    (1\neg +\neg  2  |x_*|^2 )  & \dneg \dneg \le& \dneg \dneg    2  (1  \neg  - \neg   \th) e^{\l(T_2-t_*)}        \si^2_*
     \neg + \neg  4   (1 \neg - \neg  \th) (b_0+\k) e^{\l(T_2  -  t_*\neg)} \big(1 \neg + \neg   |x_* |^2 \big) \neg + \neg 2 \k e^{\k t_*}     (1 \neg-\neg \th)    \big|  v (t_* ,x_*\neg )  \big|  \nonumber \\
&\dneg \dneg &\dneg \dneg   +\,  e^{\k t_*}   (1 \neg -\neg \th)     f   \Big(t_* ,x_* ,     0  ,      4     e^{-\k t_*    +\neg  \l(T_2  -  t_*)}  \si^T \neg (t_*, x_*)   \cd   x_*   \Big)  . \label{eq:f561}
    \eea
    Conditions \eqref{eq:uv_varpi_growth} and \eqref{cond:f_basic} imply that
     \beas
          &&   \hspace{-1.5cm}  2 \k e^{\k t_*}       \big|  v (t_* ,x_*\neg )  \big|    +\,  e^{\k t_*}         f   \Big(t_* ,x_* ,     0  ,      4     e^{-\k t_*    +\neg  \l(T_2  -  t_*)}  \si^T \neg (t_*, x_*)   \cd   x_*   \Big) \\
    &&  \le     2 \k \wt{\k} e^{\k t_*} (1+ |x_*|^\varpi) +    e^{\k t_*}  \left(\a +    8 \g  e^{-2\k t_*    +   2 \l(T_2  -  t_*)}  \si^2_*     | x_* |^2  \right) \\
    &&  \le  ( \a + 4 \k \wt{\k} )  e^{\k T} (1+ |x_*|^2) +         8 \g     e^{- \k t_*  + 1  +     \l(T_2  -  t_*)}  \si^2_*     | x_* |^2 .
     \eeas
     Plugging it back into \eqref{eq:f561} yields that
      \beas
       \l  (1 \neg-\neg \th) e^{\l(T_2-t_*)}    (1\neg +\neg  2   |x_*|^2 ) &\le &   (1 \neg-\neg \th) e^{\l(T_2-t_*)}    (1\neg +\neg     |x_*|^2 )  \Big( 4  (b_0+\k) +2  (1+ 4 \g  e      )  \si^2_* +     (\a+    4 \k  \wt{\k} )  e^{\k T} \Big)    \\
       & = &      \frac12  \l  (1 \neg-\neg \th) e^{\l(T_2-t_*)}    (1\neg +\neg  2  |x_*|^2 ) ,
      \eeas
      which results in a contradiction.   Thus we proved  claim \eqref{claim3}. Let $\{\wt{\e}_n\}_{n \in \hN} $ be the  subsequence of $\{\e_n\}_{n \in \hN} $
   as described in \eqref{claim3}. For any $n \in \hN$, since the maximum is attained at $(t_{\wt{\e}_n} \neg, x_{\wt{\e}_n}\neg, x'_{\wt{\e}_n}\neg)$,
      \bea   \label{eq:f565}
       \tilde{u} (\hat{t},\hat{x})    \neg -    \neg   \th \tilde{v}(\hat{t},\hat{x})     \neg  -    \neg  \frac{\varrho}{\hat{t}   \neg -  \neg T_1}  \neg  -  \neg  e^{\l(T_2    -    \hat{t})}    (1   \neg -   \neg  \th) (1  \neg +  \neg 2|\hat{x}|^2 ) \le  M_{\wt{\e}_n} \neg    \le \tilde{u} (t_{\wt{\e}_n} \neg, x_{\wt{\e}_n}\neg)  \neg -   \neg    \th \tilde{v}(t_{\wt{\e}_n}\neg,x'_{\wt{\e}_n}\neg)      .
      \eea
      If $ t_{\wt{\e}_n}=T_2$,  $ u (t_{\wt{\e}_n} \neg, x_{\wt{\e}_n}\neg) =u(T_2,x_{\wt{\e}_n}\neg) \le v(T_2,x_{\wt{\e}_n}\neg)   =v (t_{\wt{\e}_n} \neg, x_{\wt{\e}_n}\neg)$ by our condition. Otherwise,   $  t_{\wt{\e}_n} \in (T_1,T_2)$ and  $u(t_{\wt{\e}_n}\neg,x_{\wt{\e}_n}\neg) \le l(t_{\wt{\e}_n}\neg,x_{\wt{\e}_n}\neg)  $. As $v$ is a viscosity  supersolution  of \eqref{eq:PDE}, we always have
       $        v  (t_{\wt{\e}_n}\neg,x_{\wt{\e}_n}\neg)-     l (t_{\wt{\e}_n}\neg,x_{\wt{\e}_n}\neg)     \ge 0$. Thus we still have $u(t_{\wt{\e}_n}\neg,x_{\wt{\e}_n}\neg) \le   v  (t_{\wt{\e}_n}\neg,x_{\wt{\e}_n}\neg)$.  Then \eqref{eq:f565} reduces to
        \beas
         \tilde{u} (\hat{t},\hat{x})    \neg -    \neg   \th \tilde{v}(\hat{t},\hat{x})     \neg  -    \neg  \frac{\varrho}{\hat{t}   \neg -  \neg T_1}  \neg  -  \neg  e^{\l(T_2    -    \hat{t})}    (1   \neg -   \neg  \th) (1  \neg +  \neg 2|\hat{x}|^2 )     \le \tilde{v} (t_{\wt{\e}_n} \neg, x_{\wt{\e}_n}\neg)  \neg -   \neg    \th \tilde{v} (t_{\wt{\e}_n}\neg,x'_{\wt{\e}_n}\neg)      .
        \eeas
              As $n \to \infty$, we obtain
       \beas
         \tilde{u} (\hat{t},\hat{x})    \neg -    \neg   \th \tilde{v} (\hat{t},\hat{x})   \neg  -    \neg  \frac{\varrho}{\hat{t}  \neg -  \neg T_1}   \neg  -  \neg  e^{\l(T_2    -    \hat{t})}    (1   \neg -   \neg  \th) (1  \neg +  \neg 2|\hat{x}|^2 ) \le  (1-\th) \tilde{v} (t_*, x_*).
       \eeas
       Letting $\varrho \to 0$ and letting $\th \to 1$  yield  that $\tilde{u} (\hat{t},\hat{x})    \neg -    \neg    \tilde{v} (\hat{t},\hat{x}) \le 0$,
       thus  $  u (\hat{t},\hat{x})  \le        v(\hat{t},\hat{x})    $, which contradicts with our initial assumption. Therefore,  \eqref{eq:f505}  holds.

         \ss  \no {\bf Case 2: }   The mapping $z \to f(t,x,y,z)$ is
   concave for   all $(t, x, y) \in   [0,T] \times \hR^k    \times \hR $.

  \ss \no This case is similar to Case 1, so we only sketch out the main differences:  We  redefine
   \beas
    M_\e & \dneg \dneg \dfnn & \dneg \dneg \underset{(t,x,x')  \in (T_1,T_2] \times  \hR^k \times \hR^k   }{\sup} \big\{  \tilde{u}_\th (t,x)  -\tilde{v}_1      (t,x')    -   \F_\e(t,x,x')\big\}
   \eeas
   for any $\e>0$, and change the forms of $I^3_n$ through $I^6_n$ correspondingly. For example,
    \beas
    I^6_n \dfnn \th   f   \left(t_{\e_n}\neg, x_{\e_n}\neg,     u (t_{\e_n}\neg,x_{\e_n}\neg)  ,   \frac{1}{\th}     J_n    \right)
               \neg   -  \neg             f   \neg   \left( t_{\e_n}\neg, x_{\e_n}\neg,     u (t_{\e_n}\neg,x_{\e_n}\neg)    ,         J'_n   \right) .
               \eeas
      Then the concavity of the mapping   $z \to   f   \Big(t_{\e_n}\neg,x_{\e_n}\neg,     u (t_{\e_n}\neg,x_{\e_n}  \neg)  , z \Big)$  implies that
      \beas
       I^6_n \le - (1  -  \th) f   \left(t_{\e_n}\neg, x'_{\e_n}\neg,     v (t_{\e_n}\neg,x'_{\e_n}\neg)  ,   \frac{J'_n-J_n}{1\neg- \neg \th}     \right).
       \eeas
       All other arguments used in Case 1 still work in this case with slight adaptions.    \qed

    \begin{rem} \label{rem_gap}
     Theorem~\ref{thm:viscos_comparison} is  not a  simple extension of Theorem 3.1 of \cite{Da_Lio_Ley_2010}. In fact,  there are two gaps in the proof of Theorem 3.1 of \cite{Da_Lio_Ley_2010}:
   In their setting of $\tilde{u}^\mu  (x,t) = \mu e^{-L t} u (x,t)+C(1+|x|^p) $    \(given that $ | u (x,t)  | \le C\big(1+|x|^p\big)$ for some $p>1$\),
                  $L $ is required by their Lemma 3.4    to be dependent on $\mu \in (0,1)$,  see \eqref{eq:a215}.
         This causes a trouble when letting $\mu \to 1$ in $ \mu \tilde{u} -\tilde{v}   \le 0 $ \(or equivalently in \eqref{eq:a018}\) to obtain $  u    \le v  $,
         see    (1) of  Subsection \ref{appendix_3}.
         Another gap is still due to  Lemma 3.4 of \cite{Da_Lio_Ley_2010}, see  (2) of  Subsection \ref{appendix_3}.
         In our proof of Theorem~\ref{thm:viscos_comparison}, we set $\tilde{u}_\th (t,x) = \th e^{\k t} u  (t,x)$, $\th \in (0,1)$. Clearly, the constant  $\k$, assumed in \eqref{b_si_Lip}-\eqref{cond:f_basic}, is independent on $\th$. Moreover, since $\tilde{u}_\th $ does not contain the term $C(1+|x|^\varpi) $
         as in $\tilde{u}^\mu$, we do not need to derive and use the counterparts to Lemma 3.2-3.4 of    \cite{Da_Lio_Ley_2010}.
         \end{rem}

\appendix
\renewcommand{\thesection}{A}
\refstepcounter{section}
\makeatletter
\renewcommand{\theequation}{\thesection.\@arabic\c@equation}
\makeatother

\section{Appendix}

\subsection{Proof of  \eqref{claim1}}  \label{appendix_1}

 \begin{lemm} \label{lem_A1}
 Let $\hB$ be a generic  Banach  space with norm $|\cd|_{\hB}$ and let $p,q \in [1, \infty)$.
 If $\big\{X^n\big\}_{n \in \hN}$ is a sequence of $\hL^{p,q}_\bF\big( [0,T]; \hB \big)$ such that
 it holds \dtp ~ that
 \bea  \label{eq:x121}
 \lmt{n \to \infty} X^n_t = X_t \hb{  \; and\; }  |X^n_t|_\hB \le \cX_t, \q \fa n \in \hN
 \eea
  for some $\hB$-valued,  $\bF$-adapted process $X$  and  some $\cX \in \hL^{p,q}_\bF\big( [0,T]; \hR \big)$,
  then $X   \in \hL^{p,q}_\bF\big( [0,T]; \hB \big)$  and
      $         \|X\|_{\hL^{p,q}_\bF ( [0,T]; \hB )}    \\   = \lmt{n \to \infty} \|X^n\|_{\hL^{p,q}_\bF ( [0,T]; \hB  )}            $.
 \end{lemm}

\ss \no {\bf Proof:}   We  assume that except on a $P$-null set $\sN_1$, \eqref{eq:x121} holds for a.e. $t \in [0,T]$. Since $\cX \in \hL^{p,q}_\bF\big( [0,T]; \hR \big)$, it holds except on     another $P$-null set $\sN_2$   that
 \beas
 \left( \hb{$\int_0^T \neg \cX_t^p    dt$} \right)^{\frac{q}{p}} < \infty,  \q \hb{thus} \q \hb{$\int_0^T \neg \cX_t^p    dt$}   < \infty.
 \eeas
   For any $\o \in \sN^c_1 \cap \sN^c_2$, since it holds for a.e. $t \in [0,T]$ that
    \bea   \label{eq:x125}
  \big| X_t (\o)  \big|^p_\hB = \lmt{n \to \infty}  \big| X^n_t  (\o) \big|^p_\hB
     \q \hb{ and } \q  \big|X^n_t (\o)\big|^p_\hB \le \big(\cX_t (\o) \big)^p, \q \fa n \in \hN,
 \eea
 the Dominated Convergence Theorem implies that
  \beas
  \int_0^T\big| X_t  (\o)  \big|^p_\hB \, dt = \lmt{n \to \infty}  \int_0^T\big| X^n_t  (\o)  \big|^p_\hB \, dt  .
  \eeas
  It also follows from \eqref{eq:x125} that for any $n \in \hN$
  \beas
    \int_0^T\big| X^n_t  (\o)  \big|^p_\hB \, dt  \le \int_0^T  \big( \cX_t (\o) \big)^p    dt   ,  \q \hb{thus} \q
    \left(  \int_0^T\big| X^n_t  (\o)  \big|^p_\hB \, dt \right)^{\frac{q}{p}}   \le \left( \int_0^T  \big( \cX_t (\o) \big)^p    dt    \right)^{\frac{q}{p}}  .
  \eeas
  Applying  the Dominated Convergence Theorem  once again yields that
   \beas
  \hspace{1.2cm} E \left[  \left(  \int_0^T\big| X_t  (\o)  \big|^p_\hB \, dt \right)^{\frac{q}{p}}  \right]
  = \lmt{n \to \infty}   E \left[  \left(  \int_0^T\big| X^n_t  (\o)  \big|^p_\hB \, dt \right)^{\frac{q}{p}}  \right] \le E \left[ \left( \int_0^T  \big( \cX_t (\o) \big)^p    dt    \right)^{\frac{q}{p}} \right] <\infty .  \hspace{1.2cm} \hb{\qed}
   \eeas

  \ms Now, let us prove \eqref{claim1}.   Fix $n \in \hN$. We have seen from \eqref{eq:c261}  that
 $ \left\{\sqrt{ \big| \f '  (Y^m-Y^n  ) \big|} \big(Z^m-Z^n \big)  \right\}_{m \ge n} \neg \subset \hH^2_\bF ([0,T];\hR^d)  $.
 As $m \to \infty$ in \eqref{eq:c401}, the continuity  of function  $\f' $ implies that  \pas
  \bea   \label{eq:c269}
     \big|   \f '  (Y_t-Y^n_t) \big|       \le   e^{ \l_o \g (\sL_t+ \sY_t)  }   ,   \q   t \in [0,T].
  \eea
     Similar to \eqref{eq:c335},  applying  Young's inequality with
     $p_1=\frac{\l}{\l_o}$,  $p_2=\frac{\l'}{\l_o}$ and $p_3 = p_o$,
 we can deduce from    \eqref{eq:c255} that
 \beas
   && \hspace{-0.9cm}E \int_0^T   \big| \f'   (Y_s-Y^n_s ) \big| \,  |Z_s-Z^n_s|^2 ds   \nonumber \\
   &&  \le      c_{\l,\l'}   E \neg \left[  e^{\l_o p_1  \g    \sL_*  }   \neg +    e^{  \l_o p_2  \g      \sY_* }
    \neg + \dneg  \left( \int_0^T   \neg      | Z_s \neg - \neg Z^n_s   |^2          ds \right)^{p_o} \right]
     \le  c_{\l, \l'} \Xi      +      c_{\l,\l'}  \,   E  \neg \left[    \left( \int_0^T     \neg    | Z_s  |^2  ds \right)^{p_o} \right]     < \infty  , \qq
   \eeas
 which implies that  $\sqrt{  \big| \f '  (Y-Y^n  )  \big| } \big(Z-Z^n \big) \in \hH^2_\bF ([0,T];\hR^d)  $. \big(Note that since $Y^n$, $n \in \hN$ are $\bF$-adapted continuous processes, $Y = \lmt{n \to \infty} Y^n $   is at least an $\bF$-predictable process.\big)

  \ss   For any    $ X \in \hH^2_\bF ([0,T];\hR^d)$, we have seen from \eqref{p1_p3} that
    $    \frac{1}{p_1}+ \frac{1}{p_2}+ \frac{1}{p_o}=1$,  or equivalently $ \frac{1}{p_1}+ \frac{1}{p_2}+ 1=2- \frac{1}{p_o}$.
     Applying Young's inequality with $q_1=p_1\big(2- \frac{1}{p_o}\big)    $, $q_2=p_2\big(2- \frac{1}{p_o}\big)  $ and
  $    q_3   =  2- \frac{1}{p_o}  $,
 we see from \eqref{eq:c269} that
   \beas
   \qq &&    \hspace{-1.5cm}        E \left[ \left( \int_0^T   \big| \f '  (Y_s-Y^n_s  ) \big| \, \big|X_s\big|^2  ds  \right)^{\frac{p_o}{2p_o - 1}}\right]  \\
 &&       \le   c_{\l, \l'}   E\left[  e^{ \frac{ \l_o  p_o}{2p_o - 1} q_1  \g    \sL_*  }  +  e^{\frac{ \l_o  p_o}{2p_o - 1}  q_2 \g  \sY_* }    +
      \int_0^T    \big|X_s\big|^2  ds       \right]    \le  c_{\l, \l'}  \Xi + c_{\l, \l'}   E   \int_0^T    \big|X_s\big|^2  ds        < \infty,
 \eeas
 which means that  $X \sqrt{  \big| \f ' (Y-Y^n )  \big|} \in \hH^{2,  \frac{2 p_o}{2 p_o - 1} }_\bF ([0,T];\hR^d) $. Since
 the sequence $\{Z^m\}_{m \ge n}$
 weakly converges to    $Z$ in $\hH^{2, 2 p_o}_\bF([0,T];\hR^d)$, it follows that
  \bea   \label{eq:lmt1}
  \lmt{m \to \infty} E \int_0^T X_s   \sqrt{  \big| \f '  (Y_s-Y^n_s  )  \big| } \left(       Z_s  - Z^m_s  \right) ds  =  0  .
 \eea
    On the other hand,  for any $m \ge n$  H\"older's inequality and \eqref{eq:c255} imply that
  \bea
  \q   \;  && \hspace{-1.5cm}    \left| E \int_0^T X_s \left(    \sqrt{  \big| \f '  (Y_s-Y^n_s  )  \big|}
    - \sqrt{  \big| \f '  (Y^m_s-Y^n_s  )  \big|}  \right) \big(Z^m_s-Z^n_s \big)  ds \right|   \nonumber \\
&&  \le         \left\| \big|X_s \big|  \left(     \sqrt{  \big| \f '  (Y_s-Y^n_s  )  \big|}
     -    \sqrt{  \big|\f '  (Y^m_s-Y^n_s  )  \big|}  \right)  \right\|_{\hH^{2,   \frac{2p_o}{2 p_o - 1} }_\bF ([0,T];\hR)}
    \left\| Z^m -Z^n  \right\|_{\hH^{2, 2 p_o }_\bF ([0,T];\hR^d)} \nonumber  \\
    &&  \le       c_{\l, \l'} \,   \Xi^{\frac{1}{2p_o}}   \left\| \big|X_s \big|  \left(     \sqrt{  \big| \f '  (Y_s-Y^n_s  )  \big|}
     -    \sqrt{  \big|\f '  (Y^m_s-Y^n_s  )  \big|}  \right)  \right\|_{\hH^{2,   \frac{2p_o}{2 p_o - 1} }_\bF ([0,T];\hR)}      .   \label{eq:x131}
 \eea
    It follows from \eqref{eq:c337} that \pas
       \beas
       0 \le    \big|X_t \big|  \left(  \neg \sqrt{ \big| \f ' (Y_t-Y^n_t ) \big| }
           -    \sqrt{\f ' \big| (Y^m_t-Y^n_t  ) \big|}  \right)
        \le     \big|X_t \big|     \sqrt{ \big| \f '  (Y_t-Y^n_t  ) \big| }  ,   \q     \fa   t \in [0,T] , ~ \fa m \ge n    .
   \eeas
 Since  $|X| \sqrt{  \big| \f ' (Y-Y^n )  \big|} \in \hH^{2,  \frac{2 p_o}{2 p_o - 1} }_\bF ([0,T];\hR) $,
 one can deduce from the continuity  of function $\f' $ and Lemma \ref{lem_A1}    that
  \beas
   \lmt{m \to \infty}   \left\| \big|X_s \big|  \left(   \sqrt{ \big| \f '  (Y_s-Y^n_s ) \big| }
     -    \sqrt{ \big|  \f '  (Y^m_s-Y^n_s  ) \big| }  \right)  \right\|_{\hH^{2,   \frac{2p_o}{2 p_o - 1} }_\bF ([0,T];\hR)}   =   0  ,
 \eeas
  which together with \eqref{eq:x131} implies that
  \beas
    \lmt{m \to \infty}   E   \int_0^T X_s \left(     \sqrt{  \big|  \f '  (Y_s-Y^n_s  )  \big|  }
    - \sqrt{  \big|  \f '  (Y^m_s-Y^n_s  ) \big| }  \right) \big(Z^m_s-Z^n_s \big)  ds    = 0 .
 \eeas
 Adding this limit to that in \eqref{eq:lmt1} yields that
  \beas
   \lmt{m \to \infty}  E  \int_0^T X_s \left(    \sqrt{ \big|  \f '  (Y_s-Y^n_s  ) \big|  } \big(Z_s-Z^n_s \big)
    - \sqrt{ \big|  \f '  (Y^m_s-Y^n_s  ) \big|  }   \big(Z^m_s-Z^n_s \big)  \right) ds    = 0   .
  \eeas
Thus   \eqref{claim1} follows.   \qed

\subsection{Comparison Theorem for Quadratic RBSDEs with Bounded Obstacles}

\begin{prop}  \label{prop_compa_bdd}
    Let $(\xi_1,  f_1, L^1 ) $, $(\xi_2,  f_2, L^2 ) $   be two  parameter sets  such that

   \ss \no (\,\,i)  For $j = 1,2$,  $(\xi_j,  L^j )  \in \hL^\infty(\cF_T) \times \hC^\infty_\bF[0,T]$ and $f_i$ satisfy \eqref{cond:f_conti};

  \ss \no (\,ii)  It holds \pas ~ that  $ \xi_1 \le \xi_2  $ and that $L^1_t \le L^2_t$, $\fa t \in [0,T]$;

    \ss \no  (iii) For some  $\g > 0$ and some  function $\ell: \hR \to (0, \infty)$
 with $\int_0^\infty \frac{dx}{\ell(x)} =\infty$, it holds \dtp~ that
\bea \label{eq:x105}
-\ell( y ) - \frac{\g}{2} |z|^2 \le  f_1(t,\o,y,z) \le f_2(t,\o,y,z)  \le \ell( y ) + \frac{\g}{2} |z|^2, \q  \fa  (y,z) \in \hR \times \hR^d .
\eea
 If for $j = 1,2$,  $(Y^j, Z^j, K^j) \in \hC^\infty_\bF[0,T] \times   \hH^2_\bF([0,T];\hR^d) \times \hK_\bF[0,T]$
  be the maximal bounded solution of the RBSDE$(\xi_j, f_j, L^j)$ in the sense of Theorem 1 in \cite{KLQT_RBSDE},
       then it holds \pas ~that   $   Y^1_t \le Y^2_t$ for any  $  t \in [0,T] $.
\end{prop}

 \ss \no {\bf  Proof:} Fix $j \in \{1,2\}$.  Let us first recall the construction of the maximal bounded solution $(Y^j, Z^j, K^j) $
  of the RBSDE$(\xi_j, f_j, L^j)$ from  \cite{KLQT_RBSDE}.  Since $\int_0^\infty \frac{dx}{\ell(x)} =\infty$,
  Lemma 1 of \cite{Lep_San_98} shows that   there exists a unique solution $u^i: [0,T] \to \hR $
   to the following backward ordinary differential equation (BODE for short):
  \beas
 u^j(t) =      b_j   +   \int_t^T   \ell\left(u^j(s) \right) ds   , \q t \in [0,T] ,
\eeas
  where $b_j \dfnn \|\xi_j\|_{\hL^\infty(\cF_T)} \vee \| L^j\|_{\hC^\infty_\bF[0,T]} $.
  Correspondingly, $\wt{u}^i(t) \dfnn  e^{\g u^i(t)}$, $t \in  [0,T] $, uniquely solves the BODE:
  \beas
 \wt{u}^j(t) =     e^{\g b_j}   +   \int_t^T   \wt{\ell}\left(\wt{u}^j(s) \right) ds   , \q t \in [0,T] ,
\eeas
 where  $\wt{\ell} (y)  \dfnn    \b1_{\{y > 0\}} \g y \, \ell\left(\frac{1}{\g} \ln y \right)$, $\fa  y \in \hR $.

 \ms Let $\p: \hR \to [0,1]$ be a smooth function that equals to $1$ inside
 $[r, R]$ and vanishes outside $(r/2, 2R)$ with $r \dfnn \frac12 \exp\big\{-\g  \big(   \| L^1\|_{\hC^\infty_\bF[0,T]}  \vee  \| L^2\|_{\hC^\infty_\bF[0,T]} \big) \big\} $   and $R \dfnn 2 \,  \big(   \wt{u}^1(0) \vee  \wt{u}^2(0) \big)  $.    Clearly, the function
    \beas
    \qq  F^j_\p   (t, \o, y, z) \dfnn 
  \p(y)   \left\{ \g y f_j\left(t, \o, \frac{\ln y}{\g}, \frac{z}{\g y}\right) -\frac12 \frac{|z|^2}{y}  \right\},
   \q \fa (t, \o, y, z) \in [0,T] \times \O \times \hR \times \hR^d
 \eeas
   is  $\sP \times \sB(\hR) \times   \sB(\hR^d)/\sB(\hR)$-measurable  and  satisfies \eqref{cond:f_conti}.  By  \eqref{eq:x105},  it holds \dtp ~ that
    \beas
    - \wt{\ell} (y) -\frac{2}{\,r\,}|z|^2 \le F^j_\p   (t, \o, y, z) \le  \wt{\ell} (y) , \q  \fa ( y, z) \in   \hR \times \hR^d.
    \eeas
 Hence $F^j_\p(\cd,\cd, \rho(\cd), \cd)$ can be approximated by the following decreasing  sequence of functions: For any $n \in \hN$,
   \beas
    F^{j,n}_\p (t, \o, y, z) \dfnn  \wt{\ell}\big(\rho(y)  \big)\big(1-  \pi_n(z)  \big) + \pi_n(z) F^j_\p   \big(t,\o, \rho(y), z \big),
  \q \fa (t, \o, y, z) \in [0,T] \times \O \times \hR \times \hR^d,
  \eeas
 where $\rho: \hR \to (0, \infty)$ and  $\pi_n: \hR^d \to [0,1]$ are two smooth functions such that
    \beas
    \rho(x) = \left\{ \begin{array}{ll}
     r/2 ,  \q &  \hb{if }  x <r/2    ,\\
       x,     \q & \hb{if }  r \le x \le R, \\
     2 R,  \q &  \hb{if }  x > 2R   ,
   \end{array}
   \right.  \q \hb{ and } \q  \pi_n(z)  =   \left\{ \begin{array}{ll}
    1,   \q &   \hb{if }   |z|  \le n, \\
    0,    \q &  \hb{if }   |z|  \ge n+1.
   \end{array}
   \right.
   \eeas
 Clearly,  $F^{j,n}_\p $ is also $\sP \times \sB(\hR) \times   \sB(\hR^d)/\sB(\hR)$-measurable  and  satisfies \eqref{cond:f_conti}. Since it holds \pas ~ that
  \beas
      -\wt{\ell}\big(\rho(y)  \big) -\frac{2}{\,r\,} ( n+1 )^2  &\le&  \wt{\ell}\big(\rho(y)  \big) \big(1- 2 \pi_{n+1}(z)  \big)
      -\frac{2}{\,r\,}|z|^2 \pi_{n+1}(z) \\
      & \le  &F^{j,n+1}_\p   (t, \o, y, z)  \le  F^{j,n}_\p   (t, \o, y, z) \le   \wt{\ell}\big(\rho(y)  \big) , ~\;\,  \fa ( y, z) \in   \hR \times \hR^d,
  \eeas
  we further see that  $F^{j,n}_\p$ is a bounded function.  Thus, \cite{Matoussi_97} shows that the RBSDE$\big(e^{\g \xi_j}, F^{j,n}_\p, e^{\g L^j}\big)$ admits a
      maximal solution $\big(\wt{Y}^{j,n}, \wt{Z}^{j,n},  \wt{K}^{j,n}\big)$. We see from Remark 1 and Lemma 2.2 of \cite{KLQT_RBSDE} that
      $(\wt{u}^j(\cd), 0,0)$ is the unique solution of the RBSDE$(e^{\g b_j}, \wt{\ell} \circ \rho , e^{\g b_j} )$.
    Then Lemma 2.1 of \cite{KLQT_RBSDE} implies that \pas
       \beas
        r \le  e^{\g L^j_t}  \le \wt{Y}^{j,n+1}_t  \le \wt{Y}^{j,n}_t \le  \wt{u}^j(t) \le  \wt{u}^j(0) \le R ~\;\, \hb{and}  ~\;\, \wt{K}^{j,n}_t \le \wt{K}^{j,n+1}_t , \q  t \in [0,T].
       \eeas
       Using the fact that \dtp, $ F^{j,n}_\p(t, \o, y, z) $  converges to $F^j_\p  (t,\o, \rho(y), z ) $ for any $(y,z) \in \hR \times \hR^d$,  the proof of Theorem 2 in \cite{KLQT_RBSDE} shows that
         \bea  \label{limit_triple}
      \wt{Y}^j_t \dfnn \lmtd{n \to \infty} \wt{Y}^{j,n}_t \in [r, R],  \q  \wt{K}^j_t \dfnn  \lmtu{n \to \infty} \wt{K}^{j,n}_t, \q t \in [0,T],
 \eea
 and that the limit $\wt{Z}^j$ of $\left\{\wt{Z}^{j,n}\right\}_{n \in \hN} \subset \hH^2_\bF([0,T];\hR^d) $ constitute a maximal bounded solution of  the RBSDE$\big(e^{\g \xi_j}, \\ F^j_\p, e^{\g L^j}\big) $. Then the proof of Theorem 1 in \cite{KLQT_RBSDE}
indicates that
 \bea  \label{eq:x115}
 (Y^j, Z^j, K^j)  \dfnn \left( \frac{1}{\g} \ln(\wt{Y}^j )   , (\g \wt{Y}^j)^{-1} \wt{Z}^j,  \int_0^\cd (\g \wt{Y}^j_s)^{-1}   d \wt{K}^j_s \right)
 \eea
   is  a maximal bounded solution of the RBSDE$(\xi_j, f_j, L^j)$.

 \ms  For any $n \in \hN$,  it follows from \eqref{eq:x105}  that  \dtp
  \beas
   F^{1,n}_\p(t,\o,y,z) \le F^{2,n}_\p(t,\o,y,z) ,     \q \fa (y,z)  \in \hR \times \hR^d .
  \eeas
 Thanks to  Lemma 2.1 of \cite{KLQT_RBSDE} once again, it holds   \pas ~ that
   $  \wt{Y}^{1,n}_t \le  \wt{Y}^{2,n}_t $ for any  $t \in [0,T] $.
   As $n \to \infty$,  one can deduce from \eqref{limit_triple} and \eqref{eq:x115} that  \pas
 \beas
 \hspace{5cm}  \wt{Y}^1_t \le  \wt{Y}^2_t, \q \hb{thus} \q  Y^1_t \le  Y^2_t, \q \fa   t \in [0,T] .  \hspace{5cm}  \hb{ \qed}
 \eeas

\subsection{Two Gaps in  \cite{Da_Lio_Ley_2010}. } \label{appendix_3}

 \ss \no {\bf (1)} In \cite{Da_Lio_Ley_2010}, the authors fixed a $\mu \in (0,1)$ and chose an  
\bea  \label{eq:a215}
    L > C \big(1+ (1-\mu)^{-p'}\big)
    \eea
     for  some constant $C$
(see  line -5 of   Lemma 3.4   of  \cite{Da_Lio_Ley_2010}).
By setting $  \tilde{u}(x,t) =e^{-L t} u (x,t) +h(x)$ and  $  \tilde{v}(x,t) =e^{-L t} v(x,t) +h(x)$,
they showed  that  $ \mu  \tilde{u}-\tilde{v}   \le 0 $ step by step over each subinterval
$ \hR^N \times \big[\frac{m-1}{L}, \frac{m}{L}\big]$, $m =1, \cds, \lceil  LT \rceil$. Thus
 \bea  \label{eq:a007}
  \mu  \tilde{u} (x,t)-\tilde{v} (x,t)   \le 0  , \q \fa (x,t) \in   \hR^N \times [0,T ].
 \eea
Then they claimed that letting $\mu \to 1$   results in
 \bea   \label{eq:a211}
  u  (x,t) \le v (x,t) ,   \q \fa (x,t) \in    \hR^N \times [0,T ]
  \eea
   (see proof of Theorem 3.1 of \cite{Da_Lio_Ley_2010}).  However, this is  not true.  Actually, \eqref{eq:a007} is equivalent to
 \bea \label{eq:a018}
   \mu u (x,t) - v (x,t)    \le  e^{ L t}  (1-\mu) h(x)  , \q \fa (x,t) \in   \hR^N \times [0,T ].
 \eea
 As $\mu \to 1 $, the right-hand-side of \eqref{eq:a018} goes to $\infty $ for any $(x,t) \in   \hR^N \times (0,T ]$
  since  $L > C \big(1+ (1-\mu)^{-p'} \big)$.  Hence,  their   claim \eqref{eq:a211} fails.

  \ms \no {\bf (2)}    Moreover,  their step-by-step method in  showing that $ \mu  \tilde{u}-\tilde{v}   \le 0 $ holds over each subinterval
$ \hR^N \times \big[\frac{m-1}{L}, \frac{m}{L}\big]$ may also cause a problem. This argument
  requires that  $\cL[\F(x,t)]>0 $   holds  for any $(x,t) \in   \hR^N \times \big(\frac{m-1}{L}, \frac{m}{L}\big]$ (cf. Lemma 3.4 of \cite{Da_Lio_Ley_2010}), which entails that
\bea  \label{eq:a109}
  \dis   \frac{L}{4}-\frac{C_f e^{-Lt}}{\ol{C}} -p^{p'}C^{p'}_s\ol{C}^{p'-1}e^{L p' t} \Big( \frac{e^T}{1-\mu}+1\Big)^{p'}   \ge 1
  \eea
 holds   for any  $t \in \big(\frac{m-1}{L}, \frac{m}{L}\big]$ (see  line -8 of   Lemma 3.4   of  \cite{Da_Lio_Ley_2010}). So it is necessary to have
\bea  \label{eq:a011}
   L> \dis   \frac{4 C_f }{\ol{C}} +4 p^{p'}C^{p'}_s\ol{C}^{p'-1}e^{m p' } \Big( \frac{e^T}{1-\mu}+1\Big)^{p'} + 4 .
  \eea
  When we take $m= \lfloor LT \rfloor $, the right hand side of \eqref{eq:a011} will be much larger than $L$, a contradiction appears.

    \ss However, one might try to change the test function $\F$ in Lemma 3.4 of \cite{Da_Lio_Ley_2010} by some $\F_m$ over $ \hR^N \times \big(\frac{m-1}{L}, \frac{m}{L}\big]$, for example $\F_m(x,t) = \F(x,t -\frac{m-1}{L}) $.  Correspondingly, one has to show that
   \beas
    \cL[\F_m(x,t)]>0, \q \hb{for any } (x,t) \in \hR^N \times \Big(\frac{m-1}{L}, \frac{m}{L}\Big] .
    \eeas
     In the last term of $ \cL[\F_m (x,t)]$ (see the definition of operator $\cL$ in Lemma 3.2 of  \cite{Da_Lio_Ley_2010}), the fourth variable of function $f$  still contains $e^{Lt}$.
    Then similar to line -8 of   Lemma 3.4   of  \cite{Da_Lio_Ley_2010},  the estimation for this function $f$ still results in $e^{Lp't}$ on the right-hand-side,
    which shows that we are facing the same situation as in \eqref{eq:a109}.

\bibliographystyle{siam}
\bibliography{QRBSDE}

\end{document}